\documentclass[12pt, a4paper, reqno]{amsart}
\usepackage{color, graphicx, comment}
\usepackage{tikz}
\usepackage{framed}
\usepackage{mdframed}
\usepackage{algorithm}
\usepackage{algpseudocode}
\usepackage{booktabs}
\usepackage{array}
\usepackage[numbers,sort&compress]{natbib}
\usepackage{amssymb}

\allowdisplaybreaks

\numberwithin{equation}{section} 
\newtheorem{Theorem}{Theorem}[section]
\newtheorem{Corollary}[Theorem]{Corollary}
\newtheorem{Proposition}[Theorem]{Proposition}
\newtheorem{Conjecture}[Theorem]{Conjecture}
\newtheorem{Lemma}[Theorem]{Lemma}
\newtheorem{Remark}[Theorem]{Remark}
\newtheorem{Example}[Theorem]{Example}

\theoremstyle{definition}
\newtheorem{Definition}[Theorem]{Definition}
 
\newcommand{\rmnum}[1]{\romannumeral #1}

\font\KFracFont=cmr12 at 20pt
\def\KK{\mathop{\lower 4pt\hbox{\KFracFont K}}\limits}

\def\cFrac#1#2{%
  \begin{array}{@{}c@{}}\multicolumn{1}{c|}{#1}\\%
  \hline\multicolumn{1}{|c}{#2}\end{array}\;}

\usepackage{epsf}

\usepackage[english]{babel}

\DeclareMathOperator{\setF}{\mathbb{F}}

\title[Shifted Hankel determinants]{Shifted Hankel determinants via Hankel continued fractions}

\author{Guo-Niu Han}
\address{I.R.M.A., UMR 7501, Universit\'e de Strasbourg
et CNRS, 7 rue Ren\'e Descartes, F-67084 Strasbourg, France}
\email{guoniu.han@unistra.fr}
\author[Lihong Yang]{Lihong Yang $^*$}
\address{Institute for Advanced Study in Mathematics,
Harbin Institute of Technology,
Heilongjiang 150001, P.R. China}
\email{lihong.yang.math@gmail.com}
\address{I.R.M.A., UMR 7501, Universit\'e de Strasbourg
et CNRS, 7 rue Ren\'e Descartes, F-67084 Strasbourg, France}
\email{lihong.yang@unistra.fr}
\thanks{$^*$ Lihong Yang is the corresponding author.}


\date{November 16, 2026}

\begin{document}

\begin{abstract}
The Hankel determinants $H_n(F)$ of a power series $F(x)$, as well as 
the shifted ones $H_n^{(1)}(F)$ and $H_n^{(2)}(F)$, can be read off 
from its Jacobi continued fraction, provided that all $H_n(F)$ are 
nonzero. Without any condition on $F$, the Hankel continued 
fraction ($H$-fraction) still yields 
all $H_n(F)$. 
We show that the $H$-fraction of $F$ also determines the shifted 
Hankel determinants $H_n^{(1)}(F)$, $H_n^{(2)}(F)$ and the backward 
shifted ones $H_n^{(-1)}(F)$, and we give explicit formulas 
for them. 
As applications, we provide simple proofs of Hankel determinant 
evaluations that were previously established by means of 
complicated explicit $H$-fractions.
\end{abstract}

\maketitle

\setcounter{tocdepth}{1}
\tableofcontents


\section{Introduction}\label{sec:intro}

\subsection{Hankel determinants and continued fractions}\label{subsec:cf}

Let $\setF$ be a field and $x$ be an indeterminate. 
We identify a sequence ${\bf a}=(a_n)_{n\geq 0}$ over $\setF$ with its generating function 
$f=f(x)=a_0+a_1x+a_2x^2+\cdots\in\setF[[x]]$.
For each $n\geq 1$ and $k\geq 0$ the {\it $k$-shifted Hankel determinant} of the series $f$ (or of the sequence $\bf a$)  is defined by 
\begin{equation}\label{eq:defHnk}
H_n^{(k)} (f) := \left|
\begin{matrix}
 a_k & a_{k+1} & \ldots & a_{k+n-1} \\
a_{k+1} & a_{k+2} & \ldots & a_{k+n} \\
 \vdots  & \vdots  & \ddots &
 \vdots  \\
a_{k+n-1} & a_{k+n} & \ldots & a_{k+2n-2} 
\end{matrix} 
\right|\in\setF,
\end{equation}
with the convention $H_0^{(k)}(f)=1$.
Let $H_n(f):=H_n^{(0)}(f)$, for short; the {\it sequence of the Hankel determinants} of $f$ is
\[H(f):=(H_0(f)=1, H_1(f), H_2(f), H_3(f), \ldots).\] 

Determinant evaluations have always fascinated enumerative combinatorialists, and Hankel determinants of classical combinatorial sequences are among the most attractive ones: 
the matrix is described by a single sequence, the determinant is very 
often a ``nice'' closed formula, and yet the proof is rarely 
straightforward. 
The two surveys of Krattenthaler 
\cite{Krattenthaler1998, Krattenthaler2005} contain a long list of such 
evaluations, together with the main methods for obtaining them.
Shifted Hankel determinants appear in this list from the very 
beginning: for the Catalan numbers $C_n=\frac1{n+1}\binom{2n}{n}$, 
the Hankel determinants of the sequences 
$(C_n)_{n\ge0}$, $(C_{n+1})_{n\ge0}$ and $(C_{n+2})_{n\ge0}$ are all 
``nice'', being equal to $1$, $1$ and $n+1$, respectively 
\cite{Krattenthaler2005}, and, as shown by Cvetkovi\'c, Rajkovi\'c 
and Ivkovi\'c \cite{Cvetkovic2002RI_CatalanHDet}, the Hankel 
determinants of the sums $C_n+C_{n+1}$ of consecutive Catalan numbers 
are Fibonacci numbers; Hankel determinants of linear combinations 
of consecutive terms of more general sequences are studied in 
\cite{Cigler2021Kratt_linear, Mu2017WangYeh_LinearComb}.
Many further examples of shifted Hankel determinants, involving 
Motzkin, Schr\"oder, ternary tree numbers and their relatives, 
can be found in 
\cite{Krattenthaler2005, Gessel2006Xin, Sulanke2005Xin, Cigler2011, Cigler2013spe}.
The sequence can also be shifted in the opposite direction, by prepending one or several terms: the Hankel determinants of the {\it backward shifted} sequences $(c,a_0,a_1,\ldots)$, 
where $c$ is a constant, and more generally of $(a_{n-m})_{n\ge0}$ with prescribed values of $a_{-1},\ldots,a_{-m}$, have been studied more recently, mostly in the case where the 
prepended terms are $0$, by Mu and Wang \cite{Mu2017Wang_ShiftedCatalanLike} and by Cigler 
\cite{Cigler2023BackwardShifts}; we come back to them in Section \ref{subsec:backshift}.

\medskip

As Krattenthaler puts it \cite{Krattenthaler2005}, whenever a Hankel 
determinant evaluates nicely, the explanation is, sometimes more, 
sometimes less, hidden in the theory of continued fractions and 
orthogonal polynomials. 
The {\it continued fraction method} goes back to Heilermann 
\cite{Heilermann1846} and Stieltjes \cite{Stieltjes1894Re}; 
see also \cite{Wall1948, Flajolet1980, Krattenthaler1998, Gessel2006Xin}.
Let $\mathbf{u}=(u_1, u_2, \ldots)$ and $\mathbf{v}=(v_0, v_1, v_2, \ldots)$
be two sequences.
Recall that the {\it Jacobi continued fraction} attached to $(\mathbf{u}, \mathbf{v})$,
or {\it $J$-fraction}, for short, is a continued
fraction of the form
\begin{equation}\label{eq:Jacobi}
f(x)=
  \cFrac{v_0 }{ 1 - u_1 x} 
  - \cFrac{v_1 x^2 }{ 1 - u_2x} 
  - \cFrac{v_2 x^2 }{ 1 - u_3x} 
  - \cdots
\end{equation}
The basic properties of $J$-fractions can be found in
\cite{Krattenthaler1998, Krattenthaler2005, Flajolet1980, Wall1948, Stieltjes1894Re, Heilermann1846, Gessel2006Xin}. We emphasize the fact that
the Hankel determinants can be calculated
from the $J$-fraction by means of the following  fundamental relation, first
stated by Heilermann in 1846 \cite{Heilermann1846}:
\begin{equation}\label{eq:detH}
H_n(f)
= v_0^n v_1^{n-1} v_2^{n-2} \cdots v_{n-2}^2 v_{n-1}. 
\end{equation}
In a seminal paper \cite{Flajolet1980}, Flajolet gave a combinatorial 
interpretation of the $J$-fraction: the coefficient $a_n$ of $f(x)$ is 
the generating polynomial of the Motzkin paths of length $n$ in which 
each level step at height $h$ is weighted by $u_{h+1}$ and each down 
step from height $h$ is weighted by $v_{h}$ (with $v_0$ as an overall 
factor). 
Combined with the Lindstr\"om--Gessel--Viennot lemma on non-intersecting
lattice paths \cite{Lindstrom1973, Gessel1985Viennot}, this
interpretation yields a purely combinatorial proof of \eqref{eq:detH}:
the Hankel determinant $H_n(f)$ counts families of $n$ non-intersecting
Motzkin paths, and exactly one such family exists 
\cite{Viennot1983}. 
This path model has proved extremely fruitful; we only mention here 
its recent use in the study of Hankel-total positivity by Sokal and his
coauthors \cite{Sokal2019EM, Petreolle2023SokalZhu}, 
where the Lindstr\"om--Gessel--Viennot lemma is used to control 
{\it all} minors, and not only the leading principal minors, 
of the Hankel matrix $(a_{i+j})_{i,j\geq 0}$. 

\medskip

The same approach also applies to the shifted Hankel determinants
of a $J$-fraction.
Let $p_n(t)$ be the monic orthogonal polynomials associated with
\eqref{eq:Jacobi}, defined by the three-term recurrence, for $n \ge 1$,
\[
p_{n+1}(t)=(t-u_{n+1})p_n(t)-v_np_{n-1}(t),
\]
with initial conditions $p_0(t)=1$ and $p_1(t)=t-u_1$.
Set $q_n:=(-1)^np_n(0)$, so that $q_0=1$, $q_1=u_1$, and $q_{n+1}=u_{n+1}q_n-v_nq_{n-1}$ for $n\geq 1$.
Then
\begin{align}
H_n^{(1)}(f)&= H_n(f)\, q_n, \label{eq:H1classical}\\
H_n^{(2)}(f)&= H_n(f)\, \sum_{t=0}^{n} q_t^2\, v_{t+1}v_{t+2}\cdots v_n. \label{eq:H2classical}
\end{align}
Both formulas are stated explicitly by Krattenthaler \cite[Theorem 29]{Krattenthaler2005}, who notes that he could not find a reference for them.
The first formula is, however, essentially classical: evaluating the standard determinant representation of $p_n(t)$ in terms of the moments $a_k$ at $t=0$ gives
\[
H_n(f)\,p_n(0)=(-1)^n
\det(a_{i+j+1})_{0\leq i,j\leq n-1},
\]
which yields \eqref{eq:H1classical}.
Krattenthaler also observes that \eqref{eq:H1classical} and \eqref{eq:H2classical} follow readily from the path model.
In this model, $H_n^{(k)}(f)$ is interpreted in terms of families of $n$ non-intersecting Motzkin paths $P_i$ from $(-i,0)$ to $(i+k,0)$.
These paths are uniquely determined outside the strip $0\leq x\leq k$, and the finitely many possible configurations within this strip yield the respective factors $q_n$ and
$\sum_{t=0}^{n}q_t^2v_{t+1}\cdots v_n$. 
In the same spirit, Fulmek \cite{Fulmek2024_CiglerConjecture} recently proved conjectures of Cigler \cite{Cigler2023BackwardShifts, Cigler2023_Experimental} on shifted Hankel determinants of convolution powers of the Catalan numbers by means of non-intersecting lattice paths and the reflection principle. 
Note also that the $1$-shifted Hankel determinants play a fundamental role in the theory of moments: 
by the theorem of Stieltjes, Gantmacher and Krein \cite{Stieltjes1894Re, Gantmacher2002Krein} (see \cite{Petreolle2023SokalZhu} for a modern account), 
a real sequence $(a_n)_{n\geq 0}$ is a Stieltjes moment sequence, or equivalently its Hankel matrix $(a_{i+j})_{i,j\geq 0}$ is totally positive, if and only if both Hankel matrices 
$(a_{i+j})_{i,j\geq 0}$ and $(a_{i+j+1})_{i,j\geq 0}$ are positive semidefinite. 

\subsection{Hankel continued fractions}\label{subsec:hfrac}

The Hankel determinants of a power series $f$ can be calculated by 
the fundamental relation \eqref{eq:detH} only if the $J$-fraction 
of $f$ exists, which is equivalent to the fact that all Hankel 
determinants of $f$ are nonzero.
The same restriction applies to \eqref{eq:H1classical} and 
\eqref{eq:H2classical}.
This restriction is far from harmless: 
for many natural sequences, such as the (mixed) Euler numbers 
\cite{Han2020Trans}, 
the sequences generated by the special $C$-fractions studied by 
Scott and Wall \cite{Scott1940Wall} and by Cigler \cite{Cigler2013spe}, 
the convolution powers of the Narayana polynomials 
\cite{Cigler2018, Han2025NarayanaConv}, 
the $q$-deformed metallic numbers \cite{Han2025Pedon}, 
or the automatic sequences over finite fields considered in 
\cite{Han2016Adv},
some Hankel determinants are zero, 
the $J$-fraction expansion of $f$ does not exist, 
and the continued fraction method in its classical form breaks down.
In this situation one can use the
{\it Hankel continued fraction} ({\it $H$-fraction}, for short), introduced 
by the first author in \cite{Han2016Adv},
whose existence and uniqueness are guaranteed without any condition 
on the power series, and which also yields the Hankel determinants.
Let us recall its definition and basic properties~\cite{Han2016Adv}.

\begin{Definition}\label{def:super}
For each positive integer $\delta$, a {\it super continued fraction} associated with $\delta$, called {\it super $\delta$-fraction} for short, is defined to be a continued fraction of the following form
\begin{equation}\label{eq:super}
F(x)=
\cFrac{v_0 x^{k_0}}{ 1 - u_1(x) x} 
- \cFrac{ v_1 x^{k_0+k_1+\delta}}{ 1 - u_2(x)x}
- \cFrac{ v_2 x^{k_1+k_2+\delta}}{ 1 - u_3(x)x}
- {\cdots} 
\end{equation}
where $v_j\not=0$ are constants, $k_j$ are nonnegative integers and $u_j(x)$ are polynomials of
degree less than or equal to $k_{j-1}+\delta-2$. By convention, $0$ is of degree $-1$.
\end{Definition}

When $\delta=1$ (resp. $\delta=2$) and all $k_j=0$,
the super $\delta$-fraction \eqref{eq:super}
is the traditional $S$-fraction (resp. $J$-fraction).
A super $2$-fraction is called {\it Hankel continued fraction}:
\begin{equation}\label{eq:Hfrac}
F(x)=
\cFrac{v_0 x^{k_0}}{ 1 - u_1(x) x} 
- \cFrac{ v_1 x^{k_0+k_1+2}}{ 1 - u_2(x)x}
- \cFrac{ v_2 x^{k_1+k_2+2}}{ 1 - u_3(x)x}
- {\cdots} 
\end{equation}
where $v_j\not=0$ are constants, $k_j$ are nonnegative integers and $u_j(x)$ are polynomials of
degree less than or equal to $k_{j-1}$. By convention, $0$ is of degree $-1$.

\begin{Remark}\label{rem:sign}
In \cite{Han2016Adv} the super $\delta$-fractions are written with 
the denominators $1+u_j(x)x$ instead of $1-u_j(x)x$.
We adopt here the sign convention of \eqref{eq:Jacobi}, which is 
the traditional one for $J$-fractions \cite{Wall1948, Flajolet1980}.
Passing from one convention to the other simply amounts to replacing
each $u_j(x)$ by $-u_j(x)$; the coefficients $v_j$ and the exponents 
$k_j$ are unchanged, so that all the statements of \cite{Han2016Adv}
about $v_j$ and $k_j$, in particular Theorem \ref{th:super2} below, 
remain valid verbatim.
\end{Remark}

\begin{Theorem}[\cite{Han2016Adv}]\label{th:super2}
(i) Let $\delta$ be a positive integer.
Each super $\delta$-fraction defines
a power series, and conversely, for each power series $F(x)$,
the super $\delta$-fraction expansion of $F(x)$ exists
and is unique.
\smallskip
(ii) Let $F(x)$ be a power series such that its $H$-fraction
  is given by \eqref{eq:Hfrac}.
Then, all non-vanishing Hankel determinants of $F(x)$ are given by
  \begin{equation}\label{eq:HankelDetFundamental}
    H_{s_j}(F(x))= (-1)^{\epsilon_j} v_0^{s_j} v_1^{s_j-s_1} v_2^{s_j-s_2} \cdots v_{j-1}^{s_j-s_{j-1}},
  \end{equation}
where $\epsilon_j = \sum_{i=0}^{j-1} {k_i(k_i+1)/2}$
and
$s_j=k_0+k_1+\cdots + k_{j-1}+j$ for every $j\geq 0$.
\end{Theorem}

In other words, the $H$-fraction of $F(x)$ encodes the full sequence
$H(F)$: the non-vanishing Hankel determinants occur exactly at the
indices $s_0=0<s_1<s_2<\cdots$, and the gaps $s_{j+1}-s_j-1=k_j$ 
between two consecutive non-vanishing terms are the exponents of 
the $H$-fraction. 
When all $k_j=0$, Theorem \ref{th:super2}(ii) reduces to 
Heilermann's formula \eqref{eq:detH}. 
Part (ii) is implicit in the work of Magnus on P-fractions 
\cite{Magnus1962P}: the vanishing pattern is stated there, and the 
values follow from two of his determinant relations. 
It was obtained explicitly by Roblet 
\cite[Proposition 32]{Roblet1994}, in the language of P-fractions, 
by means of non-intersecting lattice paths, and later, by algebraic 
means, by Buslaev \cite[Theorem 1]{Buslaev2010} (up to sign) and, 
for rational functions, by Holtz and Tyaglov 
\cite[(2.18)]{Holtz2012T}; see the concluding remarks. 
The $H$-fraction has been used, for instance, to evaluate the Hankel 
determinants of various automatic sequences and of power series 
satisfying a quadratic equation over a finite field 
\cite{Han2016Adv}, 
and it is the natural tool for the sequences mentioned above whose
Hankel determinants vanish at some indices.

\subsection{Overview of the results}\label{subsec:overview}

The purpose of this paper is to derive explicit formulas for the {\it shifted} Hankel determinants $H_n^{(1)}(F)$, $H_n^{(2)}(F)$, and $H_n^{(-1)}(F)$ directly from the data
$(v_j,k_j,u_j(x))$ of the $H$-fraction of $F(x)$.
Here $H_n^{(-1)}(F)=H_n(xF(x))$, as defined in \eqref{eq:defHm1}.
In particular, we extend \eqref{eq:H1classical} and \eqref{eq:H2classical} from $J$-fractions to $H$-fractions, just as Theorem~\ref{th:super2}(ii) extends Heilermann's
formula \eqref{eq:detH}.
This extension is not immediate.
One could try to evaluate $H_n^{(1)}(F)=H_n((F(x)-a_0)/x)$ by applying Theorem~\ref{th:super2}(ii) to the shifted series.
However, the $H$-fraction of $(F(x)-a_0)/x$ generally bears no simple relation to that of $F(x)$.
The same difficulty arises for $H_n^{(2)}(F)$ and $H_n^{(-1)}(F)$.
Moreover, when some $k_j\geq 1$, the shifted determinants have vanishing patterns different from that of $H_n(F)$, and their evaluations involve additional factors and signs.
The sequence $(q_n)$ in \eqref{eq:H1classical} is replaced by sequences $(R_n)$ and $(E_n)$ built from the coefficients of highest degree in the polynomials $u_j(x)$.
The formula for $H_n^{(2)}(F)$ also takes different forms according to whether $k_j=0$ or $k_j\geq 1$.
Roblet \cite{Roblet1994} had already interpreted $H_n^{(1)}(F)$ and $H_{n+1}(c+xF(x))$ as signed sums over families of non-intersecting paths, but did not evaluate these sums.
We discuss the relation with his work in Section~\ref{sec:remarks}.

\medskip

The main results are stated in Section~\ref{sec:main}. 
Theorems~\ref{th:shiftHDet}, \ref{th:shift2HDet}, and \ref{th:backshiftHDet} give explicit formulas for $H_n^{(1)}(F)$, $H_n^{(2)}(F)$, and $H_n^{(-1)}(F)$, respectively, and show that the corresponding determinants vanish at all indices not covered by these formulas.
Theorem~\ref{th:backshiftHDet} also yields all Hankel determinants of $c+xF(x)$ for an arbitrary constant $c$ (Corollary~\ref{cor:Hc}).
As applications, we give short proofs of the evaluations of the Hankel determinants of the Euler numbers and of the convolution powers of the Narayana polynomials, as well as of a conjecture of Han and Pedon on the $(n+2)$-shifted Hankel determinants of the 
$q$-deformed metallic numbers \cite{Han2025Pedon, Han2026Pedon_nplus2}.
We also prove that the $(n+3)$-shifted ones are unbounded for $n\geq 2$, which is the first case of a related unboundedness conjecture, and we determine the $H$-fraction of an arbitrary super $1$-fraction; see Section \ref{subsec:appl}.

\subsection{Organization of the paper}\label{subsec:method}

We give two proofs of Theorems \ref{th:shiftHDet}--\ref{th:backshiftHDet}. 
The first one, given in Sections \ref{sec:shiftalg} and 
\ref{sec:doubleshalg}, is algebraic, in the spirit of the original 
proof of Theorem \ref{th:super2}(ii) in \cite{Han2016Adv}. 
The second one, given in Sections \ref{sec:flajolet}--\ref{sec:backsh}, 
is combinatorial for Theorems \ref{th:shiftHDet} and 
\ref{th:backshiftHDet}, and partly combinatorial for 
Theorem \ref{th:shift2HDet}. 
It follows the strategy sketched above for 
$J$-fractions: a path model for the coefficients, the lemma of 
Lindstr\"om, Gessel and Viennot for the determinants, 
and a careful analysis of the families of non-intersecting paths 
that survive. Each of the three steps, however, has to be 
adapted to $H$-fractions. 
A more detailed outline of both proofs is given in 
Section \ref{subsec:outline}, after the statements of the main results.

\medskip

The paper is organized as follows. 
In Section \ref{sec:main} we state the main results, specialize them 
to $J$-fractions, outline the proofs and describe the applications. 
The algebraic proof occupies Sections \ref{sec:shiftalg} and 
\ref{sec:doubleshalg}. 
The combinatorial proof begins in Section \ref{sec:flajolet}, where 
Flajolet's path interpretation is extended to $H$-fractions by means 
of quasi-Dyck paths, and in Section \ref{sec:lgv}, where this path 
model is combined with the Lindstr\"om--Gessel--Viennot lemma; 
Theorems \ref{th:shiftHDet}, \ref{th:shift2HDet} and 
\ref{th:backshiftHDet} are then proved in Sections \ref{sec:shift}, 
\ref{sec:doublesh} and \ref{sec:backsh}, respectively. 
Sections \ref{sec:euler}--\ref{sec:gfrac} are devoted to the 
applications. 
Finally, in Section \ref{sec:remarks} we discuss the relation of 
our results with the P-fractions of Magnus, the P-paths of Roblet 
and the formulas of Buslaev and of Holtz and Tyaglov, and we compare the one-step reduction lemma used in the algebraic 
proof with its earlier equivalent forms.


\section{Main results}\label{sec:main}

In this section we state our main results: the $1$-shifted, 
$2$-shifted and backward shifted Hankel determinants of a power 
series in terms of its $H$-fraction (Sections \ref{subsec:main} and 
\ref{subsec:backshift}), and their specializations to $J$-fractions 
(Section \ref{subsec:Jcase}). 
We then outline the two proofs (Section \ref{subsec:outline}) and 
describe the applications (Section \ref{subsec:appl}).

\subsection{The $1$-shifted and $2$-shifted Hankel determinants}\label{subsec:main}

Throughout the paper, $F(x)=\sum_{n\geq0}a_nx^n$ is a power series 
whose $H$-fraction is given by \eqref{eq:Hfrac}, and $s_j$, $\epsilon_j$ are defined as
in Theorem \ref{th:super2}, so that $s_0=0$ and 
$s_{j+1}=s_j+k_j+1$. 
For brevity, we write $V_j := v_0 v_1 \cdots v_j$.
For $j\geq 1$ and $0\leq i\leq k_{j-1}$, the coefficient of $x^i$ in 
$u_j(x)$ is denoted by $c_{j,i}$, so that
\[
u_j(x)=c_{j,0} +c_{j,1}x + c_{j,2}x^2 + \cdots + c_{j,k_{j-1}}x^{k_{j-1}};
\]
in particular, $c_{j,k_{j-1}}$ is the coefficient of $x^{k_{j-1}}$,
the highest admissible degree, in $u_j(x)$; it may be zero.
Define the sequence $(R_n)_{n\ge 0}$ recursively by 
	$R_0=1$, $R_1=c_{1,k_0}$, and for $n\ge 1$,
	\begin{equation}\label{def:Rn}
R_{n+1}= c_{n+1,k_{n}}R_n-v_{n}R_{n-1}.
	\end{equation}
When all $k_j=0$, we have $c_{n+1,0}=u_{n+1}$, and $(R_n)$ is 
the sequence $(q_n)$ of \eqref{eq:H1classical}, that is, 
$R_n=(-1)^np_n(0)$. 
Our first main result evaluates the $1$-shifted Hankel determinants.

\begin{Theorem}\label{th:shiftHDet}
Let $F(x)$ be a power series such that its $H$-fraction
  is given by \eqref{eq:Hfrac}.
Then, for $j \ge 1$, all non-vanishing $1$-shifted Hankel determinants of $F(x)$ are given by
	\begin{equation}\label{eq:H1sj}
H_{s_j}^{(1)}( F )
		= (-1)^{\epsilon_j-s_j+j} v_0^{s_j} v_1^{s_j-s_1} v_2^{s_j-s_2} \cdots v_{j-1}^{s_j-s_{j-1}} R_j
	\end{equation}
and for $j \ge 0$,
	\begin{equation}\label{eq:H1sjkj}
		H_{s_j+k_j}^{(1)}( F ) 
		=  (-1)^{\frac{(k_j-1)k_j}{2}} V_j^{k_j} \cdot H_{s_j}^{(1)}(F), 
	\end{equation}
where we set $H_{0}^{(1)}(F) = 1$.
\end{Theorem}

From \eqref{eq:HankelDetFundamental} and \eqref{eq:H1sj}, we have
\begin{equation}\label{eq:H1sj:Hsj}
H_{s_j}^{(1)}( F )= (-1)^{s_j+j} H_{s_j}(F)  R_j,
\end{equation}
which is the analogue of \eqref{eq:H1classical} for $H$-fractions. 
Note that, contrary to $H_{s_j}(F)$, the determinant 
$H_{s_j}^{(1)}(F)$ may vanish, namely when $R_j=0$. 

Our second main result evaluates the $2$-shifted Hankel determinants.
The formulas depend on whether $k_j=0$ or $k_j\geq 1$.

\begin{Theorem}\label{th:shift2HDet}
Let $F(x)$ be a power series such that its $H$-fraction
  is given by \eqref{eq:Hfrac}.
For $j \ge 0$, the $2$-shifted Hankel determinants of $F$ satisfy the following formulas.

\medskip

\noindent\textbf{(1)} If $k_j\geq 1$, then
\begin{equation}\label{eq:H1sj:Hsj:k1}
	H_{s_j}^{(2)}( F )=  H_{s_j}(F)  R_j^2.
\end{equation}
Moreover,
\begin{equation}\label{eq:H1sj:Hsj:k11}
	H_{s_j+k_j-1}^{(2)}( F )= (-1)^{\frac{(k_j-2)(k_j-1)}{2}} V_j^{k_j-1} \cdot  H_{s_j}^{(2)}(F),
\end{equation}
and
\begin{equation}\label{eq:H1sj:Hsj:k12}
	H_{s_j+k_j}^{(2)}( F )=  c_{j+1,k_j-1} V_j   H_{s_j+k_j-1}^{(2)}(F) + \Delta_{j},
\end{equation}
where 
\[
\Delta_j = 
\left\{
\begin{array}{ll}
 (-1)^{\frac{(k_j+1)k_j}{2}}V_j^{k_j+1} \cdot H_{s_{j-1}}^{(2)}(F) & \text{if $k_{j-1} = 0$},\\
 (-1)^{\frac{(k_j+1)k_j}{2}}V_j^{k_j+1} \cdot H_{s_{j-1}+k_{j-1}}^{(2)}(F)  & \text{if  $k_{j-1} \ge 1$},
\end{array}
\right.
\]
with $\Delta_0=0$.

\medskip

\noindent\textbf{(2)} If $k_j = 0$, then
\begin{equation}\label{eq:H1sj:Hsj:k0}
	H_{s_j}^{(2)}( F )=   H_{s_j}(F)  \sum_{t=0}^{j} d_{t+1} R_t^2 v_{t+1} v_{t+2}  \cdots v_{j},
\end{equation}
where $d_j$ is the coefficient of $x^{k_{j-1}}$ in $1-u_j(x)x$, 
that is, $d_j=1$ if $k_{j-1}=0$ and $d_j=-c_{j,k_{j-1}-1}$ if 
$k_{j-1}\geq 1$. 

\medskip

\noindent\textbf{(3)} All other $2$-shifted Hankel determinants 
vanish: $H_N^{(2)}(F)=0$ for every index $N$ which is not of the 
form $s_j$, $s_j+k_j-1$ or $s_j+k_j$ with $k_j\geq 1$, that is, 
for $N=s_j+\ell$ with $1\leq \ell\leq k_j-2$.
\end{Theorem}

When all $k_j=0$, only case (2) occurs, all $d_j$ are equal to $1$,
and \eqref{eq:H1sj:Hsj:k0} reduces to the classical formula 
\eqref{eq:H2classical}.

\subsection{Backward shifted Hankel determinants}\label{subsec:backshift}

We now consider the Hankel determinants shifted in the 
opposite direction.
The backward shifted sequence $(0, a_0, a_1, a_2, \ldots)$ of 
${\bf a}=(a_n)_{n \ge 0}$ has generating function 
$xF(x)=a_0x+a_1x^2+a_2x^3+\cdots\in\setF[[x]]$. 
The {\it $(-1)$-shifted}, or {\it backward shifted}, Hankel determinant 
of $F$ (or of the sequence ${\bf a}$) is defined by 
\[H_n^{(-1)} (F) := H_n(xF) = \det{[a_{i+j-1}]_{0 \le i,j \le n-1}},\] 
that is,
\begin{equation}\label{eq:defHm1}
H_n^{(-1)} (F) = \left|
\begin{matrix}
 0 & a_{0} & \ldots & a_{n-2} \\
a_{0} & a_{1} & \ldots & a_{n-1} \\
 \vdots  & \vdots  & \ddots &
 \vdots  \\
a_{n-2} & a_{n-1} & \ldots & a_{2n-3} 
\end{matrix} 
\right|\in\setF,
\end{equation}
with the convention $a_{-1}=0$. This is consistent with \eqref{eq:defHnk}
for $k=-1$. 
Note that multiplying \eqref{eq:Hfrac} by $x$ does not yield, 
in general, an $H$-fraction (the exponent $k_0$ is increased by one, while the 
exponents $k_{j-1}+k_j+2$ and the degrees of the $u_j(x)$ are not 
modified), so that the $H$-fraction of $xF(x)$ is not simply related 
to that of $F(x)$, and Theorem \ref{th:super2} cannot be applied 
directly to $xF(x)$.

Compared with the shifts $k\geq 1$, the backward shift has been 
studied only recently. 
Mu and Wang \cite{Mu2017Wang_ShiftedCatalanLike} evaluated the 
Hankel determinants of the shifted sequence $(0,a_0,a_1,a_2,\ldots)$
when $(a_n)_{n\ge 0}$ is a sequence of Catalan-like numbers in the 
sense of Aigner \cite{Aigner1999_Catalanlike}, that is, when $F(x)$ 
is a $J$-fraction with $v_0=1$; as applications, they settled three 
conjectures of Barry on Hankel determinants and obtained Somos-4 
sequences as Hankel determinants of shifted sequences. 
Cigler \cite{Cigler2023BackwardShifts} undertook a systematic study 
of the backward shifts: extending the sequence by $a_n=0$ for 
$n<0$, he considered the determinants $\det(a_{m+i+j})_{0\le i,j\le n-1}$ 
for negative $m$, our $H_n^{(-1)}(F)$ being the case $m=-1$. 
For the Catalan numbers he obtained closed formulas for every $m<0$,
for instance $\det(C_{i+j-1})_{0\le i,j\le n-1}=1-n$, 
together with analogous results for several related sequences and 
polynomials, and he conjectured identities relating the determinants
with negative shift $m$ to those with positive shift for the 
convolution powers of the Catalan numbers 
\cite{Cigler2023BackwardShifts, Cigler2023_Experimental}; 
these conjectures were proved in 
\cite{Fulmek2024_CiglerConjecture, Cigler2024ConvCatalan} 
and extended in \cite{Liu2025WangZhangZhang_Convolution} to 
convolution powers of power series satisfying a quadratic equation.
For the Catalan-like numbers $a_n(s)$ of Aigner with all 
$v_j=1$, Cigler \cite{Cigler2024CatalanLike_Backward} proved that 
\[
\det(a_{i+j-m}(s))_{0\le i,j\le n+m}=(-1)^{\binom{m+1}2}
\det(a_{i+j+m}(Es))_{0\le i,j\le n-1},
\]
where $Es$ is the parameter 
sequence $s$ shifted by one; the case $m=1$ is an instance of 
\eqref{eq:Hm1classical} below. 
Backward shifts of the coefficients of a partial theta function are 
treated in \cite{Cigler2024PartialTheta}. 
All these results concern specific sequences, or sequences whose 
generating function satisfies a quadratic equation. 
Our third main result evaluates $H_n^{(-1)}(F)$ for an arbitrary 
power series $F(x)$ in terms of its $H$-fraction. 

\begin{Theorem}\label{th:backshiftHDet}
Let $F(x)$ be a power series such that its $H$-fraction
  is given by \eqref{eq:Hfrac}.
Define the sequence $(E_n)_{n\ge 0}$ recursively by 
$E_0=0$, $E_1=1$, and for $n\ge 1$,
\begin{equation}\label{def:En}
    E_{n+1} = c_{n+1,k_n}E_{n} - v_n E_{n-1},
\end{equation}
Then, for $j\ge 1$, all non-vanishing backward shifted Hankel 
determinants of $F(x)$ are given by
\begin{equation}\label{eq:Hsj(-1)}
        H_{s_j +1}^{(-1)}(F) = (-1)^{\epsilon_j-s_j+j+1} v_0^{s_j +1} v_1^{s_j-s_1}v_2^{s_j-s_2}\cdots v_{j-1}^{s_{j}-s_{j-1}} E_j
    \end{equation} 
and
    \begin{equation}\label{eq:Hs_j+1(-1)}
        H_{s_{j+1}}^{(-1)}(F) = (-1)^{\frac{k_j(k_j-1)}{2}}(v_0 v_1 \cdots v_j)^{k_j}\cdot  H_{s_j +1}^{(-1)}(F).
    \end{equation}
\end{Theorem}

Note that $E_j$ is the analogue of $R_{j-1}$ for the tail
\[
  \cFrac{v_1 x^{k_1}}{ 1 - u_2(x) x} 
- \cFrac{ v_2 x^{k_1+k_2+2}}{ 1 - u_3(x)x}
- \cFrac{ v_3 x^{k_2+k_3+2}}{ 1 - u_4(x)x}
- {\cdots} 
\]
of the $H$-fraction \eqref{eq:Hfrac}: the sequence $(E_n)$ satisfies 
the same three-term recurrence as $(R_n)$, and 
$E_1=1$, $E_2=c_{2,k_1}$, $E_3=c_{2,k_1}c_{3,k_2}-v_2$, \ldots
Comparing \eqref{eq:Hsj(-1)} with \eqref{eq:HankelDetFundamental}
and \eqref{eq:H1sj}, we get
\begin{align}
H_{s_j+1}^{(-1)}(F)&= (-1)^{s_j+j+1}\, v_0\, H_{s_j}(F)\, E_j, \label{eq:Hsjm1:Hsj}\\
R_j\, H_{s_j+1}^{(-1)}(F)&= -v_0\, E_j\, H_{s_j}^{(1)}(F). \label{eq:Hsjm1:H1sj}
\end{align}
The determinant $H_2^{(-1)}(F)=-a_0^2$ is the case $j=1$ of 
\eqref{eq:Hsj(-1)} when $k_0=0$; when $k_0\geq 1$, 
the first non-vanishing backward shifted Hankel determinant is 
$H_{s_1+1}^{(-1)}(F)=(-1)^{1+k_0(k_0-1)/2}v_0^{k_0+2}$.

When all $k_j=0$, the series $F(x)=f(x)$ is the
$J$-fraction \eqref{eq:Jacobi}.
In this case, $s_j=j$, $\epsilon_j=0$, and
$c_{n+1,0}=u_{n+1}$.
Thus, for $n\ge1$, \eqref{eq:Hsjm1:Hsj} becomes
\begin{equation}\label{eq:Hm1classical}
H_n^{(-1)}(f) = -v_0H_{n-1}(f)E_{n-1},
\end{equation}
where the case $n=1$ follows from $H_1^{(-1)}(f)=E_0=0$.
Here $E_0=0$, $E_1=1$, and, for $n\ge1$, $E_{n+1}=u_{n+1}E_n-v_nE_{n-1}.$
Let $p_n^{(1)}(t)$ denote the associated orthogonal
polynomials of the $J$-fraction, also called numerator polynomials.
They are defined by $p_0^{(1)}(t)=1$,
$p_1^{(1)}(t)=t-u_2$, and, for $n\ge1$,
\[
p_{n+1}^{(1)}(t) = (t-u_{n+2})p_n^{(1)}(t)-v_{n+1}p_{n-1}^{(1)}(t).
\]
Comparing the two recurrences gives, for $n\ge1$, $E_n=(-1)^{n-1}p_{n-1}^{(1)}(0).$
For the Catalan numbers, we have $u_1=1$ and $u_n=2$ for $n\ge2$, while $v_n=1$ for $n\ge0$.
Hence $E_n=n$.
Since the ordinary Hankel determinants of the Catalan numbers are all equal to $1$, \eqref{eq:Hm1classical} gives, for $n\ge1$, $H_n^{(-1)}(C)=1-n,$ in agreement with \cite{Cigler2023BackwardShifts}.

\medskip

Theorem \ref{th:backshiftHDet} also covers the case where an 
arbitrary constant $c\in\setF$, instead of $0$, is prepended to 
the sequence, that is, the series $c+xF(x)$ with coefficient 
sequence $(c,a_0,a_1,a_2,\ldots)$. 
In the Hankel matrix of $c+xF(x)$ the constant $c$ occurs only 
in the $(0,0)$-entry, and the cofactor of this entry is 
$\det(a_{i+j+1})_{0\le i,j\le n-2}=H_{n-1}^{(1)}(F)$. 
For $n\ge1$, expanding the determinant along its first row gives
\begin{equation}\label{eq:Hc}
H_n(c+xF(x))
=
c\,H_{n-1}^{(1)}(F)+H_n^{(-1)}(F),
\end{equation}
so that Theorems \ref{th:shiftHDet} and \ref{th:backshiftHDet} 
together give all the Hankel determinants of $c+xF(x)$, 
including their vanishing pattern, in terms of the $H$-fraction 
of $F(x)$. 
Indeed, the non-vanishing $1$-shifted Hankel determinants 
$H^{(1)}_{n-1}(F)$ and the non-vanishing backward shifted 
Hankel determinants $H^{(-1)}_n(F)$ occur at the same indices, 
namely $n-1\in\{s_j,\,s_j+k_j\}$, and substituting 
\eqref{eq:H1sj:Hsj}, \eqref{eq:H1sjkj}, \eqref{eq:Hsjm1:Hsj} 
and \eqref{eq:Hs_j+1(-1)} into \eqref{eq:Hc} we obtain 
the following statement.

\begin{Corollary}\label{cor:Hc}
Let $F(x)$ be a power series whose $H$-fraction is given by 
\eqref{eq:Hfrac}, and let $c\in\setF$. 
Then, for $j\geq 0$,
\begin{align}
H_{s_j+1}(c+xF(x)) &= (-1)^{\epsilon_j+s_j+j}\, 
v_0^{s_j} v_1^{s_j-s_1} \cdots v_{j-1}^{s_j-s_{j-1}}
\,\bigl(c\,R_j - v_0\,E_j\bigr), 
\label{eq:Hc1}\\
H_{s_{j+1}}(c+xF(x)) &= (-1)^{\frac{k_j(k_j-1)}{2}}\, V_j^{k_j}\, H_{s_j+1}(c+xF(x)), 
\label{eq:Hc2}
\end{align}
and $H_n(c+xF(x))=0$ for every other index $n\geq1$, that is, 
for $n=s_j+\ell$ with $2\leq\ell\leq k_j$.
\end{Corollary}

By \eqref{eq:HankelDetFundamental}, the right-hand side of 
\eqref{eq:Hc1} is $(-1)^{s_j+j}H_{s_j}(F)\,(cR_j-v_0E_j)$, 
and the two formulas coincide when $k_j=0$. 
For $j=0$, \eqref{eq:Hc1} reads $H_1(c+xF(x))=c$, 
in accordance with $R_0=1$ and $E_0=0$. 
When all $k_j=0$, that is, when $F(x)=f(x)$ is the $J$-fraction 
\eqref{eq:Jacobi}, \eqref{eq:Hc1} becomes 
\[
H_n(c+xf(x))=H_{n-1}(f)\,(c\,R_{n-1}-v_0\,E_{n-1})
\] for $n\geq1$. 
Note that, in general, the $H$-fraction of $c+xF(x)$ is not simply 
related to that of $F(x)$, so that these determinants are not 
easily accessible through Theorem \ref{th:super2}. 

\subsection{The case of $J$-fractions}\label{subsec:Jcase}

For the convenience of the reader, we collect here the 
specializations of our results to the case where all $k_j=0$, 
that is, where $F(x)=f(x)$ is the $J$-fraction \eqref{eq:Jacobi}. 
Then $s_j=j$, $\epsilon_j=0$, $V_j=v_0v_1\cdots v_j$, 
$c_{j,0}=u_j$ and $d_j=1$, and the sequences $(R_n)$ and $(E_n)$ 
of \eqref{def:Rn} and \eqref{def:En} are defined by 
$R_0=1$, $R_1=u_1$, $E_0=0$, $E_1=1$ and the common recurrence 
$X_{n+1}=u_{n+1}X_n-v_nX_{n-1}$ $(n\geq1)$, so that 
\begin{equation}\label{eq:RdetJ}
R_n=\det\begin{pmatrix}
u_1 & 1 & \cdots & 0\\
v_1 & u_2 & \ddots & \vdots\\
\vdots & \ddots & \ddots & 1\\
0 & \cdots & v_{n-1} & u_n
\end{pmatrix},
\ \ 
E_n=\det\begin{pmatrix}
u_2 & 1 & \cdots & 0\\
v_2 & u_3 & \ddots & \vdots\\
\vdots & \ddots & \ddots & 1\\
0 & \cdots & v_{n-1} & u_n
\end{pmatrix},
\end{equation}
the determinants being of orders $n$ and $n-1$, respectively 
(expand along the last row to recover the recurrence). 
In terms of the orthogonal polynomials of Section \ref{subsec:cf}, 
$R_n=(-1)^np_n(0)$ and $E_n=(-1)^{n-1}p^{(1)}_{n-1}(0)$. 

\begin{Corollary}[$J$-fractions]\label{cor:J}
Let $f(x)$ be the $J$-fraction \eqref{eq:Jacobi}, with all 
$v_j\neq0$, let $(R_n)$, $(E_n)$ be given by \eqref{eq:RdetJ}, 
and let $c\in\setF$. Then, for $n\geq1$,
\begin{align}
H_n(f)&=v_0^{\,n}\,v_1^{\,n-1}\cdots v_{n-2}^{\,2}\,v_{n-1},\label{eq:J0}\\
H^{(1)}_n(f)&=H_n(f)\,R_n,\label{eq:J1}\\
H^{(2)}_n(f)&=H_n(f)\sum_{t=0}^{n}R_t^2\,v_{t+1}v_{t+2}\cdots v_n,\label{eq:J2}\\
H^{(-1)}_n(f)&=-\,v_0\,H_{n-1}(f)\,E_{n-1},\label{eq:Jm1}\\
H_n(c+xf(x))&=H_{n-1}(f)\,\bigl(c\,R_{n-1}-v_0\,E_{n-1}\bigr),\label{eq:Jc}\\
H_{2n}\bigl(f(x^2)\bigr)&=H_n(f)^2\,R_n,\label{eq:Jx2a}\\
H_{2n+1}\bigl(f(x^2)\bigr)&=H_n(f)\,H_{n+1}(f)\,R_n.\label{eq:Jx2b}
\end{align}
\end{Corollary}

Formula \eqref{eq:J0} is Heilermann's formula \eqref{eq:detH}, 
\eqref{eq:J1} and \eqref{eq:J2} are \eqref{eq:H1classical} 
and \eqref{eq:H2classical}, that is, 
\cite[Theorem 29]{Krattenthaler2005}, 
\eqref{eq:Jm1} is \eqref{eq:Hm1classical}, and 
\eqref{eq:Jc} is Corollary \ref{cor:Hc}. 
The last two formulas follow from \eqref{eq:J0} and \eqref{eq:J1} by the elementary identities 
\begin{align*}
    H_{2n}(f(x^2))&=H_n(f)\,H^{(1)}_n(f),\\
    H_{2n+1}(f(x^2))&=H_{n+1}(f)\,H^{(1)}_n(f),
\end{align*}
valid for every power series $f$ (Lemma \ref{gf:lem:x2} below). 

\subsection{Outline of the proofs}\label{subsec:outline}

We now describe in more detail the two proofs announced in 
Section \ref{subsec:method}.
In \cite{Han2016Adv}, Theorem \ref{th:super2}(ii) is derived from 
a one-step reduction lemma, which expresses $H_n(F)$ in terms of the 
Hankel determinants of the tail of the $H$-fraction; equivalent 
forms of this lemma had appeared earlier 
\cite{Andrews2002Wimp, Sulanke2005Xin, Holtz2012T, Cigler2013spe} 
(see the concluding remarks). 
In Section \ref{sec:shiftalg} we apply this lemma to $(F(x)-a_0)/x$ 
and to $xF(x)$, and obtain a coupled system: $H^{(1)}_n(F)$ is 
expressed in terms of the $1$-shifted and the backward shifted 
Hankel determinants of the tail, and $H^{(-1)}_n(F)$ in terms of 
the $1$-shifted ones only. 
Theorems \ref{th:shiftHDet} and \ref{th:backshiftHDet} are then 
proved together by induction; in the induction step, the sequences 
$(R_n)$ and $(E_n)$ appear through the expansion of a continuant 
along its first row, $E_j$ being the continuant $R_{j-1}$ of the tail. 
In Section \ref{sec:doubleshalg} we derive Theorem \ref{th:shift2HDet} 
from Theorems \ref{th:super2} and \ref{th:shiftHDet}, by means of 
the Desnanot--Jacobi identity and of the polynomials orthogonal with 
respect to $F$, namely the reversed denominators of the convergents of 
the $H$-fraction. 
They are only block-orthogonal, the blocks having sizes $k_j+1$, 
and the sum in \eqref{eq:H1sj:Hsj:k0} appears as a sum over the 
diagonal blocks of their Gram matrix. 

The combinatorial part begins in Section \ref{sec:flajolet}, where 
we extend Flajolet's interpretation to 
$H$-fractions. The coefficients of $F(x)$ count weighted 
{\it quasi-Dyck paths}: lattice paths with unit up steps and with 
down steps of several prescribed lengths, the lengths being 
determined by the exponents $k_j$; the polynomials $u_j(x)$ 
give the weights of the ``short'' down steps and the constants 
$v_j$ those of the ``long'' ones. 
When all $k_j=0$ these paths are in bijection with the Motzkin 
paths of Flajolet.
The proof of the correspondence relies on the first return 
decomposition of the paths, and the correspondence is then lifted 
to series in non-commutative variables in the sense of 
Sch\"utzenberger.

In Section \ref{sec:lgv} we combine this path model with the 
Lindstr\"om--Gessel--Viennot lemma. The Hankel determinant 
$H_n^{(k)}(F)$ then becomes a signed sum over families of $n$ pairwise 
non-intersecting quasi-Dyck paths, whose starting and ending points 
are shifted according to $k$.
As a first application we obtain a combinatorial proof of Theorem 
\ref{th:super2}(ii), proved by algebraic means in 
\cite{Han2016Adv}: for $k=0$ the family of non-intersecting paths 
is unique, and its weight and sign are precisely the right-hand side 
of \eqref{eq:HankelDetFundamental}. This is the $H$-fraction 
analogue of Viennot's proof of Heilermann's formula 
\cite{Viennot1983}. 
In fact, an equivalent path model, the P-paths, had been introduced 
earlier by Roblet \cite{Roblet1994} in his combinatorial study of 
Pad\'e approximants, where the same determinants are evaluated 
by essentially the same argument; see the concluding remarks.

Theorem \ref{th:shiftHDet} is proved combinatorially in 
Section \ref{sec:shift}, and Theorem \ref{th:backshiftHDet} in 
Section \ref{sec:backsh}. 
For $k=1,-1$ the families of non-intersecting paths are no longer 
unique, but the rigid structure of quasi-Dyck paths (a long down 
step has to start and end at prescribed heights) forces most of 
these families to be of a very special form: all paths are fixed, 
except one ``boundary'' path for each level, and two neighbouring 
boundary paths may only be exchanged. 
This explains both the vanishing of $H_n^{(k)}(F)$ at the indices 
not listed in the theorems and the appearance of the three-term 
recurrences defining $(R_n)$ and $(E_n)$: a boundary path either 
uses the short down step of maximal drop at its top, which produces 
the coefficient $c_{j,k_{j-1}}$ of the recurrence, or is exchanged 
with a neighbouring boundary path, which produces the 
coefficient $-v_{j-1}$. 
As in the $J$-fraction case, the only freedom left to the paths is 
thus concentrated near the origin, but the possible configurations 
now depend on the exponents $k_j$, which produces the case 
distinctions and the signs in the theorems.
The $2$-shifted case, treated in Section \ref{sec:doublesh}, 
is the most delicate one, and our proof is only partly 
combinatorial: the paths give the determinants 
$H^{(2)}_{s_j+k_j-1}(F)$ and $H^{(2)}_{s_j+k_j}(F)$ and the vanishing 
ones, and the determinants $H^{(2)}_{s_j}(F)$ are then derived from 
these and from Theorems \ref{th:super2} and \ref{th:shiftHDet} by the 
Desnanot--Jacobi identity.

\subsection{Applications}\label{subsec:appl}

The last part of the paper, Sections \ref{sec:euler}--\ref{sec:gfrac}, 
is devoted to applications. 
Theorems \ref{th:shiftHDet}--\ref{th:backshiftHDet} reduce the 
evaluation of the shifted Hankel determinants of $F(x)$ to a routine 
computation, as soon as the $H$-fraction of $F(x)$ itself is known. 
This is a real saving: so far, each shift of a sequence required 
a new $H$-fraction, to be guessed and then proved, and the 
$H$-fraction of a shifted series is often much more complicated than 
that of the series itself. 
We illustrate this on three families of sequences whose shifted 
Hankel determinants have been the subject of recent papers, and for 
which the previously known proofs are long.

In Section \ref{sec:euler} we consider the Euler numbers $E_n$, 
defined by $\tan x+\sec x=\sum_{n\ge0} E_n x^n/n!$ 
(in that section, the sequence $(E_n)$ of 
Theorem \ref{th:backshiftHDet} is renamed $(T_n)$). 
The Hankel determinants of $(E_0,E_1,E_2,\ldots)$ and of 
$(E_2,E_3,E_4,\ldots)$ were evaluated in \cite{Han2020Trans} by 
means of the $H$-fractions of their generating functions; 
both $H$-fractions have a three-term periodic pattern, and their 
proofs are long. 
We rederive both evaluations from the continued fraction of the 
intermediate sequence $(E_1,E_2,E_3,\ldots)$, which is a 
$J$-fraction and is much simpler: the first sequence is obtained 
from it by prepending $E_0=1$, and the second one is its $1$-shift, 
so that Theorem \ref{th:backshiftHDet}, in the form \eqref{eq:Hc}, 
and Theorem \ref{th:shiftHDet} apply.

In Section \ref{sec:cigler} we consider the convolution powers of 
the Narayana polynomials $\gamma_n(t)$, that is, the coefficients 
of the series $(\gamma(t,q)-1)^m/q^s$, where 
$\gamma(t,q)=\sum_{n\ge0}\gamma_n(t)q^n$. 
Their Hankel determinants were studied by Cigler \cite{Cigler2018}, 
who evaluated them for $m=2$ and conjectured their value for 
$m=3$, and were evaluated for all $m$ and $s=0,1,2,3$ in 
\cite{Han2025NarayanaConv} by establishing four separate 
$H$-fractions, three of which were obtained by a guess-and-prove 
strategy involving very large intermediate expressions. 
The $H$-fraction of the series for $s=1$ has all its exponents 
equal to $k_j=m-1$, so that, for $m\ge2$, it is a genuine 
$H$-fraction, whose Hankel determinants vanish at all indices not 
divisible by $m$, and the case distinctions of 
Theorem \ref{th:shift2HDet} fully come into play. 
We show that the evaluations for $s=2$, $s=3$ and $s=0$ all follow 
from this single $H$-fraction, which is by far the simplest of the 
four, by Theorems \ref{th:shiftHDet}, \ref{th:shift2HDet} and 
\ref{th:backshiftHDet}, respectively.

In Sections \ref{sec:metallicn2} and \ref{sec:metallicn3} we 
consider the $q$-deformed metallic numbers $\Phi_n(q)$ of 
Morier-Genoud and Ovsienko \cite{MorierGenoud2022Ovsienko}, which are $q$-analogues 
of the quadratic irrationals $(n+\sqrt{n^2+4})/2$. 
In \cite{Han2025Pedon} the first author and Pedon established the 
$H$-fraction of $\Phi_n$, which, for $n\ge2$, is periodic with 
period $6n-4$, 
and deduced that the shifted Hankel determinant sequences 
$H^{(\ell)}(\Phi_n)=(H_j^{(\ell)}(\Phi_n))_{j\ge0}$, 
$0\le\ell\le n+1$, are (anti-)periodic and take only the values 
$-1,0,1$; they conjectured that $H^{(n+2)}(\Phi_n)$ is still 
(anti-)periodic, with values in $\{-2,-1,0,1,2\}$, and that 
$H^{(\ell)}(\Phi_n)$ is unbounded for every $\ell\ge n+3$. 
Their method relies on the fact that the $H$-fractions of the 
shifts $\Phi_n^{(\ell)}$, $0\le\ell\le n+1$, are truncations of 
that of $\Phi_n$, and it does not extend to $\ell=n+2$. 
Since $H_j^{(n+2)}(\Phi_n)=H_j^{(1)}(\Phi_n^{(n+1)})$ and the 
$H$-fraction of $\Phi_n^{(n+1)}$ is known, 
Theorem \ref{th:shiftHDet} applies, and we give a short proof of 
the conjecture for $\ell=n+2$ (Theorem \ref{met:thm-main}): the sequence $(R_n)$ 
satisfies a three-term recurrence with periodic coefficients, whose 
analysis explains both the (anti-)periodicity and the appearance of 
the values $\pm2$. 
(The case $n=1$ was proved by Ovsienko and Pedon 
\cite{Ovsienko2025Pedon}; our proof covers all $n\ge2$.) 
The conjecture for $\ell=n+2$ has also been proved, for all 
$n\ge1$, in \cite{Han2026Pedon_nplus2}, by computing the 
$H$-fraction of $\Phi_n^{(n+2)}$ itself, which is eventually 
periodic; the verification of this $H$-fraction is carried out by 
a symbolic computation program. Our proof only needs the 
$H$-fraction of $\Phi_n^{(n+1)}$, which is a tail of that of 
$\Phi_n$, and involves no computer verification. 
Similarly, $H_j^{(n+3)}(\Phi_n)=H_j^{(2)}(\Phi_n^{(n+1)})$ is 
given by Theorem \ref{th:shift2HDet}, and we prove that, 
for $n\ge2$, this sequence is unbounded, with an explicit linear growth 
along an arithmetic progression of indices 
(Theorem \ref{met:thm-n+3}); together with the case $n=1$, 
proved by Pedon \cite{Pedon2026private}, this settles the first case, 
$\ell=n+3$, of the unboundedness conjecture.

Finally, in Section \ref{sec:gfrac} we turn to the super 
$1$-fractions of Definition \ref{def:super}, which generalize 
the $S$-fractions just as $H$-fractions generalize 
$J$-fractions. Classically, the even contraction of the $S$-fraction 
of a power series is its $J$-fraction. 
We determine explicitly the $H$-fraction of an arbitrary 
super $1$-fraction, hence all its Hankel determinants: 
each pair of consecutive levels of the super $1$-fraction 
is either merged into a single level of the $H$-fraction, 
as in the even contraction, or simply rewritten as two levels, 
according to whether the exponent of the second level of the 
pair is zero or not. 
The super $1$-fractions for which only the first case occurs, 
so that the $H$-fraction is exactly the even contraction, 
we call {\it $G$-fractions}; they contain the $S$-fractions. 
Theorem \ref{th:shiftHDet} yields a characterization 
in the spirit of Stieltjes: a super $1$-fraction is a $G$-fraction 
if and only if the $1$-shifted Hankel determinants 
$H^{(1)}_{s_j}(F)$ at all the indices $s_j$ where $H_{s_j}(F)\neq0$ 
are nonzero, and in that case the sequence $(R_n)$ is the 
product of the numerator constants of odd index. 
Theorems \ref{th:shiftHDet} and \ref{th:shift2HDet} then give 
the shifted and doubly shifted Hankel determinants of a 
$G$-fraction in closed form in terms of its own data. 
Specializing to $S$-fractions, we recover the classical 
evaluations of $H_n(F)$ and $H^{(1)}_n(F)$ due to Stieltjes, 
together with $H^{(2)}_n(F)$, $H^{(-1)}_n(F)$ and $H_n(F(x^2))$, 
the corresponding formulas for $J$-fractions being those of 
Corollary \ref{cor:J}.


\section{An algebraic proof of Theorems \ref{th:shiftHDet} and \ref{th:backshiftHDet}}\label{sec:shiftalg}

In \cite{Han2016Adv}, Theorem \ref{th:super2}(ii) is derived from a one-step reduction lemma, which expresses the Hankel determinants of $F$ by those of the tail of its $H$-fraction.
In this section we show that the same lemma, applied to suitable modifications of $F$, also yields Theorem \ref{th:shiftHDet}.
The point is that the $1$-shift of $F$ is expressed by the $1$-shift {\it and} the $(-1)$-shift of the tail, while the $(-1)$-shift of $F$ is expressed by the $1$-shift of the tail alone (Proposition
\ref{prop:alg:coupled} below).
Consequently, Theorems \ref{th:shiftHDet} and \ref{th:backshiftHDet} are proved together.

\subsection{Three reduction lemmas}

We first restate \cite[Lemma 2.2]{Han2016Adv} with the sign convention of Remark~\ref{rem:sign} and with a leading coefficient $v$;
replacing $u(x)$ by $-u(x)$ and $F$ by $F/v$ in \cite[Lemma 2.2]{Han2016Adv}, we get the following statement.
Note that $G$ is an arbitrary power series; it may be zero.

\begin{Lemma}[\cite{Han2016Adv}]\label{lem:alg:han}
Let $k\ge0$, let $v\neq0$ be a constant, let $u(x)$ be a polynomial of degree at most $k$, let $G(x)$ be a power series, and let
\begin{equation}\label{eq:alg:onestep}
F(x)=\frac{v\,x^k}{1-u(x)x-x^{k+2}G(x)}.
\end{equation}
Then $H_n(F)=0$ for $1\le n\le k$, and
\[
H_n(F)=(-1)^{k(k+1)/2}\,v^n\,H_{n-k-1}(G)\qquad (n\ge k+1).
\]
\end{Lemma}

The vanishing for $1\le n\le k$ is clear, since then the first row $(a_0,\ldots,a_{n-1})$ of the Hankel matrix is zero.

\begin{Remark}\label{rem:alg:hanrefs}
Equivalent forms of Lemma \ref{lem:alg:han} appeared earlier.
The case $k=0$ is \cite[Theorem 1]{Andrews2002Wimp}, written as $s(x)t(x)=s_0$.
The general case is \cite[Proposition 4.2(ii)]{Sulanke2005Xin}, stated for series satisfying a quadratic equation, but proved in general; \cite[Proposition 4.2(i)]{Sulanke2005Xin} is the variant with a $1$-shifted determinant of $G$.
It is also \cite[Lemma 2.1]{Cigler2013spe}, where $x^{-k}F$ is written as $1/(1-ax^pg(x))$ and $m=k$; the case $m=-1$ of that lemma is Lemma \ref{lem:alg:KxK} below.
Finally, it is \cite[Theorem 1.18]{Holtz2012T}, stated there for rational functions.
See the concluding remarks for a detailed comparison.
\end{Remark}

\begin{Lemma}\label{lem:alg:KxK}
For every power series $K(x)$ and every $n\ge0$ we have
\[
H_n\Bigl(\frac{K(x)}{1-xK(x)}\Bigr)=H_n\bigl(K(x)\bigr).
\]
\end{Lemma}

\begin{proof}
If $K(x)$ is the zero series, both sides are equal to $0$ for $n\ge1$.
Otherwise, let
\[
K(x)=\frac{w x^{m}}{1-u(x)x-x^{m+2}L(x)}
\]
be the first step of the $H$-fraction of $K$, where $w\neq0$, $\deg u\le m$ and $L$ is a power series.
Then
\[
\frac{K}{1-xK}=\frac{1}{1/K-x}
=\frac{w\,x^{m}}{1-\bigl(u(x)+w x^{m}\bigr)x-x^{m+2}L(x)},
\]
where $u(x)+wx^{m}$ is again a polynomial of degree at most $m$.
Thus $K$ and $K/(1-xK)$ are of the form \eqref{eq:alg:onestep} with the same $k=m$, $v=w$ and $G=L$, and Lemma \ref{lem:alg:han} gives the same value for their Hankel determinants.
\end{proof}

We shall also use identity \eqref{eq:Hc}, which holds for every power series in place of $F$ and every constant $c$.

\begin{Lemma}\label{lem:alg:shifts}
Let $F(x)$ be given by \eqref{eq:alg:onestep}, where $k\ge0$, $v\neq0$, $\deg u\le k$ and $G$ is a power series, and let $c$ be the coefficient of $x^k$ in $u(x)$. 
Then
\begin{itemize}
\item[(i)] $H^{(1)}_n(F)=0$ for $1\le n\le k-1$, and for $n\ge k$,
\begin{equation}\label{eq:alg:A}
H^{(1)}_n(F)=(-1)^{k(k-1)/2}\,v^n\,H_{n-k}\bigl(c+xG(x)\bigr);
\end{equation}
\item[(ii)] $H^{(-1)}_n(F)=0$ for $1\le n\le k+1$, and for $n\ge k+2$,
\begin{equation}\label{eq:alg:B}
H^{(-1)}_n(F)=(-1)^{(k+1)(k+2)/2}\,v^n\,H^{(1)}_{n-k-2}(G).
\end{equation}
\end{itemize}
\end{Lemma}

\begin{proof}
Write $K(x):=c+xG(x)$.

(i) Let $a_0$ be the constant term of $F$, so that $H^{(1)}_n(F)=H_n\bigl((F-a_0)/x\bigr)$.
Suppose first that $k\ge1$, so that $a_0=0$. 
Write $u(x)=u^-(x)+cx^k$ with $\deg u^-\le k-1$. 
Then
\begin{align*}
\frac{F(x)}{x}
&=\frac{v\,x^{k-1}}{1-u^-(x)x-cx^{k+1}-x^{k+2}G(x)}\\
&=\frac{v\,x^{k-1}}{1-u^-(x)x-x^{(k-1)+2}K(x)},
\end{align*}
which is of the form \eqref{eq:alg:onestep} with $k-1$ in place of $k$ and $K$ in place of $G$.
By Lemma \ref{lem:alg:han}, $H_n(F/x)=0$ for $1\le n\le k-1$ and $H_n(F/x)=(-1)^{(k-1)k/2}v^nH_{n-k}(K)$ for $n\ge k$.

Suppose now that $k=0$. Then $u(x)=c$, $a_0=v$ and $F=v/(1-xK)$, so that
\[
\frac{F(x)-v}{x}=\frac{v\,K(x)}{1-xK(x)}.
\]
By Lemma \ref{lem:alg:KxK}, $H_n\bigl((F-v)/x\bigr)=v^nH_n(K)$, which is \eqref{eq:alg:A} for $k=0$.

(ii) We have $H^{(-1)}_n(F)=H_n(xF)$. Let $g_0$ be the constant term of $G$. Then
\[
xF(x)=\frac{v\,x^{k+1}}{1-\bigl(u(x)+g_0x^{k+1}\bigr)x
-x^{(k+1)+2}\,\dfrac{G(x)-g_0}{x}},
\]
where $u(x)+g_0x^{k+1}$ is a polynomial of degree at most $k+1$.
This is of the form \eqref{eq:alg:onestep} with $k+1$ in place of $k$ and $(G-g_0)/x$ in place of $G$.
Lemma \ref{lem:alg:han} and $H_m\bigl((G-g_0)/x\bigr)=H^{(1)}_m(G)$ yield (ii).
\end{proof}

\subsection{The coupled recurrence}

For $i\ge0$, let $F_i(x)$ be the $i$-th tail of the $H$-fraction \eqref{eq:Hfrac} of $F$:
\[
F_i(x)=
\cFrac{v_i x^{k_i}}{ 1 - u_{i+1}(x) x}
- \cFrac{ v_{i+1} x^{k_i+k_{i+1}+2}}{ 1 - u_{i+2}(x)x}
- \cFrac{ v_{i+2} x^{k_{i+1}+k_{i+2}+2}}{ 1 - u_{i+3}(x)x}
- {\cdots},
\]
so that $F_0=F$ and
\[
F_i=\frac{v_ix^{k_i}}{1-u_{i+1}(x)x-x^{k_i+2}F_{i+1}(x)}.
\]
The $H$-fraction of $F_i$ has the data $(v_{i+t},k_{i+t},u_{i+t+1}(x))_{t\ge0}$; in particular, the indices $s_t(F_i)$ of Theorem \ref{th:super2} associated with $F_i$ are
\begin{equation}\label{eq:alg:sti}
s_t(F_i)=s_{i+t}-s_i\qquad(t\ge0).
\end{equation}
Similarly, we write $R_t(F_i)$ for the sequence \eqref{def:Rn} associated with $F_i$, so that
$R_0(F_i)=1$, $R_1(F_i)=c_{i+1,k_i}$ and, for $t\ge1$,
\begin{equation}\label{eq:alg:Rrec}
R_{t+1}(F_i)=c_{i+t+1,k_{i+t}}R_t(F_i)-v_{i+t}R_{t-1}(F_i).
\end{equation}
We also set $R_{-1}(F_i):=0$, which is consistent with the recurrence for $t=0$. 
Expanding along the last column shows that the determinant below satisfies \eqref{eq:alg:Rrec}, with initial values $R_0(F_i)=1$ and $R_1(F_i)=c_{i+1,k_i}$.
Hence, by induction on $t$,
\[
R_t(F_i)=\det\begin{pmatrix}
c_{i+1,k_i} & v_{i+1} & & & \\
1 & c_{i+2,k_{i+1}} & v_{i+2} & & \\
 & 1 & \ddots & \ddots & \\
 & & \ddots & \ddots & v_{i+t-1}\\
 & & & 1 & c_{i+t,k_{i+t-1}}
\end{pmatrix}.
\]
Expanding the same determinant along its first row instead, we obtain, for $t\ge1$,
\begin{equation}\label{eq:alg:Rback}
R_t(F_i)=c_{i+1,k_i}\,R_{t-1}(F_{i+1})-v_{i+1}\,R_{t-2}(F_{i+2}).
\end{equation}
Also, the sequence $(E_n)_{n \ge 1}$ of \eqref{def:En} satisfies $E_1=1=R_0(F_1)$, $E_2=c_{2,k_1}=R_1(F_1)$, and the recurrence \eqref{def:En} is \eqref{eq:alg:Rrec} for $i=1$ with $t=n-1$, so that
\begin{equation}\label{eq:alg:ER}
E_j=R_{j-1}(F_1)\qquad (j\ge1).
\end{equation}

Combining Lemmas~\ref{lem:alg:han} and~\ref{lem:alg:shifts} with \eqref{eq:Hc}, we obtain the following proposition.

\begin{Proposition}\label{prop:alg:coupled}
For every $n\ge k_0+1$ we have
\begin{equation}\label{eq:alg:coupled}
H^{(1)}_n(F_0)=(-1)^{k_0(k_0-1)/2}\,v_0^n
\Bigl(c_{1,k_0}\,H^{(1)}_{n-k_0-1}(F_1)+H^{(-1)}_{n-k_0}(F_1)\Bigr),
\end{equation}
where $H^{(-1)}_{m}(F_1)=0$ for $1\le m\le k_1+1$, and for $m\ge k_1+2$,
\[
H^{(-1)}_m(F_1)=(-1)^{(k_1+1)(k_1+2)/2}\,v_1^m\,H^{(1)}_{m-k_1-2}(F_2).
\]
Moreover, for $n\ge k_0+1$, 
\begin{equation}\label{eq:alg:hanstep}
H_n(F_0)=(-1)^{k_0(k_0+1)/2}\,v_0^n\,H_{n-k_0-1}(F_1).
\end{equation}
\end{Proposition}

\subsection{Proof of Theorem \ref{th:shiftHDet}}

By Theorem \ref{th:super2}(ii), \eqref{eq:H1sj} is equivalent to \eqref{eq:H1sj:Hsj}. 
Since every integer $n\ge0$ can be written uniquely as $n=s_j+\ell$ with $j\ge0$ and $0\le\ell\le k_j$, it suffices to prove the following three statements for every $j\ge0$:
\begin{itemize}
\item[$(\mathrm a_j)$] $H^{(1)}_{s_j}(F)=(-1)^{s_j+j}\,H_{s_j}(F)\,R_j$;
\item[$(\mathrm b_j)$] $H^{(1)}_{s_j+\ell}(F)=0$ for $1\le \ell\le k_j-1$;
\item[$(\mathrm c_j)$] $H^{(1)}_{s_j+k_j}(F)=(-1)^{k_j(k_j-1)/2}\,V_j^{k_j}\,H^{(1)}_{s_j}(F)$.
\end{itemize}
We prove $(\mathrm a_j)$, $(\mathrm b_j)$ and $(\mathrm c_j)$ by induction on $j$, simultaneously for all power series $F$.
For $j=0$, $(\mathrm a_0)$ follows from $s_0=0$ and $H_0^{(1)}(F)=H_0(F)=R_0=1$.
By Lemma~\ref{lem:alg:shifts}(i), we have $H_n^{(1)}(F)=0$ for $1\le n\le k_0-1$, and 
$H_{k_0}^{(1)}(F) =(-1)^{k_0(k_0-1)/2}v_0^{k_0}.$
Since $V_0=v_0$, these are $(\mathrm b_0)$ and $(\mathrm c_0)$, respectively.

Let $j\ge1$ and $0\le\ell\le k_j$,  and set $n=s_j+\ell$.
Then $n\ge k_0+1$ and $n-k_0=s_j-k_0+\ell$.
By \eqref{eq:alg:sti}, we have $n-k_0-1=s_{j-1}(F_1)+\ell$. 
If $j=1$, then $n-k_0=1+\ell\le k_1+1$.
If $j\ge2$, then $n-k_0\ge s_2-k_0=k_1+2$, and $n-k_0-k_1-2=s_{j-2}(F_2)+\ell$.
In both cases, Proposition \ref{prop:alg:coupled} gives, for $0\le\ell\le k_j$,
\begin{multline}\label{eq:alg:star}
H^{(1)}_{s_j+\ell}(F)=(-1)^{k_0(k_0-1)/2}\,v_0^{s_j+\ell}
\Bigl(c_{1,k_0}\,T_1(\ell)\\
+(-1)^{(k_1+1)(k_1+2)/2}\,v_1^{s_j-k_0+\ell}\,T_2(\ell)\Bigr),
\end{multline}
where $T_1(\ell)=H^{(1)}_{s_{j-1}(F_1)+\ell}(F_1)$ and
\[
T_2(\ell)=
\begin{cases}
H^{(1)}_{s_{j-2}(F_2)+\ell}(F_2), & j\ge2,\\
0, & j=1.
\end{cases}
\]
In the $H$-fraction of $F_1$, $k_j$ is the exponent at index $j-1$, and the corresponding product $V_{j-1}$ is $v_1v_2\cdots v_j$.
For $j\ge2$, $k_j$ is also the exponent at index $j-2$ in the $H$-fraction of $F_2$, where the corresponding product $V_{j-2}$ is $v_2\cdots v_j$.

For $1\le\ell\le k_j-1$, the induction hypotheses $(\mathrm b_{j-1})$ for $F_1$ and, when $j\ge2$,
$(\mathrm b_{j-2})$ for $F_2$ therefore give
\begin{equation}\label{eq:alg:Tb}
T_1(\ell)=T_2(\ell)=0.
\end{equation}
Similarly, set $\tau_j:=(-1)^{k_j(k_j-1)/2}$.
Applying $(\mathrm c_{j-1})$ to $F_1$ and, when $j\ge2$, $(\mathrm c_{j-2})$ to $F_2$, we obtain
\begin{equation}\label{eq:alg:Tbj}
\begin{aligned}
T_1(k_j)&=\tau_j\,(v_1v_2\cdots v_j)^{k_j}\,T_1(0),\\
T_2(k_j)&=\tau_j\,(v_2\cdots v_j)^{k_j}\,T_2(0).
\end{aligned}
\end{equation}
For $j=1$, the assertions involving $T_2$ also hold, since $T_2(\ell)=0$ by definition. 

To complete the induction step, we prove $(\mathrm b_j)$, $(\mathrm c_j)$, and $(\mathrm a_j)$
in this order.

\medskip

\noindent{\it{ Proof of $(\mathrm b_j)$. }}
For $1\le\ell\le k_j-1$,  substituting \eqref{eq:alg:Tb} into \eqref{eq:alg:star} gives $H^{(1)}_{s_j+\ell}(F)=0$.

\medskip

\noindent{\it{ Proof of $(\mathrm c_j)$. }}
Take $\ell=k_j$ in \eqref{eq:alg:star} and use \eqref{eq:alg:Tbj}. We distinguish two cases.

If $j=1$, then $T_2=0$ and \eqref{eq:alg:Tbj} reads $T_1(k_1)=\tau_1v_1^{k_1}T_1(0)$, so that
\begin{align*}
H^{(1)}_{s_1+k_1}(F)
&=(-1)^{k_0(k_0-1)/2}\,v_0^{s_1}v_0^{k_1}\,c_{1,k_0}\,
\tau_1v_1^{k_1}\,T_1(0)\\
&=\tau_1\,(v_0v_1)^{k_1}\,(-1)^{k_0(k_0-1)/2}\,v_0^{s_1}\,c_{1,k_0}\,T_1(0)\\
&=\tau_1\,V_1^{k_1}\,H^{(1)}_{s_1}(F),
\end{align*}
where the last equality is \eqref{eq:alg:star} for $j=1$ and $\ell=0$.

 If $j\ge2$, we use $v_0^{k_j}(v_1\cdots v_j)^{k_j}=V_j^{k_j}$ and $v_0^{k_j}v_1^{k_j}(v_2\cdots v_j)^{k_j}=V_j^{k_j}$:
\begin{align*}
H^{(1)}_{s_j+k_j}(F)
&=(-1)^{k_0(k_0-1)/2}\,v_0^{s_j}v_0^{k_j}
\Bigl(c_{1,k_0}\,\tau_j(v_1\cdots v_j)^{k_j}\,T_1(0)\\
&\qquad+(-1)^{(k_1+1)(k_1+2)/2}\,v_1^{s_j-k_0}v_1^{k_j}\,
\tau_j(v_2\cdots v_j)^{k_j}\,T_2(0)\Bigr)\\
&=\tau_j\,V_j^{k_j}\,(-1)^{k_0(k_0-1)/2}\,v_0^{s_j}
\Bigl(c_{1,k_0}\,T_1(0)\\
&\qquad+(-1)^{(k_1+1)(k_1+2)/2}\,v_1^{s_j-k_0}\,T_2(0)\Bigr)\\
&=\tau_j\,V_j^{k_j}\,H^{(1)}_{s_j}(F),
\end{align*}
where the last equality is \eqref{eq:alg:star} for $\ell=0$.

Thus $(\mathrm c_j)$ holds in both cases.

\medskip

\noindent{\it{ Proof of $(\mathrm a_j)$. }}
Take $\ell=0$. Applying \eqref{eq:alg:hanstep} to $F_0$ and, when $j\ge2$, to $F_1$, we obtain
\begin{align*}
H_{s_{j-1}(F_1)}(F_1)&=(-1)^{k_0(k_0+1)/2}\,v_0^{-s_j}\,H_{s_j}(F),\\
H_{s_{j-2}(F_2)}(F_2)&=(-1)^{k_1(k_1+1)/2}\,v_1^{-(s_j-k_0-1)}\,
H_{s_{j-1}(F_1)}(F_1)\qquad(j\ge2).
\end{align*}
Applying the induction hypothesis $(\mathrm a_{j-1})$ to $F_1$ and using $s_{j-1}(F_1)=s_j-k_0-1$, we can write the first term on the right-hand side of \eqref{eq:alg:star} as
\begin{multline*}
(-1)^{k_0(k_0-1)/2}v_0^{s_j}c_{1,k_0}(-1)^{s_j-k_0-1+j-1}
H_{s_{j-1}(F_1)}(F_1)R_{j-1}(F_1)\\
=(-1)^{s_j+j}H_{s_j}(F)\,c_{1,k_0}R_{j-1}(F_1).
\end{multline*}
Here the last equality follows from \eqref{eq:alg:hanstep} and because 
\[
\frac{k_0(k_0-1)}{2}+\frac{k_0(k_0+1)}{2}=k_0^2\equiv k_0\pmod 2.
\]
For $j\ge2$, we apply the induction hypothesis $(\mathrm a_{j-2})$ to $F_2$.
Using $s_{j-2}(F_2)=s_j-k_0-k_1-2$, we can write the second term on the right-hand side of \eqref{eq:alg:star} as
\begin{multline*}
(-1)^{e}\,v_0^{s_j}v_1^{s_j-k_0}\,v_0^{-s_j}v_1^{-(s_j-k_0-1)}\,
H_{s_j}(F)\,R_{j-2}(F_2)\\
=-(-1)^{s_j+j}H_{s_j}(F)\,v_1R_{j-2}(F_2),
\end{multline*}
where
\begin{align*}
e&=\tfrac{k_0(k_0-1)}2+\tfrac{(k_1+1)(k_1+2)}2+(s_j-k_0-k_1-2+j-2)\\
&\qquad+\tfrac{k_1(k_1+1)}2+\tfrac{k_0(k_0+1)}2\\
&\equiv k_0+(k_1+1)+(s_j-k_0-k_1+j)\equiv s_j+j+1\pmod 2.
\end{align*}
For $j=1$, the second term vanishes, consistently with the convention $R_{-1}(F_2)=0$.
Combining the two terms and using \eqref{eq:alg:Rback} with $i=0$ and $t=j$, we obtain
\begin{align*}
H^{(1)}_{s_j}(F)
&=(-1)^{s_j+j}H_{s_j}(F)
\bigl(c_{1,k_0}R_{j-1}(F_1)-v_1R_{j-2}(F_2)\bigr)\\
&=(-1)^{s_j+j}H_{s_j}(F)\,R_j.
\end{align*}
This proves $(\mathrm a_j)$ and completes the induction, and hence the proof of Theorem~\ref{th:shiftHDet}.
\qed

\subsection{Proof of Theorem \ref{th:backshiftHDet}}

Let $j\ge1$ and $0\le\ell\le k_j$.
Since $s_j\ge s_1=k_0+1$, we have $n=s_j+1+\ell\ge k_0+2$.
By \eqref{eq:alg:sti}, we have $n-k_0-2=s_j-k_0-1+\ell=s_{j-1}(F_1)+\ell.$
Applying Lemma~\ref{lem:alg:shifts}(ii) to $F_0$, we obtain
\[
H^{(-1)}_{s_j+1+\ell}(F)
=
(-1)^{(k_0+1)(k_0+2)/2}v_0^{s_j+1+\ell}
H^{(1)}_{s_{j-1}(F_1)+\ell}(F_1).
\]

For $\ell=0$, applying $(\mathrm a_{j-1})$ to $F_1$ and using \eqref{eq:alg:hanstep} and \eqref{eq:alg:ER}, we obtain
\begin{align*}
    H^{(-1)}_{s_j+1}(F)&= (-1)^{ \frac{(k_0+1)(k_0+2)}{2}
+s_j-k_0+j-2+\frac{k_0(k_0+1)}{2} } v_0 H_{s_j}(F)E_j \\
& =(-1)^{s_j+j+1}v_0H_{s_j}(F)E_j.
\end{align*}
This is \eqref{eq:Hsjm1:Hsj}, which is equivalent to \eqref{eq:Hsj(-1)} by Theorem~\ref{th:super2}(ii).

For $\ell=k_j$, applying $(\mathrm c_{j-1})$ to $F_1$ gives
\[
H^{(1)}_{s_{j-1}(F_1)+k_j}(F_1)
=
(-1)^{k_j(k_j-1)/2}(v_1\cdots v_j)^{k_j}
H^{(1)}_{s_{j-1}(F_1)}(F_1).
\]
Combining this identity with the reduction formula above, we obtain \eqref{eq:Hs_j+1(-1)}, since
$v_0^{k_j}(v_1\cdots v_j)^{k_j}=V_j^{k_j}$ and $s_j+1+k_j=s_{j+1}$.

For $1\le\ell\le k_j-1$,  applying $(\mathrm b_{j-1})$ to $F_1$ gives $H^{(1)}_{s_{j-1}(F_1)+\ell}(F_1)=0,$ and hence $H^{(-1)}_{s_j+1+\ell}(F)=0.$
Furthermore, Lemma~\ref{lem:alg:shifts}(ii) gives $H^{(-1)}_n(F)=0$ for $1\le n\le k_0+1=s_1$.
Since every $n\ge s_1+1$ can be written uniquely as $n=s_j+1+\ell$ with $j\ge1$ and $0\le\ell\le k_j$,
these are all the non-vanishing backward shifted Hankel determinants.
\qed

\begin{Remark}
If the $H$-fraction of $F$ is finite, that is, if some tail $F_{J+1}$ is zero, the same proof applies verbatim, with $H_m(0)=H^{(\pm1)}_m(0)=0$ for $m\ge1$: Lemmas \ref{lem:alg:han} and \ref{lem:alg:shifts} hold for $G=0$.
\end{Remark}


\section{An algebraic proof of Theorem \ref{th:shift2HDet}}\label{sec:doubleshalg}

In this section we give an algebraic proof of Theorem \ref{th:shift2HDet},
assuming Theorems \ref{th:super2} and \ref{th:shiftHDet}
(for an algebraic proof of the latter, see Section \ref{sec:shiftalg}).
Throughout this section, $F(x)=\sum_{n\ge0}a_nx^n$ is given by its $H$-fraction \eqref{eq:Hfrac}, which we assume to be infinite (see Remark \ref{rem:alg2:finite} for the finite case), and $\mathcal L$ denotes the linear functional on $\setF[x]$ defined by $\mathcal L(x^i)=a_i$ $(i\ge0)$, so that
$H^{(r)}_n(F)=\det\bigl(\mathcal L(x^{r+i+i'})\bigr)_{0\le i,i'\le n-1}$.

\subsection{The Desnanot--Jacobi identity}\label{subsec:doubleshalg:djiden}

Let $n\ge1$, and consider the Hankel matrix $M=(a_{i+i'})_{0\le i,i'\le n}.$
By definition, $\det M=H_{n+1}(F).$
The determinants of the relevant submatrices are listed below:
\[
\begin{array}{c|c|c}
\text{Deleted rows} & \text{Deleted columns}
& \text{Determinant}\\ \hline
\text{last} & \text{last} & H_n(F)\\
\text{first} & \text{first} & H_n^{(2)}(F)\\
\text{first} & \text{last} & H_n^{(1)}(F)\\
\text{last} & \text{first} & H_n^{(1)}(F)\\
\text{first and last} & \text{first and last}
& H_{n-1}^{(2)}(F)
\end{array}
\]
Applying the Desnanot--Jacobi identity gives \begin{equation}\label{eq:alg2:DJ}
H_{n+1}(F)\,H_{n-1}^{(2)}(F)
=
H_n(F)\,H_n^{(2)}(F)
-\bigl(H_n^{(1)}(F)\bigr)^2.
\end{equation}

\subsection{Block-orthogonal polynomials}

Let $A_t(x)$ and $B_t(x)$ be the numerators and the denominators of the convergents of \eqref{eq:Hfrac}, defined by
$A_0=0$, $A_1=v_0x^{k_0}$, $B_0=1$, $B_1=1-u_1(x)x$ and, for $t\ge1$,
\begin{align*}
A_{t+1}&=\bigl(1-u_{t+1}(x)x\bigr)A_t-v_tx^{k_{t-1}+k_t+2}A_{t-1},\\
B_{t+1}&=\bigl(1-u_{t+1}(x)x\bigr)B_t-v_tx^{k_{t-1}+k_t+2}B_{t-1}.
\end{align*}
Since $\deg u_{t+1}\le k_t$ and $s_{t+1}=s_t+k_t+1$, induction on $t$ shows that $\deg A_t\le s_t-1$, $\deg B_t\le s_t$, and $B_t(0)=1$.
For $t\ge1$, the same recurrences give
\[
A_{t+1}B_t-A_tB_{t+1}
=
v_tx^{k_{t-1}+k_t+2}
\bigl(A_tB_{t-1}-A_{t-1}B_t\bigr).
\]
Starting from $A_1B_0-A_0B_1=v_0x^{k_0}$ and using $k_0+\sum_{i=1}^{t}(k_{i-1}+k_i+2)=2s_t+k_t,$
we obtain the following identity for every $t\ge0$ by iteration:
\begin{equation}\label{eq:alg2:wronski}
A_{t+1}B_t-A_tB_{t+1}
=
V_t x^{2s_t+k_t}.
\end{equation}

As in Section~\ref{sec:shiftalg}, let $F_t$ denote the $t$-th tail of \eqref{eq:Hfrac}, with $F_0=F$.
For $t\ge0$, we have
\[
F_t(x)=
\frac{v_tx^{k_t}}{1-u_{t+1}(x)x-x^{k_t+2}F_{t+1}(x)}.
\]
Set $\xi_t:=x^{k_t+2}F_{t+1}(x)$.
Using the recurrences for the tails and for $A_t$ and $B_t$, induction on $t$ gives, for every $t\ge0$,
\[
F=\frac{A_{t+1}-\xi_tA_t}{B_{t+1}-\xi_tB_t}.
\]
Since $B_{t+1}(0)=1$ and $\xi_t(0)=0$, the denominator has constant term $1$.
Therefore, \eqref{eq:alg2:wronski} gives
\begin{equation}\label{eq:alg2:pade}
F\,B_t-A_t=\frac{A_{t+1}B_t-A_tB_{t+1}}{B_{t+1}-\xi_tB_t}
=V_t\,x^{2s_t+k_t}\bigl(1+O(x)\bigr).
\end{equation}

For $t\ge0$ let $P_t(x):=x^{s_t}B_t(1/x)$ be the reversed polynomial of $B_t$. It is a monic polynomial of degree $s_t$, and $P_0=1$, $P_1(x)=\Theta_0(x)$, and $P_{t+1}(x)=\Theta_t(x)P_t(x)-v_tP_{t-1}(x)$ for $t\ge1$, where
\[
\Theta_t(x):=x^{k_t+1}-\sum_{l=0}^{k_t}c_{t+1,l}\,x^{k_t-l}
\]
is the reversed polynomial of $1-u_{t+1}(x)x$, of degree $k_t+1$ (note that $s_{t+1}-(k_{t-1}+k_t+2)=s_{t-1}$).
Setting $x=0$ in the recurrence gives $P_{t+1}(0)=-c_{t+1,k_t}P_t(0)-v_tP_{t-1}(0)$ for $t\ge1$, with initial values $P_0(0)=1$ and $P_1(0)=-c_{1,k_0}$.
Thus the sequence $\bigl((-1)^tP_t(0)\bigr)_{t\ge0}$ satisfies the recurrence \eqref{def:Rn}, with the same initial values as $(R_t)_{t\ge0}$.
Therefore, for $t\ge0$,
\begin{equation}\label{eq:alg2:P0}
P_t(0)=(-1)^tR_t.
\end{equation}

\begin{Lemma}\label{lem:alg2:orth}
For $t\ge0$ we have $\mathcal L(x^iP_t)=0$ for $0\le i<s_t+k_t$, and $\mathcal L(x^{s_t+k_t}P_t)=V_t$.
\end{Lemma}

\begin{proof}
Write $B_t(x)=\sum_{l=0}^{s_t}\beta_lx^l$, so that  $P_t(x)=\sum_l\beta_lx^{s_t-l}$. For $i\ge0$, we have
\[
\mathcal L(x^iP_t)
=
\sum_{l=0}^{s_t}\beta_l\mathcal L(x^{i+s_t-l})
=
\sum_{l=0}^{s_t}\beta_la_{i+s_t-l}.
\]
On the other hand, the coefficient of $x^{s_t+i}$ in $F(x)B_t(x)$ is
\[
[x^{s_t+i}](FB_t)
=
\sum_{l=0}^{s_t}\beta_l[x^{s_t+i-l}]F
=
\sum_{l=0}^{s_t}\beta_la_{s_t+i-l},
\]
where $[x^r]$ denotes the coefficient of $x^r$.
Consequently,
\[
\mathcal L(x^iP_t)=[x^{s_t+i}](FB_t).
\]
Since $\deg A_t<s_t$ for $t\ge1$ and $A_0=0$, subtracting $A_t$ does not change this coefficient.
Hence
\[
\mathcal L(x^iP_t)
=
[x^{s_t+i}](FB_t-A_t).
\]
By \eqref{eq:alg2:pade},  $FB_t-A_t =V_tx^{2s_t+k_t}+O(x^{2s_t+k_t+1}).$
If $0\le i<s_t+k_t$, then $s_t+i<2s_t+k_t$, so the coefficient of $x^{s_t+i}$ is zero.
If $i=s_t+k_t$, then $s_t+i=2s_t+k_t$, and the corresponding coefficient is $V_t$.
This proves the lemma.
\end{proof}

Let $0\le t<t'$ and let $Q$ be a polynomial with $\deg Q\le k_t+k_{t'}$.
Since $s_t+k_t=s_{t+1}-1<s_{t'}$, we have  $\deg(QP_t)\le s_t+k_t+k_{t'}<s_{t'}+k_{t'}$.
By linearity, Lemma~\ref{lem:alg2:orth} therefore gives $\mathcal L(QP_tP_{t'})=0.$
For $t\ge0$ and $m\ge0$, define
\[
\mu_{t,m}:=\mathcal L(x^mP_t^2).
\]
If $0\le m<k_t$, then $\deg(x^mP_t)=s_t+m<s_t+k_t$, so Lemma~\ref{lem:alg2:orth} gives $\mu_{t,m}=0$.
For $m=k_t$, the same lemma yields
\[
\begin{aligned}
\mu_{t,k_t}
&=
\mathcal L(x^{s_t+k_t}P_t)
+\mathcal L\bigl(x^{k_t}(P_t-x^{s_t})P_t\bigr)=V_t,
\end{aligned}
\]
where the second term vanishes because $P_t$ is monic of degree $s_t$, and hence $\deg\bigl(x^{k_t}(P_t-x^{s_t})\bigr)<s_t+k_t$.
Thus, for $0\le m\le k_t$, we have
\begin{equation}\label{eq:alg2:mu}
\mu_{t,m}
=
\begin{cases}
0, & 0\le m<k_t,\\
V_t, & m=k_t.
\end{cases}
\end{equation}

\begin{Lemma}\label{lem:alg2:mu}
For $t\ge0$, let $\hat\mu_t(w):=\sum_{m\ge0}\mu_{t,k_t+m}w^m$. Then
\[
\hat\mu_t(w)\,\bigl(1-u_{t+1}(w)w\bigr)\equiv V_t\pmod{w^{k_t+2}}.
\]
\end{Lemma}

\begin{proof}
The constant term of the left-hand side is $\mu_{t,k_t}=V_t$.
It remains to show that its coefficients of $w,w^2,\ldots,w^{k_t+1}$ vanish.

Let $0\le i\le k_t$.
For $t\ge1$, the orthogonality relations established above give $\mathcal L(x^iP_tP_{t+1})=0$ and $\mathcal L(x^iP_tP_{t-1})=0$.
Using $P_{t+1}=\Theta_tP_t-v_tP_{t-1}$, we obtain
\[
\begin{aligned}
\mathcal L(x^i\Theta_tP_t^2)
&=\mathcal L(x^iP_tP_{t+1})+v_t\mathcal L(x^iP_tP_{t-1})=0.
\end{aligned}
\]
For $t=0$, this follows from $P_1=\Theta_0P_0$ and $\mathcal L(x^iP_0P_1)=0$.
Expanding
\[
\Theta_t(x)=x^{k_t+1}-\sum_{l=0}^{k_t}c_{t+1,l}x^{k_t-l}
\]
and using the definition of $\mu_{t,m}$, we get
\[
\mu_{t,k_t+i+1}-\sum_{l=0}^{k_t}c_{t+1,l}\mu_{t,k_t+i-l}=0.
\]
Whenever $i<l\le k_t$, we have $0\le k_t+i-l<k_t$, so $\mu_{t,k_t+i-l}=0$ by \eqref{eq:alg2:mu}.
Thus the left-hand side is precisely the coefficient of $w^{i+1}$ in  
$\hat\mu_t(w)\bigl(1-\sum_lc_{t+1,l}w^{l+1}\bigr) =\hat\mu_t(w)\bigl(1-u_{t+1}(w)w\bigr)$.
As $i$ ranges from $0$ to $k_t$, this proves that the coefficients of $w,w^2,\ldots,w^{k_t+1}$ vanish.
\end{proof}

\subsection{The determinants $H^{(2)}_{s_{j+1}-1}(F)$}

\begin{Proposition}\label{prop:alg2:cofactor}
For every $j\ge0$ we have
\begin{equation}\label{eq:alg2:cofactor}
H^{(2)}_{s_{j+1}-1}(F)=H_{s_{j+1}}(F)\sum_{t=0}^{j}\frac{d_{t+1}R_t^2}{V_t},
\end{equation}
where $d_{t+1}$ is the coefficient of $x^{k_t}$ in $1-u_{t+1}(x)x$, as in Theorem \ref{th:shift2HDet}.
\end{Proposition}

\begin{proof}
Let $N=s_{j+1}$ and $M=(a_{i+i'})_{0\le i,i'\le N-1}$.
By Theorem~\ref{th:super2}, $\det M=H_N(F)\neq0$, so $M$ is invertible.
Deleting the first row and the first column of $M$ gives the matrix $(a_{i+i'+2})_{0\le i,i'\le N-2}$.
The cofactor formula for the inverse therefore yields $H^{(2)}_{N-1}(F)=H_N(F)\,(M^{-1})_{0,0}$, where $(M^{-1})_{0,0}$ denotes the entry in the first row and first column of $M^{-1}$.

Let $e_0=(1,0,\ldots,0)^T$, and let $y=(y_0,\ldots,y_{N-1})^T$ be the unique solution of $My=e_0$.
Define $Y(x):=\sum_{i=0}^{N-1}y_ix^i. $
Since $y=M^{-1}e_0$, we have $(M^{-1})_{0,0}=y_0=Y(0)$.
Moreover, the equation $My=e_0$ is equivalent to
\[
\mathcal L(x^iY) =  \sum_{i'=0}^{N-1}y_{i'}\mathcal L(x^{i+i'})
=\sum_{i'=0}^{N-1}a_{i+i'}y_{i'}
=
\begin{cases}
1, & i=0,\\
0, & 1\le i\le N-1.
\end{cases}.
\]
By linearity, $Y$ is therefore the unique polynomial of degree less than $N$ such that, for every polynomial $Q$ with $\deg Q<N$,
\begin{equation}\label{eq:alg2:Y}
\mathcal L(QY)=Q(0).
\end{equation}
The polynomials $x^iP_t$, $0\le t\le j$, $0\le i\le k_t$, are monic of degrees $s_t+i$.
Since $s_t+k_t=s_{t+1}-1$, these degrees run through $0,1,\ldots,N-1$, each exactly once.
Hence these polynomials form a basis of the polynomials of degree less than $N$.
Write the expansion of $Y$ in this basis as
\[
Y(x)=
\sum_{t=0}^{j}\sum_{i'=0}^{k_t}
y_{t,i'}x^{i'}P_t(x).
\]
Taking $Q=x^iP_t$ in \eqref{eq:alg2:Y}, the terms involving $P_{t'}$ with $t'\neq t$ vanish by the
orthogonality relations established above.
Using $\mu_{t,m}=\mathcal L(x^mP_t^2)$, we obtain
\[
\sum_{i'=0}^{k_t}\mu_{t,i+i'}y_{t,i'}
=
\begin{cases}
P_t(0), & i=0,\\
0, & 1\le i\le k_t.
\end{cases}.
\]
For each $t$, define $z_t(w)=\sum_{l=0}^{k_t}y_{t,k_t-l}w^l.$
By \eqref{eq:alg2:mu}, we have $\mu_{t,m}=0$ for $0\le m<k_t$.
Consequently, for $0\le i\le k_t$,
\[
[w^i]\bigl(\hat\mu_t(w)z_t(w)\bigr)
=\sum_{l=0}^{i}
\mu_{t,k_t+i-l}y_{t,k_t-l}=\sum_{i'=0}^{k_t}
\mu_{t,i+i'}y_{t,i'}.
\]
Thus the system above is equivalent to
\[
\hat\mu_t(w)z_t(w)
\equiv P_t(0)\pmod{w^{k_t+1}}.
\]
Since $\hat\mu_t(0)=V_t\neq0$, the series $\hat\mu_t(w)$ is invertible.
Lemma~\ref{lem:alg2:mu} therefore gives
\[
z_t(w)
\equiv
\frac{P_t(0)}{\hat\mu_t(w)}
\equiv
\frac{P_t(0)}{V_t}
\bigl(1-u_{t+1}(w)w\bigr)
\pmod{w^{k_t+1}}.
\]
Comparing the coefficients of $w^{k_t}$ and using the definition of $d_{t+1}$, we obtain
\[
y_{t,0}
=
[w^{k_t}]z_t(w)
=
\frac{d_{t+1}P_t(0)}{V_t}.
\]
Evaluating the expansion of $Y$ at $x=0$ and using \eqref{eq:alg2:P0}, we conclude that
\[
\begin{aligned}
Y(0)=\sum_{t=0}^{j}y_{t,0}P_t(0)=\sum_{t=0}^{j}\frac{d_{t+1}P_t(0)^2}{V_t}=\sum_{t=0}^{j}\frac{d_{t+1}R_t^2}{V_t}.
\end{aligned}
\]
Substituting this into $H^{(2)}_{N-1}(F)=H_N(F)Y(0)$, with $N=s_{j+1}$, gives \eqref{eq:alg2:cofactor}.
\end{proof}

\begin{Remark}
    This is Lemma~\ref{lem:h2sjm1} of Section~\ref{sec:doublesh}, where it is derived from the path proofs of Theorems~\ref{th:h(2):kj11} and~\ref{th:h(2):kj12} by the Desnanot--Jacobi identity. The proof above does not use Section~\ref{sec:doublesh}.
\end{Remark}

\subsection{Proof of Theorem \ref{th:shift2HDet}}

Recall from Theorem \ref{th:super2} that, for $j\ge0$, 
\begin{equation}\label{eq:alg2:Hratio}
H_{s_{j+1}}(F)=(-1)^{k_j(k_j+1)/2}\,V_j^{k_j+1}\,H_{s_j}(F), 
\end{equation}
and that $H_n(F)=0$ if $n$ is not of the form $s_j$.

\medskip
\noindent{\it{ Proof of (2).}} 
Let $k_j=0$. Then $s_{j+1}-1=s_j$, and \eqref{eq:alg2:Hratio} gives $H_{s_{j+1}}(F)=V_jH_{s_j}(F).$
Substituting these identities into \eqref{eq:alg2:cofactor}, we obtain
\[
H^{(2)}_{s_j}(F)=H_{s_j}(F)\sum_{t=0}^{j}d_{t+1}R_t^2\frac{V_j}{V_t}.
\]
which is \eqref{eq:H1sj:Hsj:k0}.

\medskip
\noindent{\it{ Proof of (3).}} Let $k_j\ge3$, $1\le\ell\le k_j-2$ and $n=s_j+\ell$, so that 
$H^{(2)}_n(F)=\det\bigl(\mathcal L(x^{2}x^ix^{i'})\bigr)_{0\le i,i'\le n-1}$.
Since $P_j$ is monic of degree $s_j\le n-1$, we may replace the row of index $i=s_j$ by the row $\bigl(\mathcal L(x^2P_jx^{i'})\bigr)_{i'}$ without changing the determinant. 
For $0\le i'\le n-1$ we have $2+i'\le s_j+\ell+1\le s_j+k_j-1$, 
Thus every entry of the new row is zero by Lemma~\ref{lem:alg2:orth}, and consequently
$H^{(2)}_n(F)=0$.

\medskip
\noindent{\it{ Proof of (1).}}
Let $k_j\ge1$.
For $j=0$, it reads $H^{(2)}_{s_0}(F)=1=H_{s_0}(F)R_0^2$.
For $j\ge1$, apply \eqref{eq:alg2:DJ} with $n=s_j$.
Since $s_j<s_j+1<s_{j+1},$ we have $H_{s_j+1}(F)=0$. Hence, using \eqref{eq:H1sj:Hsj}, we obtain
\[
H_{s_j}(F)H^{(2)}_{s_j}(F)=\bigl(H^{(1)}_{s_j}(F)\bigr)^2=H_{s_j}(F)^2R_j^2.
\]
Dividing by $H_{s_j}(F)\neq0$ proves
\eqref{eq:H1sj:Hsj:k1}.

Next, we apply \eqref{eq:alg2:DJ} with $n=s_j+k_j=s_{j+1}-1\ge1$.
Since $k_j\ge1$, we have $s_j<n<s_{j+1}$ and hence $H_n(F)=0$.
Using \eqref{eq:H1sjkj}, \eqref{eq:alg2:Hratio}, and \eqref{eq:H1sj:Hsj:k1}, we obtain
\begin{align*}
H^{(2)}_{s_j+k_j-1}(F)&=-\frac{\bigl(H^{(1)}_{s_j+k_j}(F)\bigr)^2}{H_{s_{j+1}}(F)}
=-\frac{V_j^{2k_j}H_{s_j}(F)^2R_j^2}{(-1)^{k_j(k_j+1)/2}V_j^{k_j+1}H_{s_j}(F)}\\
&=(-1)^{\frac{(k_j-2)(k_j-1)}2}\,V_j^{k_j-1}\,H^{(2)}_{s_j}(F),
\end{align*}
since $\frac{k_j(k_j+1)}2-\frac{(k_j-2)(k_j-1)}2=2k_j-1$ is odd.
This is \eqref{eq:H1sj:Hsj:k11}.

Finally, recall that $s_j+k_j=s_{j+1}-1$. 
For $j\ge1$, applying Proposition~\ref{prop:alg2:cofactor} at indices $j$ and $j-1$ and using \eqref{eq:alg2:Hratio}, we obtain
\begin{align*}
H^{(2)}_{s_j+k_j}(F)
&=\frac{H_{s_{j+1}}(F)}{H_{s_j}(F)}\,H^{(2)}_{s_j-1}(F)
+H_{s_{j+1}}(F)\,\frac{d_{j+1}R_j^2}{V_j}\\
&=\Delta_j+(-1)^{\frac{k_j(k_j+1)}2}V_j^{k_j}\,d_{j+1}\,H_{s_j}(F)R_j^2,
\end{align*}
where the first term is absent if $j=0$ (recall $\Delta_0=0$), and where
$\Delta_j=(-1)^{k_j(k_j+1)/2}V_j^{k_j+1}H^{(2)}_{s_j-1}(F)$ in both cases of the definition of $\Delta_j$, since $s_j-1=s_{j-1}+k_{j-1}$, which is equal to $s_{j-1}$ if $k_{j-1}=0$.
Since $k_j\ge1$, we have $d_{j+1}=-c_{j+1,k_j-1}$, and by
\eqref{eq:H1sj:Hsj:k1} and \eqref{eq:H1sj:Hsj:k11},
\begin{align*}
(-1)^{\frac{k_j(k_j+1)}2}V_j^{k_j}\,d_{j+1}\,H_{s_j}(F)R_j^2
&=(-1)^{\frac{(k_j-2)(k_j-1)}2}c_{j+1,k_j-1}V_j^{k_j}H^{(2)}_{s_j}(F)\\
&=c_{j+1,k_j-1}\,V_j\,H^{(2)}_{s_j+k_j-1}(F).
\end{align*}
This is \eqref{eq:H1sj:Hsj:k12}, and completes the proof of Theorem \ref{th:shift2HDet}. \qed

\begin{Remark}\label{rem:alg2:finite}
If the $H$-fraction of $F$ is finite, that is, if $F_{J+1}=0$ for some $J$, then $F=A_{J+1}/B_{J+1}$ and $\mathcal L(x^iP_{J+1})=0$ for all $i\ge0$. The proof above applies verbatim to all indices $j\le J$.
\end{Remark}


\section{Quasi-Dyck paths and Hankel continued fractions}\label{sec:flajolet}

In this section we extend Flajolet's combinatorial interpretation of  $J$-fractions \cite{Flajolet1980} to $H$-fractions. 
The coefficients of an $H$-fraction will be expressed as weighted  sums over a family of lattice paths, called quasi-Dyck paths,  in which the down steps may have several lengths, determined by 
the exponents $k_j$ of the $H$-fraction.
An essentially equivalent model, the P-paths of Roblet  \cite[Chapter 2]{Roblet1994}, is compared with ours in the concluding remarks.

\subsection{Flajolet's theorem for $J$-fractions}\label{sec:flaj:motzkin}

Recall that a {\it Motzkin path} of length $n$ is a lattice path  from $(0,0)$ to $(n,0)$ that never goes below the $x$-axis and  consists of up steps $(1,1)$, level steps $(1,0)$ and down steps 
$(1,-1)$. 
Let $\mathbf{u}=(u_1,u_2,\ldots)$ and $\mathbf{v}=(v_0,v_1,\ldots)$ be two sequences. 
The weight of a Motzkin path $p$ is defined to be $v_0$ times the product of the weights of its steps, where an up step has weight $1$, a level step at height $h$ has weight $u_{h+1}$, and a down step from height $h$ to height $h-1$ has weight $v_h$. 

\begin{Theorem}[Flajolet \cite{Flajolet1980}]\label{th:flajolet}
Let $f(x)=\sum_{n\ge 0}a_nx^n$ be the $J$-fraction \eqref{eq:Jacobi} attached to $(\mathbf{u},\mathbf{v})$. 
Then, for each $n\geq 0$, the coefficient $a_n$ is the sum of the weights of all Motzkin paths of length $n$. 
\end{Theorem}

The proof of Theorem \ref{th:flajolet} rests on the {\it first return decomposition}: 
a Motzkin path is a sequence of primitive paths, each of which is either a level step, or an up step followed by a Motzkin path translated to height $1$ and by a down step. 
Translating this decomposition into generating functions gives 
\[
f=\frac{v_0}{1-u_1x-v_1x^2f^{(1)}},
\]
where $f^{(1)}$ is the generating function of the paths at height $1$, and iterating yields 
\eqref{eq:Jacobi}.

Let us see what has to be changed for an $H$-fraction \eqref{eq:Hfrac}. 
Compared with \eqref{eq:Jacobi}, the numerators $v_tx^2$ are replaced by $v_tx^{k_{t-1}+k_t+2}$, and the level terms $u_tx$ by 
\[
u_t(x)x=\sum_{r=0}^{k_{t-1}}c_{t,r}x^{r+1}.
\] 
In a path model in which the exponent of $x$ counts the up steps, the monomial $v_tx^{k_{t-1}+k_t+2}$ must therefore be carried by a down step that cancels $k_{t-1}+k_t+2$ up steps at once, and the monomial $c_{t,r}x^{r+1}$ by a down step that cancels $r+1$ up steps.
This leads to the following definition.

\subsection{Quasi-Dyck paths}\label{sec:flaj:qd}

Throughout this section, we fix a sequence $\mathbf{k}=(k_0,k_1,k_2,\ldots)$ of nonnegative integers.
We define two sequences of \emph{levels}, $M=(M_t)_{t\ge0}$ and $N=(N_t)_{t\ge1}$, as follows.
Set $M_0=0$ and $M_1=k_0$.
For each $t\ge1$, define
\begin{equation}\label{eq:MN}
M_{t+1}=M_t+k_{t-1}+k_t+2
\end{equation}
and
\begin{equation}\label{eq:N}
N_t=M_t+k_{t-1}+1.
\end{equation}
Thus $M_t=k_0+\sum_{i=1}^{t-1}(k_{i-1}+k_i+2)$ is the total degree in $x$ of the numerators $v_0x^{k_0}, v_1x^{k_0+k_1+2},\ldots, v_{t-1}x^{k_{t-2}+k_{t-1}+2}$ of the $H$-fraction \eqref{eq:Hfrac}, and we have
\begin{equation}\label{eq:MNorder}
0=M_0\le M_1<N_1<M_2<N_2<M_3<N_3<\cdots,
\end{equation}
since $N_t-M_t=k_{t-1}+1\geq 1$ and $M_{t+1}-N_t=k_t+1\geq 1$. 
Note that $M_1=M_0$ if and only if $k_0=0$.

All our paths are lattice paths in the plane made of steps of slope $+1$ or $-1$. 
A {\it down step of drop $\ell$} ($\ell\geq 0$) is a step from a point $(x,y)$ to the point $(x+\ell,y-\ell)$; 
thus a down step of drop $1$ is an ordinary Dyck down step, and a down step of drop $0$ is 
an empty step.

\begin{Definition}\label{def:step}
A {\it quasi-Dyck path} (associated with $\mathbf{k}$) is a lattice path starting at $(0,0)$, ending on the $x$-axis and never going below the $x$-axis, whose steps are of the following four kinds:
\begin{itemize}
\item[$(b)$] an {\it up step}, from $(x,h)$ to $(x+1,h+1)$, where $h\ge 0$;

\item[$(c_t)$] for each $t\ge 1$, a {\it short down step}, of drop $r+1$ with $0\le r\le k_{t-1}$, from height $M_t+r+1$ down to  height $M_t$; we also write $c_{t,r}$ for a short down step of 
drop $r+1$ ending at height $M_t$;

\item[$(v)$] for each $t\ge 1$, a {\it long down step}, of drop $k_{t-1}+k_t+2$, from height $M_{t+1}$ down to height $M_t$;

\item[$(d)$] the {\it terminal step}, of drop $k_0$, from height $M_1$ down to height $0$;
\end{itemize}
and in which the terminal step occurs exactly once, as the last step.
The set of all quasi-Dyck paths is denoted by $\mathcal{QD}$.
\end{Definition}

The {\it vertices} of the path are the endpoints of its steps; a down step of drop $\ell\geq2$ is a single step, and the lattice points in its interior are not vertices of the path. 

\begin{Remark}
The behaviour at height $0$ depends on $k_0$:
\begin{itemize}
\item If $k_0>0$, every internal vertex has positive height, so the path meets the $x$-axis only at its starting and ending points.

\item If $k_0=0$, internal vertices of height $0$ are allowed. 
The terminal step has drop $0$ but is still recorded as a step. 
In particular, the path consisting of the single step $d$ is a quasi-Dyck path.
\end{itemize}
Here an internal vertex is any vertex other than the starting and ending vertices.
\end{Remark}

Since the $x$-coordinate increases by $1$ along an up step and by $\ell$ along a down step of drop $\ell$, a quasi-Dyck path with $n$ up steps ends at $(2n,0)$. 
We call $|p|=2n$ the {\it length} of $p$; it is the horizontal length of the path, not its number of 
steps. 
We write $\mathcal{QD}_n$ for the set of quasi-Dyck paths of length $2n$. 
Note that $\mathcal{QD}_0$ is empty if $k_0\geq 1$ and consists of the single path $d$ if $k_0=0$.

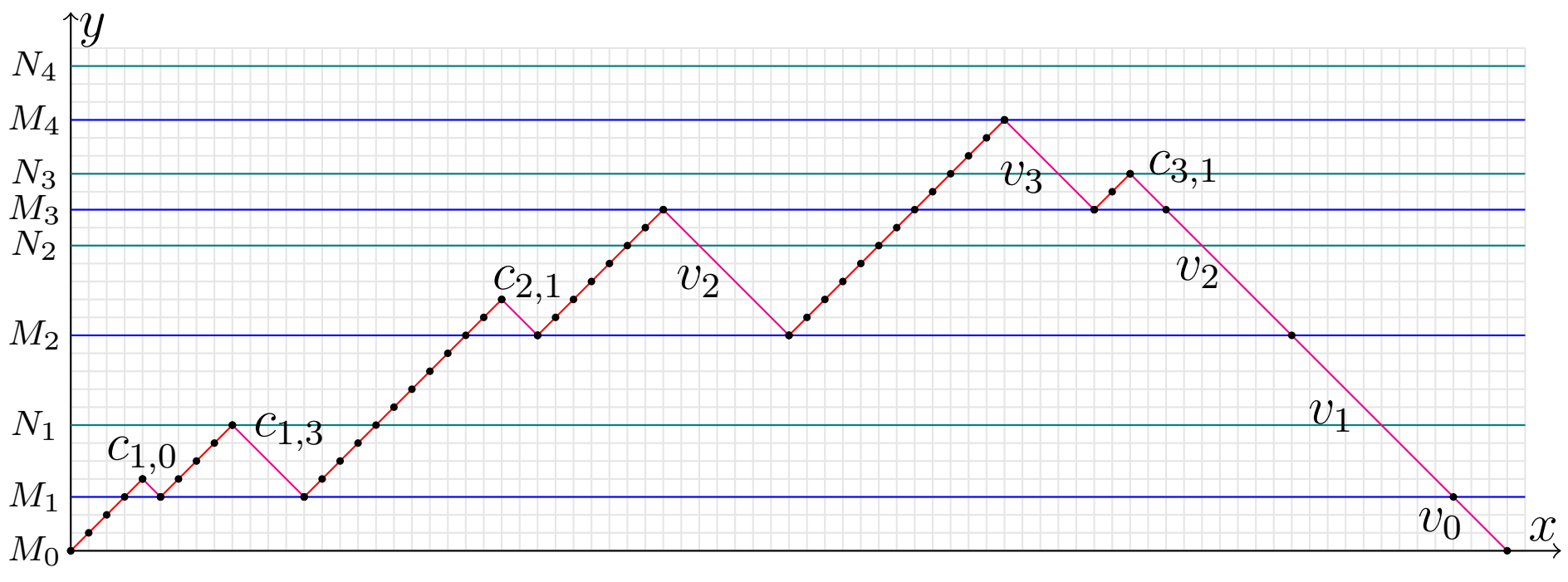
\begin{figure}
\begin{tikzpicture}
\draw [gray!20](0.0000,0.0000)--(11.5830,0.0000);
\draw [gray!20](0.0000,0.1430)--(11.5830,0.1430);
\draw [gray!20](0.0000,0.2860)--(11.5830,0.2860);
\draw [gray!20](0.0000,0.4290)--(11.5830,0.4290);
\draw [gray!20](0.0000,0.5720)--(11.5830,0.5720);
\draw [gray!20](0.0000,0.7150)--(11.5830,0.7150);
\draw [gray!20](0.0000,0.8580)--(11.5830,0.8580);
\draw [gray!20](0.0000,1.0010)--(11.5830,1.0010);
\draw [gray!20](0.0000,1.1440)--(11.5830,1.1440);
\draw [gray!20](0.0000,1.2870)--(11.5830,1.2870);
\draw [gray!20](0.0000,1.4300)--(11.5830,1.4300);
\draw [gray!20](0.0000,1.5730)--(11.5830,1.5730);
\draw [gray!20](0.0000,1.7160)--(11.5830,1.7160);
\draw [gray!20](0.0000,1.8590)--(11.5830,1.8590);
\draw [gray!20](0.0000,2.0020)--(11.5830,2.0020);
\draw [gray!20](0.0000,2.1450)--(11.5830,2.1450);
\draw [gray!20](0.0000,2.2880)--(11.5830,2.2880);
\draw [gray!20](0.0000,2.4310)--(11.5830,2.4310);
\draw [gray!20](0.0000,2.5740)--(11.5830,2.5740);
\draw [gray!20](0.0000,2.7170)--(11.5830,2.7170);
\draw [gray!20](0.0000,2.8600)--(11.5830,2.8600);
\draw [gray!20](0.0000,3.0030)--(11.5830,3.0030);
\draw [gray!20](0.0000,3.1460)--(11.5830,3.1460);
\draw [gray!20](0.0000,3.2890)--(11.5830,3.2890);
\draw [gray!20](0.0000,3.4320)--(11.5830,3.4320);
\draw [gray!20](0.0000,3.5750)--(11.5830,3.5750);
\draw [gray!20](0.0000,3.7180)--(11.5830,3.7180);
\draw [gray!20](0.0000,3.8610)--(11.5830,3.8610);
\draw [gray!20](0.0000,4.0040)--(11.5830,4.0040);
\draw [gray!20](0.0000,0.0000)--(0.0000,4.0040);
\draw [gray!20](0.1430,0.0000)--(0.1430,4.0040);
\draw [gray!20](0.2860,0.0000)--(0.2860,4.0040);
\draw [gray!20](0.4290,0.0000)--(0.4290,4.0040);
\draw [gray!20](0.5720,0.0000)--(0.5720,4.0040);
\draw [gray!20](0.7150,0.0000)--(0.7150,4.0040);
\draw [gray!20](0.8580,0.0000)--(0.8580,4.0040);
\draw [gray!20](1.0010,0.0000)--(1.0010,4.0040);
\draw [gray!20](1.1440,0.0000)--(1.1440,4.0040);
\draw [gray!20](1.2870,0.0000)--(1.2870,4.0040);
\draw [gray!20](1.4300,0.0000)--(1.4300,4.0040);
\draw [gray!20](1.5730,0.0000)--(1.5730,4.0040);
\draw [gray!20](1.7160,0.0000)--(1.7160,4.0040);
\draw [gray!20](1.8590,0.0000)--(1.8590,4.0040);
\draw [gray!20](2.0020,0.0000)--(2.0020,4.0040);
\draw [gray!20](2.1450,0.0000)--(2.1450,4.0040);
\draw [gray!20](2.2880,0.0000)--(2.2880,4.0040);
\draw [gray!20](2.4310,0.0000)--(2.4310,4.0040);
\draw [gray!20](2.5740,0.0000)--(2.5740,4.0040);
\draw [gray!20](2.7170,0.0000)--(2.7170,4.0040);
\draw [gray!20](2.8600,0.0000)--(2.8600,4.0040);
\draw [gray!20](3.0030,0.0000)--(3.0030,4.0040);
\draw [gray!20](3.1460,0.0000)--(3.1460,4.0040);
\draw [gray!20](3.2890,0.0000)--(3.2890,4.0040);
\draw [gray!20](3.4320,0.0000)--(3.4320,4.0040);
\draw [gray!20](3.5750,0.0000)--(3.5750,4.0040);
\draw [gray!20](3.7180,0.0000)--(3.7180,4.0040);
\draw [gray!20](3.8610,0.0000)--(3.8610,4.0040);
\draw [gray!20](4.0040,0.0000)--(4.0040,4.0040);
\draw [gray!20](4.1470,0.0000)--(4.1470,4.0040);
\draw [gray!20](4.2900,0.0000)--(4.2900,4.0040);
\draw [gray!20](4.4330,0.0000)--(4.4330,4.0040);
\draw [gray!20](4.5760,0.0000)--(4.5760,4.0040);
\draw [gray!20](4.7190,0.0000)--(4.7190,4.0040);
\draw [gray!20](4.8620,0.0000)--(4.8620,4.0040);
\draw [gray!20](5.0050,0.0000)--(5.0050,4.0040);
\draw [gray!20](5.1480,0.0000)--(5.1480,4.0040);
\draw [gray!20](5.2910,0.0000)--(5.2910,4.0040);
\draw [gray!20](5.4340,0.0000)--(5.4340,4.0040);
\draw [gray!20](5.5770,0.0000)--(5.5770,4.0040);
\draw [gray!20](5.7200,0.0000)--(5.7200,4.0040);
\draw [gray!20](5.8630,0.0000)--(5.8630,4.0040);
\draw [gray!20](6.0060,0.0000)--(6.0060,4.0040);
\draw [gray!20](6.1490,0.0000)--(6.1490,4.0040);
\draw [gray!20](6.2920,0.0000)--(6.2920,4.0040);
\draw [gray!20](6.4350,0.0000)--(6.4350,4.0040);
\draw [gray!20](6.5780,0.0000)--(6.5780,4.0040);
\draw [gray!20](6.7210,0.0000)--(6.7210,4.0040);
\draw [gray!20](6.8640,0.0000)--(6.8640,4.0040);
\draw [gray!20](7.0070,0.0000)--(7.0070,4.0040);
\draw [gray!20](7.1500,0.0000)--(7.1500,4.0040);
\draw [gray!20](7.2930,0.0000)--(7.2930,4.0040);
\draw [gray!20](7.4360,0.0000)--(7.4360,4.0040);
\draw [gray!20](7.5790,0.0000)--(7.5790,4.0040);
\draw [gray!20](7.7220,0.0000)--(7.7220,4.0040);
\draw [gray!20](7.8650,0.0000)--(7.8650,4.0040);
\draw [gray!20](8.0080,0.0000)--(8.0080,4.0040);
\draw [gray!20](8.1510,0.0000)--(8.1510,4.0040);
\draw [gray!20](8.2940,0.0000)--(8.2940,4.0040);
\draw [gray!20](8.4370,0.0000)--(8.4370,4.0040);
\draw [gray!20](8.5800,0.0000)--(8.5800,4.0040);
\draw [gray!20](8.7230,0.0000)--(8.7230,4.0040);
\draw [gray!20](8.8660,0.0000)--(8.8660,4.0040);
\draw [gray!20](9.0090,0.0000)--(9.0090,4.0040);
\draw [gray!20](9.1520,0.0000)--(9.1520,4.0040);
\draw [gray!20](9.2950,0.0000)--(9.2950,4.0040);
\draw [gray!20](9.4380,0.0000)--(9.4380,4.0040);
\draw [gray!20](9.5810,0.0000)--(9.5810,4.0040);
\draw [gray!20](9.7240,0.0000)--(9.7240,4.0040);
\draw [gray!20](9.8670,0.0000)--(9.8670,4.0040);
\draw [gray!20](10.0100,0.0000)--(10.0100,4.0040);
\draw [gray!20](10.1530,0.0000)--(10.1530,4.0040);
\draw [gray!20](10.2960,0.0000)--(10.2960,4.0040);
\draw [gray!20](10.4390,0.0000)--(10.4390,4.0040);
\draw [gray!20](10.5820,0.0000)--(10.5820,4.0040);
\draw [gray!20](10.7250,0.0000)--(10.7250,4.0040);
\draw [gray!20](10.8680,0.0000)--(10.8680,4.0040);
\draw [gray!20](11.0110,0.0000)--(11.0110,4.0040);
\draw [gray!20](11.1540,0.0000)--(11.1540,4.0040);
\draw [gray!20](11.2970,0.0000)--(11.2970,4.0040);
\draw [gray!20](11.4400,0.0000)--(11.4400,4.0040);
\draw [gray!20](11.5830,0.0000)--(11.5830,4.0040);
\draw [black, ->](0.0000,0.0000)--(0.0000,4.2900);
\draw [black, ->](0.0000,0.0000)--(11.8690,0.0000);
\draw [black](11.7260,0.1716) node [] {$x$};
\draw [black](0.1716,4.1470) node [] {$y$};
\draw [blue](0.0000,0.4290)--(11.5830,0.4290);
\draw [teal](0.0000,1.0010)--(11.5830,1.0010);
\draw [blue](0.0000,1.7160)--(11.5830,1.7160);
\draw [teal](0.0000,2.4310)--(11.5830,2.4310);
\draw [blue](0.0000,2.7170)--(11.5830,2.7170);
\draw [teal](0.0000,3.0030)--(11.5830,3.0030);
\draw [blue](0.0000,3.4320)--(11.5830,3.4320);
\draw [teal](0.0000,3.8610)--(11.5830,3.8610);
\draw [black](-0.2860,1.0010) node [] {{\tiny $N_1$}};
\draw [black](-0.2860,2.4310) node [] {{\tiny $N_2$}};
\draw [black](-0.2860,3.0030) node [] {{\tiny $N_3$}};
\draw [black](-0.2860,3.8610) node [] {{\tiny $N_4$}};
\draw [black](-0.2860,0.0000) node [] {{\tiny $M_0$}};
\draw [black](-0.2860,0.4290) node [] {{\tiny $M_1$}};
\draw [black](-0.2860,1.7160) node [] {{\tiny $M_2$}};
\draw [black](-0.2860,2.7170) node [] {{\tiny $M_3$}};
\draw [black](-0.2860,3.4320) node [] {{\tiny $M_4$}};
\draw [red](0.0000,0.0000)--(0.1430,0.1430)--(0.2860,0.2860)--(0.4290,0.4290)--(0.5720,0.5720);
\fill (0.0000,0.0000) circle(0.0300); 
\fill (0.5720,0.5720) circle(0.0300); 
\fill (0.0000,0.0000) circle(0.0300); 
\fill (0.1430,0.1430) circle(0.0300); 
\fill (0.2860,0.2860) circle(0.0300); 
\fill (0.4290,0.4290) circle(0.0300); 
\draw [magenta](0.5720,0.5720)--(0.7150,0.4290);
\fill (0.5720,0.5720) circle(0.0300); 
\fill (0.7150,0.4290) circle(0.0300); 
\draw [red](0.7150,0.4290)--(0.8580,0.5720)--(1.0010,0.7150)--(1.1440,0.8580)--(1.2870,1.0010);
\fill (0.7150,0.4290) circle(0.0300); 
\fill (1.2870,1.0010) circle(0.0300); 
\fill (0.7150,0.4290) circle(0.0300); 
\fill (0.8580,0.5720) circle(0.0300); 
\fill (1.0010,0.7150) circle(0.0300); 
\fill (1.1440,0.8580) circle(0.0300); 
\draw [magenta](1.2870,1.0010)--(1.4300,0.8580)--(1.5730,0.7150)--(1.7160,0.5720)--(1.8590,0.4290);
\fill (1.2870,1.0010) circle(0.0300); 
\fill (1.8590,0.4290) circle(0.0300); 
\draw [red](1.8590,0.4290)--(2.0020,0.5720)--(2.1450,0.7150)--(2.2880,0.8580)--(2.4310,1.0010)--(2.5740,1.1440)--(2.7170,1.2870)--(2.8600,1.4300)--(3.0030,1.5730)--(3.1460,1.7160)--(3.2890,1.8590)--(3.4320,2.0020);
\fill (1.8590,0.4290) circle(0.0300); 
\fill (3.4320,2.0020) circle(0.0300); 
\fill (1.8590,0.4290) circle(0.0300); 
\fill (2.0020,0.5720) circle(0.0300); 
\fill (2.1450,0.7150) circle(0.0300); 
\fill (2.2880,0.8580) circle(0.0300); 
\fill (2.4310,1.0010) circle(0.0300); 
\fill (2.5740,1.1440) circle(0.0300); 
\fill (2.7170,1.2870) circle(0.0300); 
\fill (2.8600,1.4300) circle(0.0300); 
\fill (3.0030,1.5730) circle(0.0300); 
\fill (3.1460,1.7160) circle(0.0300); 
\fill (3.2890,1.8590) circle(0.0300); 
\draw [magenta](3.4320,2.0020)--(3.5750,1.8590)--(3.7180,1.7160);
\fill (3.4320,2.0020) circle(0.0300); 
\fill (3.7180,1.7160) circle(0.0300); 
\draw [red](3.7180,1.7160)--(3.8610,1.8590)--(4.0040,2.0020)--(4.1470,2.1450)--(4.2900,2.2880)--(4.4330,2.4310)--(4.5760,2.5740)--(4.7190,2.7170);
\fill (3.7180,1.7160) circle(0.0300); 
\fill (4.7190,2.7170) circle(0.0300); 
\fill (3.7180,1.7160) circle(0.0300); 
\fill (3.8610,1.8590) circle(0.0300); 
\fill (4.0040,2.0020) circle(0.0300); 
\fill (4.1470,2.1450) circle(0.0300); 
\fill (4.2900,2.2880) circle(0.0300); 
\fill (4.4330,2.4310) circle(0.0300); 
\fill (4.5760,2.5740) circle(0.0300); 
\draw [magenta](4.7190,2.7170)--(4.8620,2.5740)--(5.0050,2.4310)--(5.1480,2.2880)--(5.2910,2.1450)--(5.4340,2.0020)--(5.5770,1.8590)--(5.7200,1.7160);
\fill (4.7190,2.7170) circle(0.0300); 
\fill (5.7200,1.7160) circle(0.0300); 
\draw [red](5.7200,1.7160)--(5.8630,1.8590)--(6.0060,2.0020)--(6.1490,2.1450)--(6.2920,2.2880)--(6.4350,2.4310)--(6.5780,2.5740)--(6.7210,2.7170)--(6.8640,2.8600)--(7.0070,3.0030)--(7.1500,3.1460)--(7.2930,3.2890)--(7.4360,3.4320);
\fill (5.7200,1.7160) circle(0.0300); 
\fill (7.4360,3.4320) circle(0.0300); 
\fill (5.7200,1.7160) circle(0.0300); 
\fill (5.8630,1.8590) circle(0.0300); 
\fill (6.0060,2.0020) circle(0.0300); 
\fill (6.1490,2.1450) circle(0.0300); 
\fill (6.2920,2.2880) circle(0.0300); 
\fill (6.4350,2.4310) circle(0.0300); 
\fill (6.5780,2.5740) circle(0.0300); 
\fill (6.7210,2.7170) circle(0.0300); 
\fill (6.8640,2.8600) circle(0.0300); 
\fill (7.0070,3.0030) circle(0.0300); 
\fill (7.1500,3.1460) circle(0.0300); 
\fill (7.2930,3.2890) circle(0.0300); 
\draw [magenta](7.4360,3.4320);
\fill (7.4360,3.4320) circle(0.0300); 
\fill (7.4360,3.4320) circle(0.0300); 
\draw [magenta](7.4360,3.4320)--(7.5790,3.2890)--(7.7220,3.1460)--(7.8650,3.0030)--(8.0080,2.8600)--(8.1510,2.7170);
\fill (7.4360,3.4320) circle(0.0300); 
\fill (8.1510,2.7170) circle(0.0300); 
\draw [red](8.1510,2.7170)--(8.2940,2.8600)--(8.4370,3.0030);
\fill (8.1510,2.7170) circle(0.0300); 
\fill (8.4370,3.0030) circle(0.0300); 
\fill (8.1510,2.7170) circle(0.0300); 
\fill (8.2940,2.8600) circle(0.0300); 
\draw [magenta](8.4370,3.0030)--(8.5800,2.8600)--(8.7230,2.7170)--(8.8660,2.5740)--(9.0090,2.4310)--(9.1520,2.2880)--(9.2950,2.1450)--(9.4380,2.0020)--(9.5810,1.8590)--(9.7240,1.7160)--(9.8670,1.5730)--(10.0100,1.4300)--(10.1530,1.2870)--(10.2960,1.1440)--(10.4390,1.0010)--(10.5820,0.8580)--(10.7250,0.7150)--(10.8680,0.5720)--(11.0110,0.4290)--(11.1540,0.2860)--(11.2970,0.1430)--(11.4400,0.0000);
\fill (8.4370,3.0030) circle(0.0300); 
\fill (11.4400,0.0000) circle(0.0300); 
\draw [black](0.5720,0.7579) node [] {$c_{1,0}$};
\draw [black](1.7446,0.9438) node [] {$c_{1,3}$};
\draw [black](3.6465,2.1164) node [] {$c_{2,1}$};
\draw [black](5.0050,2.1450) node [] {$v_2$};
\draw [black](7.5790,2.9744) node [] {$v_3$};
\draw [black](8.8660,3.0316) node [] {$c_{3,1}$};
\draw [black](8.9804,2.2165) node [] {$v_2$};
\draw [black](10.0386,1.0725) node [] {$v_1$};
\draw [black](10.9109,0.2145) node [] {$v_0$};
\fill (8.7230,2.7170) circle(0.0300); 
\fill (9.7240,1.7160) circle(0.0300); 
\fill (11.0110,0.4290) circle(0.0300);

\end{tikzpicture}
\caption{A quasi-Dyck path $p$ for $(k_0,k_1,k_2,k_3)=(3,4,1,2)$}\label{fig:path}
\end{figure}

A quasi-Dyck path $p$ is encoded by the word $p=p_1p_2\cdots p_m$ over the alphabet $\{b,\,c_t\ (t\ge 1),\,v,\,d\}$ recording the types of its successive steps; 
as usual, $x^r$ denotes the word consisting of $r$ copies of the letter $x$. 
For example, let $(k_0,k_1,k_2,k_3)=(3,4,1,2)$. Then
\begin{align*}
(M_0,M_1,M_2,M_3,M_4)&=(0,3,12,19,24), \\
(N_1,N_2,N_3,N_4)&=(7,17,21,27),
\end{align*}
and the path 
\begin{equation}\label{eq:examplepath}
p=b^4c_1b^4c_1b^{11}c_2b^7vb^{12}vb^2c_3v^2d
\end{equation}
shown in Figure~\ref{fig:path} is a quasi-Dyck path with $40$ up steps and $9$ down steps; it ends at $(80,0)$. 
Its two short down steps $c_1$ have drops $1$ and $4$, its short down steps $c_2$ and $c_3$ both have drop $2$, and its long down steps have drops $7$ ($M_3\to M_2$), $5$ ($M_4\to M_3$), $7$ 
($M_3\to M_2$) and $9$ ($M_2\to M_1$). 

We now assign weights to the steps. 
Let $\mathbf{v}=(v_0,v_1,v_2,\ldots)$ and $\mathbf{c}=(c_{t,r})_{t\geq 1,\,0\le r\le k_{t-1}}$ 
be given; 
in the applications $v_t$ and $c_{t,r}$ will be the coefficients of the $H$-fraction \eqref{eq:Hfrac}, $c_{t,r}$ being the coefficient of $x^r$ in $u_t(x)$. 

\begin{Definition}\label{def:weight}
The {\it weight} $w(\xi)$ of a step $\xi$ of a quasi-Dyck path is 
\[
w(\xi)=
\begin{cases}
1, & \text{if $\xi$ is an up step $b$}, \\
v_0, & \text{if $\xi$ is the terminal step $d$}, \\
v_t, & \text{if $\xi$ is the long down step from $M_{t+1}$ to $M_t$}, \\
c_{t,r}, & \text{if $\xi$ is a short down step $c_{t,r}$}.
\end{cases}
\]
\end{Definition}
The weight of a quasi-Dyck path $p=p_1p_2\cdots p_m$ is the product of the weights of its steps, 
\[
w(p)=\prod_{i=1}^{m}w(p_i).
\]
Since every quasi-Dyck path contains exactly one terminal step, of weight $v_0$, we have 
\begin{equation}\label{eq:what}
w(p)=v_0\,\widehat{w}(p),
\end{equation}
where $\widehat{w}(p)$ denotes the product of the weights of the non-terminal steps of $p$.

Since up steps have weight $1$ and the terminal step has weight $v_0$, the weight of a quasi-Dyck path is the product of the labels of its down steps. 
It is convenient to record these labels in their order of appearance: the {\it labeled word} of $p$ is the word $\lambda(p)$ over the infinite alphabet 
\begin{equation}\label{eq:alphabetX}
X=\{v_t : t\ge 0\}\cup\{c_{t,r}: t\ge 1,\ 0\le r\le k_{t-1}\}
\end{equation}
obtained by writing the label of each down step of $p$ (the terminal step being labeled $v_0$), in the order in which the down steps occur. 
For the path \eqref{eq:examplepath} of Figure~\ref{fig:path}, we obtain
\begin{align*}
\lambda(p)&=c_{1,0}\,c_{1,3}\,c_{2,1}\,v_2\,v_3\,c_{3,1}\,v_2\,v_1\,v_0,\\
w(p)&=c_{1,0}c_{1,3}c_{2,1}c_{3,1}\,v_0v_1v_2^2v_3 .
\end{align*}
Between two consecutive down steps a quasi-Dyck path only climbs, from the ending height of the first down step to the starting height of the second one, and both heights are determined by the labels. 
Hence $p$ can be recovered from $\lambda(p)$, and $\lambda$ is a bijection from $\mathcal{QD}$ onto the set $\mathcal{P}:=\lambda(\mathcal{QD})$ of {\it labeled quasi-Dyck paths}, which is a subset of the free monoid $X^*$. 

\subsection{The first return decomposition}\label{sec:flaj:decomp}

For $t\ge1$, a \emph{path at level $t$} is a lattice path that starts and ends at height $M_t$, stays at or above $M_t$, and uses only up steps and short or long down steps.
Let $\mathcal{QD}^{(t)}$ denote the set of these paths, including the empty path, with weights $w(p)$ and labeled words $\lambda(p)$ defined as before.

A path at level $t$ is \emph{primitive} if it is nonempty and returns to $M_t$ only at its last step.
Cutting at each return to $M_t$ gives a unique decomposition into primitive paths.
Every quasi-Dyck path $p$ can be written uniquely as
\begin{equation}\label{eq:factor:qd}
p=b^{M_1}q\,d,
\end{equation}
where $q\in\mathcal{QD}^{(1)}$. The following lemma describes the primitive paths.

\begin{Lemma}\label{lem:primitive}
Let $t\geq 1$. A path $e$ at level $t$ is primitive if and only if it is of one of the following two forms:
\begin{itemize}
\item[(i)] $e=b^{r+1}c_{t,r}$ for some $0\le r\le k_{t-1}$;
\item[(ii)] $e=b^{k_{t-1}+k_t+2}\,q\,v$, where $q\in\mathcal{QD}^{(t+1)}$ and $v$ is the long down step from $M_{t+1}$ to $M_t$.
\end{itemize}
Moreover, the form and the corresponding $r$ or $q$ are uniquely determined by $e$.
\end{Lemma}

\begin{proof}
Let $e$ be primitive.
It begins with an up step, and  it contains a down step, since it ends at height $M_t$. 
Let $y$ be the starting height of its first down step.
Thus $e$ begins with $b^{y-M_t}$, and this down step ends at some level $M_s$ with $s\ge t$.
We distinguish two cases.

\medskip

\noindent{\bf{ Case 1.}} Suppose that $y\le N_t$.
Since $M_s<y\le N_t<M_{t+1}$, we must have $s=t$.
The first down step is therefore $c_{t,r}$, where $r=y-M_t-1$ and $0\le r\le k_{t-1}$.
It returns to $M_t$ and hence, by primitivity, must be the last step. Thus $e=b^{r+1}c_{t,r}$.

\medskip

\noindent{\bf{ Case 2.}} Suppose that $y>N_t$.
The first down step is either the long down step from $M_{t+1}$ to $M_t$, or a down step ending at a level at least $M_{t+1}$.
In either case, $y\ge M_{t+1}$.
After its first $M_{t+1}-M_t=k_{t-1}+k_t+2$ up steps, the path is at height $M_{t+1}$.
Since no down step ends strictly between $M_t$ and $M_{t+1}$, descending below $M_{t+1}$ would return the path to $M_t$.
By primitivity, this can occur only at the final step.
The final step must therefore be the long down step $v$ from $M_{t+1}$ to $M_t$.
The portion $q$ between the prescribed initial climb and this final step starts and ends at $M_{t+1}$ and stays at or above that height.
Hence $q\in\mathcal{QD}^{(t+1)}$ and $e=b^{k_{t-1}+k_t+2}qv.$ 

\medskip

The two cases are disjoint. In Case 1, $r$ is determined by $y$.
In Case 2, deleting the first $M_{t+1}-M_t$ up steps and the final step determines $q$ uniquely.
\end{proof}

\subsection{The path interpretation of $H$-fractions}\label{sec:flaj:main}

As usual, the infinite continued fraction \eqref{eq:Hfrac} is interpreted as the $x$-adic limit in $\setF[[x]]$ of its \emph{convergents} $C_h(x)$, defined for $h\ge1$ by
\begin{equation}\label{eq:convergent}
C_h(x)=
\cFrac{v_0 x^{k_0}}{ 1 - u_1(x) x} 
- \cFrac{ v_1 x^{k_0+k_1+2}}{ 1 - u_2(x)x}
- {\cdots} 
- \cFrac{ v_{h-1} x^{k_{h-2}+k_{h-1}+2}}{ 1 - u_h(x)x}.
\end{equation}
The existence of this limit is part of Theorem \ref{th:super2}(i). 
It will also follow from the proof of the next theorem, which is the $H$-fraction analogue of Theorem \ref{th:flajolet}.

\begin{Theorem}\label{th:charhankel}
Let $F(x)=\sum_{n\geq 0}a_nx^n$ be the $H$-fraction \eqref{eq:Hfrac}, and let the weights of Definition \ref{def:weight} be defined by the coefficients $v_t$ of \eqref{eq:Hfrac} and by the coefficients $c_{t,r}$ of $x^r$ in $u_t(x)$. 
Then, for each $n\geq 0$, 
\begin{equation}\label{eq:an=sumw}
a_n=\sum_{p\in\mathcal{QD}_n}w(p),
\end{equation}
that is, $a_n$ is the sum of the weights of all quasi-Dyck paths of length $2n$. 
\end{Theorem}

\begin{proof}
For a path $p$, let $\mathrm{up}(p)$ denote its number of up steps, so that $|p|=2\,\mathrm{up}(p)$ for $p\in\mathcal{QD}$. 
For $t\geq 1$, consider the generating function
\[
G_t(x)=\sum_{q\in\mathcal{QD}^{(t)}}w(q)\,x^{\mathrm{up}(q)}\in\setF[[x]].
\]
It is a well-defined formal power series, since a path at level $t$ with $n$ up steps has at most $n$ down steps, and hence there are finitely many of them. 
By the unique factorization into primitive paths and by Lemma 
\ref{lem:primitive}, we have
\begin{align*}
G_t(x)&=\sum_{m\geq 0}\Bigl(\sum_{r=0}^{k_{t-1}}c_{t,r}x^{r+1}
+v_tx^{k_{t-1}+k_t+2}G_{t+1}(x)\Bigr)^m\\
&=\frac{1}{1-u_t(x)x-v_tx^{k_{t-1}+k_t+2}G_{t+1}(x)},
\end{align*}
where the geometric series converges because every primitive path has at least one up step. Similarly, by \eqref{eq:factor:qd},
\[
G(x):=\sum_{p\in\mathcal{QD}}w(p)\,x^{\mathrm{up}(p)}
=v_0x^{k_0}G_1(x).
\]
Iterating the relation between $G_t$ and $G_{t+1}$ $h$ times, we get
\begin{equation}\label{eq:Gfinite}
G(x)=
\cFrac{v_0 x^{k_0}}{ 1 - u_1(x) x} 
- \cFrac{ v_1 x^{k_0+k_1+2}}{ 1 - u_2(x)x}
- {\cdots} 
- \cFrac{ v_{h-1} x^{k_{h-2}+k_{h-1}+2}}{ D_h(x)},
\end{equation}
where $D_h(x)=1 - u_h(x)x-v_hx^{e_h}G_{h+1}(x)$, and where we have put $e_0=k_0$ and $e_t=k_{t-1}+k_t+2$ for $t\geq 1$, so that $M_{h+1}=e_0+e_1+\cdots+e_h$ by \eqref{eq:MN}. 
The right-hand side of \eqref{eq:Gfinite} differs from the convergent $C_h(x)$ only in the last denominator, and $D_h(x)\equiv 1-u_h(x)x \pmod{x^{e_h}}$. 
Going up one level, the two fractions $v_{h-1}x^{e_{h-1}}/D_h(x)$ and $v_{h-1}x^{e_{h-1}}/(1-u_h(x)x)$ are congruent modulo $x^{e_{h-1}+e_h}$, since the numerator contributes the factor 
$x^{e_{h-1}}$ and inverting a series with constant term $1$ preserves congruences. 
Repeating this argument at each of the $h$ levels, and taking the final numerator $v_0x^{e_0}$ into account, we obtain
\[
G(x)\equiv C_h(x)\pmod{x^{M_{h+1}}}.
\]
Since $M_{h+1}\geq 2h$ by \eqref{eq:MN}, the convergents $C_h(x)$ converge to $G(x)$ as $h\to\infty$, that is, $F(x)=G(x)$. 
Extracting the coefficient of $x^n$ gives \eqref{eq:an=sumw}.
\end{proof}

The congruence $G\equiv C_h \pmod{x^{M_{h+1}}}$ has a transparent combinatorial meaning: $C_h(x)$ is the generating function of the quasi-Dyck paths that never reach level $M_{h+1}$, and a path 
reaching level $M_{h+1}$ has at least $M_{h+1}$ up steps; the shortest one is $b^{M_{h+1}}v^{h}d$.

When all $k_t=0$, the $H$-fraction reduces to the $J$-fraction \eqref{eq:Jacobi}.
For $t\ge1$, we have $M_t=2t-2$ and $N_t=2t-1$.
The short down steps have drop $1$ and start at odd heights, while the long down steps have drop $2$
and connect consecutive even levels.
The terminal step has drop $0$ but is still recorded.
Apart from the terminal step, every such path decomposes uniquely into factors $bc_{t,0}$, $bb$, and long down steps.
Replace $bc_{t,0}$ by a level step at height $t-1$, $bb$ from height $2t-2$ to $2t$ by an up step
from height $t-1$ to $t$, and the long down step from height $2t$ to $2t-2$ by a down step from height $t$ to $t-1$.
The terminal step is omitted from this geometric compression, with its weight $v_0$ retained as an
overall factor.
We thus obtain a weight-preserving bijection between quasi-Dyck paths of horizontal length $2n$ and
Motzkin paths of length $n$.
Under this correspondence, Theorem~\ref{th:charhankel} reduces to Theorem~\ref{th:flajolet}.

\subsection{Non-commutative version}\label{sec:flaj:noncomm}

Following Flajolet \cite{Flajolet1980}, the identity of Theorem~\ref{th:charhankel} can be lifted to the algebra of formal power series in non-commutative variables, in the sense of Sch\"utzenberger~\cite{Schutzenberger1963}. 
We briefly recall the setting. 

Let $X$ be an alphabet, let $X^*$ be the free monoid of all words over $X$ (including the empty word $\varepsilon$), and let $K$ be a field. 
A {\it non-commutative formal power series} over $X$ with coefficients in $K$ is a mapping $S : X^* \to K$, written in the form
\[
S = \sum_{w \in X^*} S_w \, w.
\]
The set of all such series is denoted by $K\langle\!\langle X \rangle\!\rangle$; it is a $K$-algebra for the addition 
\[
S + T = \sum_{w} (S_w + T_w)\, w\]
and the Cauchy product
\[
S \cdot T = \sum_{w \in X^*} \Bigl( \sum_{uv = w} S_u T_v \Bigr) w.
\]
The valuation of a series $S\neq 0$ is $\mathrm{val}(S) = \min \{|w| : S_w \ne 0\}$, where $|w|$ denotes the length of the word $w$, and $\mathrm{val}(0) = +\infty$. 
A sequence $(S_n)_{n \ge 0}$ converges to $S$ if $\mathrm{val}(S - S_n)\to +\infty$; in other words, the coefficients of $S_n$ stabilize progressively, starting from the shortest words. 
A series $S$ with $\mathrm{val}(S) > 0$ has a {\it quasi-inverse}
\[
(1-S)^{-1} = \sum_{m \ge 0}S^m,
\]
which is the inverse of $1-S$. 

For a subset (a {\it language}) $L \subset X^*$, the {\it characteristic series} of $L$ is 
\[
\mathrm{char}(L) = \sum_{u \in L} u.
\] 
For languages $E, F \subset X^*$, we write $E + F$ for $E \cup F$, 
$E \cdot F = \{ uv : u \in E,\ v \in F \}$, and 
$E^* = \varepsilon + E + E^2 + E^3 + \cdots$. 
The following classical lemma is immediate from the definitions.

\begin{Lemma}\label{lem:char}
Let $E$ and $F$ be languages over $X$. 
\begin{itemize}
    \item[(1)] If $E$ and $F$ are disjoint, then
    \[
    \mathrm{char}(E + F) = \mathrm{char}(E) + \mathrm{char}(F).
    \]
    \item[(2)] If the product $E \cdot F$ is unambiguous, i.e., every word $w \in E \cdot F$ admits a unique factorization $w = uv$ with $u \in E$ and $v \in F$, then
    \[
    \mathrm{char}(E \cdot F) = \mathrm{char}(E)\,\mathrm{char}(F).
    \]
    \item[(3)] If $\varepsilon\notin E$ and the star $E^*$ is unambiguous, that is, every word in $E^*$ admits a unique factorization into elements of $E$, then
    \[
    \mathrm{char}(E^*) = \bigl(1 - \mathrm{char}(E)\bigr)^{-1}.
    \]
\end{itemize}
\end{Lemma}

Let now $X$ be the alphabet \eqref{eq:alphabetX} and $\mathcal{P}=\lambda(\mathcal{QD})\subset X^*$ be the language of labeled quasi-Dyck paths. 
For $t\geq 1$, put $\mathcal{P}^{(t)}=\lambda(\mathcal{QD}^{(t)})$ and 
\[
u_t = c_{t,0} + c_{t,1} + \cdots + c_{t,k_{t-1}}\in K\langle\!\langle X \rangle\!\rangle .
\]
After deleting the up steps, we can rewrite \eqref{eq:factor:qd} as $\mathcal{P}=\mathcal{P}^{(1)}\cdot v_0.$
For every $t\ge1$, Lemma~\ref{lem:primitive} similarly gives
\[
\mathcal{P}^{(t)}
=
\bigl(
\{c_{t,0},\ldots,c_{t,k_{t-1}}\}
+\mathcal{P}^{(t+1)}\cdot v_t
\bigr)^*.
\]
The uniqueness statements in Section~\ref{sec:flaj:decomp} ensure that all products and stars in these expressions are unambiguous.
Let $\mathcal{P}^{[h]}$ (resp. $\mathcal{P}^{(t),[h]}$) be the set of labeled paths in $\mathcal{P}$ (resp. $\mathcal{P}^{(t)}$) whose underlying path never reaches level $M_{h+1}$. 
Since a path can only come down from level $M_{h+1}$ by the long down step $v_h$, these are the labeled paths using no letter $v_s$ with $s\geq h$ and no letter $c_{s,r}$ with $s\geq h+1$. 
A path at level $h+1$ starts at height $M_{h+1}$, so $\mathcal{P}^{(h+1),[h]}=\emptyset$.
For every $1\le t\le h$, Lemma~\ref{lem:char} gives
\[
\mathrm{char}(\mathcal{P}^{(t),[h]})
=
\bigl(
1-u_t-\mathrm{char}(\mathcal{P}^{(t+1),[h]})\,v_t
\bigr)^{-1}.
\]
The same lemma also yields
\[
\mathrm{char}(\mathcal{P}^{[h]})
=
\mathrm{char}(\mathcal{P}^{(1),[h]})\,v_0.
\]
Iterating these relations yields the following noncommutative analogue of the convergents defined in \eqref{eq:convergent}.

\begin{Theorem}\label{th:lim=char}
For $h\ge1$, set $Q_h^{[h]}=(1-u_h)^{-1}$.
For $1\le t<h$, define $Q_t^{[h]}$ recursively by
\begin{equation}\label{eq:setlabel}
Q_t^{[h]}
=\bigl(1-u_t-Q_{t+1}^{[h]}v_t\bigr)^{-1}.
\end{equation}
Set $Q^{[h]}=Q_1^{[h]}v_0$.
Then $Q^{[h]}=\mathrm{char}(\mathcal{P}^{[h]})$.
Moreover, the sequence $(Q^{[h]})_{h\ge1}$ converges in $K\langle\!\langle X\rangle\!\rangle$, and
\begin{equation}\label{eq:char}
\lim_{h\to\infty}Q^{[h]}
=\mathrm{char}(\mathcal{P}).
\end{equation}
\end{Theorem}

\begin{proof}

By Lemma~\ref{lem:char}, the characteristic series $\mathrm{char}(\mathcal{P}^{(t),[h]})$ satisfy
the same recurrence \eqref{eq:setlabel} as $Q_t^{[h]}$ and have the same value $(1-u_h)^{-1}$ at $t=h$.
For every $1\le t\le h$, backward induction therefore gives
\[
Q_t^{[h]}
=
\mathrm{char}(\mathcal{P}^{(t),[h]}).
\]
Taking $t=1$ and multiplying on the right by $v_0$, we obtain
\[
Q^{[h]}=\mathrm{char}(\mathcal{P}^{[h]}).
\]
By definition, $\mathcal{P}^{[h]}$ consists of the labeled paths that never reach height $M_{h+1}$.
Thus, $\mathrm{char}(\mathcal{P})-Q^{[h]}$ is the characteristic series of the labeled paths that
reach height $M_{h+1}$.

Every such path must eventually return to height $0$.
For each $1\le t\le h$, the only allowed step from height at least $M_{t+1}$ to a height below $M_{t+1}$ is the long down step from $M_{t+1}$ to $M_t$.
Hence its label contains at least one occurrence of each of the letters $v_h,\ldots,v_1$, together with the letter $v_0$ corresponding to the terminal step.
Its length is therefore at least $h+1$, so
\[
\mathrm{val}\bigl(\mathrm{char}(\mathcal{P})-Q^{[h]}\bigr)
\ge h+1.
\]
Since this lower bound tends to infinity as $h\to\infty$, we obtain \eqref{eq:char}.
\end{proof}

We recover Theorem~\ref{th:charhankel} from Theorem~\ref{th:lim=char} by the following substitution.
Let $X'=X\setminus\{v_0\}$.
For $t\ge1$ and $0\le r\le k_{t-1}$, set
\[
\begin{aligned}
\phi(v_t)&=v_tx^{k_{t-1}+k_t+2},\\
\phi(c_{t,r})&=c_{t,r}x^{r+1},
\end{aligned}
\]
and extend $\phi$ multiplicatively to words in $X'$.
Every word occurring in $Q^{[h]}$ or $\mathrm{char}(\mathcal{P})$ has the form $w\,v_0$, where $w\in X'^*$, and we map
\[
w\,v_0\longmapsto\phi(w)\,v_0x^{k_0}.
\]
This substitution sends $u_t$ to $u_t(x)x$ and, by \eqref{eq:setlabel}, sends $Q^{[h]}$ to $C_h(x)$.
It also sends $\lambda(p)$ to $w(p)x^{\mathrm{up}(p)}$, since the total drop of a quasi-Dyck path equals its number of up steps.
The substitution is well defined on these path series, since only finitely many paths have a given number of up steps.
Every letter of $X'$ contributes a positive power of $x$, so the substitution preserves convergence.
Applying it to \eqref{eq:char} therefore gives
\[
\sum_{p\in\mathcal{QD}}w(p)x^{\mathrm{up}(p)}
=
\lim_{h\to\infty}C_h(x)
=
F(x),
\]
as asserted in Theorem~\ref{th:charhankel}.


\section{Non-intersecting lattice paths}\label{sec:lgv}
In \cite{Han2016Adv}, Han gave an algebraic proof of the Hankel determinants in Theorem \ref{th:super2}; a combinatorial evaluation of these determinants, in terms of Roblet's P-paths, is already contained in \cite[Proposition 32]{Roblet1994} (see the concluding remarks). 
The purpose of this section is to explain how these Hankel determinants can be interpreted combinatorially in terms of non-intersecting lattice paths, and to set up the framework used in Sections~\ref{sec:shift}--\ref{sec:backsh} for the $1$-shifted, $2$-shifted and backward shifted Hankel determinants of $F(x)$.

\subsection{The directed graph model}
\label{subsec:lgv:graph}

To apply the Lindstr\"om--Gessel--Viennot theorem, we regard quasi-Dyck paths as directed paths
in a weighted graph $\Gamma$.

The vertices of $\Gamma$ are the ordinary vertices $(x,h)\in\mathbb Z\times\mathbb Z_{\ge0}$, together with a \emph{source} $(x,0)_{\mathrm{src}}$ and a \emph{sink} $(x,0)_{\mathrm{snk}}$ for each $x\in\mathbb Z$.
The source $(x,0)_{\mathrm{src}}$ and the sink $(x,0)_{\mathrm{snk}}$ always lie on the boundary level $M_0=0$.
The ordinary vertex $(x,0)$ lies on level $M_1$ if and only if $k_0=0$, in which case $M_0=M_1=0$.
Nevertheless, the three vertices $(x,0)_{\mathrm{src}}$, $(x,0)$, and $(x,0)_{\mathrm{snk}}$ are pairwise distinct, although they are drawn at the same point in the plane.

The edges between ordinary vertices are the up steps and the short and long down steps of Definition~\ref{def:step}, translated horizontally.
Each source has exactly one outgoing edge, and each sink has exactly one incoming edge. These edges are defined as follows.
\begin{itemize}
\item If $k_0=0$, they are
\[
(x,0)_{\mathrm{src}}
\longrightarrow (x,0)
\xrightarrow{d}(x,0)_{\mathrm{snk}},
\]
where the first edge is an auxiliary \emph{entrance edge} and the second is a terminal step.

\item If $k_0>0$, they are
\[
(x,0)_{\mathrm{src}}\xrightarrow{b}(x+1,1), \qquad (x-k_0,M_1)\xrightarrow{d}(x,0)_{\mathrm{snk}},
\]
namely, an initial up step and a terminal step, respectively.
\end{itemize}
Sources have no incoming edges, and sinks have no outgoing edges.
Each auxiliary entrance edge has weight $1$.
All other edges carry the weights of Definition~\ref{def:weight}; in particular, every terminal step has weight $v_0$.
The weight of a path is the product of the weights of its edges.

Under this convention, quasi-Dyck paths are identified with directed paths from sources to sinks in $\Gamma$.
The auxiliary entrance edge, when present, is not counted among the steps of the underlying quasi-Dyck path.
A path starting at $(x_0,0)_{\mathrm{src}}$ with $n$ up steps ends at $(x_0+2n,0)_{\mathrm{snk}}$.

The graph $\Gamma$ is acyclic: every edge strictly increases the $x$-coordinate, except the entrance and terminal edges when $k_0=0$, which respectively leave a source and enter a sink.
Moreover, only finitely many directed paths join any fixed source and sink.

For each $i\ge0$, define
\[
A_i:=(-2i,0)_{\mathrm{src}}, \qquad B_i:=(2i,0)_{\mathrm{snk}}.
\]
The vertices of a path are the endpoints of its edges.
For $n\ge1$, a family $(P_0,\ldots,P_{n-1})$ of paths in $\Gamma$ is called \emph{non-intersecting} if no two distinct paths share a vertex.

\begin{Remark}\label{rem:path-intersections}
Intersection is understood in terms of graph vertices, with the following consequences.
\begin{enumerate}
\item[(i)]
A down step of drop $\ell\ge2$ is a single edge, and the lattice points in its interior are not vertices of the path.
Thus, non-intersecting paths may cross geometrically at non-lattice points, or at a lattice point that is a vertex of one path and lies in the interior of a down step of another.
For example, in Figure~\ref{fig:h(0)}, the long down step of $P_4$ from $(-1,7)$ to $(4,2)$ passes through the vertex $(0,6)$ of $P_3$ without creating a common graph vertex.

\item[(ii)]
Suppose that $k_0=0$, and let $C_0=(0,0)$ be the ordinary vertex at the origin.
The quasi-Dyck path consisting of the single terminal step $d$ is represented by
\[
A_0\longrightarrow C_0\xrightarrow{d}B_0
\]
and has weight $v_0$.
Any path from $A_s$ to $B_0$ and any path from $A_0$ to $B_t$, with $s,t\ge1$, both pass through $C_0$.
They therefore cannot belong to the same non-intersecting family.

\item[(iii)]
Suppose that $k_0>0$.
No path from $A_0$ to $B_0$ exists, since every edge has positive horizontal displacement.
Any existing path from $A_s$ to $B_0$ lies strictly to the left of the vertical axis except at its sink, whereas any existing path from $A_0$ to $B_t$ lies strictly to the right except at its source.
Since $A_0$ and $B_0$ are distinct graph vertices, these two paths share no vertex and are therefore non-intersecting.
\end{enumerate}
\end{Remark}

\subsection{Non-intersecting path families}
\label{subsec:lgv:families}

Let $\mathfrak{S}_{[m,n]}$ denote the set of all permutations of $\{m, m+1, \dots, n\}$,  that is, all rearrangements of the integers $m, m+1, \dots, n$. 
In particular, we write $[0,n-1] = \{0,1,2,\dots,n-1\}.$
For a permutation $\sigma = \sigma_1 \sigma_2 \cdots \sigma_{n-m+1} \in \mathfrak{S}_{[m,n]}$, 
a pair $(i,j)$ is called \emph{inversion pair} of $\sigma$ if $i < j$ and $\sigma_i > \sigma_j$. 
We denote by $\mathrm{inv}(\sigma)$  the number of inversion pairs of $\sigma$.

Combining Theorem~\ref{th:charhankel} with the Lindstr\"om--Gessel--Viennot theorem gives the following theorem.

\begin{Theorem}\label{th:lGV}
Let $A_i=(-2i,0)_{\mathrm{src}}$ and $B_i=(2i,0)_{\mathrm{snk}}$ for $0\le i\le n-1$.
Then
\begin{equation}
\det[a_{i+j}]_{0\le i,j\le n-1}
=
\sum_{(\sigma;P_0,\ldots,P_{n-1})}
(-1)^{\mathrm{inv}(\sigma)}
\prod_{i=0}^{n-1}w(P_i),
\end{equation}
where the sum runs over all $\sigma\in\mathfrak{S}_{[0,n-1]}$ and all non-intersecting families $(P_0,\ldots,P_{n-1})$ of quasi-Dyck paths such that, for each $0\le i\le n-1$, the path $P_i$ goes from $A_i$ to $B_{\sigma(i)}$.
\end{Theorem}

To give a combinatorial interpretation of \eqref{eq:HankelDetFundamental}, we introduce the following notation.
Let $\Omega_n$ denote the set of all non-intersecting families $(P_0,\ldots,P_{n-1})$ of quasi-Dyck paths such that
\[
P_i:A_i\longrightarrow B_{\sigma(i)}
\qquad (0\le i\le n-1)
\]
for some permutation $\sigma$ of $\{0,\ldots,n-1\}$. For convenience, we write $P_i=(A_i,B_{\sigma(i)})$ to indicate that $P_i$ is a path from $A_i$ to $B_{\sigma(i)}$.
The permutation $\sigma$ is uniquely determined by the family and is called its associated permutation.
Writing $\sigma=\sigma_1\sigma_2\cdots\sigma_n$ as a word, we have $\sigma(i)=\sigma_{i+1}$,  so that $P_i$ ends at $B_{\sigma_{i+1}}$.
For such a family and $0\le m\le r\le n-1$, we write
\[
\mathcal{Q}_{m,r}:=\prod_{i=m}^{r}w(P_i).
\]

The following lemma determines the initial and final segments forced by non-intersection.

\begin{Lemma}\label{lem:forced-segments}
Let $j\ge1$, and let $(P_0,\ldots,P_{s_j-1})\in\Omega_{s_j}$ have associated permutation $\sigma$.
For each $1\le i\le j$ and $0\le t\le s_j-1$, the following hold:
\begin{enumerate}
\item[(i)]
If $t\ge s_{i-1}$, then $P_t$ begins with $b^{M_i}$.
\item[(ii)]
If $\sigma_{t+1}\ge s_{i-1}$, then $P_t$ ends with $v^{i-1}d$.
\end{enumerate}
Entrance edges, when present, are omitted from the step words.
\end{Lemma}

\begin{proof}
We prove both assertions simultaneously by induction on $i$.
The case $i=1$ follows directly from the allowed steps.
Suppose that both assertions hold through level $i$, where $i<j$.

\medskip

\noindent{ \it{Proof of (i).}}
Consider $P_t$ with $t\ge s_i$.
By the induction hypothesis, it begins with $b^{M_i}$.
Suppose that its first down step starts below $M_{i+1}$.
By Definition~\ref{def:step}, we distinguish three cases.

\medskip
\noindent\textbf{Case 1.}
Suppose that the first down step is a short down step of drop $r$, where $1\le r\le k_{i-1}+1$.
It lands at $(M_i-2(t-r),M_i).$
Since $t-r\ge s_i-(k_{i-1}+1)=s_{i-1}$, the induction hypothesis shows that this vertex is the endpoint of the initial climb of $P_{t-r}$.
This contradicts non-intersection.

\medskip
\noindent\textbf{Case 2.}
Suppose that the first down step is the long down step from $M_i$ to $M_{i-1}$, necessarily with $i\ge2$.
It lands at the endpoint, at height $M_{i-1}$, of the initial climb of $P_{t-(M_i-M_{i-1})}$.
Since
\[
t-(M_i-M_{i-1})
\ge s_i-(M_i-M_{i-1})
=s_{i-2},
\]
this climb is forced by the induction hypothesis.
This again contradicts non-intersection.

\medskip
\noindent\textbf{Case 3.}
Suppose that the first down step is the terminal step, necessarily with $i=1$.
The path then has step word $b^{M_1}d$, which requires a sink with index $M_1-t<0,$ since $t\ge s_1=M_1+1$.
This is impossible.

\medskip
Thus the first down step cannot start below $M_{i+1}$, and $P_t$ begins with $b^{M_{i+1}}$.

\medskip

\noindent{\it{Proof of (ii).}} Consider $P_t$ with $x:=\sigma_{t+1}\ge s_i$.
By the induction hypothesis, it ends with $v^{i-1}d$.
We show that the step immediately preceding this final descent can be neither a short down step
nor an up step.

\medskip
\noindent\textbf{Case 1.}
Suppose that the preceding step is a short down step of drop $r$, where $1\le r\le k_{i-1}+1$.
Since no down step ends strictly between $M_i$ and $M_{i+1}$, it must be preceded by $r$ consecutive up steps starting at $(2(x-r)-M_i,M_i).$
Since
\[
x-r\ge s_i-(k_{i-1}+1)=s_{i-1},
\]
the induction hypothesis shows that this vertex is also the starting point of the forced final descent
of the path ending at $B_{x-r}$.
The two paths have different sinks, so this contradicts non-intersection.

\medskip
\noindent\textbf{Case 2.}
Suppose that the preceding step is an up step.
We distinguish two cases.
\begin{itemize}
\item[(1)]
If $i\ge2$, this step belongs to a consecutive climb passing through height $M_{i-1}$, since no down step ends strictly between $M_{i-1}$ and $M_i$.
The vertex of this climb at height $M_{i-1}$ is $(2x-2M_i+M_{i-1},M_{i-1}).$
Since
\[
x-(M_i-M_{i-1})
\ge s_i-(M_i-M_{i-1})
=s_{i-2},
\]
the induction hypothesis shows that the forced final descent
of the path ending at $B_{x-(M_i-M_{i-1})}$ also starts at this vertex.
This contradicts non-intersection.

\item[(2)]
If $i=1$, the allowed steps force the path to have step word $b^{M_1}d$.
Its endpoints would then satisfy $t=M_1-x<0,$ since $x\ge s_1=M_1+1$.
This contradicts $t\ge0$.
\end{itemize}

The same endpoint calculation excludes the path consisting only of $d$ when $k_0=0$.
Thus the preceding step must be the long down step from $M_{i+1}$ to $M_i$.
Hence $P_t$ ends with $v^id$, completing the induction.
\end{proof}

\subsection{Proof of Theorem~\ref{th:super2}}\label{sec:lgv:proofthm}
Let $j\ge1$. For each $0\le i\le s_j-1$, let $A_i=(-2i,0)_{\mathrm{src}}$ and $B_i=(2i,0)_{\mathrm{snk}}$
denote the corresponding source and sink.
Let $\sigma=\sigma_1\cdots\sigma_{s_j}$ be a permutation in $\mathfrak{S}_{[0,s_j-1]}$.
We decompose $\sigma$ into consecutive blocks
\[
\sigma=\sigma^1\sigma^2\cdots\sigma^j.
\]
For each $1\le i\le j$, the block $\sigma^i$ is defined by
\[
\sigma^i
=
\sigma_{s_{i-1}+1}\sigma_{s_{i-1}+2}\cdots\sigma_{s_i}.
\]
Using Lemma~\ref{lem:forced-segments}, we now determine the permutation and the paths in each block.

\begin{Lemma}\label{lem:psigma}
Let $j\ge1$, and let $(P_0,\ldots,P_{s_j-1})\in\Omega_{s_j}$ have associated permutation $\sigma$.
For each $1\le i\le j$, the block $\sigma^i=\sigma_{s_{i-1}+1}\cdots\sigma_{s_i}$ satisfies
\[
\sigma^i=(s_i-1)\,(s_i-2)\cdots s_{i-1}.
\]
For $0\le\ell\le k_{i-1}$, the corresponding paths are
\[
P_{s_{i-1}+\ell}
=
(A_{s_{i-1}+\ell},B_{s_i-1-\ell})
=
b^{M_i}v^{i-1}d.
\]
Moreover, the paths specified above form the unique family in $\Omega_{s_j}$.
\end{Lemma}

\begin{proof}
We first prove the formula for $\sigma^i$.
Suppose otherwise.
Let $i$ be the smallest index for which the identity fails, and let $\ell\in\{1,\ldots,k_{i-1}+1\}$ be the smallest integer such that $\sigma_{s_{i-1}+\ell}\ne s_i-\ell.$
Set $x:=\sigma_{s_{i-1}+\ell}$.
The corresponding path is $P_{s_{i-1}+\ell-1} = (A_{s_{i-1}+\ell-1},B_x).$
By the minimality of $i$, all values smaller than $s_{i-1}$ already occur in the preceding blocks.
Since $\sigma$ is a permutation, we have $x\ge s_{i-1}$.

\medskip
\noindent\textbf{Case 1.}
Suppose that $s_{i-1}\le x\le s_i-1$.
By the minimality of $\ell$ and the assumption $x\ne s_i-\ell$, we have $x\le s_i-\ell-1$.
Thus $P_{s_{i-1}+\ell-1}$ has
\[
s_{i-1}+\ell-1+x
\le s_{i-1}+s_i-2
=M_i-1
\]
up steps.
However, Lemma~\ref{lem:forced-segments}(i) shows that this path begins with $b^{M_i}$, a contradiction.

\medskip
\noindent\textbf{Case 2.}
Suppose that $x\ge s_i$.
Since $\sigma$ is a permutation, there exists an index $y$ such that $P_y=(A_y,B_{s_i-\ell}).$
By the minimality of $i$ and $\ell$, we have $y\ge s_{i-1}+\ell$.
By Lemma~\ref{lem:forced-segments}(i), the initial climb of $P_{s_{i-1}+\ell-1}$ ends at
\[
(M_i-2(s_{i-1}+\ell-1),M_i)
=
(2(s_i-\ell)-M_i,M_i).
\]
Since $s_i-\ell\ge s_{i-1}$, Lemma~\ref{lem:forced-segments}(ii) shows that $P_y$ ends with $v^{i-1}d$, whose starting vertex is the same vertex.
These two distinct paths therefore share an ordinary vertex, contradicting non-intersection.
Hence the formula for $\sigma^i$ holds.

\medskip
We now determine the paths.
Fix $1\le i\le j$ and $0\le\ell\le k_{i-1}$.
The formula for $\sigma^i$ gives $P_{s_{i-1}+\ell} = (A_{s_{i-1}+\ell},B_{s_i-1-\ell}).$
This path has
\[
s_{i-1}+\ell+s_i-1-\ell
=
s_{i-1}+s_i-1
=
M_i
\]
up steps.
By Lemma~\ref{lem:forced-segments}, it begins with $b^{M_i}$ and ends with $v^{i-1}d$.
These two segments account for its entire horizontal length, so
$
P_{s_{i-1}+\ell}
=
b^{M_i}v^{i-1}d.
$
Thus $\Omega_{s_j}$ contains at most one family.

Conversely, these step words define quasi-Dyck paths with the required endpoints.
We verify that they are pairwise non-intersecting.
Distinct initial climbs lie on distinct parallel lines, as do the final descents to distinct sinks.
A common ordinary vertex must therefore lie on the initial climb from some $A_t$ and the final descent
to some $B_x$, at a height $M_h$.
Comparing their horizontal coordinates gives
\[
t+x=M_h=s_{h-1}+s_h-1.
\]
By construction, paths in earlier blocks stay below $M_h$, so $t,x\ge s_{h-1}$.
Consequently, $s_{h-1}\le t,x\le s_h-1.$
Since $x=M_h-t$, the prescribed pairing in block $h$ connects $A_t$ to $B_x$.
Thus the two segments belong to the same path.

Hence the constructed paths are pairwise non-intersecting, and $\Omega_{s_j}$ consists of exactly this family.
\end{proof}

We now prove Theorem~\ref{th:super2}.

\begin{proof}[Proof of Theorem~\ref{th:super2}]

The case $n=0$ follows from the convention $H_0(F)=1$.
Let $j\ge1$.
By Lemma~\ref{lem:psigma}, $\Omega_{s_j}$ consists of a unique family, with associated permutation $\sigma$.
We compute its permutation sign and its weight.

Each block $\sigma^i$ is decreasing, and every entry in an earlier block is smaller than every entry in a later block. 
Therefore,
\[
\mathrm{inv}(\sigma)
=
\sum_{i=1}^{j}\mathrm{inv}(\sigma^i)
=
\sum_{i=0}^{j-1}\frac{k_i(k_i+1)}{2}
=
\epsilon_j.
\]
For each $0\le i\le j-1$, every path in block $i+1$ has weight $V_i=v_0v_1\cdots v_i$.
Since this block contains $k_i+1$ paths, we have $\mathcal{Q}_{s_i,s_{i+1}-1}=V_i^{k_i+1}.$
By Theorem~\ref{th:lGV}, we obtain
\[
\begin{aligned}
H_{s_j}(F)
&=(-1)^{\mathrm{inv}(\sigma)}
  \prod_{i=0}^{j-1}\mathcal{Q}_{s_i,s_{i+1}-1}
  =(-1)^{\epsilon_j}
  \prod_{r=0}^{j-1}v_r^{s_j-s_r}.
\end{aligned}
\]
Here the last equality follows because, for each $0\le r\le j-1$, the exponent of $v_r$ is
\[
\sum_{i=r}^{j-1}(k_i+1)=s_j-s_r.
\]
This proves \eqref{eq:HankelDetFundamental}.

Finally, let $j\ge0$ and $s_j<n<s_{j+1}$.
Suppose, for a contradiction, that $\Omega_n$ is nonempty.
The proof of Lemma~\ref{lem:forced-segments} applies to families in $\Omega_n$ for arbitrary $n$.
Indeed, the induction for initial climbs uses only smaller source indices, and that for final descents uses only smaller sink indices.

Let $P_t$ be the path ending at $B_{s_j}$.
Applying the final-descent assertion with $i=j+1$, we see that $P_t$ ends with $v^jd$ and therefore reaches height $M_{j+1}$.
Since a path from $A_t$ to $B_{s_j}$ has $t+s_j$ up steps, we obtain
\[
t+s_j\ge M_{j+1}=s_j+s_{j+1}-1.
\]
Thus $t\ge s_{j+1}-1>n-1,$ contradicting $0\le t\le n-1$.
Hence $\Omega_n=\emptyset$, and Theorem~\ref{th:lGV} gives $H_n(F)=0$.
This completes the proof.

\end{proof}

\begin{figure}
\begin{tikzpicture}
\draw [gray!20](0.0000,0.0000)--(11.8800,0.0000);
\draw [gray!20](0.0000,0.3300)--(11.8800,0.3300);
\draw [gray!20](0.0000,0.6600)--(11.8800,0.6600);
\draw [gray!20](0.0000,0.9900)--(11.8800,0.9900);
\draw [gray!20](0.0000,1.3200)--(11.8800,1.3200);
\draw [gray!20](0.0000,1.6500)--(11.8800,1.6500);
\draw [gray!20](0.0000,1.9800)--(11.8800,1.9800);
\draw [gray!20](0.0000,2.3100)--(11.8800,2.3100);
\draw [gray!20](0.0000,2.6400)--(11.8800,2.6400);
\draw [gray!20](0.0000,2.9700)--(11.8800,2.9700);
\draw [gray!20](0.0000,3.3000)--(11.8800,3.3000);
\draw [gray!20](0.0000,3.6300)--(11.8800,3.6300);
\draw [gray!20](0.0000,3.9600)--(11.8800,3.9600);
\draw [gray!20](0.0000,4.2900)--(11.8800,4.2900);
\draw [gray!20](0.0000,4.6200)--(11.8800,4.6200);
\draw [gray!20](0.0000,4.9500)--(11.8800,4.9500);
\draw [gray!20](0.0000,5.2800)--(11.8800,5.2800);
\draw [gray!20](0.0000,5.6100)--(11.8800,5.6100);
\draw [gray!20](0.0000,5.9400)--(11.8800,5.9400);
\draw [gray!20](0.0000,0.0000)--(0.0000,5.9400);
\draw [gray!20](0.3300,0.0000)--(0.3300,5.9400);
\draw [gray!20](0.6600,0.0000)--(0.6600,5.9400);
\draw [gray!20](0.9900,0.0000)--(0.9900,5.9400);
\draw [gray!20](1.3200,0.0000)--(1.3200,5.9400);
\draw [gray!20](1.6500,0.0000)--(1.6500,5.9400);
\draw [gray!20](1.9800,0.0000)--(1.9800,5.9400);
\draw [gray!20](2.3100,0.0000)--(2.3100,5.9400);
\draw [gray!20](2.6400,0.0000)--(2.6400,5.9400);
\draw [gray!20](2.9700,0.0000)--(2.9700,5.9400);
\draw [gray!20](3.3000,0.0000)--(3.3000,5.9400);
\draw [gray!20](3.6300,0.0000)--(3.6300,5.9400);
\draw [gray!20](3.9600,0.0000)--(3.9600,5.9400);
\draw [gray!20](4.2900,0.0000)--(4.2900,5.9400);
\draw [gray!20](4.6200,0.0000)--(4.6200,5.9400);
\draw [gray!20](4.9500,0.0000)--(4.9500,5.9400);
\draw [gray!20](5.2800,0.0000)--(5.2800,5.9400);
\draw [gray!20](5.6100,0.0000)--(5.6100,5.9400);
\draw [gray!20](5.9400,0.0000)--(5.9400,5.9400);
\draw [gray!20](6.2700,0.0000)--(6.2700,5.9400);
\draw [gray!20](6.6000,0.0000)--(6.6000,5.9400);
\draw [gray!20](6.9300,0.0000)--(6.9300,5.9400);
\draw [gray!20](7.2600,0.0000)--(7.2600,5.9400);
\draw [gray!20](7.5900,0.0000)--(7.5900,5.9400);
\draw [gray!20](7.9200,0.0000)--(7.9200,5.9400);
\draw [gray!20](8.2500,0.0000)--(8.2500,5.9400);
\draw [gray!20](8.5800,0.0000)--(8.5800,5.9400);
\draw [gray!20](8.9100,0.0000)--(8.9100,5.9400);
\draw [gray!20](9.2400,0.0000)--(9.2400,5.9400);
\draw [gray!20](9.5700,0.0000)--(9.5700,5.9400);
\draw [gray!20](9.9000,0.0000)--(9.9000,5.9400);
\draw [gray!20](10.2300,0.0000)--(10.2300,5.9400);
\draw [gray!20](10.5600,0.0000)--(10.5600,5.9400);
\draw [gray!20](10.8900,0.0000)--(10.8900,5.9400);
\draw [gray!20](11.2200,0.0000)--(11.2200,5.9400);
\draw [gray!20](11.5500,0.0000)--(11.5500,5.9400);
\draw [gray!20](11.8800,0.0000)--(11.8800,5.9400);
\draw [blue](0.0000,0.6600)--(11.8800,0.6600);
\draw [teal](0.0000,1.6500)--(11.8800,1.6500);
\draw [blue](0.0000,2.3100)--(11.8800,2.3100);
\draw [teal](0.0000,2.9700)--(11.8800,2.9700);
\draw [blue](0.0000,4.2900)--(11.8800,4.2900);
\draw [teal](0.0000,5.6100)--(11.8800,5.6100);
\draw [black](-0.3300,1.6500) node [] {{\tiny $N_1$}};
\draw [black](-0.3300,2.9700) node [] {{\tiny $N_2$}};
\draw [black](-0.3300,5.6100) node [] {{\tiny $N_3$}};
\draw [black](-0.3300,0.0000) node [] {{\tiny $M_0$}};
\draw [black](-0.3300,0.6600) node [] {{\tiny $M_1$}};
\draw [black](-0.3300,2.3100) node [] {{\tiny $M_2$}};
\draw [black](-0.3300,4.2900) node [] {{\tiny $M_3$}};
\draw [black](0.6600,0.0000)--(0.9900,0.3300)--(1.3200,0.6600)--(1.6500,0.9900)--(1.9800,1.3200)--(2.3100,1.6500)--(2.6400,1.9800)--(2.9700,2.3100)--(3.3000,2.6400)--(3.6300,2.9700)--(3.9600,3.3000)--(4.2900,3.6300)--(4.6200,3.9600)--(4.9500,4.2900);
\fill (0.6600,0.0000) circle(0.0300); 
\fill (4.9500,4.2900) circle(0.0300); 
\fill (0.6600,0.0000) circle(0.0300); 
\fill (0.9900,0.3300) circle(0.0300); 
\fill (1.3200,0.6600) circle(0.0300); 
\fill (1.6500,0.9900) circle(0.0300); 
\fill (1.9800,1.3200) circle(0.0300); 
\fill (2.3100,1.6500) circle(0.0300); 
\fill (2.6400,1.9800) circle(0.0300); 
\fill (2.9700,2.3100) circle(0.0300); 
\fill (3.3000,2.6400) circle(0.0300); 
\fill (3.6300,2.9700) circle(0.0300); 
\fill (3.9600,3.3000) circle(0.0300); 
\fill (4.2900,3.6300) circle(0.0300); 
\fill (4.6200,3.9600) circle(0.0300); 
\draw [black](4.9500,4.2900)--(5.2800,3.9600)--(5.6100,3.6300)--(5.9400,3.3000)--(6.2700,2.9700)--(6.6000,2.6400)--(6.9300,2.3100)--(7.2600,1.9800)--(7.5900,1.6500)--(7.9200,1.3200)--(8.2500,0.9900)--(8.5800,0.6600)--(8.9100,0.3300)--(9.2400,0.0000);
\fill (4.9500,4.2900) circle(0.0300); 
\fill (9.2400,0.0000) circle(0.0300); 
\draw [black](1.3200,0.0000)--(1.6500,0.3300)--(1.9800,0.6600)--(2.3100,0.9900)--(2.6400,1.3200)--(2.9700,1.6500)--(3.3000,1.9800)--(3.6300,2.3100)--(3.9600,2.6400)--(4.2900,2.9700)--(4.6200,3.3000)--(4.9500,3.6300)--(5.2800,3.9600)--(5.6100,4.2900);
\fill (1.3200,0.0000) circle(0.0300); 
\fill (5.6100,4.2900) circle(0.0300); 
\fill (1.3200,0.0000) circle(0.0300); 
\fill (1.6500,0.3300) circle(0.0300); 
\fill (1.9800,0.6600) circle(0.0300); 
\fill (2.3100,0.9900) circle(0.0300); 
\fill (2.6400,1.3200) circle(0.0300); 
\fill (2.9700,1.6500) circle(0.0300); 
\fill (3.3000,1.9800) circle(0.0300); 
\fill (3.6300,2.3100) circle(0.0300); 
\fill (3.9600,2.6400) circle(0.0300); 
\fill (4.2900,2.9700) circle(0.0300); 
\fill (4.6200,3.3000) circle(0.0300); 
\fill (4.9500,3.6300) circle(0.0300); 
\fill (5.2800,3.9600) circle(0.0300); 
\draw [black](5.6100,4.2900)--(5.9400,3.9600)--(6.2700,3.6300)--(6.6000,3.3000)--(6.9300,2.9700)--(7.2600,2.6400)--(7.5900,2.3100)--(7.9200,1.9800)--(8.2500,1.6500)--(8.5800,1.3200)--(8.9100,0.9900)--(9.2400,0.6600)--(9.5700,0.3300)--(9.9000,0.0000);
\fill (5.6100,4.2900) circle(0.0300); 
\fill (9.9000,0.0000) circle(0.0300); 
\draw [black](1.9800,0.0000)--(2.3100,0.3300)--(2.6400,0.6600)--(2.9700,0.9900)--(3.3000,1.3200)--(3.6300,1.6500)--(3.9600,1.9800)--(4.2900,2.3100)--(4.6200,2.6400)--(4.9500,2.9700)--(5.2800,3.3000)--(5.6100,3.6300)--(5.9400,3.9600)--(6.2700,4.2900);
\fill (1.9800,0.0000) circle(0.0300); 
\fill (6.2700,4.2900) circle(0.0300); 
\fill (1.9800,0.0000) circle(0.0300); 
\fill (2.3100,0.3300) circle(0.0300); 
\fill (2.6400,0.6600) circle(0.0300); 
\fill (2.9700,0.9900) circle(0.0300); 
\fill (3.3000,1.3200) circle(0.0300); 
\fill (3.6300,1.6500) circle(0.0300); 
\fill (3.9600,1.9800) circle(0.0300); 
\fill (4.2900,2.3100) circle(0.0300); 
\fill (4.6200,2.6400) circle(0.0300); 
\fill (4.9500,2.9700) circle(0.0300); 
\fill (5.2800,3.3000) circle(0.0300); 
\fill (5.6100,3.6300) circle(0.0300); 
\fill (5.9400,3.9600) circle(0.0300); 
\draw [black](6.2700,4.2900)--(6.6000,3.9600)--(6.9300,3.6300)--(7.2600,3.3000)--(7.5900,2.9700)--(7.9200,2.6400)--(8.2500,2.3100)--(8.5800,1.9800)--(8.9100,1.6500)--(9.2400,1.3200)--(9.5700,0.9900)--(9.9000,0.6600)--(10.2300,0.3300)--(10.5600,0.0000);
\fill (6.2700,4.2900) circle(0.0300); 
\fill (10.5600,0.0000) circle(0.0300); 
\draw [black](2.6400,0.0000)--(2.9700,0.3300)--(3.3000,0.6600)--(3.6300,0.9900)--(3.9600,1.3200)--(4.2900,1.6500)--(4.6200,1.9800)--(4.9500,2.3100)--(5.2800,2.6400)--(5.6100,2.9700)--(5.9400,3.3000)--(6.2700,3.6300)--(6.6000,3.9600)--(6.9300,4.2900);
\fill (2.6400,0.0000) circle(0.0300); 
\fill (6.9300,4.2900) circle(0.0300); 
\fill (2.6400,0.0000) circle(0.0300); 
\fill (2.9700,0.3300) circle(0.0300); 
\fill (3.3000,0.6600) circle(0.0300); 
\fill (3.6300,0.9900) circle(0.0300); 
\fill (3.9600,1.3200) circle(0.0300); 
\fill (4.2900,1.6500) circle(0.0300); 
\fill (4.6200,1.9800) circle(0.0300); 
\fill (4.9500,2.3100) circle(0.0300); 
\fill (5.2800,2.6400) circle(0.0300); 
\fill (5.6100,2.9700) circle(0.0300); 
\fill (5.9400,3.3000) circle(0.0300); 
\fill (6.2700,3.6300) circle(0.0300); 
\fill (6.6000,3.9600) circle(0.0300); 
\draw [black](6.9300,4.2900)--(7.2600,3.9600)--(7.5900,3.6300)--(7.9200,3.3000)--(8.2500,2.9700)--(8.5800,2.6400)--(8.9100,2.3100)--(9.2400,1.9800)--(9.5700,1.6500)--(9.9000,1.3200)--(10.2300,0.9900)--(10.5600,0.6600)--(10.8900,0.3300)--(11.2200,0.0000);
\fill (6.9300,4.2900) circle(0.0300); 
\fill (11.2200,0.0000) circle(0.0300); 
\draw [black](3.3000,0.0000)--(3.6300,0.3300)--(3.9600,0.6600)--(4.2900,0.9900)--(4.6200,1.3200)--(4.9500,1.6500)--(5.2800,1.9800)--(5.6100,2.3100);
\fill (3.3000,0.0000) circle(0.0300); 
\fill (5.6100,2.3100) circle(0.0300); 
\fill (3.3000,0.0000) circle(0.0300); 
\fill (3.6300,0.3300) circle(0.0300); 
\fill (3.9600,0.6600) circle(0.0300); 
\fill (4.2900,0.9900) circle(0.0300); 
\fill (4.6200,1.3200) circle(0.0300); 
\fill (4.9500,1.6500) circle(0.0300); 
\fill (5.2800,1.9800) circle(0.0300); 
\draw [black](5.6100,2.3100)--(5.9400,1.9800)--(6.2700,1.6500)--(6.6000,1.3200)--(6.9300,0.9900)--(7.2600,0.6600)--(7.5900,0.3300)--(7.9200,0.0000);
\fill (5.6100,2.3100) circle(0.0300); 
\fill (7.9200,0.0000) circle(0.0300); 
\draw [black](3.9600,0.0000)--(4.2900,0.3300)--(4.6200,0.6600)--(4.9500,0.9900)--(5.2800,1.3200)--(5.6100,1.6500)--(5.9400,1.9800)--(6.2700,2.3100);
\fill (3.9600,0.0000) circle(0.0300); 
\fill (6.2700,2.3100) circle(0.0300); 
\fill (3.9600,0.0000) circle(0.0300); 
\fill (4.2900,0.3300) circle(0.0300); 
\fill (4.6200,0.6600) circle(0.0300); 
\fill (4.9500,0.9900) circle(0.0300); 
\fill (5.2800,1.3200) circle(0.0300); 
\fill (5.6100,1.6500) circle(0.0300); 
\fill (5.9400,1.9800) circle(0.0300); 
\draw [black](6.2700,2.3100)--(6.6000,1.9800)--(6.9300,1.6500)--(7.2600,1.3200)--(7.5900,0.9900)--(7.9200,0.6600)--(8.2500,0.3300)--(8.5800,0.0000);
\fill (6.2700,2.3100) circle(0.0300); 
\fill (8.5800,0.0000) circle(0.0300); 
\draw [black](4.6200,0.0000)--(4.9500,0.3300)--(5.2800,0.6600);
\fill (4.6200,0.0000) circle(0.0300); 
\fill (5.2800,0.6600) circle(0.0300); 
\fill (4.6200,0.0000) circle(0.0300); 
\fill (4.9500,0.3300) circle(0.0300); 
\draw [black](5.2800,0.6600)--(5.6100,0.3300)--(5.9400,0.0000);
\fill (5.2800,0.6600) circle(0.0300); 
\fill (5.9400,0.0000) circle(0.0300); 
\draw [black](5.2800,0.0000)--(5.6100,0.3300)--(5.9400,0.6600);
\fill (5.2800,0.0000) circle(0.0300); 
\fill (5.9400,0.6600) circle(0.0300); 
\fill (5.2800,0.0000) circle(0.0300); 
\fill (5.6100,0.3300) circle(0.0300); 
\draw [black](5.9400,0.6600)--(6.2700,0.3300)--(6.6000,0.0000);
\fill (5.9400,0.6600) circle(0.0300); 
\fill (6.6000,0.0000) circle(0.0300); 
\draw [black](5.9400,0.0000)--(6.2700,0.3300)--(6.6000,0.6600);
\fill (5.9400,0.0000) circle(0.0300); 
\fill (6.6000,0.6600) circle(0.0300); 
\fill (5.9400,0.0000) circle(0.0300); 
\fill (6.2700,0.3300) circle(0.0300); 
\draw [black](6.6000,0.6600)--(6.9300,0.3300)--(7.2600,0.0000);
\fill (6.6000,0.6600) circle(0.0300); 
\fill (7.2600,0.0000) circle(0.0300); 
\draw [black](8.5800,3.2010) node [] {$v_2$};
\draw [black](10.2300,1.4520) node [] {$v_1$};
\draw [black](11.4180,0.3300) node [] {$v_0$};
\draw [black](0.6600,-0.3300) node [] {{\tiny $A_8$}};
\draw [black](1.3200,-0.3300) node [] {{\tiny $A_7$}};
\draw [black](1.9800,-0.3300) node [] {{\tiny $A_6$}};
\draw [black](2.6400,-0.3300) node [] {{\tiny $A_5$}};
\draw [black](3.3000,-0.3300) node [] {{\tiny $A_4$}};
\draw [black](3.9600,-0.3300) node [] {{\tiny $A_3$}};
\draw [black](4.6200,-0.3300) node [] {{\tiny $A_2$}};
\draw [black](5.2800,-0.3300) node [] {{\tiny $A_1$}};
\draw [black](5.9400,-0.3300) node [] {{\tiny $A_0B_0$}};
\draw [black](6.6000,-0.3300) node [] {{\tiny $B_1$}};
\draw [black](7.2600,-0.3300) node [] {{\tiny $B_2$}};
\draw [black](7.9200,-0.3300) node [] {{\tiny $B_3$}};
\draw [black](8.5800,-0.3300) node [] {{\tiny $B_4$}};
\draw [black](9.2400,-0.3300) node [] {{\tiny $B_5$}};
\draw [black](9.9000,-0.3300) node [] {{\tiny $B_6$}};
\draw [black](10.5600,-0.3300) node [] {{\tiny $B_7$}};
\draw [black](11.2200,-0.3300) node [] {{\tiny $B_8$}};
\fill (8.9100,2.3100) circle(0.0300); 
\fill (8.2500,2.3100) circle(0.0300); 
\fill (7.5900,2.3100) circle(0.0300); 
\fill (6.9300,2.3100) circle(0.0300); 
\fill (10.5600,0.6600) circle(0.0300); 
\fill (9.9000,0.6600) circle(0.0300); 
\fill (9.2400,0.6600) circle(0.0300); 
\fill (8.5800,0.6600) circle(0.0300); 
\fill (7.9200,0.6600) circle(0.0300); 
\fill (7.2600,0.6600) circle(0.0300); 
\end{tikzpicture}
\caption{The non-intersecting paths from $A_0,A_1,\dots,A_8$ to $B_0,B_1,\dots,B_8$ }\label{fig:h(0)}
\end{figure}
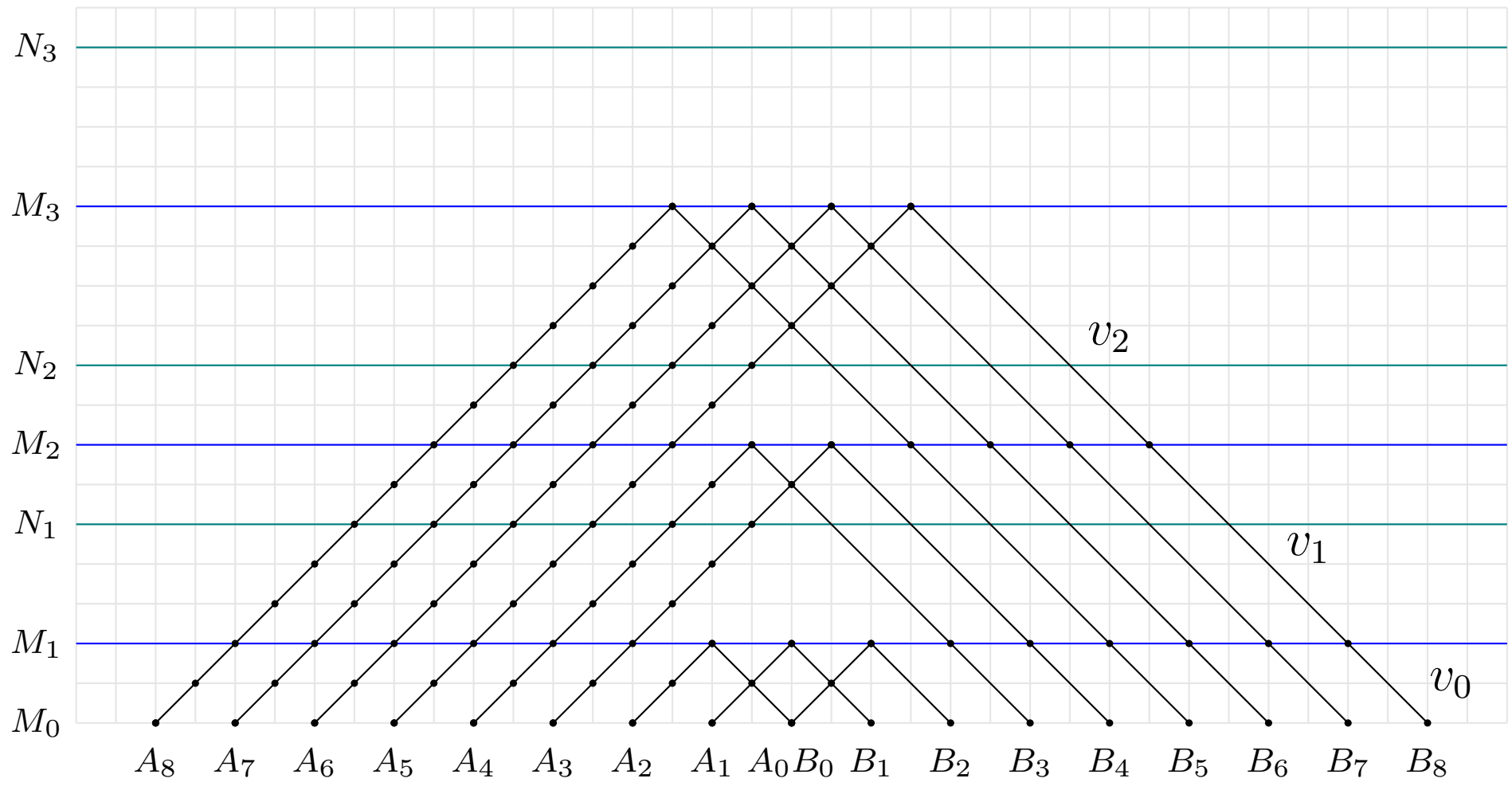

For example, let $(k_0,k_1,k_2)=(2,1,3)$ and $(s_1,s_2,s_3)=(3,5,9)$. 
Then $(M_0,M_1,M_2,M_3)= (0,2,7,13)$.
As illustrated in Figure~\ref{fig:h(0)}, the unique family of non-intersecting paths from the sources $A_0,A_1,\dots,A_8$ to the sinks $B_0,B_1,\dots,B_8$ is described as follows:
\begin{itemize}
    \item[(1)]  For $i=0,1,2$, $A_i$ is connected to $B_{2-i}$, and $P_i= (A_i,B_{2-i})=b^{2}d$.
    \item[(2)]  For $i=3,4$, writing $i=3+\ell$ with $\ell=0,1$, we have $A_i$ connected to $B_{4-\ell}$, and $P_i=(A_i, B_{4-\ell})=b^{7}v d$.
    \item[(3)]  For $i=5,6,7,8$, writing $i=5+\ell$ with $\ell=0,1,2,3$, we have $A_i$ connected to $B_{8-\ell}$, and $P_i= (A_i,B_{8-\ell})=b^{13}v^{2}d.$
\end{itemize}
Hence, the total weight is 
\[
\mathcal{Q}_{0,8} =  v_0^3 \cdot (v_0v_1)^2 \cdot (v_0v_1v_2)^4 = v_0^9 v_1^6 v_2^4.
\]


\section{A combinatorial proof of Theorem \ref{th:shiftHDet}}\label{sec:shift}

In this section, we study the $1$-shifted Hankel determinants of $F(x)$ and give a combinatorial proof of Theorem~\ref{th:shiftHDet}.
Building on the quasi-Dyck path model and the combinatorial proof of Theorem~\ref{th:super2}, we determine the structure of the relevant non-intersecting path families.
Evaluating their signed weights gives the formulas for $H_{s_j}^{(1)}(F)$ and $H_{s_j+k_j}^{(1)}(F)$, 
while the absence of such families at the remaining indices proves the vanishing assertions.
An algebraic proof of Theorem~\ref{th:shiftHDet} is given in Section~\ref{sec:shiftalg}.

\subsection{The path model for $1$-shifted Hankel determinants}
\label{subsec:shift:paths}

Let $n\ge1$.
For $0\le i\le n-1$, define the sources and sinks in the graph $\Gamma$ of Section~\ref{sec:lgv} by
\[
A_i:=(-2i,0)_{\mathrm{src}},
\qquad
B_i:=(2(i+1),0)_{\mathrm{snk}}
\]
Let $\Omega_n^{(1)}$ denote the set of all non-intersecting families $(P_0,\ldots,P_{n-1})$ of quasi-Dyck paths in $\Gamma$ such that, for some permutation $\sigma\in\mathfrak{S}_{[0,n-1]}$ and every $0\le i\le n-1$,
\[
P_i=(A_i,B_{\sigma(i)}).
\]

For $i,t\ge0$, a quasi-Dyck path from $A_i$ to $B_t$ has $i+t+1$ up steps.
The path interpretation of the coefficients of $F$ therefore gives
\[
a_{i+t+1}
=
\sum_{P:\,A_i\to B_t}w(P),
\]
where the sum runs over all paths of $\Gamma$ from $A_i$ to $B_t$.
For $n\ge1$, the corresponding matrix of path-weight sums is the $1$-shifted Hankel matrix of $F$,
whose determinant is
\[
H_n^{(1)}(F)
=
\det[a_{i+t+1}]_{0\le i,t\le n-1}.
\]
Since $B_i$ is the sink denoted $B_{i+1}$ in Section~\ref{sec:lgv}, and the proof of Theorem~\ref{th:lGV} applies verbatim to any sources and sinks of the acyclic graph $\Gamma$, we can express $H_n^{(1)}(F)$ as a signed sum over non-intersecting path families.

Applying the same argument as in Theorem~\ref{th:lGV} to these sources and sinks gives the following formula for $1$-shifted Hankel determinants.

  \begin{Theorem}\label{th:hanke(1)}
      For every $n \ge 1$, 
      \[
      H_{n}^{(1)}(F) = \sum_{(P_0, P_1, \dots, P_{n-1}) \in \Omega_{n}^{(1)}} (-1)^{\mathrm{inv}(\sigma)} \prod_{i=0}^{n-1} w(P_i),
      \]
     where the sum ranges over all non-intersecting families of paths $P_i$ going from $A_i$ to $B_{\sigma(i)}$ for $0 \le i \le n-1$, $\sigma$ being some permutation of $\mathfrak{S}_{[0,n-1]}$.
  \end{Theorem}

As in Section~\ref{sec:lgv}, we write $\sigma=\sigma_1\cdots\sigma_n$ in one-line notation, so that $\sigma(i)=\sigma_{i+1}$ for $0\le i\le n-1$.
Auxiliary entrance edges, when present, are included in the paths but omitted from their step words.

For $i\ge1$, an initial climb $b^{M_i}$ of a path from $A_t$ ends at the ordinary vertex $(M_i-2t,M_i)$, while a final descent $v^{i-1}d$ to $B_x$ starts at $(2x+2-M_i,M_i)$.
These two vertices coincide if and only if $t+x=M_i-1.$
A path from $A_t$ to $B_x$ has $t+x+1$ up steps.
Thus every height $h$ reached by the path satisfies $h \le t+x+1$.
Moreover, any vertex at height $h$ has abscissa at least $-2t+h$.

The following lemma is an analogue of Lemma~\ref{lem:forced-segments}.
Its proof is similar and is omitted.

\begin{Lemma}\label{lem:forced-segments1}
Let $n\ge1$, and let $(P_0,\ldots,P_{n-1})\in\Omega_n^{(1)}$ have associated permutation $\sigma$.
For each $i\ge1$ and $0\le t\le n-1$, the following hold:
\begin{enumerate}
\item[(i)]
If $t\ge s_{i-1}$, then $P_t$ begins with $b^{M_i}$.
\item[(ii)]
If $\sigma_{t+1}\ge s_{i-1}$, then $P_t$ ends with $v^{i-1}d$.
\end{enumerate}
\end{Lemma}

We next use these forced segments to determine the block structure and the boundary paths.

\subsection{Block structure and boundary paths}
\label{subsec:shift:structure}

Let $(P_0,\ldots,P_{n-1})\in\Omega_n^{(1)}$ have associated permutation $\sigma$.
For each $i\ge1$ with $s_i-1\le n$, define the consecutive block
\[
\sigma^i
=
\sigma_{s_{i-1}+1}\sigma_{s_{i-1}+2}\cdots
\sigma_{s_i-1},
\]
which has $k_{i-1}$ entries and is empty when $k_{i-1}=0$.
In particular, when $n=s_j$, we have the decomposition
\[
\sigma
=
\sigma^1\,\sigma_{s_1}\,
\sigma^2\,\sigma_{s_2}\cdots
\sigma^j\,\sigma_{s_j}.
\]
The following lemma determines the paths corresponding to the entries of each block $\sigma^i$.
For $n=s_j$, this leaves only the paths $P_{s_1-1},\ldots,P_{s_j-1}$ to be determined later.

\begin{Lemma}\label{lem:psigma1}    
Let $n\ge1$, and let $(P_0,\ldots,P_{n-1})\in\Omega_n^{(1)}$ have associated permutation $\sigma$.
For each $i\ge1$ with $k_{i-1}\ge1$ and $s_i-2\le n-1$, we have
\[
\sigma^i = (s_i - 2) \; (s_i - 3) \; \cdots \; s_{i-1},
\]
that is,  for $0 \le \ell \le k_{i-1}-1$,
\[
P_{s_{i-1} + \ell} = (A_{s_{i-1} + \ell},\, B_{s_i - 2 - \ell}) = b^{M_i} v^{i-1} d.
\]
Moreover, for any $i\ge1$, if $s_{i-1}\le n-1<s_i-2$, then $\Omega^{(1)}_n$ is empty.
\end{Lemma}

\begin{proof}
Fix $i\ge1$ with $k_{i-1}\ge1$ and $s_i-2\le n-1$.
Let $0\le\ell\le k_{i-1}-1$, and set $x:=\sigma_{s_{i-1}+1+\ell}$, so that $P_{s_{i-1}+\ell}=(A_{s_{i-1}+\ell},B_x).$
Since $s_{i-1}+\ell\ge s_{i-1}$, Lemma~\ref{lem:forced-segments1}(i) shows that this path begins with $b^{M_i}$.
In particular, it has at least $M_i$ up steps.
Using $M_i=s_{i-1}+s_i-1$, we see that its initial climb ends at
\[
(M_i-2(s_{i-1}+\ell),M_i)
=
(2(s_i-2-\ell)+2-M_i,M_i).
\]
We show that $x=s_i-2-\ell$ by excluding the following two cases.

\medskip
\noindent\textbf{Case 1.}
Suppose that $x<s_i-2-\ell$.
Then $P_{s_{i-1}+\ell}$ has
\[
s_{i-1}+\ell+x+1 < s_{i-1}+s_i-1 = M_i
\]
up steps, contradicting its forced initial climb.

\medskip
\noindent\textbf{Case 2.}
Suppose that $x>s_i-2-\ell$.
Since $\sigma$ is a permutation and $s_{i-1}\le s_i-2-\ell\le n-1,$ there exists $y\ne s_{i-1}+\ell$ such that $P_y=(A_y,B_{s_i-2-\ell}).$
By Lemma~\ref{lem:forced-segments1}(ii), $P_y$ ends with $v^{i-1}d$, whose starting vertex is $(2(s_i-2-\ell)+2-M_i,M_i).$
This is the endpoint of the initial climb of $P_{s_{i-1}+\ell}$.
Thus two distinct paths share an ordinary vertex, contradicting non-intersection.

\medskip

  Hence, we have  $\sigma_{s_{i-1}+1+\ell} = s_i - 2-\ell$ for $\ell \in [0,k_{i-1}-1]$. 
  The path $P_{s_{i-1}+\ell}=(A_{s_{i-1} + \ell},\, B_{s_i - 2 - \ell})$ has $M_i$ up steps. 
  It begins with $b^{M_i}$ and, by Lemma~\ref{lem:forced-segments1}(ii), ends with $v^{i-1}d$. 
  These two segments account for its entire horizontal length $2M_i$, so that
  \[
  P_{s_{i-1} + \ell} = (A_{s_{i-1} + \ell},\, B_{s_i - 2 - \ell}) = b^{M_i} v^{i-1} d
  \]
  for $0 \le \ell \le k_{i-1}-1$.

Finally, suppose that $s_{i-1}\le n-1<s_i-2$ for some $i\ge1$ and, for a contradiction, that
$\Omega_n^{(1)}$ is nonempty.
Let $P_t$ be the path ending at $B_{s_{i-1}}$.
By Lemma~\ref{lem:forced-segments1}(ii), it ends with $v^{i-1}d$ and therefore reaches height $M_i$.
Since it has $t+s_{i-1}+1$ up steps, we obtain $t+s_{i-1}+1\ge M_i.$
Thus
\[
t\ge M_i-s_{i-1}-1=s_i-2>n-1,
\]
contradicting $0\le t\le n-1$.
Hence $\Omega_n^{(1)}$ is empty.
    
\end{proof}

As shown in Lemma~\ref{lem:psigma1}, when $n=s_j$, the paths $P_{s_{i-1}+\ell}$ for $ 1 \le i \le j $ and $0 \le \ell \le k_{i-1}-1$ are completely determined. 
It remains to describe the boundary case $P_{s_i-1}$. We first determine the path ending at the sink $B_{s_i-1}$.

\begin{Lemma}\label{lem:free1}
Let $n\ge1$, and let $(P_0,\ldots,P_{n-1})\in\Omega_n^{(1)}$.
Let $i\ge1$ with $s_i\le n$.
Then the path ending at $B_{s_i-1}$ is described by exactly one of the following cases:
\begin{enumerate}
\item[(i)]
The path is
\[
P_{s_i-1}=b^{N_i}c_iv^{i-1}d.
\]

\item[(ii)]
For $i\ge2$, the path is
\[
P_{s_{i-1}-1}=b^{M_i}v^{i-1}d.
\]

\item[(iii)]
If $s_{i+1}\le n$, then the path is
\[
P_{s_{i+1}-1}=b^{M_{i+1}}v^id.
\]
\end{enumerate}
Here $c_i=c_{i,k_{i-1}}$ denotes the short down step of drop $k_{i-1}+1$ ending at height $M_i$.

\end{Lemma}

\begin{proof}
Let $P_t$ be the path ending at $B_{s_i-1}$.
By Lemma~\ref{lem:forced-segments1}(ii), it ends with $v^{i-1}d$, whose starting vertex is
\[
D=(2s_i-M_i,M_i)=(k_{i-1}+2,M_i).
\]
Since $P_t$ has $t+s_i$ up steps and reaches height $M_i$, we have 
\[
t\ge M_i-s_i=s_{i-1}-1.
\]
When $k_{i-1}\ge1$, Lemma~\ref{lem:psigma1} shows that the paths with source indices in $[s_{i-1},s_i-2]$ end at sinks with indices in the same interval.
Thus none of these paths ends at $B_{s_i-1}$. 
When $k_{i-1}=0$, this interval is empty.
We therefore distinguish the following three cases.

\medskip

\noindent\textbf{Case 1.}
Suppose that $t=s_{i-1}-1$, necessarily with $i\ge2$.
Then $P_t$ has exactly $t+s_i=M_i$ up steps.
Since it reaches height $M_i$ at $D$, all these up steps occur before $D$, with no intervening down step.
Thus, $ P_t=b^{M_i}v^{i-1}d.$

\medskip

\noindent\textbf{Case 2.}
Suppose that $t=s_i-1$.
Then $P_t$ has $2s_i-1=N_i$ up steps.
By Lemma~\ref{lem:forced-segments1}(i), it begins with $b^{M_i}$, whose endpoint is
\[
C=(M_i-2(s_i-1),M_i)=(-k_{i-1},M_i).
\]
The portion from $C$ to $D$ has $k_{i-1}+1$ up steps and the same total drop.
This portion cannot go below $M_i$.
For $i\ge2$, a long down step from $M_i$ to $M_{i-1}$ would require at least $M_i-M_{i-1}>k_{i-1}+1$ subsequent up steps to return to $M_i$.
For $i=1$, the terminal step would end the path before $D$.
Moreover, this portion cannot reach $M_{i+1}$, since $M_{i+1}-M_i>k_{i-1}+1$.
By Lemma~\ref{lem:primitive}, its decomposition into primitive paths consists only of factors
$b^r c_{i,r-1}$ with $1\le r\le k_{i-1}+1$.
If the first factor has $r\le k_{i-1}$, it ends at
\[
(-k_{i-1}+2r,M_i)
=
(M_i-2(s_i-1-r),M_i).
\]
Since $s_{i-1}\le s_i-1-r\le s_i-2,$ Lemma~\ref{lem:forced-segments1}(i) shows that this vertex is also the endpoint of the initial climb of $P_{s_i-1-r}$.
This contradicts non-intersection.
Hence $r=k_{i-1}+1$, which uses all the available up steps.
Consequently, $P_t=b^{N_i}c_i v^{i-1}d.$

\medskip
\noindent\textbf{Case 3.}
Suppose that $t\ge s_i$.
By Lemma~\ref{lem:forced-segments1}(i), $P_t$ begins with $b^{M_{i+1}}$.
Thus $t+s_i\ge M_{i+1}$, giving $t\ge s_{i+1}-1.$
In particular, $P_{s_{i+1}-1}$ belongs to the family.
Consider the step entering $D$ immediately before the final descent $v^{i-1}d$.
We first exclude the following two possibilities.

\medskip
\noindent\textbf{Subcase 3.1.}
Suppose that this step is a short down step of drop $r$, where $1\le r\le k_{i-1}+1$.
Since no down step ends strictly between $M_i$ and $M_{i+1}$, it must be preceded by $r$ consecutive up steps starting at
\[
Z = (2s_i-M_i-2r,M_i) = (M_i-2(s_{i-1}+r-1),M_i).
\]
By Lemma~\ref{lem:forced-segments1}(i), this is also the endpoint of the initial climb of $P_{s_{i-1}+r-1}$.
Since
\[
s_{i-1}\le s_{i-1}+r-1\le s_i-1<t,
\]
the two paths are distinct, contradicting non-intersection.

\medskip
\noindent\textbf{Subcase 3.2.}
Suppose that this step is an up step.
If $i\ge2$, it belongs to a consecutive climb passing through height $M_{i-1}$, since no down step ends strictly between $M_{i-1}$ and $M_i$.
The vertex of this climb at height $M_{i-1}$ is $(M_{i-1}-2(s_{i-1}-1),M_{i-1}).$
Since $s_{i-1}-1\ge s_{i-2}$, Lemma~\ref{lem:forced-segments1}(i) shows that this is also the endpoint of the initial climb of $P_{s_{i-1}-1}$.
As $s_{i-1}-1<t$, this contradicts non-intersection.
If $i=1$, the allowed steps force the path to have step word $b^{M_1}d$.
Its endpoints would then satisfy $t+s_1=M_1$, giving $t=-1$, which is impossible.

\medskip
The same endpoint calculation excludes the word $d$ alone when $k_0=0$.
The step entering $D$ must therefore be the long down step from $M_{i+1}$ to $M_i$.
Hence $P_t$ ends with $v^id$, whose starting vertex is
\[
D' = (2s_i-M_{i+1},M_{i+1}) = (M_{i+1}-2(s_{i+1}-1),M_{i+1}).
\]
This is also the endpoint of the forced initial climb of $P_{s_{i+1}-1}$.
If $t>s_{i+1}-1$, the two distinct paths would share the vertex $D'$, contradicting non-intersection.
Therefore $t=s_{i+1}-1$.
The path has exactly $M_{i+1}$ up steps, all used in its initial climb, so $P_t=b^{M_{i+1}}v^id.$
This completes the proof.
\end{proof}

We now describe the boundary paths $P_{s_i-1}$ in the following lemma.

\begin{Lemma}\label{lem:psi1}
Let $j\ge1$, and let $(P_0,\ldots,P_{s_j-1})\in\Omega_{s_j}^{(1)}$ have associated permutation $\sigma$.
For each $1\le i\le j$, exactly one of the following cases occurs:
\begin{itemize}
\item[(i)]
We have $\sigma_{s_i}=s_i-1$, and
\[
P_{s_i-1} =(A_{s_i-1},B_{s_i-1}) =b^{N_i}c_i v^{i-1}d.
\]

\item[(ii)]
For $i\ge2$, we have $\sigma_{s_i}=s_{i-1}-1,$ $ \sigma_{s_{i-1}}=s_i-1,$ and
\[
\begin{aligned}
P_{s_{i-1}-1}
&=(A_{s_{i-1}-1},B_{s_i-1})
=b^{M_i}v^{i-1}d,\\
P_{s_i-1}
&=(A_{s_i-1},B_{s_{i-1}-1})
=b^{M_i}v^{i-1}d.
\end{aligned}
\]

\item[(iii)]
For $i\le j-1$, we have $\sigma_{s_i}=s_{i+1}-1,$ $\sigma_{s_{i+1}}=s_i-1.$ This is case (ii) with $i$ replaced by $i+1$.
\end{itemize}
 Conversely, every collection of disjoint adjacent transpositions of $(s_1-1,\ldots,s_j-1)$ determines a unique family in $\Omega_{s_j}^{(1)}$, with the remaining paths given by Lemma~\ref{lem:psigma1}.
\end{Lemma}

\begin{proof}
By Lemma~\ref{lem:psigma1}, the entries $\sigma_{s_1},\ldots,\sigma_{s_j}$ form a permutation of $s_1-1,\ldots,s_j-1$.
By Lemma~\ref{lem:free1}, a connection from $A_{s_m-1}$ to $B_{s_i-1}$ is possible only if
$|m-i|\le1$.
Thus each boundary source can be connected only to its corresponding boundary sink or to a neighboring boundary sink.
If not all boundary entries are fixed, let $i$ be the smallest index for which $\sigma_{s_i}\ne s_i-1$.
Since all earlier boundary positions are fixed, $A_{s_i-1}$ must be connected to $B_{s_{i+1}-1}$.
The only available source for $B_{s_i-1}$ is then $A_{s_{i+1}-1}$.
Thus, for $i<j$, we have $\sigma_{s_i}=s_{i+1}-1$, and $\sigma_{s_{i+1}}=s_i-1$.
Removing this pair and repeating the argument for the remaining indices gives disjoint adjacent transpositions.
The corresponding path forms are determined by Lemma~\ref{lem:free1}, proving (i)--(iii).

Conversely, choose any collection of disjoint adjacent transpositions of $(s_1-1,\ldots,s_j-1)$.
Construct the boundary paths using (i)--(iii) and the interior paths using Lemma~\ref{lem:psigma1}.
All prescribed steps are allowed.
For an interior path in block $i$ or a path in the transposed pair $(s_{i-1}-1,s_i-1)$, its source and sink indices $t,x$ satisfy
\[
t+x+1=s_{i-1}+s_i-1=M_i.
\]
For a fixed boundary path $P_{s_i-1}$, we have
$t=x=s_i-1$, so
\[
t+x+1=2s_i-1=N_i.
\]
Thus all constructed paths have the required endpoints.
Each step word and its source determine a unique path.

It remains to verify non-intersection.
Every ordinary vertex lies on an initial climb or a final descent.
The endpoints of a short down step lie on these two segments, and its interior contains no path vertices.
Initial climbs from $A_t$ lie on the lines $u-h=-2t,$ while final descents to $B_x$ lie on the lines $u+h=2(x+1).$
Thus distinct climbs are disjoint, as are distinct final descents.

Suppose that two distinct paths share a vertex.
Since sources and sinks are distinct graph vertices, the common vertex lies on a climb from $A_t$ and a final descent to $B_x$.
Its height is $M_m$ for some $1\le m\le j$, so
\[
M_m-2t=2(x+1)-M_m.
\]
Consequently,
\[
t+x=M_m-1=s_{m-1}+s_m-2.
\]

For $m\ge2$, the prescribed path forms give, for every ordinary vertex $(u,h)$ of a constructed path from $A_a$ to $B_b$,
\[
\min\{a,b\}<s_{m-1}-1
\quad\Longrightarrow\quad
h\le M_{m-1}<M_m.
\]
Applying this to the two paths through the common vertex yields $t,x\ge s_{m-1}-1.$
For $m=1$, the same bound follows from $t,x\ge0$.
Together with the formula for $t+x$, this gives
\[
s_{m-1}-1\le t,x\le s_m-1.
\]
We distinguish two cases.

\medskip

\noindent\textbf{Case 1.}
One of $t$ and $x$ equals $s_{m-1}-1$.
Then the other equals $s_m-1$, and necessarily $m\ge2$.
Since the path involving the lower boundary index reaches height $M_m$, the transposition in case (ii) for the index $m$ must have been chosen.
Its two paths connect $A_{s_{m-1}-1}$ to $B_{s_m-1}$ and $A_{s_m-1}$ to $B_{s_{m-1}-1}$.
Hence the construction connects $A_t$ to $B_x$.

\medskip

\noindent\textbf{Case 2.}
Both $t$ and $x$ are greater than $s_{m-1}-1$.
The equality $t+x=s_{m-1}+s_m-2$ then implies $s_{m-1}\le t,x\le s_m-2.$
These are interior indices, for which the construction in Lemma~\ref{lem:psigma1} gives 
\[
P_t=(A_t,B_{M_m-1-t})=(A_t,B_x).
\]

In both cases, the initial climb from $A_t$ and the final descent to $B_x$ belong to the same constructed path, contradicting the assumption that the common vertex belongs to two distinct paths.
The constructed family is therefore non-intersecting.
This proves the converse and completes the proof.
\end{proof}

\subsection{Proof of Theorem~\ref{th:shiftHDet}}
We first compute the signed weight sum of the boundary paths.
For each family $(P_0,\ldots,P_{s_j-1})\in\Omega_{s_j}^{(1)}$ with associated permutation $\sigma$, let
\[
\tau=(\sigma_{s_1},\ldots,\sigma_{s_j}).
\]
We define the signed weight sum
\[
\mathcal{W}_j
:=
\sum_{(P_0,\ldots,P_{s_j-1})\in\Omega_{s_j}^{(1)}}
(-1)^{\mathrm{inv}(\tau)}
\prod_{i=1}^{j}w(P_{s_i-1}).
\]
We also take $\Omega_0^{(1)}$ to consist of the empty family and set $\mathcal{W}_0=1$.

\begin{Lemma}\label{lem:Wj1}
For each $j\ge1$,
\[
\mathcal{W}_j
=
v_0^jv_1^{j-1}\cdots v_{j-1}R_j,
\]
where $(R_n)_{n\ge0}$ is defined by $R_0=1,$ and $R_1=c_{1,k_0},$ and for $n\ge1$,
\[
R_{n+1} = c_{n+1,k_n}R_n-v_nR_{n-1}.
\]
\end{Lemma}

\begin{proof}
Recall that $V_h=v_0v_1\cdots v_h$.
For $j=1$, Lemma~\ref{lem:psi1} allows only a fixed boundary entry, so  $\mathcal{W}_1=V_0c_{1,k_0}=v_0R_1.$
Let $j\ge2$.
We denote by $\mathcal{W}_j^{\mathrm{fix}}$ the contribution to $\mathcal{W}_j$ from families satisfying $\sigma_{s_j}=s_j-1$.
We denote by $\mathcal{W}_j^{\mathrm{swap}}$ the contribution from families satisfying
$\sigma_{s_j}=s_{j-1}-1$ and $\sigma_{s_{j-1}}=s_j-1$.
By Lemma~\ref{lem:psi1}, these are the only possibilities, and hence
\[
\mathcal{W}_j
=
\mathcal{W}_j^{\mathrm{fix}}
+
\mathcal{W}_j^{\mathrm{swap}}.
\]
We evaluate these two contributions separately.

\medskip

\noindent\textbf{Case 1.}
Suppose that $\sigma_{s_j}=s_j-1$.
By Lemma~\ref{lem:psi1}, deleting the last block gives a bijection with $\Omega_{s_{j-1}}^{(1)}$.
The deleted boundary path is $P_{s_j-1}=b^{N_j}c_jv^{j-1}d,$ whose weight is $V_{j-1}c_{j,k_{j-1}}$.
Since the last entry of $\tau$ is its largest entry, deleting it does not change the inversion number.
Thus
\begin{equation}\label{eq:Wj:fix}
\mathcal{W}_j^{\mathrm{fix}}
=
V_{j-1}c_{j,k_{j-1}}\mathcal{W}_{j-1}.
\end{equation}

\medskip

\noindent\textbf{Case 2.}
Suppose that $\sigma_{s_j}=s_{j-1}-1$ and $\sigma_{s_{j-1}}=s_j-1.$
By Lemma~\ref{lem:psi1}, deleting the last two blocks gives a bijection with $\Omega_{s_{j-2}}^{(1)}$.
The two deleted boundary paths both have step word $b^{M_j}v^{j-1}d$ and weight $V_{j-1}$.
The last two entries of $\tau$ are its two largest entries in reverse order, so they contribute exactly one additional inversion.
Hence
\begin{equation}\label{eq:Wj:swap}
\mathcal{W}_j^{\mathrm{swap}}
=
-V_{j-1}^{\,2}\mathcal{W}_{j-2}.
\end{equation}

Combining \eqref{eq:Wj:fix} and \eqref{eq:Wj:swap},
we obtain
\[
\mathcal{W}_j =\mathcal{W}_j^{\mathrm{fix}}
+
\mathcal{W}_j^{\mathrm{swap}}
=
V_{j-1}c_{j,k_{j-1}}\mathcal{W}_{j-1}
-
V_{j-1}^{\,2}\mathcal{W}_{j-2}.
\]
Together with $\mathcal{W}_0=1$ and $\mathcal{W}_1=v_0R_1$, this proves the desired formula by induction.
Indeed, using $V_{j-1}=V_{j-2}v_{j-1}$, the induction hypothesis gives
\begin{align*}
\mathcal{W}_j
&=
\left(\prod_{h=0}^{j-1}V_h\right)
\left(
c_{j,k_{j-1}}R_{j-1}-v_{j-1}R_{j-2}
\right)\\
&=\left(\prod_{h=0}^{j-1}V_h\right)R_j\\
&=
v_0^jv_1^{j-1}\cdots v_{j-1}R_j.
\end{align*}
This completes the proof.
\end{proof}

Now we are ready to prove Theorem \ref{th:shiftHDet}.

\begin{proof}[Proof of Theorem~\ref{th:shiftHDet}]
\noindent\textit{Proof of \eqref{eq:H1sj}.}
Fix $j\ge1$.
Lemmas~\ref{lem:psigma1} and~\ref{lem:psi1} describe all families in $\Omega_{s_j}^{(1)}$.
For each $1\le i\le j$, Lemma~\ref{lem:psigma1} gives
\begin{equation}\label{eq:invsi}
\mathrm{inv}(\sigma^i)
=
\frac{k_{i-1}(k_{i-1}-1)}{2}
\qquad \text{and} \qquad
\mathcal{Q}_{s_{i-1},s_i-2}
=
V_{i-1}^{k_{i-1}}.
\end{equation}
When $k_{i-1}=0$, the block $\sigma^i$ is empty, its inversion number is $0$, and the corresponding weight product is $1$.

Let $\tau=(\sigma_{s_1},\ldots,\sigma_{s_j})$.
Entries in different interior blocks form no inversions.
By Lemma~\ref{lem:psi1}, a fixed boundary entry $\sigma_{s_i}=s_i-1$ forms no inversions with the interior entries.
For an exchanged pair $ \sigma_{s_{i-1}}=s_i-1$ and $\sigma_{s_i}=s_{i-1}-1, $ the intervening entries are precisely those of $\sigma^i$.
Each of these $k_{i-1}$ entries forms two inversions with the exchanged pair, giving $2k_{i-1}$ inversions in total.
The other interior entries form no inversions with the exchanged pair, since those before position $s_{i-1}$ are at most $s_{i-1}-2$ and those after position $s_i$ are at least $s_i$.
The inversion between the two boundary entries is counted in $\mathrm{inv}(\tau)$.
Consequently,
\[
(-1)^{\mathrm{inv}(\sigma)}
=
(-1)^{\mathrm{inv}(\tau)}
\prod_{i=1}^{j}(-1)^{\mathrm{inv}(\sigma^i)}.
\]
Thus Theorem~\ref{th:hanke(1)}, the definition of $\mathcal{W}_j$, \eqref{eq:invsi} and Lemma~\ref{lem:Wj1} yield
\begin{align*}
H_{s_j}^{(1)}(F)
&=
\sum_{(P_0,\ldots,P_{s_j-1})\in\Omega_{s_j}^{(1)}}
(-1)^{\mathrm{inv}(\sigma)}
\prod_{t=0}^{s_j-1}w(P_t)\\
&=
\left(
\prod_{i=1}^{j}
(-1)^{\mathrm{inv}(\sigma^i)}
\mathcal{Q}_{s_{i-1},s_i-2}
\right)\mathcal{W}_j\\
&=
(-1)^{\sum_{i=0}^{j-1}\frac{k_i(k_i-1)}{2}}
\left(\prod_{i=0}^{j-1}V_i^{k_i+1}\right)R_j\\
&=
(-1)^{\epsilon_j-s_j+j}
v_0^{s_j}v_1^{s_j-s_1}\cdots
v_{j-1}^{s_j-s_{j-1}}R_j.
\end{align*}
Here $\epsilon_j=\sum_{i=0}^{j-1}\frac{k_i(k_i+1)}{2}$. This proves \eqref{eq:H1sj}.

\medskip

\noindent\textit{Proof of \eqref{eq:H1sjkj}.}
If $k_j=0$, identity \eqref{eq:H1sjkj} is immediate.
Now assume that $k_j\ge1$, and set $n=s_j+k_j=s_{j+1}-1$.
For any family in $\Omega_n^{(1)}$ and each
$0\le\ell\le k_j-1$, Lemma~\ref{lem:psigma1} gives
\[
P_{s_j+\ell}
=
(A_{s_j+\ell},B_{s_{j+1}-2-\ell})
=
b^{M_{j+1}}v^jd.
\]
These last $k_j$ paths use exactly the sinks $B_{s_j},\ldots,B_{s_{j+1}-2}$.
Deleting them therefore leaves a family in $\Omega_{s_j}^{(1)}$.

Conversely, these paths can be adjoined to every family in $\Omega_{s_j}^{(1)}$.
To verify non-intersection, we use the vertex comparison from the proof of Lemma~\ref{lem:psi1}.
A common vertex on the initial climb from $A_t$ and the final descent to $B_x$, at height $M_m$, would require
\[
t+x=M_m-1=s_{m-1}+s_m-2.
\]
For $m\le j$, the prescribed path forms give $t,x\ge s_{m-1}-1$ (for the added paths this is clear, since their indices are at least $s_j$), hence $t,x\le s_m-1<s_j$.
Thus such a vertex cannot involve an added path, whose source and sink indices are both at least $s_j$.
At height $M_{j+1}$, only the added paths occur, and the equality $x=M_{j+1}-1-t$ pairs the source and sink of the same added path.
Hence no intersection is introduced.
This establishes a bijection between $\Omega_{s_j}^{(1)}$ and $\Omega_n^{(1)}$.

The appended block
\[
\sigma^{j+1}
=
(s_{j+1}-2)(s_{j+1}-3)\cdots s_j
\]
has $k_j$ entries and inversion number $k_j(k_j-1)/2$.
Its entries exceed all preceding entries, so there are no inversions between the two parts.
Each added path has weight $V_j$.
Therefore,
\begin{align*}
H_{s_j+k_j}^{(1)}(F)
&=
H_{s_j}^{(1)}(F)\,
(-1)^{\mathrm{inv}(\sigma^{j+1})}
\mathcal{Q}_{s_j,s_{j+1}-2}\\
&=
(-1)^{\frac{k_j(k_j-1)}{2}}
V_j^{k_j}H_{s_j}^{(1)}(F),
\end{align*}
where $H_0^{(1)}(F)=1$ corresponds to the empty family.
This proves \eqref{eq:H1sjkj}.

\medskip
\noindent\textit{Proof of the vanishing assertion.}
Let $j\ge0$.
For $1\le\ell\le k_j-1$, set $n=s_j+\ell$.
Then $s_j\le n-1<s_{j+1}-2.$
The last assertion of Lemma~\ref{lem:psigma1}, applied with $i=j+1$, gives $\Omega_n^{(1)}=\emptyset$.
By Theorem~\ref{th:hanke(1)}, we obtain $ H_{s_j+\ell}^{(1)}(F)=0.$

Every positive integer $n$ can be written uniquely as $n=s_j+\ell$, where $j\ge0$ and $0\le\ell\le k_j$.
The remaining indices are $n=s_j$ and $n=s_j+k_j$, which are covered by \eqref{eq:H1sj} and \eqref{eq:H1sjkj}.
This completes the proof.
\end{proof}


\section{A partial combinatorial proof of Theorem \ref{th:shift2HDet}}\label{sec:doublesh}

In this section, we study the $2$–shifted Hankel determinants of $F(x)$ and, building on the previous results, prove Theorem \ref{th:shift2HDet}.
The non-intersecting paths give the determinants $H^{(2)}_{s_j+k_j-1}(F)$ and $H^{(2)}_{s_j+k_j}(F)$, and the vanishing of all the other ones except $H^{(2)}_{s_j}(F)$.
The determinants $H^{(2)}_{s_j}(F)$ are then obtained from these results and from Theorems~\ref{th:super2} and~\ref{th:shiftHDet} by means of the Desnanot--Jacobi identity.
A purely algebraic proof of Theorem \ref{th:shift2HDet} is given in Section~\ref{sec:doubleshalg}.

\subsection{Path families and preliminary lemmas}
\label{subsec:doubleshift:paths}

Let $n\ge1$.
For $0\le i\le n-1$, define the sources and sinks in the graph $\Gamma$ of Section~\ref{sec:lgv} by
\[
A_i:=(-2(i+1),0)_{\mathrm{src}},
\qquad
B_i:=(2(i+1),0)_{\mathrm{snk}}.
\]
Let $\Omega_n^{(2)}$ denote the set of all non-intersecting families $(P_0,\ldots,P_{n-1})$ of quasi-Dyck paths in $\Gamma$ such that, for some permutation $\sigma\in\mathfrak{S}_{[0,n-1]}$ and every $0\le i\le n-1$, $P_i=(A_i,B_{\sigma(i)}).$
When $k_0=0$, each path includes the auxiliary entrance edge of weight $1$.
As before, path vertices are the endpoints of edges, and two paths are non-intersecting if they share no vertex of $\Gamma$.
We also take $\Omega_0^{(2)}$ to consist of the empty family and set $H_0^{(2)}(F)=1$.

For $i,t\ge0$, a quasi-Dyck path from $A_i$ to $B_t$ has $i+t+2$ up steps.
The path interpretation of the coefficients of $F$ therefore gives
\[
a_{i+t+2}
=
\sum_{P:\,A_i\to B_t}w(P),
\]
where the sum runs over all paths of $\Gamma$ from $A_i$ to $B_t$.
For $n\ge1$, the corresponding matrix of path-weight sums is the $2$-shifted Hankel matrix of $F$,
whose determinant is
\[
H_n^{(2)}(F)
=
\det[a_{i+t+2}]_{0\le i,t\le n-1}.
\]
Since the proof of Theorem~\ref{th:lGV} applies verbatim to any sources and sinks of the acyclic graph $\Gamma$, we obtain the following analogue of Theorem~\ref{th:hanke(1)}.

 \begin{Theorem}\label{th:hankel(2)}
      For every $n \ge 1$, 
      \[
      H_{n}^{(2)}(F) = \sum_{(P_0, P_1, \dots, P_{n-1}) \in \Omega_{n}^{(2)}} (-1)^{\mathrm{inv}(\sigma)} \prod_{i=0}^{n-1} w(P_i),
      \]
     where the sum ranges over all non-intersecting families of paths $P_i$ going from $A_i$ to $B_{\sigma(i)}$ for $0 \le i \le n-1$, $\sigma$ being some permutation of $\mathfrak{S}_{[0,n-1]}$.
  \end{Theorem}

We introduce the following notation for signed weights. 
Let $n\ge1$, and let $P=(P_0,\ldots,P_{n-1})\in\Omega_n^{(2)}$ have associated permutation $\sigma=\sigma_1\cdots\sigma_n$.
We define the signed weight of $P$ by
\[
\mathcal Z_n(P)
:=
(-1)^{\mathrm{inv}(\sigma)}
\prod_{i=0}^{n-1}w(P_i).
\]
For any subset $\mathcal A\subseteq\Omega_n^{(2)}$, its total signed weight is
\[
\mathcal Z_n(\mathcal A)
:=
\sum_{P\in\mathcal A}\mathcal Z_n(P).
\]
In particular, Theorem~\ref{th:hankel(2)} gives
\begin{equation}\label{eq:znomega(2)}
    H_n^{(2)}(F)
=
\mathcal Z_n\bigl(\Omega_n^{(2)}\bigr).
\end{equation}
For a subfamily $Q=(P_a,\ldots,P_b)$, where $0\le a\le b\le n-1$, set
\[
\sigma(Q)
:=
(\sigma_{a+1},\ldots,\sigma_{b+1}).
\]
We define its partial signed weight by
\[
\mathcal Z(Q)
:=
(-1)^{\mathrm{inv}(\sigma(Q))}
\prod_{i=a}^{b}w(P_i),
\]
where $\mathrm{inv}(\sigma(Q))$ counts inversions within this subsequence.

To evaluate these signed weight sums, we first identify the initial climbs and final descents forced by non-intersection.
The following lemma is proved by the same induction as Lemma~\ref{lem:forced-segments1}.
We omit the details.

\begin{Lemma}\label{lem:forced-segments2}
Let $n\ge1$, and let $(P_0,\ldots,P_{n-1})\in\Omega_n^{(2)}$ have associated permutation
$\sigma=\sigma_1\cdots\sigma_n$.
For each $i\ge1$ and $0\le t\le n-1$, the following hold:
\begin{enumerate}
\item[(i)]
If $t\ge s_{i-1}$, then $P_t$ begins with $b^{M_i}$.

\item[(ii)]
If $\sigma_{t+1}\ge s_{i-1}$, then $P_t$ ends
with $v^{i-1}d$.
\end{enumerate}
Entrance edges, when present, are omitted from the step words.
\end{Lemma}

Let $(P_0,\ldots,P_{n-1})\in\Omega_n^{(2)}$ have associated permutation $\sigma$.
For each $i\ge1$ with $s_i-2\le n$, define the interior block
\[
\sigma^i
=
\sigma_{s_{i-1}+1}\sigma_{s_{i-1}+2}\cdots
\sigma_{s_i-2}
\qquad(k_{i-1}\ge2),
\]
and let $\sigma^i$ be empty when $k_{i-1}\le1$.
Thus $\sigma^i$ has $\max\{k_{i-1}-1,0\}$ entries.
The following lemma determines the paths corresponding to the entries of each interior block $\sigma^i$.
It is proved in the same way as Lemma~\ref{lem:psigma1}, using Lemma~\ref{lem:forced-segments2}, we omit the details.

\begin{Lemma}\label{lem:psigma2}
Let $n\ge1$, and let $(P_0,\ldots,P_{n-1})\in\Omega_n^{(2)}$ be a family of pairwise non-intersecting paths, with associated permutation $\sigma$.
For each $i\ge1$ with $k_{i-1}\ge2$ and $s_i-3\le n-1$, we have
\[
\sigma^i
=
(s_i-3)\,(s_i-4)\cdots s_{i-1}.
\]
More precisely, for $0\le\ell\le k_{i-1}-2$,
\[
P_{s_{i-1}+\ell}
=
(A_{s_{i-1}+\ell},B_{s_i-3-\ell})
=
b^{M_i}v^{i-1}d.
\]
These paths are uniquely determined.
Moreover, if there exists an index $i\ge1$ such that $ s_{i-1}\le n-1<s_i-3,$ then $\Omega_n^{(2)}$ is empty.
\end{Lemma}

Lemma~\ref{lem:psigma2} determines the paths in the interior blocks.
For $n=s_j$ and $1\le i\le j$, the remaining paths in block $i$ are $P_{s_i-1}$ if $k_{i-1}=0$, and $P_{s_i-2}$ and $P_{s_i-1}$ if $k_{i-1}\ge1$.
Instead of describing these remaining paths completely, we combine the path results above with the Desnanot--Jacobi identity from
Subsection~\ref{subsec:doubleshalg:djiden} to evaluate the determinants.
In the next two subsections, we treat the cases $k_j\ge1$ and $k_j=0$, respectively.

\subsection{The case $k_j\ge1$.}\label{subsec:doubleshift:kge1}
Recall from Subsection~\ref{subsec:doubleshalg:djiden} that, for every $n\ge1$, the Desnanot--Jacobi identity gives
\begin{equation}\label{eq:DJ2}
H_{n+1}(F)\,H^{(2)}_{n-1}(F)
=
H_n(F)\,H^{(2)}_n(F)
-
\bigl(H^{(1)}_n(F)\bigr)^2.
\end{equation}
Since all $v_i$ are nonzero, Theorem~\ref{th:super2} gives $H_{s_j}(F)\neq0$ for every $j\ge0$.
Moreover, for $j\ge0$, \eqref{eq:HankelDetFundamental} yields
\begin{equation}\label{eq:Hratio2}
H_{s_{j+1}}(F)
=
(-1)^{\frac{k_j(k_j+1)}{2}}
V_j^{k_j+1}H_{s_j}(F).
\end{equation}

We first apply the Desnanot--Jacobi identity to evaluate $H_n^{(2)}(F)$ at the indices $n=s_j$ with $j\ge1$ and $k_j\ge1$.

\begin{Theorem}\label{th:h(2)r}
For every $j\ge1$ with $k_j\ge1$, we have
\begin{equation}\label{eq:hsj2}
    H_{s_j}^{(2)}(F) = H_{s_j}(F) R_j^2,
\end{equation}
where $(R_n)_{n \ge 0}$ is the sequence \eqref{def:Rn}.
\end{Theorem}

\begin{proof}
Set $n=s_j$.
Since $s_{j+1}=s_j+k_j+1$ and $k_j\ge1$, we have $s_j<n+1<s_{j+1}$.
Thus $H_{n+1}(F)=0$ by Theorem~\ref{th:super2}.
Moreover, \eqref{eq:H1sj:Hsj} gives
\[
\bigl(H_n^{(1)}(F)\bigr)^2=H_n(F)^2R_j^2.
\]
Substituting these relations into \eqref{eq:DJ2}, we obtain
\[
H_n(F)H_n^{(2)}(F)=H_n(F)^2R_j^2.
\]
Since $H_n(F)\neq0$ by Theorem~\ref{th:super2}, dividing by $H_n(F)$ gives \eqref{eq:hsj2}.
\end{proof}

We next evaluate $H_n^{(2)}(F)$ at $n=s_j+k_j-1$ for $j\ge0$ with $k_j\ge1$.

\begin{Theorem}\label{th:h(2):kj11}
For every $j\ge0$ with $k_j\ge1$, we have
\begin{equation}\label{eq:thh2k1}
H_{s_j+k_j-1}^{(2)}(F)
=
(-1)^{\frac{(k_j-1)(k_j-2)}{2}}
V_j^{k_j-1}H_{s_j}^{(2)}(F).
\end{equation}
\end{Theorem}

\begin{proof}
If $k_j=1$, both sides of \eqref{eq:thh2k1} are equal to $H_{s_j}^{(2)}(F)$.
Assume now that $k_j\ge2$, and set $n=s_j+k_j-1=s_{j+1}-2$.
For a family $P=(P_0,\ldots,P_{n-1})\in\Omega_n^{(2)}$, write
\[
P^{\mathrm L}=(P_0,\ldots,P_{s_j-1})
\qquad
\text{and}
\qquad
P^{\mathrm R}=(P_{s_j},\ldots,P_{n-1}).
\]
When $j=0$, the family $P^{\mathrm L}$ is empty, and we use the convention $\mathcal Z_0(\Omega_0^{(2)})=1$.
By Lemma~\ref{lem:psigma2}, applied with $i=j+1$, for each $0\le\ell\le k_j-2$ we have
\[
P_{s_j+\ell}
=
(A_{s_j+\ell},B_{s_j+k_j-2-\ell})
=
b^{M_{j+1}}v^jd.
\]
Thus $P^{\mathrm R}$ is uniquely determined, and it uses exactly the sinks $B_{s_j}$, \ldots, $B_{n-1}$.
Hence the remaining paths use precisely the sinks $B_0$, \ldots, $B_{s_j-1}$, and $P^{\mathrm L}\in\Omega_{s_j}^{(2)}$.

Conversely, take any family $P^{\mathrm L}\in\Omega_{s_j}^{(2)}$ and append the prescribed paths above.
At height $h$, the abscissa $u$ of any ordinary vertex of $P^{\mathrm L}$ satisfies
$-2s_j+h\le u\le 2s_j-h.$
Indeed, reaching height $h$ requires a horizontal displacement of at least $h$, as does descending from that height to a sink.
For an appended path, the abscissa of a vertex on its initial climb satisfies $u\le -2s_j-2+h,$
whereas the abscissa of a vertex on its final descent satisfies $u\ge 2s_j+2-h.$
Thus the appended paths share no ordinary vertex with $P^{\mathrm L}$.

The appended paths are themselves pairwise non-intersecting.
Distinct climbs lie on distinct parallel lines, as do distinct final descents.
Below height $M_{j+1}$, the descent vertices occur only at heights at most $M_j<2s_j+2$, where the two bounds above exclude a common vertex between a climb and a descent.
At height $M_{j+1}$, each descent starts at the endpoint of its own climb, and these endpoints are distinct.
Since all sources and sinks are distinct, the extended family belongs to $\Omega_n^{(2)}$.

Deletion and extension are therefore inverse operations.
All sink indices in $P^{\mathrm R}$ are larger than those in $P^{\mathrm L}$, so there are no inversions between the two parts.
Their signed weights therefore multiply, giving
\begin{equation}\label{eq:zpskj1}
\mathcal Z_{s_j+k_j-1}
\bigl(\Omega_{s_j+k_j-1}^{(2)}\bigr)
=
\mathcal Z_{s_j}
\bigl(\Omega_{s_j}^{(2)}\bigr)
\mathcal Z(P^{\mathrm R}).
\end{equation}
The sink indices of $P^{\mathrm R}$ form the decreasing sequence
\[
\sigma^{j+1}=(n-1)(n-2)\cdots s_j,
\]
which has $k_j-1$ entries.
Each appended path has weight $V_j$.
Therefore,
\begin{equation*}
\mathrm{inv}(\sigma^{j+1})
=
\frac{(k_j-1)(k_j-2)}{2}
\qquad \text{and} \qquad
\mathcal{Q}_{s_j,s_{j+1}-3}
=
V_j^{k_j-1}.
\end{equation*}
Consequently,
\begin{equation}\label{eq:zpkjpr}
\mathcal Z(P^{\mathrm R})
=
(-1)^{\mathrm{inv}(\sigma^{j+1})}
\mathcal Q_{s_j,s_{j+1}-3}
=
(-1)^{\frac{(k_j-1)(k_j-2)}{2}}V_j^{k_j-1}.
\end{equation}
Substituting \eqref{eq:zpkjpr} into \eqref{eq:zpskj1} and combining \eqref{eq:znomega(2)} gives \eqref{eq:thh2k1}.
\end{proof}

We then derive a recurrence for $H_n^{(2)}(F)$ at $n=s_j+k_j$, again assuming $j\ge0$ and $k_j\ge1$.

\begin{Theorem}\label{th:h(2):kj12}
For $j \ge 0$ with $k_j \ge 1$, we have
    \begin{equation}\label{eq:thh2k}
	H_{s_j+k_j}^{(2)}( F )=  c_{j+1,k_j-1} V_j   H_{s_j+k_j-1}^{(2)}(F) + \Delta_{j},
\end{equation}
where 
\[
\Delta_j = 
\left\{
\begin{array}{ll}
 (-1)^{\frac{(k_j+1)k_j}{2}}V_j^{k_j+1} \cdot H_{s_{j-1}}^{(2)}(F) & \text{if $k_{j-1} = 0$},\\
 (-1)^{\frac{(k_j+1)k_j}{2}}V_j^{k_j+1} \cdot H_{s_{j-1}+k_{j-1}}^{(2)}(F)  & \text{if  $k_{j-1} \ge 1$},
\end{array}
\right.
\]
with $\Delta_0=0$, where $c_{j+1,k_j-1}$ denotes the coefficient of $x^{k_j-1}$ in $u_{j+1}(x)$.
\end{Theorem}

\begin{proof}
Set $n=s_j+k_j$, so that $M_{j+1}=2s_j+k_j=s_j+n. $
Let $P=(P_0,\ldots,P_{n-1})\in\Omega_n^{(2)}$.
For $k_j\ge2$ and $0\le\ell\le k_j-2$, Lemma~\ref{lem:psigma2}, applied to block $j+1$, gives
\[
P_{s_j+\ell} = (A_{s_j+\ell},B_{n-2-\ell}) = b^{M_{j+1}}v^jd.
\]
This list is empty when $k_j=1$.
These paths use precisely the sinks $B_{s_j},\ldots,B_{n-2}$.
By Lemma~\ref{lem:forced-segments2}(i), $P_{n-1}$ begins with $b^{M_{j+1}}$.
If its endpoint is $B_x$, it has $n+x+1$ up steps.
Hence
\[
n+x+1\ge M_{j+1}=s_j+n,
\]
which gives $x\ge s_j-1$.
The remaining sink indices are $0,\ldots,s_j-1$ and $n-1$, so only $x=s_j-1$ and $x=n-1$ are possible.
We distinguish these two cases. 
When $j=0$, the sink $B_{s_j-1}$ does not exist, and only Case~(1) occurs.

\begin{itemize}
\item[(1)]
If $P_{n-1}$ ends at $B_{n-1}$, then the total signed weight of the families in this case is
\[
c_{j+1,k_j-1}V_jH_{n-1}^{(2)}(F).
\]

\item[(2)]
If $P_{n-1}$ ends at $B_{s_j-1}$, then the total signed weight of the families in this case is
\[
(-1)^{\frac{k_j(k_j+1)}{2}}
V_j^{k_j+1}H_{s_j-1}^{(2)}(F).
\]
\end{itemize}
The evaluations in (1) and (2) are proved in Lemmas~\ref{lem:h2:last-fixed} and~\ref{lem:h2:last-paired} below, respectively.

For $j\ge1$, combining the two cases, we obtain 
\[
H_{s_j+k_j}^{(2)}(F) = 
c_{j+1,k_j-1}V_j
H_{s_j+k_j-1}^{(2)}(F) 
+
(-1)^{\frac{k_j(k_j+1)}{2}}
V_j^{k_j+1}H_{s_j-1}^{(2)}(F),
\]
Since $s_j-1=s_{j-1}+k_{j-1}$, the second term is precisely $\Delta_j$.
When $j=0$, only Case~(1) occurs, and $\Delta_0=0$.
This proves \eqref{eq:thh2k}.
\end{proof}

We now prove the two signed-weight evaluations used above.
The empty family is understood to have signed weight $1$.

\begin{Lemma}\label{lem:h2:last-fixed}
Let $j\ge0$ with $k_j\ge1$, and set $n=s_j+k_j$.
Let $\mathcal A_j$ denote the set of all families $P=(P_0,\ldots,P_{n-1})\in\Omega_n^{(2)}$ such that $P_{n-1}$ ends at $B_{n-1}$.
Then
\[
\mathcal Z_n(\mathcal A_j)
=
c_{j+1,k_j-1}V_jH_{n-1}^{(2)}(F).
\]
\end{Lemma}

\begin{proof}
Let $P\in\mathcal A_j$.
By Lemma~\ref{lem:forced-segments2}, $P_{n-1}$ begins with $b^{M_{j+1}}$ and ends with $v^jd$.
The portion between these segments runs from $(-k_j,M_{j+1})$ to $(k_j,M_{j+1})$ and has exactly $k_j$ up steps.
This portion cannot go below $M_{j+1}$.
For $j\ge1$, doing so would require a long down step to $M_j$, after which returning to $M_{j+1}$ would require at least
\[
M_{j+1}-M_j=k_{j-1}+k_j+2>k_j
\]
up steps.
For $j=0$, the only step descending below $M_1$ is the terminal step, which ends the path.
The portion also cannot reach $M_{j+2}$, since
\[
M_{j+2}-M_{j+1}=k_j+k_{j+1}+2>k_j.
\]
Its first down step must therefore be a short down step ending at $M_{j+1}$.
If its drop is $r$, then $1\le r\le k_j$, and its landing vertex is $(-k_j+2r,M_{j+1}).$
If $r<k_j$, then $s_j\le n-1-r\le n-2$, and this vertex is the endpoint of the initial climb of the fixed path $P_{n-1-r}$.
This contradicts non-intersection.
Thus $r=k_j$, using all the remaining up steps, and consequently
\[
P_{n-1}
=
(A_{n-1},B_{n-1})
=
b^{M_{j+1}+k_j}c_{j+1,k_j-1}v^jd.
\]
Its weight is $w(P_{n-1})=c_{j+1,k_j-1}V_j$.

Deleting this path leaves a family in $\Omega_{n-1}^{(2)}$.
Conversely, it can be appended to any family in $\Omega_{n-1}^{(2)}$ without creating an intersection.
Indeed, at height $h$, every ordinary vertex of such a family has abscissa between $-2n+2+h$ and $2n-2-h$, whereas the vertices of the appended path have abscissa $-2n+h$ on its initial climb or $2n-h$ on its final descent.
The endpoints of its short down step lie on these two segments, and the interior of that step contains no path vertices.
Thus no intersection is introduced.

Deletion and extension are therefore inverse operations.
The appended path creates no new inversion, since it joins the largest source index to the largest sink index.
Multiplying by its weight and applying Theorem~\ref{th:hankel(2)} gives
\[
\mathcal Z_n(\mathcal A_j)
=
w(P_{n-1})H_{n-1}^{(2)}(F)
=
c_{j+1,k_j-1}V_jH_{n-1}^{(2)}(F).
\]
\end{proof}

\begin{Lemma}\label{lem:h2:last-paired}
Let $j\ge1$ with $k_j\ge1$, and set $n=s_j+k_j$.
Let $\mathcal B_j$ denote the set of all families $P=(P_0,\ldots,P_{n-1})\in\Omega_n^{(2)}$ such that $P_{n-1}$ ends at $B_{s_j-1}$.
Then
\[
\mathcal Z_n(\mathcal B_j)
=
(-1)^{\frac{k_j(k_j+1)}{2}}
V_j^{k_j+1}H_{s_j-1}^{(2)}(F).
\]
\end{Lemma}

\begin{proof}
Let $P\in\mathcal B_j$.
The path $P_{n-1}$ has exactly $(n-1)+(s_j-1)+2=M_{j+1}$ up steps.
By Lemma~\ref{lem:forced-segments2}(i), it begins with $b^{M_{j+1}}$, which uses all these up steps.
The remaining steps must therefore be the successive long down steps followed by the terminal step.
Hence
\[
P_{n-1}
=
(A_{n-1},B_{s_j-1})
=
b^{M_{j+1}}v^jd.
\]
For $k_j\ge2$ and $0\le\ell\le k_j-2$, Lemma~\ref{lem:psigma2} gives
\[
P_{s_j+\ell}
=
(A_{s_j+\ell},B_{n-2-\ell})
=
b^{M_{j+1}}v^jd.
\]
This list is empty when $k_j=1$.

Let $P_t$ be the path ending at $B_{n-1}$.
By Lemma~\ref{lem:forced-segments2}(ii), it ends with $v^jd$ and therefore reaches height $M_{j+1}$.
Since it has $t+n+1$ up steps,
\[
t+n+1\ge M_{j+1}=s_j+n,
\]
so $t\ge s_j-1$.
The fixed paths and $P_{n-1}$ already account for all sources with indices at least $s_j$.
Thus $t=s_j-1$.
This path also has exactly $M_{j+1}$ up steps.
Since it reaches height $M_{j+1}$, all these up steps occur before the final descent, with no intervening down step.
Therefore,
\[
P_{s_j-1} = (A_{s_j-1},B_{n-1}) = b^{M_{j+1}}v^jd.
\]

Together with the fixed interior paths, these are the $k_j+1$ paths with source indices $s_j-1,\ldots,n-1$.
Their sink indices form the decreasing sequence $(n-1)(n-2)\cdots(s_j-1)$, and all have step word $b^{M_{j+1}}v^jd$.
Deleting them leaves a family in $\Omega_{s_j-1}^{(2)}$.

Conversely, these prescribed paths can be appended to any family in $\Omega_{s_j-1}^{(2)}$.
At height $h$, the ordinary vertices of the remaining family have abscissae between $-2s_j+2+h$ and $ 2s_j-2-h$,  while those of the appended paths have abscissae at most $-2s_j+h$ on their climbs and at least $2s_j-h$ on their descents.
Thus the two parts have no common ordinary vertex.

The appended paths are also pairwise non-intersecting.
Their climbs lie on distinct parallel lines, as do their descents.
Below height $M_{j+1}$, descent vertices occur only at heights at most $M_j=2s_j-k_{j-1}-2<2s_j.$
At these heights, $-2s_j+h<2s_j-h$, so the preceding bounds exclude a common vertex between a climb and a descent.
At height $M_{j+1}$, each descent starts at the endpoint of its own climb, and these endpoints are distinct.
Since all sources and sinks are distinct graph vertices, the extension is non-intersecting.

Deletion and extension therefore give a bijection between $\mathcal B_j$ and $\Omega_{s_j-1}^{(2)}$.
The deleted block is decreasing and has $k_j+1$ entries, so it contributes $\frac{k_j(k_j+1)}{2}$ inversions.
There are no inversions between this block and the remaining entries, whose positions and values are all smaller.
Each path in the block has weight $V_j$.
Thus Theorem~\ref{th:hankel(2)} gives
\[
\mathcal Z_n(\mathcal B_j)
=
(-1)^{\frac{k_j(k_j+1)}{2}}
V_j^{k_j+1}H_{s_j-1}^{(2)}(F).
\]
\end{proof}

\subsection{The case $k_j=0$ and completion of the proof}\label{subsec:doubleshift:k=0}
We now treat the indices $s_{j+1}-1$ for every $j\ge0$, and then the case $k_j=0$.
For $t\ge0$, let $d_{t+1}$ denote the coefficient of $x^{k_t}$ in $1-u_{t+1}(x)x$. Thus
\[
d_{t+1}
=
\begin{cases}
1, & \text{if } k_t=0,\\
-c_{t+1,k_t-1}, & \text{if } k_t\ge1.
\end{cases}
\]

\begin{Lemma}\label{lem:h2sjm1}
For every $j\ge0$, we have
\begin{equation}\label{eq:h2sjm1}
H^{(2)}_{s_{j+1}-1}(F)=H_{s_{j+1}}(F)\sum_{t=0}^{j}\frac{d_{t+1}R_t^2}{V_t}.
\end{equation}
\end{Lemma}

\begin{proof}
Set $S_{-1}=0$ and, for $j\ge0$, define
\[
S_j=\sum_{t=0}^{j}\frac{d_{t+1}R_t^2}{V_t}.
\]
Thus, for every $j\ge0$,
\[
S_j=S_{j-1}+\frac{d_{j+1}R_j^2}{V_j}.
\]
We prove \eqref{eq:h2sjm1} by induction on $j$.
For $j\ge1$, the induction hypothesis is
\[
H^{(2)}_{s_j-1}(F)=H_{s_j}(F)S_{j-1}.
\]

\medskip

\noindent\textbf{Case 1.} 
Suppose that $k_j=0$.
Then $s_{j+1}-1=s_j$ and $d_{j+1}=1$.
By \eqref{eq:Hratio2}, we have
\[
H_{s_{j+1}}(F)=V_jH_{s_j}(F).
\]
If $j=0$, then $s_1=1$, and
\[
H_{s_1}(F)S_0
=
V_0\frac{R_0^2}{V_0}
=
1
=
H_0^{(2)}(F).
\]
Hence \eqref{eq:h2sjm1} holds in this case.
Now let $j\ge1$.
Since $H_{s_j}(F)\neq0$ by Theorem~\ref{th:super2}, applying \eqref{eq:DJ2} with $n=s_j$ and using \eqref{eq:H1sj:Hsj} gives
\begin{align*}
H^{(2)}_{s_j}(F)
&=
H_{s_j}(F)R_j^2
+
V_jH^{(2)}_{s_j-1}(F)\\
&=
V_jH_{s_j}(F)
\left(S_{j-1}+\frac{R_j^2}{V_j}\right)\\
&=
H_{s_{j+1}}(F)S_j,
\end{align*}
where the second equality uses the induction hypothesis and the last uses $d_{j+1}=1$.

\medskip
\noindent\textbf{Case 2.} 
Suppose that $k_j\ge1$.
Then $s_{j+1}-1=s_j+k_j$ and $d_{j+1}=-c_{j+1,k_j-1}$.
Set
\[
\rho_j
=
(-1)^{\frac{k_j(k_j+1)}{2}}V_j^{k_j+1},
\]
so that \eqref{eq:Hratio2} becomes $H_{s_{j+1}}(F)=\rho_jH_{s_j}(F).$
For $j\ge1$, Theorem~\ref{th:h(2)r} gives $H_{s_j}^{(2)}(F)=H_{s_j}(F)R_j^2$.
This identity also holds for $j=0$, since both sides equal $1$.
Combining it with Theorem~\ref{th:h(2):kj11}, we obtain
\begin{align*}
c_{j+1,k_j-1}V_jH^{(2)}_{s_j+k_j-1}(F)
&=
(-1)^{\frac{(k_j-1)(k_j-2)}{2}}
c_{j+1,k_j-1}V_j^{k_j}H_{s_j}(F)R_j^2\\
&=
H_{s_{j+1}}(F)\frac{d_{j+1}R_j^2}{V_j}.
\end{align*}
The last equality follows from \eqref{eq:Hratio2} and $d_{j+1}=-c_{j+1,k_j-1}$, since
\[
\frac{k_j(k_j+1)}{2}
-
\frac{(k_j-1)(k_j-2)}{2}
=
2k_j-1
\]
is odd.

For $j\ge1$, the identity $s_j-1=s_{j-1}+k_{j-1}$ shows that both cases in the definition of $\Delta_j$ give $\Delta_j=\rho_jH^{(2)}_{s_j-1}(F)$.
The induction hypothesis therefore yields
\[
\Delta_j=\rho_jH^{(2)}_{s_j-1}(F)=\rho_jH_{s_j}(F)S_{j-1}=H_{s_{j+1}}(F)S_{j-1}.
\]
For $j=0$, we also have $\Delta_0=H_{s_1}(F)S_{-1}=0$.
Theorem~\ref{th:h(2):kj12} now gives
\[
H^{(2)}_{s_{j+1}-1}(F) 
=  H_{s_{j+1}}(F)
\left(S_{j-1}+\frac{d_{j+1}R_j^2}{V_j}\right)
= H_{s_{j+1}}(F)S_j.
\]
This completes the induction.
\end{proof}

\begin{Remark}
Formula~\eqref{eq:h2sjm1} is also proved in Proposition~\ref{prop:alg2:cofactor} by a cofactor argument using orthogonal polynomials.
\end{Remark}

Applying Lemma~\ref{lem:h2sjm1} with $k_j=0$ gives the following evaluation of $H_{s_j}^{(2)}(F)$.

\begin{Theorem}\label{th:h(2):kj0}
For $j \ge 0$ with $k_j = 0$, we have 
    \begin{equation*}
	H_{s_j}^{(2)}( F )=   H_{s_j}(F)  \sum_{t=0}^{j} d_{t+1} R_t^2 v_{t+1} v_{t+2}  \cdots v_{j},
\end{equation*}
where $d_j$ is the coefficient of $x^{k_{j-1}}$ in $1-u_j(x)x$. 
\end{Theorem}

\begin{proof}
Since $k_j=0$, we have $s_{j+1}-1=s_j$.
By \eqref{eq:Hratio2}, we also have $H_{s_{j+1}}(F)=V_jH_{s_j}(F)$.
Lemma~\ref{lem:h2sjm1} therefore gives
\[
H^{(2)}_{s_j}(F)
=
H_{s_j}(F)\sum_{t=0}^{j}
d_{t+1}R_t^2\,\frac{V_j}{V_t}.
\]
Since $V_j/V_t=v_{t+1}v_{t+2}\cdots v_j$, the desired formula follows.
\end{proof}

We now combine the preceding results to complete the proof of Theorem~\ref{th:shift2HDet}.

\begin{proof}[Proof of Theorem~\ref{th:shift2HDet}]
Equation~\eqref{eq:H1sj:Hsj:k1} follows from Theorem~\ref{th:h(2)r} for $j\ge1$.
For $j=0$, it reads $H_0^{(2)}(F)=1=H_0(F)R_0^2$.
Equations~\eqref{eq:H1sj:Hsj:k11}, \eqref{eq:H1sj:Hsj:k12}, and \eqref{eq:H1sj:Hsj:k0} follow from Theorems~\ref{th:h(2):kj11}, \ref{th:h(2):kj12}, and \ref{th:h(2):kj0}, respectively.
It remains to prove the vanishing assertion in part~(3).

Let $j\ge0$ and $1\le\ell\le k_j-2$, and set $n=s_j+\ell$.
Assume, for contradiction, that $\Omega_n^{(2)}$ contains a family $(P_0,\ldots,P_{n-1})$.
Since $s_j\le n-1$, this family contains a path $P_t$ ending at $B_{s_j}$.
By Lemma~\ref{lem:forced-segments2}(ii), applied with $i=j+1$, the path $P_t$ ends with $v^jd$ and therefore reaches height $M_{j+1}=2s_j+k_j.$
However, a path from $A_t$ to $B_{s_j}$ has $t+s_j+2$ up steps.
Since $t\le n-1$, we have
\[
t+s_j+2 \le 2s_j+\ell+1 \le 2s_j+k_j-1 = M_{j+1}-1.
\]
Thus $P_t$ cannot reach height $M_{j+1}$, a contradiction.
Therefore $\Omega_n^{(2)}=\emptyset$, and Theorem~\ref{th:hankel(2)} gives $H_{s_j+\ell}^{(2)}(F)=0.$
This proves part~(3) and completes the proof.
\end{proof}


\section{A combinatorial proof of Theorem~\ref{th:backshiftHDet}}\label{sec:backsh}

In this section, we study the $(-1)$-shifted Hankel determinants of $F(x)$ and give a combinatorial proof of Theorem~\ref{th:backshiftHDet}.
Following the approach of Section~\ref{sec:shift}, we determine the structure of the relevant non-intersecting families of quasi-Dyck paths.
Evaluating their signed weights gives the formulas for $H_{s_j+1}^{(-1)}(F)$ and $H_{s_{j+1}}^{(-1)}(F)$,
while the absence of such families at the remaining indices proves the vanishing assertions.
An algebraic proof of Theorem~\ref{th:backshiftHDet} is given in Section~\ref{sec:shiftalg}.

\subsection{The path model for $(-1)$-shifted Hankel determinants}
\label{subsec:backshift:paths}

For $i\ge0$, define the sources and sinks in the graph $\Gamma$ of Section~\ref{sec:lgv} by
\[
A_i:=(-2i,0)_{\mathrm{src}},
\qquad
B_i:=(2(i-1),0)_{\mathrm{snk}}.
\]
For each $n\ge1$, let $\Omega_n^{(-1)}$ denote the set of all non-intersecting families
$(P_0,\ldots,P_{n-1})$ of quasi-Dyck paths in $\Gamma$ such that, for some permutation
$\sigma\in\mathfrak{S}_{[0,n-1]}$ and every $0\le i\le n-1$,
\[
P_i=(A_i,B_{\sigma(i)}).
\]

For $i,t\ge0$ with $i+t\ge1$, a quasi-Dyck path from $A_i$ to $B_t$ has $i+t-1$ up steps.
By Theorem~\ref{th:charhankel}, we therefore have
\[
a_{i+t-1}
=
\sum_{P:\,A_i\to B_t}w(P).
\]
There is no path from $A_0$ to $B_0$, in agreement with the convention $a_{-1}=0$.
Thus the corresponding matrix of path-weight sums is the matrix in \eqref{eq:defHm1}, whose determinant is $H_n^{(-1)}(F)$.
Applying the Lindstr\"om--Gessel--Viennot theorem to these sources and sinks gives the following formula.

\begin{Theorem}\label{th:hankel(-1)}
For every $n\ge1$,
\[
H_{n}^{(-1)}(F) = \sum_{(P_0, P_1, \dots, P_{n-1}) \in \Omega_{n}^{(-1)}} (-1)^{\mathrm{inv}(\sigma)}\prod_{i=0}^{n-1}w(P_i),
\]
where the sum ranges over all non-intersecting families of paths $P_i$ going from $A_i$ to $B_{\sigma(i)}$ for $0 \le i \le n-1$, $\sigma$ being some permutation of $\mathfrak{S}_{[0,n-1]}$.
\end{Theorem}

As in Section~\ref{sec:shift}, we use the one-line notation $\sigma=\sigma_1\cdots\sigma_n$, where $\sigma_{i+1}=\sigma(i)$.
Auxiliary entrance edges, when present, are included in the paths but omitted from their step words.
When $k_0=0$, the paths from $A_0$ to $B_1$ and from $A_1$ to $B_0$ have step word $d$, and each has weight $v_0$.
When $k_0\ge1$, neither path exists.

We record the vertex coordinates used below.
For $i\ge1$, an initial climb $b^{M_i}$ from $A_t$ ends at the ordinary vertex $(M_i-2t,M_i),$ while a final descent $v^{i-1}d$ to $B_x$ starts at $(2x-2-M_i,M_i).$
Since $M_i=s_{i-1}+s_i-1$, these vertices coincide if and only if
\[
t+x=M_i+1=s_{i-1}+s_i.
\]
Moreover, any ordinary vertex $(u,h)$ of a path from $A_t$ to $B_x$ satisfies $u\ge-2t+h.$
Indeed, reaching height $h$ requires a horizontal displacement of at least $h$.

The following analogue of Lemma~\ref{lem:forced-segments} describes the initial climbs and final descents forced by non-intersection.
Its proof uses the same induction and is omitted.

\begin{Lemma}\label{lem:forced-segments-1}
Let $n\ge1$, and let $(P_0,\ldots,P_{n-1})\in\Omega_n^{(-1)}$ have associated permutation $\sigma$.
For each $i\ge1$ and $0\le t\le n-1$, the following hold:
\begin{enumerate}
\item[(i)]
If $i=1$ or $t\ge s_{i-1}+1$, then $P_t$ begins with $b^{M_i}$.
\item[(ii)]
If $i=1$ or $\sigma_{t+1}\ge s_{i-1}+1$, then $P_t$ ends
with $v^{i-1}d$.
\end{enumerate}
\end{Lemma}

\subsection{Block structure and boundary paths}
\label{subsec:backshift:structure}
Let $n\ge1$, and let $(P_0,\ldots,P_{n-1})\in\Omega_n^{(-1)}$ have associated permutation $\sigma$.
If $n\ge s_1+1$, define the first block
\[
\sigma^1=\sigma_1\sigma_2\cdots\sigma_{s_1+1},
\]
which has $s_1+1=k_0+2$ entries.
For each $i\ge2$ with $s_i\le n$, define the consecutive block
\[
\sigma^i
=
\sigma_{s_{i-1}+2}\,\sigma_{s_{i-1}+3}\cdots
\sigma_{s_i},
\]
which has $k_{i-1}$ entries and is empty when $k_{i-1}=0$.
In particular, 
When $n=s_1+1$, we have $\sigma=\sigma^1$.
For $n=s_j+1$ with $j\ge2$, we have the decomposition
\[
\sigma
=
\sigma^1\,
\sigma^2\,\sigma_{s_2+1}\,
\sigma^3\,\sigma_{s_3+1}\cdots
\sigma^j\,\sigma_{s_j+1}.
\]
The following lemma determines the paths corresponding to the entries of the blocks $\sigma^i$.
Its proof follows that of Lemma~\ref{lem:psigma1}, using Lemma~\ref{lem:forced-segments-1}, and is omitted.

\begin{Lemma}\label{lem:psigma-1}
Let $n\ge1$, and let $(P_0, P_1, \dots, P_{n-1}) \in \Omega_{n}^{(-1)}$ have associated permutation $\sigma$.
\begin{enumerate}
\item[(i)]
If $s_1\le n-1$, then
\[
\sigma^1=s_1\,(s_1-1)\cdots 1\,0,
\]
that is, for $0\le\ell\le s_1$,
\[
P_{\ell}
= (A_{\ell},B_{s_1-\ell})
= b^{M_1}d.
\]

\item[(ii)]
For each $i\ge2$ with $k_{i-1}\ge 1$ and $s_i-1\le n-1$, we have
\[
\sigma^i = (s_i-1)\,(s_i-2)\cdots (s_{i-1}+1),
\]
that is, for $1\le \ell\le k_{i-1}$,
\[
P_{s_{i-1}+\ell}
= (A_{s_{i-1}+\ell},B_{s_i-\ell})
= b^{M_i}v^{i-1}d.
\]

\item[(iii)]
If $n\le s_1$, or if $s_{i-1}+1\le n-1<s_i-1$ for some $i\ge2$, then $\Omega_n^{(-1)}$ is empty.
\end{enumerate}
\end{Lemma}

For $n=s_j+1$, Lemma~\ref{lem:psigma-1} leaves only the boundary paths $P_{s_i}$, with $2\le i\le j$, to be determined.
We first classify the paths ending at $B_{s_i}$.
The following lemma is proved as in Lemma~\ref{lem:free1}, using Lemmas~\ref{lem:primitive}, \ref{lem:forced-segments-1}, and~\ref{lem:psigma-1}.
We omit the details.

\begin{Lemma}\label{lem:free-1}
Let $n\ge1$, and let $(P_0,\ldots,P_{n-1})\in\Omega_n^{(-1)}$.
Let $i\ge2$ with $s_i\le n-1$.
Then the path ending at $B_{s_i}$ is described by exactly one of the following cases:
\begin{itemize}
\item[(i)] The path is
\[
P_{s_i}
=
(A_{s_i},B_{s_i})
=
b^{N_i}c_iv^{i-1}d.
\]

\item[(ii)]
For $i\ge3$, the path is 
\[
P_{s_{i-1}}
=
(A_{s_{i-1}},B_{s_i})
=
b^{M_i}v^{i-1}d.
\]

\item[(iii)]
If $s_{i+1}\le n-1$, then the path is 
\[
P_{s_{i+1}}
=
(A_{s_{i+1}},B_{s_i})
=
b^{M_{i+1}}v^id.
\]
\end{itemize}
Here $c_i=c_{i,k_{i-1}}$ denotes the short down step of drop $k_{i-1}+1$ ending at height $M_i$.
\end{Lemma}

We now use Lemma~\ref{lem:free-1} to describe the possible connections between the remaining boundary sources and sinks in families of $\Omega_{s_j+1}^{(-1)}$.

\begin{Lemma}\label{lem:psi-11}
Let $j\ge1$, and let $(P_0,\ldots,P_{s_j})\in\Omega_{s_j+1}^{(-1)}$ have associated permutation $\sigma$.
For each $2\le i\le j$, exactly one of the following cases occurs:
\begin{itemize}
\item[(i)]
We have $\sigma_{s_i+1}=s_i$, and
\[
P_{s_i}=(A_{s_i},B_{s_i})=b^{N_i}c_iv^{i-1}d.
\]

\item[(ii)]
For $i\ge3$, we have $\sigma_{s_i+1}=s_{i-1}$, $\sigma_{s_{i-1}+1}=s_i$, and
\[
\begin{aligned}
P_{s_{i-1}}
&=(A_{s_{i-1}},B_{s_i})
=b^{M_i}v^{i-1}d,\\
P_{s_i}
&=(A_{s_i},B_{s_{i-1}})
=b^{M_i}v^{i-1}d.
\end{aligned}
\]

\item[(iii)]
For $i\le j-1$, we have $\sigma_{s_i+1}=s_{i+1}$ and $\sigma_{s_{i+1}+1}=s_i$. This is case (ii) with $i$ replaced by $i+1$.
\end{itemize}

Conversely, every collection of disjoint adjacent transpositions of $(s_2,\ldots,s_j)$ determines a unique family in $\Omega_{s_j+1}^{(-1)}$, with the remaining paths given by Lemma~\ref{lem:psigma-1}.
\end{Lemma}

\begin{proof}
By Lemma~\ref{lem:psigma-1}, the entries $\sigma_{s_2+1},\ldots,\sigma_{s_j+1}$ form a permutation of $s_2,\ldots,s_j$.
By Lemma~\ref{lem:free-1}, a connection from $A_{s_m}$ to $B_{s_i}$ is possible only if $|m-i|\le1$.
Thus each boundary source can be connected only to its corresponding boundary sink or to a neighboring boundary sink in this list.
If not all boundary entries are fixed, let $i$ be the smallest index for which $\sigma_{s_i+1}\ne s_i$.
Since all earlier boundary positions are fixed, $A_{s_i}$ must be connected to $B_{s_{i+1}}$.
The only available source for $B_{s_i}$ is then $A_{s_{i+1}}$.
Thus, for $i<j$, we have $\sigma_{s_i+1}=s_{i+1}$ and $\sigma_{s_{i+1}+1}=s_i$.
Removing this pair and repeating the argument for the remaining indices gives disjoint adjacent transpositions.
The corresponding path forms are determined by Lemma~\ref{lem:free-1}, proving (i)--(iii).

Conversely, choose any collection of disjoint adjacent transpositions of the ordered list $(s_2,\ldots,s_j)$.
Construct the boundary paths using (i)--(iii) and the remaining paths using Lemma~\ref{lem:psigma-1}.
All prescribed steps are allowed.
A path from $A_t$ to $B_x$ requires $t+x-1$ up steps.
For a path in the first block,
\[
t+x-1=s_1-1=M_1.
\]
For $i\ge2$, the interior pairing and any chosen transposition $(s_{i-1},s_i)$ give
\[
t+x-1=s_{i-1}+s_i-1=M_i.
\]
For a fixed boundary path $P_{s_i}$, we have
$t=x=s_i$, so
\[
t+x-1=2s_i-1=N_i.
\]
These counts agree with the prescribed initial climbs.
Thus all constructed paths have the required endpoints, and each is uniquely determined by its source and step word.

It remains to verify non-intersection.
Every ordinary vertex lies on an initial climb or a final descent.
The endpoints of each short down step lie on these two segments, and its interior contains no path vertices.
The supporting lines of the climb from $A_t$ and the descent to $B_x$ are, respectively, $u-h=-2t $ and $ u+h=2x-2.$
Hence distinct climbs are disjoint, as are distinct final descents.

Suppose that two distinct paths share a vertex.
Since sources and sinks are distinct graph vertices, the common vertex lies on a climb from $A_t$ and a final descent to $B_x$.
The ordinary vertices of final descents occur only at the levels $M_1,\ldots,M_j$.
Thus, for some $1\le m\le j$,
\[
M_m-2t=2x-2-M_m,
\]
which gives
\[
t+x=M_m+1=s_{m-1}+s_m.
\]
We distinguish two cases.

\medskip
\noindent\textbf{Case 1.} Suppose that $m=1$.
Then $t+x=s_1$, so $0\le t,x\le s_1$.
By Lemma~\ref{lem:psigma-1}(i), the first block connects $A_t$ to $B_{s_1-t}=B_x$.

\medskip
\noindent\textbf{Case 2.} Suppose that $m\ge2$.
The prescribed path forms and $N_{m-1}<M_m$ show that any path reaching height $M_m$ has both its source and sink indices at least $s_{m-1}$.
Equality can occur only in case (ii) with $i=m$.
Together with $t+x=s_{m-1}+s_m$, this gives
\[
s_{m-1}\le t,x\le s_m.
\]
If either index equals $s_{m-1}$, the other equals $s_m$, and the chosen transposition connects $A_t$ to $B_x$.
Otherwise,
\[
s_{m-1}+1\le t,x\le s_m-1,
\]
and Lemma~\ref{lem:psigma-1}(ii) gives
\[
P_t=(A_t,B_{s_{m-1}+s_m-t})=(A_t,B_x).
\]

In both cases, the initial climb from $A_t$ and the final descent to $B_x$ belong to the same constructed path.
This contradicts the assumption that the common vertex belongs to two distinct paths.
The constructed family is therefore non-intersecting, which proves the converse.
\end{proof}

\subsection{Proof of Theorem~\ref{th:backshiftHDet}}

For each family $(P_0,\ldots,P_{s_j})\in\Omega_{s_j+1}^{(-1)}$ with associated permutation $\sigma$, let
\[
\zeta=(\sigma_{s_2+1},\ldots,\sigma_{s_j+1}),
\]
which is empty when $j=1$.
We define the signed weight sum
\[
\mathcal{X}_j
:=
\sum_{(P_0,\ldots,P_{s_j})\in\Omega_{s_j+1}^{(-1)}}
(-1)^{\mathrm{inv}(\zeta)}
\prod_{i=2}^{j}w(P_{s_i}),
\]
so that $\mathcal{X}_1=1$, and we set $\mathcal{X}_0=0$.

\begin{Lemma}\label{lem:Xj-1}
For each $j\ge1$,
\[
\mathcal{X}_j
=
\left(\prod_{h=1}^{j-1}V_h\right)E_j
=
v_0^{j-1}v_1^{j-1}v_2^{j-2}\cdots v_{j-1}\,E_j,
\]
where $(E_n)_{n\ge0}$ is the sequence of Theorem~\ref{th:backshiftHDet}.
\end{Lemma}

\begin{proof}
For $j=1$, Lemma~\ref{lem:psi-11} gives a unique family in $\Omega_{s_1+1}^{(-1)}$.
Since $\zeta$ is empty, we have $\mathcal{X}_1=1=E_1$.
Let $j\ge2$.
We denote by $\mathcal{X}_j^{\mathrm{fix}}$ the contribution to $\mathcal{X}_j$ from families satisfying $\sigma_{s_j+1}=s_j$.
For $j\ge3$, we denote by $\mathcal{X}_j^{\mathrm{swap}}$ the contribution from families satisfying $\sigma_{s_j+1}=s_{j-1}$ and $\sigma_{s_{j-1}+1}=s_j$.
We set $\mathcal{X}_2^{\mathrm{swap}}=0$.
By Lemma~\ref{lem:psi-11}, these are the only possibilities, and hence
\[
\mathcal{X}_j
=
\mathcal{X}_j^{\mathrm{fix}}
+
\mathcal{X}_j^{\mathrm{swap}}.
\]
We evaluate these two contributions separately.

\medskip

\noindent\textbf{Case 1.}
Suppose that $\sigma_{s_j+1}=s_j$. 
By Lemma~\ref{lem:psi-11}, deleting the paths with source indices $s_{j-1}+1,\ldots,s_j$ defines a bijection between the families in this case and $\Omega_{s_{j-1}+1}^{(-1)}$.
Among the deleted paths, only $P_{s_j}$ contributes to the boundary weight product.
Its step word is $b^{N_j}c_jv^{j-1}d$, so $w(P_{s_j})=c_{j,k_{j-1}}V_{j-1}.$
The last entry of $\zeta$ is its largest entry, so deleting it does not change the inversion number.
Thus the contribution of this case is
\begin{equation}\label{eq:xp1}
\mathcal{X}_j^{\mathrm{fix}}
=
V_{j-1}c_{j,k_{j-1}}\mathcal{X}_{j-1}.
\end{equation}

\medskip

\noindent\textbf{Case 2.}
Suppose that $j\ge3$, $\sigma_{s_j+1}=s_{j-1}$ and $\sigma_{s_{j-1}+1}=s_j$.
By Lemma~\ref{lem:psi-11}, deleting the paths with source indices $s_{j-2}+1,\ldots,s_j$ defines a bijection between the families in this case and $\Omega_{s_{j-2}+1}^{(-1)}$.
The two deleted boundary paths both have step word $b^{M_j}v^{j-1}d$, and hence
\[
w(P_{s_{j-1}})=w(P_{s_j})=V_{j-1}.
\]
The last two entries of $\zeta$ are $s_j,s_{j-1}$.
They form one inversion and are both larger than all preceding entries.
Thus deleting them removes exactly one inversion, and the contribution of this case is
\begin{equation}\label{eq:xp2}
\mathcal{X}_j^{\mathrm{swap}}
=
-V_{j-1}^{\,2}\mathcal{X}_{j-2}.
\end{equation}

\medskip

For $j=2$, only Case~1 occurs, and \eqref{eq:xp2} is also valid since $\mathcal{X}_0=0$.
Combining \eqref{eq:xp1} and \eqref{eq:xp2}, we obtain, for every $j\ge2$,
\[
\mathcal{X}_j = \mathcal{X}_j^{\mathrm{fix}}
+
\mathcal{X}_j^{\mathrm{swap}}
=
V_{j-1}c_{j,k_{j-1}}\mathcal{X}_{j-1}
-
V_{j-1}^{\,2}\mathcal{X}_{j-2}.
\]
Together with $\mathcal{X}_0=0=E_0$ and $\mathcal{X}_1=1=E_1$, this proves the desired formula by induction.
Indeed, using $V_{j-1}=V_{j-2}v_{j-1}$, the induction hypothesis gives
\[
\mathcal{X}_j
=
\left(\prod_{h=1}^{j-1}V_h\right)
\left(c_{j,k_{j-1}}E_{j-1}-v_{j-1}E_{j-2}\right)
=
\left(\prod_{h=1}^{j-1}V_h\right)E_j,
\]
where the last equality follows from \eqref{def:En}. 
\end{proof}

Now we are ready to prove Theorem~\ref{th:backshiftHDet}.

\begin{proof}[Proof of Theorem~\ref{th:backshiftHDet}]

\noindent{\it{Proof of  \eqref{eq:Hsj(-1)}. }}
Fix $j\ge1$ and let $n=s_j+1$.
By Lemma~\ref{lem:psigma-1}(i), the first block is the reversal of $0,1,\ldots,s_1$, with
\[
\mathrm{inv}(\sigma^1)=\binom{k_0+2}{2}
\qquad \text{and} \qquad
\mathcal{Q}_{0,s_1}=v_0^{k_0+2},
\]
since each of its $k_0+2$ paths has step word $b^{M_1}d$ and weight $v_0$.
For $2\le i\le j$, Lemma~\ref{lem:psigma-1}(ii) gives
\begin{equation}\label{eq:invsi-1}
\mathrm{inv}(\sigma^i)
=
\frac{k_{i-1}(k_{i-1}-1)}{2}
\qquad \text{and} \qquad
\mathcal{Q}_{s_{i-1}+1,s_i-1}
=
V_{i-1}^{k_{i-1}}.
\end{equation}
These formulas also hold when $k_{i-1}=0$, with the empty block having no inversions and the empty product equal to $1$.

The value ranges of the blocks $\sigma^1,\ldots,\sigma^j$ are ordered as their positions, so entries in different blocks form no inversions.
By Lemma~\ref{lem:psi-11}, a fixed boundary entry $\sigma_{s_i+1}=s_i$ forms no inversions with the block entries.
For an exchanged pair $\sigma_{s_{i-1}+1}=s_i$, $\sigma_{s_i+1}=s_{i-1}$, the intervening entries are precisely those of $\sigma^i$, each of which forms two inversions with the pair; the other block entries form no inversions with the pair, since those before position $s_{i-1}+1$ are smaller than $s_{i-1}$ and those after position $s_i+1$ are larger than $s_i$.
Consequently,
\[
(-1)^{\mathrm{inv}(\sigma)}
=
(-1)^{\mathrm{inv}(\zeta)}
\prod_{i=1}^{j}(-1)^{\mathrm{inv}(\sigma^i)}.
\]
Moreover,
\[
\binom{k_0+2}{2}=\frac{k_0(k_0-1)}{2}+2k_0+1
\equiv 1+\frac{k_0(k_0-1)}{2}\pmod 2 .
\]
Thus Theorem~\ref{th:hankel(-1)}, \eqref{eq:invsi-1}, the definition of $\mathcal{X}_j$, and Lemma~\ref{lem:Xj-1} yield
\begin{align*}
H_{s_j+1}^{(-1)}(F)
&=
\sum_{(P_0,\ldots,P_{s_j})\in\Omega_{s_j+1}^{(-1)}}
(-1)^{\mathrm{inv}(\sigma)}
\prod_{t=0}^{s_j}w(P_t)\\
&=
(-1)^{\mathrm{inv}(\sigma^1)}\mathcal{Q}_{0,s_1}
\left(
\prod_{i=2}^{j}
(-1)^{\mathrm{inv}(\sigma^i)}
\mathcal{Q}_{s_{i-1}+1,s_i-1}
\right)\mathcal{X}_j\\
&=
(-1)^{1+\sum_{i=0}^{j-1}\frac{k_i(k_i-1)}{2}}\,
v_0^{k_0+2}
\left(\prod_{i=1}^{j-1}V_i^{k_i+1}\right)E_j\\
&=
(-1)^{\epsilon_j-s_j+j+1}
v_0^{s_j+1}v_1^{s_j-s_1}\cdots
v_{j-1}^{s_j-s_{j-1}}E_j.
\end{align*}
Here we used
\[
\sum_{i=0}^{j-1}\frac{k_i(k_i-1)}{2}
=\epsilon_j-s_j+j.
\]
The exponent of $v_0$ is
$k_0+2+\sum_{i=1}^{j-1}(k_i+1)=s_j+1$,
and for $1\le m\le j-1$, the exponent of $v_m$ is
$\sum_{i=m}^{j-1}(k_i+1)=s_j-s_m$.
This proves \eqref{eq:Hsj(-1)}.

\medskip

\noindent{\it{Proof of  \eqref{eq:Hs_j+1(-1)}. }}
Let $j\ge1$.
If $k_j=0$, then $s_{j+1}=s_j+1$ and \eqref{eq:Hs_j+1(-1)} is immediate.
Suppose that $k_j\ge1$, and set $n=s_{j+1}$.
By Lemma~\ref{lem:psigma-1}(ii), applied with $i=j+1$, the last $k_j$ paths of any family in $\Omega_n^{(-1)}$ are prescribed.
For $1\le\ell\le k_j$, they are
\[
P_{s_j+\ell}=(A_{s_j+\ell},B_{s_{j+1}-\ell})=b^{M_{j+1}}v^jd.
\]
Their sink indices are precisely $s_j+1,\ldots,s_{j+1}-1$, so deleting these paths leaves a family in $\Omega_{s_j+1}^{(-1)}$.

Conversely, append these paths to any family in $\Omega_{s_j+1}^{(-1)}$.
As in the proof of Lemma~\ref{lem:psi-11}, a common ordinary vertex of two distinct paths
must lie on an initial climb from some $A_t$ and a final descent to some $B_x$, at a height $M_m$ with $1\le m\le j+1$.
Comparing its abscissae gives
\[
t+x=M_m+1=s_{m-1}+s_m.
\]
For $m\le j$, the prescribed path forms imply $t,x\ge s_{m-1}$, and hence $t,x\le s_m\le s_j. $
Thus such a vertex cannot involve an added path, whose source and sink indices are both at least $s_j+1$.

Every original path has height at most $N_j<M_{j+1}$.
At height $M_{j+1}$, only the added paths occur, and the equality $x=s_j+s_{j+1}-t$ pairs the source and sink of the same added path.
Hence no intersection is introduced.
Deletion and extension are therefore inverse operations, giving a bijection between
$\Omega_{s_j+1}^{(-1)}$ and $\Omega_{s_{j+1}}^{(-1)}$.

The appended block
\[
\sigma^{j+1}
=
(s_{j+1}-1)(s_{j+1}-2)\cdots(s_j+1)
\]
has $k_j$ entries and inversion number $\frac{k_j(k_j-1)}{2}$.
Its entries exceed all preceding entries, so there are no inversions between the two parts.
Each added path has weight $V_j$.
By Theorem~\ref{th:hankel(-1)}, we therefore obtain
\[
H_{s_{j+1}}^{(-1)}(F)
=
(-1)^{k_j(k_j-1)/2}
V_j^{k_j}H_{s_j+1}^{(-1)}(F),
\]
which proves \eqref{eq:Hs_j+1(-1)}.

\medskip
\noindent\textit{Proof of the vanishing assertion.}
For $1\le n\le s_1$,
Lemma~\ref{lem:psigma-1}(iii) and Theorem~\ref{th:hankel(-1)} give $H_n^{(-1)}(F)=0$.
Every integer $n>s_1$ can be written uniquely as $n=s_j+\ell$, where $j\ge1$ and $1\le\ell\le k_j+1$.
If $2\le\ell\le k_j$, then
\[
s_j+1\le n-1<s_{j+1}-1.
\]
Lemma~\ref{lem:psigma-1}(iii), applied with $i=j+1$, gives $\Omega_n^{(-1)}=\emptyset$, and hence $H_n^{(-1)}(F)=0$ by Theorem~\ref{th:hankel(-1)}.
The remaining indices are $n=s_j+1$ and $n=s_{j+1}$, which are covered by \eqref{eq:Hsj(-1)} and \eqref{eq:Hs_j+1(-1)}.
This completes the proof.
\end{proof}


\section{Application: Hankel determinants of the Euler numbers}\label{sec:euler}

In \cite{Han2020Trans} the first author evaluated the Hankel  determinants of the Euler numbers $(E_0,E_1,\ldots)$ and of  the shifted Euler numbers $(E_2,E_3,\ldots)$ by establishing  explicitly the $H$-fractions of their generating functions. 
Neither of these two $H$-fractions reduces to a $J$-fraction: 
both sequences have vanishing Hankel determinants, both expansions have a three-term periodic pattern, and their proofs are long. 

In this section we give short proofs of the same two evaluations. 
The only continued fraction we use is the $H$-fraction of the  sequence $(E_1,E_2,E_3,\ldots)$ lying in between, which is in  fact a $J$-fraction and is much simpler. 
Since $(E_2,E_3,E_4,\ldots)$ is the $1$-shift of $(E_1,E_2,E_3,\ldots)$, and $(E_0,E_1,E_2,\ldots)$ is obtained  from it by prepending $E_0=1$, the two evaluations follow from our main theorems, namely the shift Theorem \ref{th:shiftHDet} and the backward shift Theorem \ref{th:backshiftHDet},  by routine computations.

To avoid confusion with the Euler numbers, throughout this section we denote the auxiliary sequence $(E_n)_{n\ge0}$ of Theorem~\ref{th:backshiftHDet} by $(T_n)_{n\ge0}$.
Thus $T_0=0$, $T_1=1$, and, for $n\ge1$,
\[
T_{n+1}=c_{n+1,k_n}T_n-v_nT_{n-1}.
\]

\subsection{The Euler numbers and the continued fraction (F7)}\label{subsec:eu:F7}

The {\it Euler numbers} $E_n\ (n\geq 0)$ are defined by the exponential generating function
\begin{equation*}
\tan(x)+\sec(x)=\sum_{n\geq 0} E_n \frac{x^n}{n!}.
\end{equation*}
The first values are
\begin{equation*}
(E_0, E_1, E_2, E_3, E_4,\ldots)=(1, 1, 1, 2, 5, 16, 61, 272, 1385, \ldots).
\end{equation*}
As shown by Andr\'e~\cite{Andre1879} as early as 1879, the Euler numbers count the alternating permutations, and they appear in many combinatorial contexts; 
see, for example, \cite{Andre1879, Euler1755, Foata2008Han, Viennot1982}.
Continued fractions and Hankel determinants of the even or of the odd Euler numbers have a long history going back to Stieltjes \cite{Stieltjes1894Re, Stieltjes1889Red}; 
see also \cite[p.~206 (53.11), p.~369]{Wall1948} and \cite{Flajolet1980}. 
In \cite{Han2020Trans} the first author collected a list of continued fractions (F1)--(F11), and of the corresponding Hankel determinant evaluations (H1)--(H11), for various sequences built from the Euler numbers. 
Among them, the following {\it Jacobi} continued fraction expansion of the ordinary generating function of the sequence $(E_{n+1})_{n\geq 0}$, labeled (F7) in \cite{Han2020Trans}, is particularly simple: 
\begin{equation}\label{eq:eu:F7}
F(x):=\sum_{n\geq 0} E_{n+1} x^n =
\cFrac{1}{1-x}
- \cFrac{x^2}{1-2x}
- \cFrac{3x^2}{1-3x}
- \cFrac{6x^2}{1-4x}
-\cdots
\end{equation}
with general patterns
\begin{equation*}
a_1=1, \quad
a_k =-\binom{k}{2}x^2; \qquad
b_0=0, \quad
b_k=1-kx.
\end{equation*}
Comparing \eqref{eq:eu:F7} with the general form~\eqref{eq:Hfrac} of an $H$-fraction, we read off the following data term by term:
\begin{align*}
k_j &= 0 \qquad \text{for every $j\geq 0$,} \\
v_0 &= 1, \qquad v_j = \binom{j+1}{2} = \frac{j(j+1)}{2}
\qquad \text{for every $j\geq 1$,} \\
u_j(x) &= j \qquad (\text{a constant polynomial})
\qquad \text{for every $j\geq 1$.}
\end{align*}
In particular, since all $k_j=0$, the continued fraction~\eqref{eq:eu:F7} is a $J$-fraction, so that
$s_j=j$ and $\epsilon_j=0$ for every $j\geq 0$, and $c_{j,k_{j-1}}=c_{j,0}=j$ for every $j\geq 1$. 
By Theorem \ref{th:super2}, or directly by Heilermann's formula 
\eqref{eq:detH}, 
\begin{equation}\label{eq:eu:HnF}
H_n(F)
= v_0^n v_1^{n-1} \cdots v_{n-1}
= \prod_{\ell=1}^{n-1} \binom{\ell+1}{2}^{n-\ell}
= \frac{n!}{2^{n(n-1)/2}}\prod_{\ell=2}^{n-1}(\ell!)^2,
\end{equation}
which is precisely the evaluation (H7) of \cite{Han2020Trans} of the Hankel determinants of $(E_1,E_2,E_3,\ldots)$. 
The sequences $(R_n)$ and $(T_n)$ of Theorems \ref{th:shiftHDet} and \ref{th:backshiftHDet} are defined by $R_0=1$, $R_1=c_{1,0}=1$, $T_0=0$, $T_1=1$, 
together with the common recurrence, for $X=R,T$ and $n\ge1$,
\begin{equation}\label{eq:eu:RnRec}
X_{n+1} = c_{n+1,0}\, X_n - v_n\, X_{n-1}
= (n+1)\, X_n - \binom{n+1}{2}\, X_{n-1}.
\end{equation}
Their first values are
\begin{equation}\label{eq:eu:RTvalues}
\begin{array}{*{13}c}
n   &=& 0 & 1 & 2 & 3 & 4 & 5 & 6 & 7 & 8 & 9 & 10\\
R_n &=& 1 & 1 & 1 & 0 & -6 & -30 & -90 & 0 & 2520 & 22680 & 113400\\
T_n &=& 0 & 1 & 2 & 3 & 0 & -30 & -180 & -630 & 0 & 22680 & 226800.
\end{array}
\end{equation}

\smallskip

In the computations below we shall use the following elementary identities, in which empty products are equal to $1$. 
For every $m\geq 0$, 
\begin{equation}\label{eq:eu:prodfact}
\begin{aligned}
\prod_{\ell=2}^{2m+1}(\ell!)^2
&= \prod_{i=1}^{m}(2i+1)^2 \prod_{i=1}^{m}((2i)!)^4,\\
\prod_{j=1}^{m}(2j+1)!^4 &= \prod_{j=1}^{m}(2j+1)^4 \prod_{j=1}^{m}((2j)!)^4,
\end{aligned}
\end{equation}
the first one being obtained by splitting the product into even and odd factorials and writing $(2i+1)!^2=(2i+1)^2(2i)!^2$. 
Moreover, for every $m\ge0$, the classical identity relating factorials and double factorials is
\begin{equation*}\label{eq:eu:doublefact}
\prod_{i=1}^{m}(2i+1) = (2m+1)!! = \frac{(2m+1)!}{2^m\, m!}.
\end{equation*}
Taking $m=2k-1$ for $k\ge1$ and $m=2k$ for $k\ge0$, respectively, gives
\begin{align}
(4k)! &= 2^{2k}\, (2k)! \prod_{i=1}^{2k-1}(2i+1), \label{eq:eu:key4k}\\
(4k+2)! &= 2^{2k+1}\,(2k+1)! \prod_{i=1}^{2k}(2i+1), \label{eq:eu:key4k2}
\end{align}
since $(4k)!=4k\cdot(4k-1)!$, $2^{2k}(2k)!=2\cdot 2k\cdot 2^{2k-1}(2k-1)!$, 
$(4k+2)!=(4k+2)\,(4k+1)!$ and $2^{2k+1}(2k+1)!=(4k+2)\cdot 2^{2k}(2k)!$.

\subsection{The shifted sequence $(E_2,E_3,E_4,\ldots)$}\label{subsec:eu:H10}

The Hankel determinants of the shifted sequence $(E_2,E_3,E_4,\ldots)$ are exactly the $1$-shifted Hankel determinants $H_n^{(1)}(F)$ of the series $F(x)$ of \eqref{eq:eu:F7}. 
In \cite{Han2020Trans} they were derived from the $H$-fraction of $\sum_{n\geq 0}E_{n+2}x^n$, labeled (F10) there, namely
\begin{multline}\label{eq:eu:F10}
\sum_{n\geq 0} E_{n+2} x^n
= \cFrac{1}{1-2x} - \cFrac{x^2}{1-4x}
- \cFrac{18x^3}{1-4x-16x^2} - \cFrac{50x^3}{1-8x} \\
- \cFrac{9x^2}{1-10x}
- \cFrac{196x^3}{1-8x-64x^2} - \cFrac{324x^3}{1-14x}
- \cdots
\end{multline}
with general patterns
\begin{align*}
a_{1}&=1,                    &  b_{0}&=0,                            \\
a_{3k+0} &=-2(4k - 1)^2kx^3, &  b_{3k+0} &= -16k^2x^2 - 4kx + 1,    \\
a_{3k+1} &=-2(4k + 1)^2kx^3, &  b_{3k+1} &= -2(3k + 1)x + 1,        \\
a_{3k+2} &=-(2k + 1)^2x^2;   &  b_{3k+2} &= -2(3k + 2)x + 1,
\end{align*}
where $a_{3k}$, $a_{3k+1}$ and $b_{3k}$ are given for $k\geq1$, and $a_{3k+2}$, $b_{3k+1}$ and $b_{3k+2}$ for $k\geq0$.
Unlike \eqref{eq:eu:F7}, the continued fraction \eqref{eq:eu:F10} is a genuine $H$-fraction (it does not reduce to a $J$-fraction), with a three-term periodic pattern, and its proof in \cite{Han2020Trans} is much longer than that of \eqref{eq:eu:F7}.

\begin{Theorem}[{\cite[(H10)]{Han2020Trans}}]\label{th:eu:H10}
The Hankel determinants $H_n^{(1)}(F)$ of the shifted Euler numbers $(E_2, E_3, E_4, E_5, \ldots)$ are given, for every $k\geq 0$, by
\begin{equation}\label{eq:eu:H10}
\left\{
\begin{aligned}
H_{4k}^{(1)}(F)&= \frac{ (-1)^{k}k^2(2k-1)!^2 }{ 2^{8k^2-4k-2} }
\prod_{j=1}^{2k-1}(2j+1)!^4 \qquad (k\geq 1), \\
H_{4k+1}^{(1)}(F)&= \frac{ (-1)^{k}(2k)!^2 }{ 2^{8k^2} (4k+1)!^2 }
\prod_{j=1}^{2k}(2j+1)!^4, \\
H_{4k+2}^{(1)}(F)&= \frac{ (-1)^{k}(2k+1)!^2 }{ 2^{8k^2+4k} }
\prod_{j=1}^{2k}(2j+1)!^4, \\
H_{4k+3}^{(1)}(F)&=0,
\end{aligned}
\right.
\end{equation}
together with $H_0^{(1)}(F)=1$.
\end{Theorem}

The first values are $(1, 1, 1, 0, -324, -291600, -2361960000, 0, \ldots)$.
We now derive Theorem~\ref{th:eu:H10} directly from \eqref{eq:eu:F7} and Theorem~\ref{th:shiftHDet}, without using \eqref{eq:eu:F10}.
Since \eqref{eq:eu:F7} is a $J$-fraction, the $J$-fraction form \eqref{eq:H1classical} of Theorem \ref{th:shiftHDet} reads
\begin{equation}\label{eq:eu:MasterR}
H_n^{(1)}(F) = H_n(F)\, R_n \qquad (n\geq 0).
\end{equation}
Since $H_n(F)$ is given by \eqref{eq:eu:HnF}, it remains to find  a closed form for $R_n$. 

\begin{Lemma}\label{lem:eu:Rn}
For every $k\geq 0$,
\begin{equation}\label{eq:eu:Rn-closed}
\begin{aligned}
R_{4k}   &= \frac{(-1)^k (4k)!}{4^k}, &
R_{4k+1} &= \frac{(-1)^k (4k+1)!}{4^k}, \\
R_{4k+2} &= \frac{(-1)^k (4k+2)!}{2\cdot 4^k}, &
R_{4k+3} &= 0.
\end{aligned}
\end{equation}
\end{Lemma}

\begin{proof}
We proceed by induction on $k$. The base case $k=0$ follows from $R_0=R_1=R_2=1$ and $R_3=0$. Assume that~\eqref{eq:eu:Rn-closed} holds for some $k\geq 0$. Using the recurrence~\eqref{eq:eu:RnRec} and the induction hypothesis, we compute
\begin{align*}
R_{4k+4}
&= (4k+4)\,R_{4k+3}-\binom{4k+4}{2} R_{4k+2} \\
&= -\frac{(4k+4)(4k+3)}{2}\cdot\frac{(-1)^k(4k+2)!}{2\cdot 4^k}
=  \frac{(-1)^{k+1}(4k+4)!}{4^{k+1}}, \\[4pt]
R_{4k+5}
&= (4k+5)\,R_{4k+4}-\binom{4k+5}{2} R_{4k+3} \\
&= (4k+5)\cdot\frac{(-1)^{k+1}(4k+4)!}{4^{k+1}}
=  \frac{(-1)^{k+1}(4k+5)!}{4^{k+1}}, \\[4pt]
R_{4k+6}
&= (4k+6)\,R_{4k+5}-\binom{4k+6}{2} R_{4k+4} \\
&= \frac{(-1)^{k+1}(4k+6)!}{4^{k+1}}
   -\frac{(-1)^{k+1}(4k+6)!}{2\cdot 4^{k+1}}
=  \frac{(-1)^{k+1}(4k+6)!}{2\cdot 4^{k+1}}, \\[4pt]
R_{4k+7}
&= (4k+7)\,R_{4k+6}-\binom{4k+7}{2} R_{4k+5} \\
&= \frac{(-1)^{k+1}(4k+7)!}{2\cdot 4^{k+1}}
   - \frac{(-1)^{k+1}(4k+7)!}{2\cdot 4^{k+1}}
= 0,
\end{align*}
which is \eqref{eq:eu:Rn-closed} for $k+1$. This completes the induction.
\end{proof}

With these preparations, we now prove Theorem~\ref{th:eu:H10}.

\begin{proof}[Proof of Theorem \ref{th:eu:H10}]
We substitute \eqref{eq:eu:Rn-closed} and \eqref{eq:eu:HnF} into \eqref{eq:eu:MasterR} and compare the resulting expressions with \eqref{eq:eu:H10}.
The case $n=0$ follows from the convention $H_0^{(1)}(F)=1$.
For $n\ge1$, we distinguish four cases according to the residue of $n$ modulo $4$.
All empty products below are understood to be $1$.

\medskip

\noindent\textbf{Case 1.}
Suppose that $n=4k+3$ with $k\ge0$.
Since $R_{4k+3}=0$, we obtain $H_{4k+3}^{(1)}(F)=0$, in agreement with~\eqref{eq:eu:H10}.

\medskip
\noindent\textbf{Case 2.}
Suppose that $n=4k$ with $k\ge1$.
By \eqref{eq:eu:HnF} and \eqref{eq:eu:Rn-closed}, we obtain
\begin{equation}\label{eq:eu:H4k-LHS}
H_{4k}^{(1)}(F)
= \frac{(4k)!}{2^{2k(4k-1)}}\prod_{\ell=2}^{4k-1}(\ell!)^2
\cdot \frac{(-1)^k (4k)!}{4^k}
= \frac{(-1)^k (4k)!^2}{2^{8k^2}} \prod_{\ell=2}^{4k-1}(\ell!)^2.
\end{equation}
Applying \eqref{eq:eu:prodfact} with $m=2k-1$ gives
\[
\prod_{\ell=2}^{4k-1}(\ell!)^2
=
\left(\prod_{i=1}^{2k-1}((2i)!)^4\right)
\left(\prod_{i=1}^{2k-1}(2i+1)^2\right).
\]
The odd-factorial product in the corresponding expression in \eqref{eq:eu:H10} satisfies
\[
\prod_{i=1}^{2k-1}((2i+1)!)^4
=
\left(\prod_{i=1}^{2k-1}((2i)!)^4\right)
\left(\prod_{i=1}^{2k-1}(2i+1)^4\right).
\]
Thus, after cancelling $(-1)^k$ and $\prod_{i=1}^{2k-1}((2i)!)^4$, the equality between
\eqref{eq:eu:H4k-LHS} and the corresponding expression in \eqref{eq:eu:H10} reduces to
\[
\frac{(4k)!^2}{2^{8k^2}}
\prod_{i=1}^{2k-1}(2i+1)^2
=
\frac{k^2(2k-1)!^2}{2^{8k^2-4k-2}}
\prod_{i=1}^{2k-1}(2i+1)^4.
\]
Multiplying by $2^{8k^2}$ and cancelling
$\prod_{i=1}^{2k-1}(2i+1)^2$, we obtain
\[
(4k)!^2
=
2^{4k+2}k^2(2k-1)!^2
\prod_{i=1}^{2k-1}(2i+1)^2.
\]
Since $2k(2k-1)!=(2k)!$, this is precisely the square of \eqref{eq:eu:key4k}.

\medskip
\noindent\textbf{Case 3.}
Suppose that $n=4k+1$ with $k\ge0$.
By \eqref{eq:eu:HnF} and \eqref{eq:eu:Rn-closed}, we obtain
\begin{equation}\label{eq:eu:H4k1-LHS}
\begin{aligned}
H_{4k+1}^{(1)}(F)
&=
\frac{(4k+1)!}{2^{2k(4k+1)}}
\left(\prod_{\ell=2}^{4k}(\ell!)^2\right)
\frac{(-1)^k(4k+1)!}{4^k}\\
&=
\frac{(-1)^k(4k+1)!^2(4k)!^2}
{2^{8k^2+4k}}
\prod_{\ell=2}^{4k-1}(\ell!)^2.
\end{aligned}
\end{equation}
For $k=0$, both \eqref{eq:eu:H4k1-LHS} and \eqref{eq:eu:H10} give $H_1^{(1)}(F)=1=E_2$.
For $k \ge 1$, we apply \eqref{eq:eu:prodfact} with $m=2k-1$ to the product in \eqref{eq:eu:H4k1-LHS}.
For the corresponding expression in \eqref{eq:eu:H10}, we separate the last factor of the odd-factorial product and obtain
\begin{align*}
\prod_{i=1}^{2k}((2i+1)!)^4
&=
(4k+1)!^4
\prod_{i=1}^{2k-1}((2i+1)!)^4\\
&=
(4k+1)^4(4k)!^4
\left(\prod_{i=1}^{2k-1}((2i)!)^4\right)
\left(\prod_{i=1}^{2k-1}(2i+1)^4\right).
\end{align*}
After cancelling $(-1)^k$ and $\prod_{i=1}^{2k-1}((2i)!)^4$, the desired equality reduces to
\begin{multline*}
\frac{(4k+1)!^2(4k)!^2}{2^{8k^2+4k}}
\prod_{i=1}^{2k-1}(2i+1)^2
=
\frac{(2k)!^2(4k+1)^4(4k)!^4}
{2^{8k^2}(4k+1)!^2}
\prod_{i=1}^{2k-1}(2i+1)^4,
\end{multline*}
that is, since $(4k+1)!=(4k+1)\,(4k)!$, to
\[
(4k)!^2 = 2^{4k}\,(2k)!^2\prod_{i=1}^{2k-1}(2i+1)^2,
\]
which is the square of \eqref{eq:eu:key4k}.

\medskip
\noindent\textbf{Case 4.}
Suppose that $n=4k+2$ with $k\ge0$.
By \eqref{eq:eu:HnF} and \eqref{eq:eu:Rn-closed}, we obtain
\begin{equation}\label{eq:eu:H4k2-LHS}
\begin{aligned}
H_{4k+2}^{(1)}(F)
&=
\frac{(4k+2)!}{2^{(2k+1)(4k+1)}}
\left(\prod_{\ell=2}^{4k+1}(\ell!)^2\right)
\frac{(-1)^k(4k+2)!}{2\cdot4^k}\\
&=
\frac{(-1)^k(4k+2)!^2}{2^{8k^2+8k+2}}
\prod_{\ell=2}^{4k+1}(\ell!)^2.
\end{aligned}
\end{equation}
Applying \eqref{eq:eu:prodfact} with $m=2k$ gives
\[
\prod_{\ell=2}^{4k+1}(\ell!)^2
=
\left(\prod_{i=1}^{2k}((2i)!)^4\right)
\left(\prod_{i=1}^{2k}(2i+1)^2\right).
\]
The odd-factorial product in the corresponding expression in \eqref{eq:eu:H10} is expanded in the same way.
After cancelling $(-1)^k$ and $\prod_{i=1}^{2k}((2i)!)^4$, the equality between \eqref{eq:eu:H4k2-LHS} and that expression reduces to
\[
\frac{(4k+2)!^2}{2^{8k^2+8k+2}}
\prod_{i=1}^{2k}(2i+1)^2
=
\frac{(2k+1)!^2}{2^{8k^2+4k}}
\prod_{i=1}^{2k}(2i+1)^4.
\]
Multiplying by $2^{8k^2+8k+2}$ and cancelling
$\prod_{i=1}^{2k}(2i+1)^2$, we obtain
\[
(4k+2)!^2
=
2^{4k+2}(2k+1)!^2
\prod_{i=1}^{2k}(2i+1)^2.
\]
Since $2^{4k+2}=(2^{2k+1})^2$, this is the square of \eqref{eq:eu:key4k2}.
\medskip

Together with the case $n=0$, these four cases prove \eqref{eq:eu:H10}.
\end{proof}

\subsection{The sequence $(E_0,E_1,E_2,\ldots)$ of all Euler numbers}\label{subsec:eu:En}

Some Hankel determinants of the {\it mixed} Euler numbers $(E_0, E_1, E_2, \ldots)$ vanish, so that their $J$-fraction does  not exist. 
In~\cite{Han2020Trans} the first author derived the  $H$-fraction of the ordinary generating function of the Euler numbers \cite[Theorem 1.1]{Han2020Trans}, namely
\begin{align}
	\sum_{n\geq 0} E_n x^n
	&=
	 \cFrac{ 1}{ 1-x} 
	- \cFrac{ x^3}{1-2x-4x^2}
	- \cFrac{ 9x^3}{ 1-5x} \label{eq:eu:En_H}  \\
	&\qquad
	- \cFrac{ 4x^2}{ 1-7x} 
	- \cFrac{ 75x^3}{ 1-6x-36x^2}
	- \cFrac{ 147x^3}{1-11x} \nonumber \\
	&\qquad
	- \cFrac{ 16x^2}{1-13x} - \cdots \nonumber
\end{align}
with general patterns
\begin{alignat*}{2}
a_1 &= 1,
&\qquad a_{3k} &= -(4k-1)^2(2k-1)x^3,\\
a_{3k+1} &= -4k^2x^2,
&\qquad a_{3k+2} &= -(4k+1)^2(2k+1)x^3,\\[8pt]
b_0 &= 0, &\qquad b_{3k} &= 1-(6k-1)x,\\
b_{3k+1} &= 1-(6k+1)x,
&\qquad b_{3k+2} &= 1-2(2k+1)x-4(2k+1)^2x^2.
\end{alignat*}
where $a_{3k}$, $a_{3k+1}$ and $b_{3k}$ are given for $k\geq1$, and  $a_{3k+2}$, $b_{3k+1}$ and $b_{3k+2}$ for $k\geq0$, and deduced from it, by Theorem \ref{th:super2}, the following 
evaluation of the Hankel determinants of the Euler numbers.

\begin{Theorem}[{\cite[Theorem 1.2]{Han2020Trans}}]\label{th:eu:En_Det}
Let $H_n=H_n(E_0,E_1,E_2,\ldots)$ be the Hankel determinants of the Euler numbers. Then $H_0=1$ and, for every $k\geq 0$,
\begin{equation}\label{eq:eu:En_Det}
\left\{
\begin{aligned}
	H_{4k+1} &= (-1)^{k} \frac{(2k)!^2}{2^{4k(2k-1)}} \prod_{j=1}^{2k-1} (2j+1)!^4, \\
	H_{4k+2}&=0,\\
	H_{4k+3} &= (-1)^{k+1}\frac{(2k+1)!^2}{2^{4k(2k+1)}} \prod_{j=1}^{2k} (2j+1)!^4,\\
	H_{4k+4} &= (-1)^{k+1}\frac{(2k+1)!^2(4k+3)!^2}{2^{2(2k+1)^2}} \prod_{j=1}^{2k} (2j+1)!^4.
\end{aligned}
\right.
\end{equation}
\end{Theorem}

The first values are 
\begin{align*}
(H_0, H_1, H_2, \ldots)=(1, 1, 0, -1, -9, -324, 0,\ &2361960000, \\
&937461924000000,\ldots).
\end{align*}
The proof of \eqref{eq:eu:En_H} in \cite{Han2020Trans} is rather  long: it makes use of Flajolet's continued fraction for a  permutation statistic and of several combinatorial models for the Euler numbers, and the expansion itself, with its three-term periodic pattern, is considerably more involved than the classical  continued fractions of Stieltjes for the tangent and secant numbers.
We now derive Theorem \ref{th:eu:En_Det} from \eqref{eq:eu:F7},  bypassing \eqref{eq:eu:En_H} altogether.
The key observation is that the sequence $(E_0, E_1, E_2, \ldots)$ is obtained from the sequence $(E_1, E_2, E_3, \ldots)$ by prepending the term $E_0=1$; in other words, 
\begin{equation}\label{eq:eu:1plusxF}
\sum_{n\geq 0} E_n x^n = 1 + xF(x).
\end{equation}
By \eqref{eq:eu:1plusxF} and \eqref{eq:Hc} with $c=1$, for $n \ge 1$,
\begin{equation*}
H_n = H_n(1+xF(x)) 
= H_{n-1}^{(1)}(F) + H_{n}^{(-1)}(F).
\end{equation*}
Since \eqref{eq:eu:F7} is a $J$-fraction with $v_0=1$, \eqref{eq:eu:MasterR} and the $J$-fraction form \eqref{eq:Hm1classical} of Theorem \ref{th:backshiftHDet} give
\begin{align*}
    H_{n-1}^{(1)}(F)&=H_{n-1}(F)\,R_{n-1},\\
    H_{n}^{(-1)}(F)&=-H_{n-1}(F)\,T_{n-1}.
\end{align*}
For $n \ge 0$, setting 
\begin{equation}\label{eq:eu:Dn}
D_n:=R_n-T_n,
\end{equation}
we therefore obtain
\begin{equation}\label{eq:eu:Master}
H_n = H_{n-1}(F)\, D_{n-1} \qquad (n\geq 1).
\end{equation}
Since $(R_n)$ and $(T_n)$ satisfy the same linear recurrence \eqref{eq:eu:RnRec}, so does $(D_n)$, with the initial values $D_0=R_0-T_0=1$ and $D_1=R_1-T_1=0$; by \eqref{eq:eu:RTvalues}, its first values are
\begin{equation*}
\begin{array}{*{13}c}
n   &=& 0 & 1 & 2 & 3 & 4 & 5 & 6 & 7 & 8 & 9 & 10\\
D_n &=& 1 & 0 & -1 & -3 & -6 & 0 & 90 & 630 & 2520 & 0 & -113400.
\end{array}
\end{equation*}
Again, it remains to find a closed form for $D_n$.

\begin{Lemma}\label{lem:eu:Dn}
For every $k\geq 0$,
\begin{equation}\label{eq:eu:Dn-closed}
\begin{aligned}
D_{4k}   &= \frac{(-1)^k (4k)!}{4^k}, &
D_{4k+1} &= 0, \\
D_{4k+2} &= \frac{(-1)^{k+1} (4k+2)!}{2\cdot 4^k}, &
D_{4k+3} &= \frac{(-1)^{k+1} (4k+3)!}{2\cdot 4^k}.
\end{aligned}
\end{equation}
\end{Lemma}

\begin{proof}
We proceed by induction on $k$. The base case $k=0$ follows from $D_0=1$, $D_1=0$, $D_2=-1$ and $D_3=-3$. 
Assume that~\eqref{eq:eu:Dn-closed} holds for some $k\geq 0$. Using the recurrence~\eqref{eq:eu:RnRec} and the induction hypothesis, we compute
\begin{align*}
D_{4k+4}
&= (4k+4)\,D_{4k+3}-\binom{4k+4}{2} D_{4k+2} \\
&= \frac{(-1)^{k+1}(4k+4)!}{2\cdot 4^k}
   -\frac{(-1)^{k+1}(4k+4)!}{4\cdot 4^{k}}
=  \frac{(-1)^{k+1}(4k+4)!}{4^{k+1}}, \\[4pt]
D_{4k+5}
&= (4k+5)\,D_{4k+4}-\binom{4k+5}{2} D_{4k+3} \\
&= \frac{(-1)^{k+1}(4k+5)!}{4^{k+1}}
   -\frac{(-1)^{k+1}(4k+5)!}{4^{k+1}}
= 0, \\[4pt]
D_{4k+6}
&= (4k+6)\,D_{4k+5}-\binom{4k+6}{2} D_{4k+4} \\
&= -\frac{(4k+6)(4k+5)}{2}\cdot\frac{(-1)^{k+1}(4k+4)!}{4^{k+1}}
=  \frac{(-1)^{k+2}(4k+6)!}{2\cdot 4^{k+1}}, \\[4pt]
D_{4k+7}
&= (4k+7)\,D_{4k+6}-\binom{4k+7}{2} D_{4k+5} \\
&= (4k+7)\cdot\frac{(-1)^{k+2}(4k+6)!}{2\cdot 4^{k+1}}
=  \frac{(-1)^{k+2}(4k+7)!}{2\cdot 4^{k+1}},
\end{align*}
which is \eqref{eq:eu:Dn-closed} for $k+1$. This completes the induction.
\end{proof}

Since $T_n=R_n-D_n$, Lemmas \ref{lem:eu:Rn} and \ref{lem:eu:Dn} also give a closed form for $T_n$: 
$T_{4k}=0$, $T_{4k+1}=R_{4k+1}$, $T_{4k+2}=2R_{4k+2}$ and $T_{4k+3}=-D_{4k+3}$, in accordance with \eqref{eq:eu:RTvalues}.

\medskip
Now we are ready to prove Theorem~\ref{th:eu:En_Det}.

\begin{proof}[Proof of Theorem \ref{th:eu:En_Det}]
We insert the closed form~\eqref{eq:eu:Dn-closed} of $D_n$ and the formula~\eqref{eq:eu:HnF} for $H_n(F)$ into \eqref{eq:eu:Master}, and verify that the result agrees with \eqref{eq:eu:En_Det}. 
The case $n=0$ follows from the convention $H_0=1$.
For $n=1$, we have $H_1=E_0=H_0(F)D_0=1$, which agrees with the formula for $H_{4k+1}$ in \eqref{eq:eu:En_Det} at $k=0$.
All empty products below are understood to be $1$.
For $n\ge2$, we distinguish four cases according to the residue of $n$ modulo $4$.

\medskip

\noindent\textbf{Case 1.}
Suppose that $n=4k+2$ with $k\ge0$.
Since $D_{4k+1}=0$, \eqref{eq:eu:Master} gives $H_{4k+2}=0$, in agreement with \eqref{eq:eu:En_Det}.

\medskip
\noindent\textbf{Case 2.}
Suppose that $n=4k+1$ with $k\ge1$.
By \eqref{eq:eu:Master}, \eqref{eq:eu:HnF}, and \eqref{eq:eu:Dn-closed}, we obtain
\begin{equation}\label{eq:eu:H4k1E-LHS}
H_{4k+1}
= \frac{(4k)!}{2^{2k(4k-1)}}\prod_{\ell=2}^{4k-1}(\ell!)^2
\cdot \frac{(-1)^k (4k)!}{4^k}
= \frac{(-1)^k (4k)!^2}{2^{8k^2}} \prod_{\ell=2}^{4k-1}(\ell!)^2.
\end{equation}
On the other hand, the desired formula in~\eqref{eq:eu:En_Det} reads
\begin{equation}\label{eq:eu:H4k1E-RHS}
H_{4k+1} = \frac{(-1)^k (2k)!^2}{2^{8k^2-4k}}
\prod_{j=1}^{2k-1}(2j+1)!^4.
\end{equation}
By \eqref{eq:eu:prodfact} with $m=2k-1$, the equality between~\eqref{eq:eu:H4k1E-LHS} and~\eqref{eq:eu:H4k1E-RHS} reduces, after cancellation of $(-1)^k$ and of $\prod_{i=1}^{2k-1}((2i)!)^4$, to
\begin{equation*}
\frac{(4k)!^2}{2^{8k^2}}\prod_{i=1}^{2k-1}(2i+1)^2
= \frac{(2k)!^2}{2^{8k^2-4k}}\prod_{i=1}^{2k-1}(2i+1)^4.
\end{equation*}
Equivalently,
\[
(4k)!^2
=
2^{4k}(2k)!^2
\prod_{i=1}^{2k-1}(2i+1)^2,
\]
which is the square of \eqref{eq:eu:key4k}.

\medskip

\noindent\textbf{Case 3.}
Suppose that $n=4k+3$ with $k\ge0$.
By \eqref{eq:eu:Master}, \eqref{eq:eu:HnF}, and \eqref{eq:eu:Dn-closed}, we obtain
\begin{equation}\label{eq:eu:H4k3E-LHS}
\begin{aligned}
H_{4k+3}
&= \frac{(4k+2)!}{2^{(2k+1)(4k+1)}}\prod_{\ell=2}^{4k+1}(\ell!)^2
\cdot \frac{(-1)^{k+1}(4k+2)!}{2\cdot 4^k} \\
&= \frac{(-1)^{k+1}(4k+2)!^2}{2^{8k^2+8k+2}}\prod_{\ell=2}^{4k+1}(\ell!)^2,
\end{aligned}
\end{equation}
while the desired formula in~\eqref{eq:eu:En_Det} reads
\begin{equation}\label{eq:eu:H4k3E-RHS}
H_{4k+3} = \frac{(-1)^{k+1}(2k+1)!^2}{2^{8k^2+4k}}
\prod_{j=1}^{2k}(2j+1)!^4.
\end{equation}
By \eqref{eq:eu:prodfact} with $m=2k$, the equality between~\eqref{eq:eu:H4k3E-LHS} and~\eqref{eq:eu:H4k3E-RHS} reduces, after cancellation of $(-1)^{k+1}$ and of $\prod_{i=1}^{2k}((2i)!)^4$, to
\[
(4k+2)!^2 = 2^{4k+2}\,(2k+1)!^2 \prod_{i=1}^{2k}(2i+1)^2,
\]
that is, to the square of \eqref{eq:eu:key4k2}.

\medskip

\noindent\textbf{Case 4.}
Suppose that $n=4k+4$ with $k\ge0$.
By \eqref{eq:eu:Master}, \eqref{eq:eu:HnF}, and \eqref{eq:eu:Dn-closed}, we obtain
\begin{equation}\label{eq:eu:H4k4E-LHS}
\begin{aligned}
H_{4k+4}
&= \frac{(4k+3)!}{2^{(4k+3)(2k+1)}}\prod_{\ell=2}^{4k+2}(\ell!)^2
\cdot \frac{(-1)^{k+1}(4k+3)!}{2\cdot 4^k} \\
&= \frac{(-1)^{k+1}(4k+3)!^2\,(4k+2)!^2}{2^{8k^2+12k+4}}
\prod_{\ell=2}^{4k+1}(\ell!)^2,
\end{aligned}
\end{equation}
while the desired formula in~\eqref{eq:eu:En_Det} reads
\begin{equation}\label{eq:eu:H4k4E-RHS}
H_{4k+4} = \frac{(-1)^{k+1}(2k+1)!^2(4k+3)!^2}{2^{8k^2+8k+2}}
\prod_{j=1}^{2k}(2j+1)!^4.
\end{equation}
After cancelling $(-1)^{k+1}(4k+3)!^2$ and applying \eqref{eq:eu:prodfact} with $m=2k$, the equality between \eqref{eq:eu:H4k4E-LHS} and \eqref{eq:eu:H4k4E-RHS} reduces to the same factorial identity as in Case~3.
It therefore also follows from \eqref{eq:eu:key4k2}.

\medskip

Together with the cases $n=0$ and $n=1$, these four cases prove \eqref{eq:eu:En_Det}.
\end{proof}

\begin{Remark}\label{rem:eu:remarks}
The proofs above use the Jacobi continued fraction  \eqref{eq:eu:F7} of $\sum_{n\geq 0} E_{n+1}x^n$ only, together  with Theorems \ref{th:shiftHDet} and \ref{th:backshiftHDet} in  their $J$-fraction forms \eqref{eq:H1classical} and \eqref{eq:Hm1classical}. 
Note that, conversely, Theorems \ref{th:eu:H10} and  \ref{th:eu:En_Det} determine the exponents $k_j$ of the  $H$-fractions \eqref{eq:eu:F10} and \eqref{eq:eu:En_H},  namely the positions of the vanishing Hankel determinants, but not their coefficients $u_j(x)$, so that the present approach 
does not yield the $H$-fractions themselves.

More generally, for any constant $c$, the series $c+xF(x)$ has coefficient sequence $(c,E_1,E_2,\ldots)$.
Since $v_0=F(0)=E_1=1$, equations~\eqref{eq:Hc}, \eqref{eq:H1classical}, and~\eqref{eq:Hm1classical}
give, for every $n\ge1$,
\[
H_n(c+xF(x))
=cH_{n-1}^{(1)}(F)+H_n^{(-1)}(F)
=H_{n-1}(F)\bigl(cR_{n-1}-T_{n-1}\bigr).
\]
Since $E_0=1$, the choice $c=1$ gives the determinants evaluated in Theorem~\ref{th:eu:En_Det}.
Taking $c=0$ gives the backward shifted determinants $H_n^{(-1)}(F)=-H_{n-1}(F)T_{n-1}$.
\end{Remark}


\section{Application: Hankel determinants of convolution powers of the Narayana polynomials}\label{sec:cigler}

In \cite{Han2025NarayanaConv} the first author evaluated the Hankel determinants of the convolution powers of the Narayana polynomials, that is, of the power series $(\gamma(t,q)-1)^m/q^s$ for $s=0,1,2,3$, by establishing explicitly the $H$-fractions of these four series.
For $m\geq2$, all four are genuine $H$-fractions, with vanishing Hankel determinants and a periodic pattern of period $1$, $2$ or $3$.
The derivations for $s=0$, $2$, and $3$ are lengthy and rely on a guess-and-prove strategy involving large intermediate expressions, many of which are verified by computer algebra.

In this section we give short proofs of the same evaluations.
We use only the simplest of the four $H$-fractions, namely that of
\[
F(q)=\frac{(\gamma(t,q)-1)^m}{q}.
\]
For $m\ge2$, the series with $s=2$ is the $1$-shift of $F$.
For $m\ge3$, the series with $s=3$ is its $2$-shift.
The series with $s=0$ is $qF(q)$, whose coefficient sequence is obtained from that of $F$ by prepending a zero.
The three remaining evaluations therefore follow by straightforward computations from
Theorems~\ref{th:shiftHDet}, \ref{th:shift2HDet}, and~\ref{th:backshiftHDet}.
Two of these evaluations generalize a theorem and a conjecture of Cigler \cite{Cigler2018}.
Following \cite{Cigler2018, Han2025NarayanaConv}, the variable of the power series is denoted by $q$ instead of $x$ throughout this section, and $t$ is a second indeterminate.

\subsection{Narayana polynomials and their convolution powers}\label{subsec:cig:intro}

The {\it Narayana polynomials} are defined by
\begin{equation*}
	\gamma_n(t) = \sum_{k=0}^n \binom{n}{k} \binom{n-1}{k} \frac{1}{k+1}\, t^k
	\qquad (n\geq 0),
\end{equation*}
so that $\gamma_0(t)=\gamma_1(t)=1$, $\gamma_2(t)=1+t$, $\gamma_3(t)=1+3t+t^2$, \ldots, and $\gamma_n(1)$ is the Catalan number $C_n$. Their generating function
$\gamma(t,q)=\sum_{n\ge 0}\gamma_n(t)\, q^n$ satisfies the quadratic relation
\begin{equation*}
	-1 + (1 - q + tq)\,\gamma(t,q) - tq\,\gamma(t,q)^2 = 0.
\end{equation*}
Since $\gamma(t,q)-1=q+(1+t)q^2+\cdots$ has valuation $1$, the series $G(t,q)=(\gamma(t,q)-1)/q$ is a power series with constant term $1$.
The {\it convolution powers of the Narayana polynomials}, denoted $\gamma^{(\tau)}_n(t)$, are defined by
\begin{align*}
	\sum_{n\ge0} \gamma^{(2\tau)}_n(t)\, q^n &= G(t,q)^\tau,\\
	\sum_{n\ge0} \gamma^{(2\tau+1)}_n(t)\, q^n &= \gamma(t,q)\, G(t,q)^\tau,
\end{align*}
with the convention $\gamma^{(\tau)}_n(t)=0$ for $n<0$.
At $t=1$, $\gamma^{(\tau)}_n(1)=\frac{\tau}{2n+\tau}\binom{2n+\tau}{n}$ is the $\tau$-fold convolution power of the Catalan numbers.
Motivated by the evaluation  of the Hankel determinants of the sums of consecutive Catalan numbers in \cite{Cvetkovic2002RI_CatalanHDet},
Cigler \cite{Cigler2018} studied the Hankel determinants of order $N\ge1$ given by
\begin{equation}\label{eq:cig:Delta}
	\Delta^{(\tau)}_{\varepsilon} =
	\det \bigl({\gamma}^{(\tau)}_{i+j-\varepsilon} (t) \bigr)_{i,j=0}^{N-1}.
\end{equation}
He obtained explicit formulas for $\tau=3$ and $\tau=4$ and conjectured a formula for $\tau=6$.
For related results on the convolution powers of the Catalan numbers, see
\cite{Wang2018Xin_Convolution, Cigler2023_Experimental, Cigler2024ConvCatalan, Fulmek2024_CiglerConjecture}.
In \cite{Han2025NarayanaConv} the first author evaluated $\Delta^{(2m)}_{m-s}$ for $s=0,1,2,3$ and every $m\geq s$, thereby generalizing Cigler's theorem for $\tau=4$
(the case $m=2$, $s=2$) and proving his conjecture for $\tau=6$ (the case $m=3$, $s=3$).

\smallskip

For the rest of this section, $m\geq 1$ is a fixed integer, $[k]_q=1+q+\cdots+q^{k-1}$ is the usual $q$-integer, and we write
\begin{equation*}
	\xi := (-1)^{\binom{m}{2}} = (-1)^{m(m-1)/2},
\end{equation*}
so that $\xi^n$ is the sign denoted by $\xi_1$ in \cite{Han2025NarayanaConv}.
For $n\ge0$, we also define the polynomials
\begin{equation}\label{eq:cig:R}
	R(m; t,n)
	= m[m]_t
	\sum_{i=0}^{2n} \binom{\min(i,\,2n-i)+2}{2}\, t^{mi}.
\end{equation}
The sum in \eqref{eq:cig:R} is a palindromic polynomial of degree $2n$ in $t^m$, with coefficients $1,3,6,\ldots,\binom{n+2}{2},\ldots,6,3,1$.
For instance, 
\begin{align*}
R(m;t,0) &= m[m]_t,\\
R(m;t,1) &= m[m]_t(1+3t^m+t^{2m}).
\end{align*}
The three results of \cite{Han2025NarayanaConv} which we reprove here are the cases $s=0$, $s=2$ and $s=3$.

\begin{Theorem}[{\cite[Theorem~4]{Han2025NarayanaConv}}]\label{th:cig:Det0}
	For $m\ge 1$ we have
	\begin{align*}
		\Delta^{(2m)}_{m}
		&=\begin{cases}
			1, & \text{if $N=0$},\\[2pt]
			-\xi^n\, t^{m(n-1)(mn-2)/2}\, [n-1]_{t^m},
			& \text{if $N=mn$}, \\[2pt]
			(-1)^{m}\,\xi^n\, t^{m^2 n(n-1)/2}\, [n]_{t^m},
			& \text{if $N=mn+1$}, \\[2pt]
			0, & \text{otherwise}.
		\end{cases}
	\end{align*}
\end{Theorem}

\begin{Theorem}[{\cite[Theorem~6]{Han2025NarayanaConv}}]\label{th:cig:Det2}
	For $m\ge 2$ we have
	\begin{align*}
		\Delta^{(2m)}_{m-2}
		&=\begin{cases}
			\xi^n\, t^{m^2 n(n-1)/2}\, [n+1]_{t^m},
			& \text{if $N=mn$}, \\[2pt]
			(-1)^{m-1}\,\xi^n\, t^{m(n-1)(mn-2)/2}\, [n]_{t^m},
			& \text{if $N=mn-1$}, \\[2pt]
			0, & \text{otherwise}.
		\end{cases}
	\end{align*}
\end{Theorem}

\begin{Theorem}[{\cite[Theorem~7]{Han2025NarayanaConv}}]\label{th:cig:Det3}
	For $m\geq 3$ we have
	\begin{align*}
		\Delta^{(2m)}_{m-3}
		&=\begin{cases}
			\xi^n\, t^{m^2 n(n-1)/2 }\, [n+1]^2_{t^m},
			& \text{if $N=mn$}, \\[2pt]
			(-1)^{m-1}\xi^n\, t^{m(n-1)(mn-2)/2}\, R(m; t, n-1),
			& \text{if $N=mn-1$} ,\\[2pt]
			-\xi^n\,  t^{m(n-1)(mn-4)/2}\,  [n]^2_{t^m},
			& \text{if $N=mn-2$}, \\[2pt]
			0, & \text{otherwise}.
		\end{cases}
	\end{align*}
\end{Theorem}

Fix $m\ge1$.
The determinants $\Delta^{(2m)}_{m-s}$ \eqref{eq:cig:Delta} are shifted Hankel determinants of a single power series. 
By definition, $\gamma_r^{(2m)}(t)$ is the coefficient of $q^{r+m}$ in $(\gamma(t,q)-1)^m$,
which has valuation $m$.
For $0\le s\le m$ and $i,j\ge0$, we therefore have
\[
\gamma_{i+j-m+s}^{(2m)}(t)
=
[q^{i+j}]\frac{(\gamma(t,q)-1)^m}{q^s},
\]
with the convention $\gamma_r^{(2m)}(t)=0$ for $r<0$.
Taking determinants gives
\[
\Delta^{(2m)}_{m-s}
=
H_N\!\left(\frac{(\gamma(t,q)-1)^m}{q^s}\right).
\]
Let
\begin{equation}\label{eq:cig:F}
	F(q) := \frac{(\gamma(t,q)-1)^m}{q} = q^{m-1}G(t,q)^m .
\end{equation}
The correspondence between the generating functions and their Hankel determinants is summarized below.
Each row applies when $s\le m$.
\begin{equation}\label{eq:cig:DeltaH}
\renewcommand{\arraystretch}{1.4}
\begin{array}{cccc}
\hline
s & \text{Generating function}
  & \text{Cigler's notation}
  & \text{Our notation} \\
\hline
0 & qF(q)
  & \Delta^{(2m)}_{m}
  & H_N^{(-1)}(F) \\
1 & F(q)
  & \Delta^{(2m)}_{m-1}
  & H_N(F) \\
2 & F(q)/q
  & \Delta^{(2m)}_{m-2}
  & H_N^{(1)}(F) \\
3 & F(q)/q^2
  & \Delta^{(2m)}_{m-3}
  & H_N^{(2)}(F) \\
\hline
\end{array}
\end{equation}
Here $N$ is the order of the determinant in \eqref{eq:cig:Delta}.
Consequently, Theorems \ref{th:cig:Det0}, \ref{th:cig:Det2} and \ref{th:cig:Det3} are evaluations of the backward shifted, $1$-shifted and $2$-shifted Hankel determinants of $F(q)$, and they will follow from Theorems \ref{th:backshiftHDet}, \ref{th:shiftHDet} and \ref{th:shift2HDet} once the $H$-fraction of $F(q)$ is known.

\subsection{The $H$-fraction of $(\gamma(t,q)-1)^m/q$}\label{subsec:cig:hfrac}

The $H$-fraction of $F(q)$ was determined in \cite{Han2025NarayanaConv}. To state it, define,
for $d=0,1,\ldots,m-1$,
\begin{equation*}
	\rho(m; t, d)=\sum_{i=0}^{d} \frac{m(m-d)}{(m-i)(m-d+i)}
	\binom{m-d+i}{i} \binom{ m-i}{d-i}\, t^i ,
\end{equation*}
together with $\rho(m;t,m)=1+t^m$, and let
\begin{equation}\label{eq:cig:beta}
	\beta(m; t, q) = \sum_{d=0}^{m} \rho(m; t, d)\,(-q)^d .
\end{equation}
Thus $\beta(m;t,q)$ is a polynomial of degree $m$ in $q$ with constant term $1$; explicitly,
\begin{align}
	\beta(m;t,q)
	&= 1 - m(1+t)\,q
	+ \Bigl(\binom{m}{2} + m(m-2)t + \binom{m}{2}t^2\Bigr) q^2
	+ \cdots \notag\\
	&\qquad + m[m]_t\, (-q)^{m-1} + (1+t^m)\,(-q)^m,
	\label{eq:cig:betaexp}
\end{align}
where the two top coefficients come from $\rho(m;t,m-1)=m(1-t^m)/(1-t)=m[m]_t$ and $\rho(m;t,m)=1+t^m$.
For the first three values of $m$, we have
\begin{align*}
\beta(1;t,q)
&=1-(1+t)q,\\
\beta(2;t,q)
&=1-2(1+t)q+(1+t^2)q^2,\\
\beta(3;t,q)
&=1-3(1+t)q+3(1+t+t^2)q^2-(1+t^3)q^3.
\end{align*}
It is shown in \cite{Han2025NarayanaConv} that $\beta(m;t,q)=L_m(1-q-tq,\,-tq^2)$, where $L_m(x,s)$ is the Lucas polynomial, but we shall not need this fact.

\begin{Lemma}[{\cite[Lemma~13]{Han2025NarayanaConv}}]\label{lem:cig:HFrac}
	For $m\ge 1$, the $H$-fraction expansion of the power series
	$F(q)=(\gamma(t,q)-1)^m/q$ is
	\begin{equation}\label{eq:cig:HFrac}
		F(q)
		=
		\cFrac{v_0 q^{k_0}}{ 1-u_1(q) q}
		- \cFrac{ v_1 q^{k_0+k_1+2}}{ 1-u_2(q)q}
		- \cFrac{ v_2 q^{k_1+k_2+2}}{ 1-u_3(q)q}
		- \cdots,
	\end{equation}
	where, for every $j\ge 0$,
	\begin{equation}\label{eq:cig:HFracdata}
		k_j = m-1, \quad
		v_0 = 1, \quad
		v_{j+1} = t^m , \quad
		1-u_{j+1}(q)\,q = \beta(m;t,q).
	\end{equation}
\end{Lemma}

In \cite{Han2025NarayanaConv}, the denominators of the $H$-fraction are written as $1+u_j(q)q$.
By Remark~\ref{rem:sign}, passing to our convention \eqref{eq:Hfrac} simply changes the sign of each $u_j(q)$.
The data in \eqref{eq:cig:HFracdata} have already been converted to this convention.
These data also satisfy the required degree condition.
Indeed, since $\beta(m;t,q)$ has constant term $1$ and degree $m$ in $q$, the expression
\[
u_{j+1}(q)=\frac{1-\beta(m;t,q)}{q}
\]
is a polynomial of degree $m-1=k_j$.
For $m=1$, we have $F(q)=G(t,q)$ and $k_j=0$ for all $j$.
Thus \eqref{eq:cig:HFrac} reduces to the classical $J$-fraction of the shifted Narayana polynomials,
\[
G(t,q)
=
\cFrac{1}{1-(1+t)q}
- \cFrac{tq^2}{1-(1+t)q}
- \cFrac{tq^2}{1-(1+t)q}
- \cdots.
\]
For $m\ge2$, we have $k_j=m-1\ge1$ and $s_j=jm$.
Theorem~\ref{th:super2}(ii) therefore gives $H_N(F)=0$ whenever $m\nmid N$.
In this case, $F(q)$ does not admit a $J$-fraction of the form \eqref{eq:Jacobi}.

\subsection{The data of the $H$-fraction}\label{subsec:cig:data}

We now extract from Lemma \ref{lem:cig:HFrac} all the quantities entering Theorems \ref{th:super2}, \ref{th:shiftHDet}, \ref{th:shift2HDet} and \ref{th:backshiftHDet}.

By \eqref{eq:cig:HFracdata}, we have $k_j=m-1$ for every $j\ge0$, together with $v_0=1$
and $v_j=t^m$ for $j\ge1$.
Hence, for every $j\ge0$,
\begin{equation}\label{eq:cig:sjepsj}
	s_j = (m-1)j + j = mj, \quad
	\epsilon_j = j\binom{m}{2}, \quad
	V_j = v_0v_1\cdots v_j = t^{mj}.
\end{equation}

We next determine the highest-degree coefficients of $u_j(q)$.
For every $j\ge1$, equations \eqref{eq:cig:HFracdata} and \eqref{eq:cig:beta} give
\[
u_j(q)\, q = 1-\beta(m;t,q) = \sum_{d=1}^m (-1)^{d-1}\rho(m;t,d)\,q^d,
\]
so that the coefficient $c_{j,i}$ of $q^i$ in $u_j(q)$ is $(-1)^{i}\rho(m;t,i+1)$.
Reading off the two top coefficients with \eqref{eq:cig:betaexp}, we get, for every $j\ge 1$,
\begin{align}
	c_{j,m-1} &= (-1)^{m-1}(1+t^m), \label{eq:cig:cjm1}\\
	c_{j,m-2} &= (-1)^{m}\, m[m]_t  \qquad (m\geq 2).   \label{eq:cig:cjm2}
\end{align}

Substituting \eqref{eq:cig:HFracdata} and \eqref{eq:cig:sjepsj} into \eqref{eq:HankelDetFundamental}, we obtain, for every $j\ge 0$,
\begin{equation}\label{eq:cig:HsjF}
	H_{mj}(F)
	= (-1)^{\epsilon_j}\,
	v_0^{s_j}\prod_{i=1}^{j-1}v_i^{s_j-s_i}
	= \xi^j \prod_{i=1}^{j-1}(t^m)^{m(j-i)}
	= \xi^j\, t^{m^2 j(j-1)/2},
\end{equation}
and $H_N(F)=0$ if $m$ does not divide $N$.
By \eqref{eq:cig:DeltaH}, this recovers \cite[Theorem~5]{Han2025NarayanaConv}, namely
\[
\Delta^{(2m)}_{m-1}
=
\begin{cases}
\xi^n t^{m^2n(n-1)/2}, & \text{if } N=mn \text{ for some } n\ge0,\\
0, & \text{otherwise}.
\end{cases}
\]

We now consider the sequence $(R_n)_{n\ge0}$ associated with the $H$-fraction \eqref{eq:cig:HFrac}.
By \eqref{def:Rn} and \eqref{eq:cig:cjm1}, its initial values are $R_0=1$ and
$R_1=c_{1,m-1}=(-1)^{m-1}(1+t^m)$.
For $n\ge1$, we have $v_n=t^m$, so the recurrence becomes
\begin{equation}\label{eq:cig:Rnrec}
	R_{n+1} = (-1)^{m-1}(1+t^m)\, R_n - t^m R_{n-1}.
\end{equation}
This sequence should not be confused with the polynomials $R(m;t,n)$ defined in \eqref{eq:cig:R}. 
The following lemma gives a closed form for the sequence $(R_n)_{n\ge0}$.

\begin{Lemma}\label{lem:cig:Rn}
	For every $n\ge 0$, we have $R_n = (-1)^{(m-1)n}\, [n+1]_{t^m}$.
\end{Lemma}

\begin{proof}
	The formula holds for $n=0,1$. Let $n\geq 1$ and assume that it holds for $n-1$ and $n$. 
    Since $(-1)^{(m-1)(n-1)}=(-1)^{(m-1)(n+1)}$, the recurrence \eqref{eq:cig:Rnrec} gives
	\begin{align*}
		R_{n+1}
		&= (-1)^{m-1}(1+t^m)\,(-1)^{(m-1)n} [n+1]_{t^m}
		   - t^m\, (-1)^{(m-1)(n-1)} [n]_{t^m}\\
		&= (-1)^{(m-1)(n+1)}\bigl[(1+t^m) [n+1]_{t^m} - t^m [n]_{t^m}\bigr]\\
		&= (-1)^{(m-1)(n+1)} [n+2]_{t^m},
	\end{align*}
	because $(1+Q)[n+1]_Q-Q[n]_Q=[n+1]_Q+Q^{n+1}=[n+2]_Q$ for
	any $Q$.
\end{proof}

Similarly, the sequence $(E_n)_{n\ge0}$ of Theorem~\ref{th:backshiftHDet} has initial values
$E_0=0$ and $E_1=1$.
For $n\ge1$, \eqref{def:En} gives the same recurrence as for $(R_n)$,
\begin{equation}\label{eq:cig:Enrec}
	E_{n+1} = (-1)^{m-1}(1+t^m)\, E_n - t^m E_{n-1}.
\end{equation}
The following lemma gives a closed form for the sequence $(E_n)_{n\ge0}$.

\begin{Lemma}\label{lem:cig:En}
	For every $n\ge 1$, we have $E_n = (-1)^{(m-1)(n-1)}\, [n]_{t^m}$.
\end{Lemma}

\begin{proof}
	The formula holds for $n=1$ and $n=2$, since $E_1=1=[1]_{t^m}$ and
$E_2=(-1)^{m-1}(1+t^m)=(-1)^{m-1}[2]_{t^m}$.
Let $n\ge2$, and assume that the formula holds for $n-1$ and $n$.
By \eqref{eq:cig:Enrec}, we obtain
\[
		E_{n+1}
		= (-1)^{(m-1)n}\bigl[(1+t^m) [n]_{t^m} - t^m [n-1]_{t^m}\bigr]
		= (-1)^{(m-1)n} [n+1]_{t^m},
	\]
as in the proof of Lemma~\ref{lem:cig:Rn}.
This completes the induction.
\end{proof}

Note that $E_n=R_{n-1}$ for $n\geq 1$, in accordance with the remark following Theorem \ref{th:backshiftHDet}: the tail of the $H$-fraction \eqref{eq:cig:HFrac} obtained by deleting its first term has the same data $k_j=m-1$, $v_j=t^m$, $1-u_j(q)q=\beta(m;t,q)$ as the $H$-fraction itself, except for $v_0$, which does not enter the recurrence.

We also record the sign identities needed in \eqref{eq:H1sjkj} and \eqref{eq:H1sj:Hsj:k11}.
For $m\ge2$, the identities
$\binom{m-1}{2}=\binom{m}{2}-(m-1)$ and
$\binom{m-2}{2}=\binom{m}{2}-(2m-3)$ give
\begin{equation}\label{eq:cig:signs}
	(-1)^{\binom{m-1}{2}} = (-1)^{m-1}\,\xi,
	\qquad
	(-1)^{\binom{m-2}{2}} = -\,\xi .
\end{equation}
The first identity also holds for $m=1$.

The vanishing assertions in Theorems~\ref{th:cig:Det2} and~\ref{th:cig:Det3} follow from Theorems~\ref{th:shiftHDet} and~\ref{th:shift2HDet}.
The following lemma gives a more general result for arbitrary positive shifts.
Its proof uses only the vanishing pattern of the ordinary Hankel determinants and elementary
linear algebra.

\begin{Lemma}\label{lem:cig:vanish}
Let $F(x)=\sum_{n\ge0}a_nx^n$ be a formal power series over a field, and let $s\ge0$ and $k\ge1$ be integers.
Suppose that $H_s(F)\ne0$ and
\[
H_{s+1}(F)=H_{s+2}(F)=\cdots=H_{s+k}(F)=0.
\]
Then, for all integers $r,\ell\ge1$ satisfying $r+\ell\le k$, we have $ H_{s+\ell}^{(r)}(F)=0.$
\end{Lemma}

\begin{proof}
	For $i\geq 0$ let $\mathbf{a}_i=(a_{i+q})_{q\geq 0}$ be the $i$-th row of the infinite Hankel matrix $(a_{i+q})_{i,q\geq 0}$, so that $H^{(r)}_N(F)$ is the determinant of the matrix formed by the first $N$ entries of the rows
	$\mathbf{a}_r,\mathbf{a}_{r+1},\ldots,\mathbf{a}_{r+N-1}$.

	\smallskip
    
	\noindent{\bf Step 1.}\ \ We claim that there exist	$c_0,c_1,\ldots,c_{s-1}\in\setF$ such that the sequence
	\[
	e_n := a_{s+n}-\sum_{p=0}^{s-1}c_p\, a_{p+n} \qquad (n\geq 0)
	\]
	satisfies $e_n=0$ for $0\leq n\leq s+k-1$.
Since $H_s(F)\ne0$, the first $s$ rows of the $(s+1)\times(s+1)$ Hankel matrix $(a_{p+q})_{0\le p,q\le s}$ are linearly independent.
Since $H_{s+1}(F)=0$, its last row is a linear combination of the first $s$ rows.
We may therefore choose $c_0,\ldots,c_{s-1}$ so that $e_n=0$ for $0\le n\le s$.
For $r\ge0$, set
\[
\mathbf e_r
:=
\mathbf a_{s+r}
-\sum_{p=0}^{s-1}c_p\mathbf a_{p+r}
=
(e_{r+q})_{q\ge0}.
\]
The row $\mathbf e_r$ is obtained from $\mathbf a_{s+r}$ by subtracting a linear combination of the rows $\mathbf a_r,\ldots,\mathbf a_{s+r-1}$ of smaller index.
Hence, for every $i\ge1$, the matrix
\[
\widetilde A_{s+i}
=
\bigl(
\mathbf a_0,\ldots,\mathbf a_{s-1},
\mathbf e_0,\ldots,\mathbf e_{i-1}
\bigr)^{\mathsf T}
\]
is obtained from the Hankel matrix
\[
A_{s+i}
=
\bigl(
\mathbf a_0,\ldots,\mathbf a_{s-1},
\mathbf a_s,\ldots,\mathbf a_{s+i-1}
\bigr)^{\mathsf T}
\]
by left multiplication by a lower triangular matrix $T$ with all diagonal entries equal to $1$.
Here each row is restricted to its first $s+i$ entries.
Thus $\widetilde A_{s+i}=T A_{s+i}$, and since $\det T=1$, we obtain
\[
\det \widetilde A_{s+i}
=
\det A_{s+i}
=
H_{s+i}(F).
\]
We now prove by induction on $i\in\{1,\ldots,k\}$ that $e_n=0$ for $n\leq s+i-1$. The case $i=1$ has just been established. 
Let $2\leq i\leq k$ and assume that $e_n=0$ for $n\leq s+i-2$. 
Then the entry $e_{r+q}$ of $\widetilde A_{s+i}$ in row $\mathbf{e}_r$ $(0\le r\le i-1)$ and column $q$ vanishes unless $q\geq s+i-1-r\geq s$. 
Hence
	\[
	\widetilde A_{s+i}=
	\begin{pmatrix}
		(a_{i'+q})_{0\le i',q\le s-1} & * \\
		0 & B
	\end{pmatrix},
	\qquad
	B=(e_{s+r+q})_{0\le r,q\le i-1}.
	\]
By the induction hypothesis, the entries of $B$ vanish when $r+q<i-1$, and its antidiagonal entries are all $e_{s+i-1}$.
Reversing its columns gives a lower triangular matrix, so $\det B=(-1)^{ i(i-1)/2}e_{s+i-1}^{\,i}.$
Consequently,
\[
0=H_{s+i}(F)
=
(-1)^{i(i-1)/2 }H_s(F)e_{s+i-1}^{\,i}.
\]
Since $H_s(F)\ne0$ and $\setF$ is a field, we obtain $e_{s+i-1}=0$.
This completes the induction and proves the claim.

	\smallskip
	\noindent{\bf Step 2.}\ \ Let $r,\ell\ge1$ satisfy $r+\ell\le k$, and set $N=s+\ell$.
Consider the matrix
\[
M=
\bigl(
\mathbf a_r,\ldots,\mathbf a_{s+r-1},
\mathbf a_{s+r},\ldots,\mathbf a_{r+N-1}
\bigr)^{\mathsf T},
\]
where each row is restricted to its first $N$ entries, so that $\det M=H_N^{(r)}(F)$.
For $r\le p\le r+\ell-1$, the row $\mathbf e_p$ is obtained from $\mathbf a_{s+p}$ by subtracting a linear combination of the original rows $\mathbf a_p,\ldots,\mathbf a_{s+p-1}$.
These are all rows of $M$ of smaller index.
Thus, as in Step 1, the matrix
\[
M'=
\bigl(
\mathbf a_r,\ldots,\mathbf a_{s+r-1},
\mathbf e_r,\ldots,\mathbf e_{r+\ell-1}
\bigr)^{\mathsf T}.
\]
By the same argument as in Step 1, we obtain $\det M'=\det M.$
For $r\le p\le r+\ell-1$ and $0\le q\le s$, we have
\[
p+q\le s+r+\ell-1\le s+k-1.
\]
Step 1 therefore gives $e_{p+q}=0$.
Thus the last $\ell$ rows of $M'$ vanish in the first $s+1$ columns.
Since $M'$ has $s+\ell$ columns, these rows are supported on at most $\ell-1$ columns and are therefore linearly dependent.
It follows that $H_N^{(r)}(F)=\det M=\det M'=0$, which completes this proof.    
\end{proof}

For $m\ge2$, consider the series $F(q)$ in \eqref{eq:cig:F}.
By Lemma~\ref{lem:cig:HFrac} and Theorem~\ref{th:super2}(ii), we have $H_{mj}(F)\ne0$ and $H_{mj+1}(F)=\cdots=H_{mj+m-1}(F)=0$ for every $j\ge0$.
Applying Lemma~\ref{lem:cig:vanish} with $s=mj$ and $k=m-1$, we obtain, for $r\in\{1,2\}$
and $1\le\ell\le m-r-1$,
\begin{equation}\label{eq:cig:vanish}
H_{mj+\ell}^{(r)}(F)=0.
\end{equation}
Consequently, $H_N^{(1)}(F)=0$ unless $N\equiv0,-1\pmod m$, and $H_N^{(2)}(F)=0$ unless $N\equiv0,-1,-2\pmod m$.
These vanishing statements are also contained in Theorems~\ref{th:shiftHDet} and~\ref{th:shift2HDet}(3), respectively.

\subsection{Proofs of the Hankel determinant formulas}\label{subsec:cig:proofs} 
We now apply the preceding computations to prove Theorems~\ref{th:cig:Det2}, \ref{th:cig:Det3}, and \ref{th:cig:Det0}.
We first prove Theorem~\ref{th:cig:Det2} by applying the formula for $1$-shifted Hankel determinants.

\begin{proof}[Proof of Theorem~\ref{th:cig:Det2}]
Let $m\geq 2$. By \eqref{eq:cig:DeltaH}, we have to evaluate $H^{(1)}_N(F)$, and we apply Theorem \ref{th:shiftHDet} to the $H$-fraction \eqref{eq:cig:HFrac}, whose exponents $k_j=m-1$ are all positive.

\medskip
\noindent\textbf{Case 1.}
Suppose that $N=mn$.
By \eqref{eq:H1sj:Hsj} with $j=n$, then \eqref{eq:cig:HsjF} and Lemma \ref{lem:cig:Rn},
\begin{align*}
H^{(1)}_{mn}(F)
&= (-1)^{s_n+n} H_{s_n}(F)\, R_n\\
&= (-1)^{mn+n}\, \xi^n t^{m^2n(n-1)/2}\, (-1)^{(m-1)n}[n+1]_{t^m}\\
&= \xi^n\, t^{m^2 n(n-1)/2}\, [n+1]_{t^m},
\end{align*}
since $(mn+n)+(m-1)n=2mn$ is even.

\medskip
\noindent\textbf{Case 2.}
Suppose that $N=mn-1$ with $n\ge1$.
We apply \eqref{eq:H1sjkj} with $j=n-1\geq 0$, for which $s_j+k_j=m(n-1)+m-1=mn-1$. Using \eqref{eq:cig:signs}, $V_{n-1}=t^{m(n-1)}$ and the case just proved
(with $n-1$ in place of $n$; for $n=1$ we use $H^{(1)}_0(F)=1$),
we get
\begin{align*}
	H^{(1)}_{mn-1}(F)
	&= (-1)^{\binom{m-1}{2}}\, V_{n-1}^{m-1}\, H^{(1)}_{m(n-1)}(F)\\
	&= (-1)^{m-1}\xi \cdot t^{m(m-1)(n-1)} \cdot
	   \xi^{n-1}\, t^{m^2(n-1)(n-2)/2}\, [n]_{t^m}\\
	&= (-1)^{m-1}\xi^n\, t^{m(n-1)(mn-2)/2} \, [n]_{t^m},
\end{align*}
because $m(m-1)(n-1)+m^2(n-1)(n-2)/2 = m(n-1)(mn-2)/2.$

\medskip
\noindent\textbf{Case 3.}
Suppose that $N\not\equiv0,-1\pmod m$. Then \eqref{eq:cig:vanish} gives $H_N^{(1)}(F)=0$.

\medskip

Together with \eqref{eq:cig:DeltaH}, these formulas prove Theorem~\ref{th:cig:Det2}.
This completes the proof. 
\end{proof}

We next prove Theorem~\ref{th:cig:Det3} using the formula for $2$-shifted Hankel determinants.

\begin{proof}[Proof of Theorem~\ref{th:cig:Det3}]

Let $m\geq 3$. By \eqref{eq:cig:DeltaH}, we have to evaluate $H^{(2)}_N(F)$, and we apply Theorem \ref{th:shift2HDet} to the $H$-fraction \eqref{eq:cig:HFrac}. 
Since all exponents $k_j=m-1\geq 2$ are positive, only part (1) of Theorem~\ref{th:shift2HDet} is needed, and in the definition of $\Delta_j$ only the case $k_{j-1}\geq 1$ occurs.

\medskip
\noindent\textbf{Case 1.} Suppose that  $N=mn$.
By \eqref{eq:H1sj:Hsj:k1} with $j=n$, \eqref{eq:cig:HsjF} and Lemma \ref{lem:cig:Rn},
\[
H^{(2)}_{mn}(F)
= H_{s_n}(F)\, R_n^2
= \xi^n\, t^{m^2 n(n-1)/2}\, [n+1]_{t^m}^2.
\]

\medskip
\noindent\textbf{Case 2.} Suppose that  $N=mn-2$ with $n\geq 1$.
We apply \eqref{eq:H1sj:Hsj:k11} with $j=n-1\geq 0$, for which $s_j+k_j-1 = m(n-1)+(m-1)-1 = mn-2$.
By \eqref{eq:cig:signs}, $V_{n-1}=t^{m(n-1)}$ and the case just proved (with $H^{(2)}_0(F)=1$ when $n=1$),
\begin{align*}
	H^{(2)}_{mn-2}(F)
	&= (-1)^{\binom{m-2}{2}}\, V_{n-1}^{m-2}\,
	   H^{(2)}_{m(n-1)}(F) \\
	&= -\xi\cdot t^{m(m-2)(n-1)}\cdot
	   \xi^{n-1}\, t^{m^2 (n-1)(n-2)/2}\, [n]_{t^m}^2\\
	&= -\xi^n\, t^{m(n-1)(mn-4)/2}\, [n]_{t^m}^2,
\end{align*}
because
$m(m-2)(n-1) + m^2(n-1)(n-2)/2 =  m(n-1)(mn-4)/2$.

\medskip
\noindent\textbf{Case 3.}
Suppose that $N=mn-1$ with $n\ge1$.
Set $D_0=0$ and, for $n\ge1$, define
\[
D_n := H^{(2)}_{mn-1}(F),
\qquad
\eta_n := \frac{m(n-1)(mn-2)}{2},
\]
We shall prove that, for every $n\ge1$,
\begin{equation}\label{eq:cig:Dn:target}
D_n
=
(-1)^{m-1}\xi^n t^{\eta_n}R(m;t,n-1).
\end{equation}
We apply \eqref{eq:H1sj:Hsj:k12} with $j=n-1\geq 0$, for which $s_j+k_j = mn-1$ and $c_{j+1,k_j-1}=c_{n,m-2}$:
\begin{equation}\label{eq:cig:Dn:rec0}
	D_n
	= c_{n,m-2}\, V_{n-1}\, H^{(2)}_{mn-2}(F) + \Delta_{n-1}.
\end{equation}
By \eqref{eq:cig:cjm2} and Case~2, the first term is
\begin{align*}
c_{n,m-2}V_{n-1}H_{mn-2}^{(2)}(F)
&=
-(-1)^m m[m]_t\,t^{m(n-1)}
\xi^n t^{m(n-1)(mn-4)/2}[n]_{t^m}^2\\
&=
(-1)^{m-1}\xi^n m[m]_t
t^{\eta_n}[n]_{t^m}^2,
\end{align*}
since
\[
m(n-1)+\frac{m(n-1)(mn-4)}{2}=\eta_n.
\]
For $n\ge2$, the definition of $\Delta_{n-1}$ and $s_{n-2}+k_{n-2}=m(n-1)-1$ give
\[
\Delta_{n-1}
=
(-1)^{\binom{m}{2}}V_{n-1}^{m}
H_{m(n-1)-1}^{(2)}(F)
=
\xi t^{m^2(n-1)}D_{n-1}.
\]
This also holds for $n=1$, since $\Delta_0=D_0=0$.
Substituting into \eqref{eq:cig:Dn:rec0}, we obtain, for every $n\ge1$,
\begin{equation}\label{eq:cig:Dn:rec}
D_n
=
(-1)^{m-1}\xi^n m[m]_t
t^{\eta_n}[n]_{t^m}^2
+
\xi t^{m^2(n-1)}D_{n-1}.
\end{equation}
Set $W_0=0$ and, for $n\ge1$, define the Laurent polynomial
\[
W_n
:=
\frac{D_n}{(-1)^{m-1}\xi^n t^{\eta_n}}.
\]
For $n\ge2$, we have $\eta_n-\eta_{n-1}=m^2(n-1)-m.$
Thus, after normalization, \eqref{eq:cig:Dn:rec} gives, for every $n\ge1$,
\begin{equation}\label{eq:cig:Wn:rec}
W_n
=
m[m]_t[n]_{t^m}^2+t^mW_{n-1}.
\end{equation}
The case $n=1$ follows directly from $D_0=0$ and $\eta_1=0$.
It remains to show that $R(m;t,n-1)$ satisfies the same recurrence.
More precisely, we claim that, for every $n\ge2$,
\begin{equation}\label{eq:cig:Rrec}
R(m;t,n-1)
=
m[m]_t[n]_{t^m}^2+t^mR(m;t,n-2).
\end{equation}
For $h\ge0$, write $R(m;t,h)=m[m]_tT_h(t^m)$, where
\[
T_h(Q)
=
\sum_{i=0}^{h}\binom{i+2}{2}Q^i
+
\sum_{i=h+1}^{2h}\binom{2h-i+2}{2}Q^i.
\]
For $h\ge1$, comparing coefficients and using $\binom{i+2}{2}-\binom{i+1}{2}=i+1$, we obtain
\begin{equation}\label{eq:cig:Trec}
\begin{aligned}
T_h(Q)-QT_{h-1}(Q)
=
\sum_{i=0}^{h}(i+1)Q^i
+
\sum_{i=h+1}^{2h}(2h-i+1)Q^i=[h+1]_Q^2.
\end{aligned}
\end{equation}
Indeed, the coefficient of $Q^i$ in $[h+1]_Q^2$ counts the pairs $(a,b)\in\{0,\ldots,h\}^2$
with $a+b=i$.
Taking $h=n-1$ and $Q=t^m$ in \eqref{eq:cig:Trec} and multiplying both sides by $m[m]_t$, we obtain \eqref{eq:cig:Rrec}.
Since $W_1=m[m]_t=R(m;t,0),$ the recurrences \eqref{eq:cig:Wn:rec} and \eqref{eq:cig:Rrec} imply $W_n=R(m;t,n-1)$ for every $n\ge1$.
This proves \eqref{eq:cig:Dn:target}.

\medskip
\noindent\textbf{Case 4.}
Suppose that $N\not\equiv0,-1,-2\pmod m$.
Then \eqref{eq:cig:vanish} gives $H_N^{(2)}(F)=0$.

\medskip

Together with \eqref{eq:cig:DeltaH}, these formulas prove Theorem~\ref{th:cig:Det3}.
\end{proof}

Finally, we prove Theorem~\ref{th:cig:Det0} using the formula for backward-shifted Hankel determinants.

\begin{proof}[Proof of Theorem~\ref{th:cig:Det0}]
Let $m\geq 1$. By \eqref{eq:cig:DeltaH}, it suffices to evaluate $H_N^{(-1)}(F)$.
We apply Theorem~\ref{th:backshiftHDet} to the $H$-fraction \eqref{eq:cig:HFrac}.
The case $N=0$ follows from the convention $H_0^{(-1)}(F)=1$.

\medskip
\noindent\textbf{Case 1.} Suppose that $N=mn+1$ with $n\geq 1$.
By \eqref{eq:Hsjm1:Hsj} with $j=n$, then \eqref{eq:cig:HsjF}, $v_0=1$ and Lemma \ref{lem:cig:En}, we obtain
\begin{align*}
H^{(-1)}_{mn+1}(F)
&= (-1)^{s_n+n+1}\, v_0\, H_{s_n}(F)\, E_n\\
&= (-1)^{mn+n+1}\, \xi^n t^{m^2n(n-1)/2}\, (-1)^{(m-1)(n-1)}[n]_{t^m}\\
&= (-1)^{m}\xi^n\, t^{m^2 n(n-1)/2}\, [n]_{t^m},
\end{align*}
since $mn+n+1+(m-1)(n-1)=2mn-m+2\equiv m \pmod 2$.

\medskip
\noindent\textbf{Case 2.} Suppose that $N=mn$ with $n\geq 2$.
We apply \eqref{eq:Hs_j+1(-1)} with $j=n-1\geq 1$, for which $s_{j+1}=mn$ and $s_j+1=m(n-1)+1$. 
Using \eqref{eq:cig:signs}, $V_{n-1}=t^{m(n-1)}$, and Case~1, we obtain
\begin{align*}
	H^{(-1)}_{mn}(F)
	&= (-1)^{\binom{m-1}{2}}\, V_{n-1}^{m-1}\, H^{(-1)}_{m(n-1)+1}(F)\\
	&= (-1)^{m-1}\xi \, t^{m(m-1)(n-1)} \,
	   (-1)^m\xi^{n-1}\, t^{m^2(n-1)(n-2)/2}\, [n-1]_{t^m}\\
	&= -\xi^n\, t^{m(n-1)(mn-2)/2}\, [n-1]_{t^m}.
\end{align*}
The exponent of $t$ is simplified as in the proof of Theorem~\ref{th:cig:Det2}.
When $m=1$, the two index families in Cases~1 and~2 coincide, and the formulas agree after reindexing.

\medskip
\noindent\textbf{Case 3.}
Suppose that $N\ge1$ is not among the indices considered in Cases~1 and~2.
Since $s_j=mj$, the index $N$ is neither $s_j+1$ nor $s_{j+1}$ for any $j\ge1$.
Hence Theorem~\ref{th:backshiftHDet} gives $H_N^{(-1)}(F)=0$.
In particular, this holds for $1\le N\le m$.
It remains to check the stated formulas at $N=1$ and $N=m$.
For $N=1$, setting $n=0$ in the expression for $N=mn+1$ gives $(-1)^m[0]_{t^m}=0$.
For $N=m$, setting $n=1$ in the expression for $N=mn$ gives $-\xi[0]_{t^m}=0$.
Both agree with the vanishing established above.

\medskip

Together with \eqref{eq:cig:DeltaH}, these evaluations prove Theorem~\ref{th:cig:Det0}.
\end{proof}


\section{Application: the $(n+2)$-shifted Hankel determinants of $q$-deformed metallic numbers}\label{sec:metallicn2}

\makeatletter
\providecommand{\CF}{}
\DeclareRobustCommand{\CF}{\DOTSB\gaussk@\slimits@}
\newcommand{\gaussk@}{\mathop{\vphantom{\sum}\mathpalette\bigcal@{K}}}
\newcommand{\bigcal@}[2]{%
  \vcenter{\m@th
    \sbox\z@{$#1\sum$}%
    \dimen@=\dimexpr\ht\z@+\dp\z@
    \hbox{\resizebox{!}{0.8\dimen@}{$\mathcal{K}$}}%
  }%
}
\providecommand{\cfm}{\,\mathrel{\cfracminus@}\,}
\newcommand{\cfracminus@}{%
  \sbox\z@{$\dfrac{1}{1}$}%
  \sbox\tw@{$-$}%
  \raisebox{\dimexpr\dp\tw@-\dp\z@\relax}{$-$}%
}
\providecommand{\cfp}{\,\mathrel{\cfracplus@}\,}
\newcommand{\cfracplus@}{%
  \sbox\z@{$\dfrac{1}{1}$}%
  \sbox\tw@{$+$}%
  \raisebox{\dimexpr\dp\tw@-\dp\z@\relax}{$+$}%
}
\providecommand{\cfd}{\mathord{\cfracdots@}}
\newcommand{\cfracdots@}{%
  \sbox\z@{$\dfrac{1}{1}$}%
  \sbox\tw@{$+$}%
  \raisebox{\dimexpr\dp\tw@-\dp\z@\relax}{$\cdots$}%
}
\makeatother
\providecommand{\cron}{\langle n\rangle_q}

In this section we study the $(n+2)$-shifted Hankel determinants of the $q$-metallic numbers $\Phi_n$.
For $n\ge2$, we give an alternative proof of the (anti-)periodicity and the bounds stated in
Theorem~\ref{met:thm-main} by applying Theorem~\ref{th:shiftHDet} to the shifted series $\Phi_n^{(n+1)}$, whose $H$-fraction is known.

\subsection{The $q$-metallic numbers and the conjecture of Han and Pedon}\label{met:sec-intro}

In this section the variable of the power series is denoted by $q$ instead of $x$, and the letter $n$ is reserved for the index of the metallic number. 
Accordingly, for a power series $F(q)=\sum_{i\geq 0} f_i q^i$, we write
\begin{equation}\label{met:HankDet}
H_j^{(\ell)}(F):=\det\bigl(f_{a+b+\ell}\bigr)_{0 \le a,b \le j-1},\qquad \ell, j\in\mathbb{Z}_{\geq 0},
\end{equation}
for the $\ell$-shifted Hankel determinant of order $j$ of \eqref{eq:defHnk}, with the convention $H_0^{(\ell)}(F):=1$, $H_j(F):=H_j^{(0)}(F)$, and $H^{(\ell)}(F):=\bigl(H_j^{(\ell)}(F)\bigr)_{j\geq 0}$ for the whole shifted sequence. 
For $\ell\geq 0$, we also write $F^{(\ell)}(q):=\sum_{i\geq 0} f_{i+\ell}\,q^i$ for the $\ell$-shift of $F$, so that $H_j^{(\ell)}(F)=H_j\bigl(F^{(\ell)}\bigr)$. 
Finally, the level index of the $H$-fraction data $(v_p,k_p,u_p)$, $s_p$, $\epsilon_p$, $R_p$, $V_p$, $d_p$ of Theorems \ref{th:super2}, \ref{th:shiftHDet} and \ref{th:shift2HDet} is denoted by $p$ instead of $j$.

We consider the Hankel determinants of a family of power series $\Phi_n(q)$ ($n\in\mathbb{Z}_{\geq 1}$) introduced by S.~Morier-Genoud and V.~Ovsienko as $q$-deformations of the so-called \emph{metallic numbers}
\begin{equation*}
\phi_n=\frac{n+\sqrt{n^2+4}}{2},\qquad n\geq 1.
\end{equation*}
By definition (see \cite{MorierGenoud2022Ovsienko}),
\begin{equation}\label{met:defPhin}
\Phi_n(q)=[n]_q+\left(\frac{q^n}{[n]_{q^{-1}}}\cfp\frac{q^{-n}}{[n]_q}\cfp\right)^*,
\end{equation}
where $[n]_q:=(1-q^n)/(1-q)$ is the classical $q$-deformation of the integer~$n$. 
Throughout, the symbols $\cfp$ and $\cfm$ denote continued-fraction concatenation with addition and subtraction in the denominators, respectively.
The superscript $*$ denotes infinite repetition of the block in parentheses.
The function $\Phi_n$ expands as a Taylor series around $q=0$ with integer coefficients (see \cite{MorierGenoud2022Ovsienko}), and the family of its Hankel determinants is therefore well defined. 
We will write
\[
H_j^{(\ell)}:=H_j^{(\ell)}(\Phi_n)
\qquad\text{and}\qquad 
H^{(\ell)}:=\bigl(H_j^{(\ell)}\bigr)_{j\geq 0}
\]
for short.

In a recent article, the first author and E.~Pedon \cite{Han2025Pedon} proved that the first $n+2$ sequences $H^{(0)},H^{(1)}, \ldots, H^{(n+1)}$ enjoy several remarkable common properties: 
they are (anti-)periodic, take only the values $-1,0,1$, satisfy a three-term Gale--Robinson recurrence, and are connected to each other by an explicit contiguity formula. 
The main tool is the explicit Hankel continued fraction expansion of $\Phi_n$, the central result of~\cite{Han2025Pedon}, which we recall here in a streamlined form.

\begin{Theorem}[{\cite[Theorem~A]{Han2025Pedon}}]\label{met:thm-HF}
Suppose $n\geq 2$. 
Define the polynomial
\begin{equation}\label{met:defnq}
\cron:=q[n]_q+(1+q^n)(1-q).
\end{equation}
The $H$-fraction expansion of the $q$-deformed metallic number $\Phi_n$ is $(6n-4)$-periodic with offset~$1$ and reads
\begin{equation}\label{met:Hfrac-Phin}
\Phi_n(q)=\frac{1}{1-q}\cfm\left(U_n(q)\cfm V_n(q)\cfm W_n(q)\cfm\right)^*,
\end{equation}
where
\begin{align*}
U_n(q)&=\CF_{i=0}^{n-3}\Bigl(\frac{-q^{n-i}}{[n-i]_q}\cfm \frac{-q^n}{[i+2]_q-q}\cfm \frac{-q^{i+2}}{1-q}\Bigr)\\
&\qquad\cfm\frac{-q^2}{[2]_q}\cfm \frac{-q^n}{[n]_q-q}\cfm \frac{-q^n}{1},\\
V_n(q)&=\frac{q^{n+1}}{\cron+q^{n+1}}\cfm\frac{-q^{2n+1}}{\cron}\cfm\frac{-q^{2n+1}}{\cron+q^{n+1}}\cfm\frac{q^{n+1}}{1},\\
W_n(q)&=\CF_{i=0}^{n-3}\Bigl(\frac{-q^{n-i}}{[n-i]_q-q}\cfm \frac{-q^n}{[i+2]_q}\cfm \frac{-q^{i+2}}{1-q}\Bigr)
\cfm\frac{-q^2}{1}.
\end{align*}
Here $\CF$ denotes the concatenation of the indicated blocks in increasing order of $i$, with $\cfm$ between successive blocks.
\end{Theorem}

In \cite[Theorem~A]{Han2025Pedon}, this expansion is stated with the convention $1+u(q)q$ for the denominators and the operator $+$ between consecutive blocks, whereas here we use $1-u(q)q$ and the operator $-$, as in Remark~\ref{rem:sign}. 
The two conventions are equivalent.
Passing from one to the other reverses the sign of every non-leading numerator and leaves the exponents unchanged.

The properties of the sequences $H^{(\ell)}$ for $\ell\leq n+1$ follow from Theorem~\ref{th:super2}, applied to the explicit $H$-fraction expansions of $\Phi_n$ and its shifts. See \cite[Section~5]{Han2025Pedon} for $\ell=0$ and \cite[Section~7]{Han2025Pedon} for $1\leq \ell\leq n+1$.
For the shifts $\ell\geq n+2$, computer experiments led the authors of \cite{Han2025Pedon} to formulate the following conjecture \cite[Conjecture~E]{Han2025Pedon}.

\begin{Conjecture}\label{met:conj-HKP}
{\rm (a)} The sequence $H^{(n+2)}$ is $2n(n+1)$-periodic when $n$ is even, and $2n(n+1)$-antiperiodic (hence $4n(n+1)$-periodic) when $n$ is odd. 
It takes only the values $-2,-1,0,1,2$, except in the case $n=1$ where its values lie in $\{-1,0,1\}$.

{\rm (b)} All shifted sequences $H^{(\ell)}$ with $\ell\geq n+3$ are unbounded.
\end{Conjecture}

Part~(a) of Conjecture~\ref{met:conj-HKP} was proved for $n=1$ in \cite[Theorem~1.6]{Ovsienko2025Pedon} and recently for all $n\geq 1$ in \cite{Han2026Pedon_nplus2}. 
The latter proof establishes the eventual periodicity of the $H$-fraction of $\Phi_n^{(n+2)}$ using an algorithm for quadratic power series, with its steps verified by symbolic computation. 
The Hankel determinants are then evaluated through telescoping products of the partial numerators of this $H$-fraction.

In this section, we give a short proof of part~(a) of Conjecture~\ref{met:conj-HKP} that does not require the $H$-fraction of $\Phi_n^{(n+2)}$ or any computer verification.
In Section~\ref{sec:metallicn3}, we prove the unboundedness asserted in part~(b) for the first relevant shift $\ell=n+3$ and obtain an explicit linear growth rate.
We reprove below the following result, established by Ovsienko and Pedon~\cite{Ovsienko2025Pedon} for $n=1$ and by Han and Pedon~\cite{Han2026Pedon_nplus2} for all $n\ge1$.

\begin{Theorem}[{\cite{Ovsienko2025Pedon, Han2026Pedon_nplus2}}]\label{met:thm-main}
For every $n\ge1$, the $(n+2)$-shifted Hankel determinant sequence $H^{(n+2)}=\bigl(H_j^{(n+2)}\bigr)_{j\ge0}$ of the $q$-deformed metallic number $\Phi_n$ is
$2n(n+1)$-periodic if $n$ is even and $2n(n+1)$-antiperiodic if $n$ is odd.
More precisely, for every $j\ge0$,
\begin{equation}\label{met:eq-antiper-n2}
H_{j+2n(n+1)}^{(n+2)} = (-1)^n H_j^{(n+2)}.
\end{equation}
All terms of this sequence belong to $\{-2,-1,0,1,2\}$.
When $n=1$, they belong to $\{-1,0,1\}$.
\end{Theorem}

\begin{Example}
For $n=1,\ldots,5$, Theorem~\ref{met:thm-main} gives the following periods and antiperiods.

\begin{center}
\begin{tabular}{ccl}
$n$ & Shift $\ell=n+2$ & Periodicity of $H^{(\ell)}(\Phi_n)$ \\
\hline
$1$ & $3$ & $4$-antiperiodic \\
$2$ & $4$ & $12$-periodic \\
$3$ & $5$ & $24$-antiperiodic \\
$4$ & $6$ & $40$-periodic \\
$5$ & $7$ & $60$-antiperiodic
\end{tabular}
\end{center}

For $n=1$, the sequence is given explicitly in \cite[Theorem~1.6]{Ovsienko2025Pedon} and takes values in $\{-1,0,1\}$. 
For $2\le n\le5$, all terms lie in $\{-2,-1,0,1,2\}$.
For $n=2$, one full period is given by:
\[
\bigl(H_j^{(4)}(\Phi_2)\bigr)_{j=0}^{11}
=
(1,1,-2,-1,2,-1,-2,1,1,0,0,0).
\]
These values follow from the data in Lemma~\ref{met:lem-n2} below and formulas~\eqref{met:key-formula} and~\eqref{met:key-formula-bis}. 
In particular, $H_2^{(4)}(\Phi_2)=-2$ and $H_4^{(4)}(\Phi_2)=2$, so both extremal values are attained.
\end{Example}

The case $n=1$ being known, we assume $n\ge2$ in the proof below.
 The explicit computations are carried out separately for $n\ge3$ and $n=2$. In the latter case, the indexed triples in $U_n$ and $W_n$ are absent, and the corresponding computations are given in Section~\ref{met:subsec-n2}.

\subsection{Reduction to $1$-shifted Hankel determinants}\label{met:sec-bg}

In~\cite{Han2025Pedon}, the $H$-fraction expansion of every shift $\Phi_n^{(\ell)}$, $1\leq \ell\leq n+1$, is obtained from that of $\Phi_n=\Phi_n^{(0)}$ given in Theorem~\ref{met:thm-HF} by deleting the first blocks and adjusting the power of $q$ in the leading numerator. 
Throughout this subsection, $n\geq 2$.

Let us write the $H$-fraction \eqref{met:Hfrac-Phin} of $\Phi_n$ in the compact form
\begin{equation}\label{met:Hfrac-Phin-compact}
\Phi_n(q)=\frac{1}{\beta_0}\cfm\CF_{i=1}^{+\infty}\frac{\alpha_i}{\beta_i},
\end{equation}
so that $\beta_0=1-q$, and the sequences $(\alpha_i)_{i\geq 1}, (\beta_i)_{i\geq 1}$ are $(6n-4)$-periodic.
The following result, proved in \cite[\S7.1]{Han2025Pedon} (modulo the convention change discussed in Theorem~\ref{met:thm-HF} and Remark~\ref{rem:sign}), expresses the $H$-fraction of each shift $\Phi_n^{(\ell)}$ as a tail of \eqref{met:Hfrac-Phin-compact}.

\begin{Theorem}[{\cite[\S7.1]{Han2025Pedon}}]\label{met:thm-HF-shifts}
For every $\ell\in\{1,\ldots, n+1\}$, the $H$-fraction expansion of $\Phi_n^{(\ell)}$ reads
\begin{align*}
\Phi_n^{(\ell)}(q)&=\frac{1}{\beta_{3\ell}}\cfm\CF_{i=3\ell+1}^{+\infty}\frac{\alpha_i}{\beta_i}\qquad\text{for } 1\leq\ell\leq n-1,\\
\Phi_n^{(n)}(q)&=\frac{q^n}{\beta_{3n-1}}\cfm\CF_{i=3n}^{+\infty}\frac{\alpha_i}{\beta_i},\\
\Phi_n^{(n+1)}(q)&=\frac{q^{n-1}}{\beta_{3n}}\cfm\CF_{i=3n+1}^{+\infty}\frac{\alpha_i}{\beta_i}.
\end{align*}
\end{Theorem}

In particular, the $H$-fraction of $\Phi_n^{(n+1)}$ is also $(6n-4)$-periodic.
We describe its periodic structure explicitly in Section~\ref{met:sec-proof}.
The following lemma expresses the $(n+2)$-shifted Hankel determinants of $\Phi_n$ as the $1$-shifted Hankel determinants of $\Phi_n^{(n+1)}$.

\begin{Lemma}\label{met:lem-reduction}
For every $n\geq 1$ and $j\geq 0$, one has
\begin{equation}\label{met:eq-reduction}
H_j^{(n+2)}=H_j^{(1)}\bigl(\Phi_n^{(n+1)}\bigr).
\end{equation}
\end{Lemma}

\begin{proof}
Write $\Phi_n(q)=\sum_{i\ge0}f_iq^i$.
By definition, the coefficient of $q^i$ in $\Phi_n^{(n+1)}$ is $f_{i+n+1}$.
Thus, for every $j\ge1$, \eqref{met:HankDet} gives
\[
\begin{aligned}
H_j^{(1)}\bigl(\Phi_n^{(n+1)}\bigr) = \det\bigl(f_{a+b+n+2}\bigr)_{0\le a,b\le j-1}=H_j^{(n+2)}.
\end{aligned}
\]
For $j=0$, both sides equal $1$ by the empty determinant convention.
\end{proof}

To obtain the $H$-fraction of $\Phi_n^{(n+1)}$ needed on the right-hand side of \eqref{met:eq-reduction} in Lemma~\ref{met:lem-reduction}, we specialize Theorem~\ref{met:thm-HF-shifts} to $\ell=n+1$.
We use the polynomial $\cron$ and the continued-fraction blocks $U_n$, $V_n$, and $W_n$ defined in Theorem~\ref{met:thm-HF}.

\begin{Proposition}\label{met:prop-HF-n+1}
For every $n\geq 2$, the $H$-fraction expansion of $\Phi_n^{(n+1)}$ is $(6n-4)$-periodic with offset~$2$ and reads
\begin{equation}\label{met:Hfrac-Phin-n+1}
\Phi_n^{(n+1)}(q)=\frac{q^{n-1}}{\cron+q^{n+1}}\cfm\frac{q^{n+1}}{1}\cfm\bigl(W_n(q)\cfm U_n(q)\cfm V_n(q)\cfm\bigr)^*.
\end{equation}
\end{Proposition}

\begin{proof}
We apply Theorem~\ref{met:thm-HF-shifts} with $\ell=n+1$:
\begin{equation}\label{met:tail-3n}
\Phi_n^{(n+1)}(q)=\frac{q^{n-1}}{\beta_{3n}}\cfm\CF_{i=3n+1}^{+\infty}\frac{\alpha_i}{\beta_i},
\end{equation}
where $\alpha_i$ and $\beta_i$ are the partial numerators and denominators in \eqref{met:Hfrac-Phin-compact} 
 By the explicit formula \eqref{met:Hfrac-Phin} of Theorem~\ref{met:thm-HF}, the blocks of $\Phi_n$ are arranged in one period as follows (the periodic part starting at block index~$1$):
\begin{equation*}
\underbrace{1,2,\ldots,3n-3}_{U_n}\ \mid\ \underbrace{3n-2,3n-1,3n,3n+1}_{V_n}\ \mid\ \underbrace{3n+2,\ldots,6n-4}_{W_n}.
\end{equation*}
Thus block~$3n$ is the third block of $V_n$, namely $\alpha_{3n}=-q^{2n+1}$ and $\beta_{3n}=\cron+q^{n+1}$, while block~$3n+1$ is the last (fourth) block of $V_n$, namely $\alpha_{3n+1}=q^{n+1}$ and $\beta_{3n+1}=1$. 
Inserting these into~\eqref{met:tail-3n} we obtain
\[
\Phi_n^{(n+1)}(q)=\frac{q^{n-1}}{\cron+q^{n+1}}\cfm\frac{q^{n+1}}{1}\cfm\CF_{i=3n+2}^{+\infty}\frac{\alpha_i}{\beta_i}.
\]
The remaining tail $\CF_{i=3n+2}^{+\infty}$ starts at the first block of $W_n$ and inherits the $(6n-4)$-periodicity of $(\alpha_i,\beta_i)_{i\geq 1}$. 
The cyclic rotation of the period $U_n V_n W_n$ that begins with $W_n$ is the period $W_n U_n V_n$, which gives the announced expression~\eqref{met:Hfrac-Phin-n+1}.
\end{proof}

\subsection{Proof of Theorem~\ref{met:thm-main}}\label{met:sec-proof} 

Fix an integer $n\ge2$.
To apply Theorem~\ref{th:shiftHDet}, we write the $H$-fraction of $\Phi_n^{(n+1)}$ given in Proposition~\ref{met:prop-HF-n+1} in standard form:
\begin{equation}\label{met:Hfrac-Phin-shift}
\Phi_n^{(n+1)}(q)=\frac{v_0\,q^{k_0}}{1-u_1(q)\,q}\cfm\frac{v_1\,q^{k_0+k_1+2}}{1-u_2(q)\,q}\cfm\frac{v_2\,q^{k_1+k_2+2}}{1-u_3(q)\,q}\cfm\cfd
\end{equation}
By Proposition~\ref{met:prop-HF-n+1}, this expansion consists of two initial terms followed by a repeating block of length $6n-4$.
Let $(s_p)_{p\ge0}$ and $(\epsilon_p)_{p\ge0}$ be the associated sequences of Theorem~\ref{th:super2}, and let $(R_p)_{p\ge0}$ be the auxiliary sequence defined by \eqref{def:Rn}.
We also set
\[
\mathcal S:=\{s_p:p\ge0\}, \qquad 
\mathcal S+k:=\{s_p+k_p:p\ge0\}.
\]

By Lemma~\ref{met:lem-reduction} and Theorem~\ref{th:shiftHDet}, in the form \eqref{eq:H1sj:Hsj} and \eqref{eq:H1sjkj},  we obtain, for every $p\ge0$,
\begin{equation}\label{met:key-formula}
H_{s_p}^{(n+2)}=H_{s_p}^{(1)}\bigl(\Phi_n^{(n+1)}\bigr)=(-1)^{s_p+p}\,H_{s_p}\bigl(\Phi_n^{(n+1)}\bigr)\,R_p,
\end{equation}
and
\begin{equation}\label{met:key-formula-bis}
H_{s_p+k_p}^{(n+2)}=(-1)^{\frac{(k_p-1)k_p}{2}}\bigl(v_0 v_1\cdots v_p\bigr)^{k_p}\,H_{s_p}^{(n+2)}.
\end{equation}
Here the first identity at $p=0$ follows from $s_0=0$, $R_0=1$, and the empty determinant convention.
For every $j\ge0$ outside $\mathcal S\cup(\mathcal S+k)$, the same theorem gives $H_j^{(n+2)}(\Phi_n)=0.$ 
Recall that, for each $j\ge1$, $c_{j,k_{j-1}}$ denotes the coefficient of $q^{k_{j-1}}$ in $u_j(q)$.

The proof of Theorem~\ref{met:thm-main} will proceed in three steps:
\begin{enumerate}
\item compute the explicit data $\bigl((v_i, k_i, c_{i,k_{i-1}})\bigr)$ for the $H$-fraction \eqref{met:Hfrac-Phin-shift};
\item deduce the values of $H_{s_p}\bigl(\Phi_n^{(n+1)}\bigr)$ and the auxiliary numbers $R_p$ for $p$ in one full period $0\leq p\leq 6n-4$;
\item combine these in \eqref{met:key-formula}--\eqref{met:key-formula-bis} to obtain the values of $H_j^{(n+2)}$ for $0\leq j\leq 2n(n+1)-1$ (one (anti-)period), and check that they belong to $\{-2,-1,0,1,2\}$.
\end{enumerate}
The (anti-)periodicity \eqref{met:eq-antiper-n2} will then follow at once from the periodicity of the $H$-fraction \eqref{met:Hfrac-Phin-shift} and from the way the prefactors $(-1)^{s_p+p}$, $(-1)^{\epsilon_p}$ behave under the shift $p\mapsto p+(6n-4)$.

\subsubsection{Computations for $n\ge3$}
\label{met:subsec-pf}
To carry out Step~1, we compare the explicit expansion \eqref{met:Hfrac-Phin-n+1} with its standard form \eqref{met:Hfrac-Phin-shift}. 
This gives the parameters in the following lemma.

 \begin{Lemma}\label{met:lem-step1}
Let $n\ge3$.
Consider the $H$-fraction \eqref{met:Hfrac-Phin-shift} of $\Phi_n^{(n+1)}$.
The triples $(v_j,k_j,c_{j+1,k_j})$ are $(6n-4)$-periodic for $j\ge2$.
It therefore suffices to specify $(v_j,k_j)$ for $0\le j\le6n-3$ and $c_{j,k_{j-1}}$ for $1\le j\le6n-2$. These values are given as follows.

\medskip

\noindent{(\rmnum{1})} The initial terms are
 \begin{align*}
(v_0,k_0)&=(1,n-1), & c_{1,n-1}&=-2,\\
(v_1,k_1)&=(1,0),   & c_{2,0}&=0.
\end{align*}

\noindent{(\rmnum{2})} The block $W_n$ occupies the indices $2\le j\le3n-4$.
For $0\le i\le n-3$, the three consecutive terms beginning at $j=3i+2$ give
 \begin{align*}
 (v_{3i+2},k_{3i+2})&=(-1,n-2-i), & c_{3i+3,n-2-i}&=-1,\\
 (v_{3i+3},k_{3i+3})&=(-1,i), & c_{3i+4,i}&=-1,\\
 (v_{3i+4},k_{3i+4})&=(-1,0), & c_{3i+5,0}&=1.
 \end{align*}
 At $j=3n-4$, the final term $-q^2/1$ of $W_n$ gives
\[
(v_{3n-4},k_{3n-4})=(-1,0),
\qquad c_{3n-3,0}=0.
\]

\noindent{(\rmnum{3})} The block $U_n$ occupies the indices $3n-3\le j\le6n-7$.
For $0\le i\le n-3$, the three consecutive terms beginning at $j=3n-3+3i$ give
 \begin{align*}
 (v_{3n-3+3i},k_{3n-3+3i})&=(-1, n-2-i), &   c_{3n-2+3i, n-2-i}& =-1,\\
 (v_{3n-2+3i},k_{3n-2+3i})&=(-1,i),   & c_{3n-1+3i, i} &=-1+\pi(i=0),\\
 (v_{3n-1+3i},k_{3n-1+3i})&=(-1,0), &  c_{3n+3i,0} &= 1,
 \end{align*}
 where $\pi(P) = 1$ if $P$ is true and $0$ otherwise.
 For the final three indices of $U_n$, namely $j=6n-9,\,6n-8,\,6n-7$, we have
 \begin{align*}
 (v_{6n-9},k_{6n-9})&=(-1,0), & c_{6n-9,k_{6n-10}}&=1,\\
  (v_{6n-8},k_{6n-8})&=(-1,n-2), & c_{6n-8,k_{6n-9}}&=-1,\\
 (v_{6n-7},k_{6n-7})&=(-1,0), & c_{6n-7,k_{6n-8}}&=-1.
 \end{align*}

\noindent{(\rmnum{4})} The block $V_n$ occupies the indices $6n-6\le j\le6n-3$.
 \begin{align*}
 (v_{6n-6},k_{6n-6})&=(1,n-1), & c_{6n-6,k_{6n-7}}&=0,\\
 (v_{6n-5},k_{6n-5})&=(-1,n), & c_{6n-5,k_{6n-6}}&=-2,\\
 (v_{6n-4},k_{6n-4})&=(-1,n-1), & c_{6n-4,k_{6n-5}}&=1,\\
 (v_{6n-3},k_{6n-3})&=(1,0), & c_{6n-3,k_{6n-4}}&=-2.
 \end{align*}
 Finally, since the last term of $V_n$ has denominator $1$ and $k_{6n-3}=0$, we have $c_{6n-2,k_{6n-3}}=c_{6n-2,0}=0$, in agreement with $c_{2,0}=0$ in part~(\rmnum{1}).

\medskip

 In particular, every $k_j$ belongs to $\{0,1,\ldots,n\}$, and $v_j\in\{-1,1\}$.
 \end{Lemma}

\begin{proof}
We compare the explicit expansion \eqref{met:Hfrac-Phin-n+1} with the standard form \eqref{met:Hfrac-Phin-shift}, using the definitions of $W_n$, $U_n$, and $V_n$ in
Theorem~\ref{met:thm-HF}.

\medskip
\noindent{\it{ Proof of the values $(v_j,k_j)$.}}
The first two terms of \eqref{met:Hfrac-Phin-n+1} are
\[
\frac{q^{n-1}}{\cron+q^{n+1}}
\qquad\text{and}\qquad
\frac{q^{n+1}}{1}.
\]
Comparing the first numerator with $v_0q^{k_0}$ in \eqref{met:Hfrac-Phin-shift} gives $v_0=1$ and $k_0=n-1$.
Comparing the second numerator with $v_1q^{k_0+k_1+2}$ gives $v_1=1$ and $k_0+k_1+2=n+1$, so $k_1=0$.
This proves the formulas for $(v_j,k_j)$ in part~(\rmnum{1}).

For $0\le i\le n-3$, the three terms of $W_n$ beginning at $j=3i+2$ are
\[
\frac{-q^{n-i}}{[n-i]_q-q},\qquad
\frac{-q^n}{[i+2]_q},\qquad
\frac{-q^{i+2}}{1-q}.
\]
These triples are followed by the final term $-q^2/1$ at $j=3n-4$.
Since every numerator in $W_n$ has coefficient $-1$, we have $v_j=-1$ for every $2\le j\le3n-4$.
If $k_{3i+1}=0$, applying the recurrence from \eqref{met:Hfrac-Phin-shift} to the numerator degrees $n-i$, $n$, and $i+2$ gives the corresponding values $n-2-i$, $i$, and $0$ for $k_j$.
Since $k_1=0$ and the last equality gives $k_{3(i+1)+1}=0$ for the next triple, these formulas hold by induction on $i$.
Taking $i=n-3$ gives $k_{3n-5}=0$, so the final numerator $-q^2$ gives $k_{3n-4}=2-k_{3n-5}-2=0$.
This proves the formulas for $(v_j,k_j)$ in part~(\rmnum{2}).

For $0\le i\le n-3$, the three terms of $U_n$ beginning at $j=3n-3+3i$ are
\[
\frac{-q^{n-i}}{[n-i]_q},\qquad
\frac{-q^n}{[i+2]_q-q},\qquad
\frac{-q^{i+2}}{1-q}.
\]
Their numerators agree with those of the corresponding triple in $W_n$, and $k_{3n-4}=k_1=0$. 
Hence the same recurrence gives the successive $k$-values $n-2-i$, $i$, and $0$.
The final three terms of $U_n$ are 
\[
\frac{-q^2}{[2]_q}, \qquad
\frac{-q^n}{[n]_q-q}, \qquad
\frac{-q^n}{1}.
\]
The last triple ends with $k_{6n-10}=0$. Applying the same recurrence to the numerator degrees $2$, $n$, and $n$ then gives the corresponding values $0$, $n-2$, and $0$ for $k_j$.
Every numerator in $U_n$ has coefficient $-1$, so $v_j=-1$ for $3n-3\le j\le6n-7$. This proves the formulas for $(v_j,k_j)$ in part~(\rmnum{3}).

For $6n-6\le j\le6n-3$, the four terms of $V_n$ are
\[
\frac{q^{n+1}}{\cron+q^{n+1}},\qquad
\frac{-q^{2n+1}}{\cron},\qquad
\frac{-q^{2n+1}}{\cron+q^{n+1}},\qquad
\frac{q^{n+1}}{1}.
\]
Their numerator coefficients give the successive values $1,-1,-1,1$ for $v_j$. 
Applying the recurrence with initial value $k_{6n-7}=0$ gives the corresponding values $n-1,n,n-1,0$ for $k_j$.
This proves the formulas for $(v_j,k_j)$ in part~(\rmnum{4}).

\medskip

\noindent{\it{Proof of the coefficients $c_{j,k_{j-1}}$.}}
Let $B_j(q)$ denote the denominator at index $j$ in \eqref{met:Hfrac-Phin-n+1}. 
For $j\ge1$, the standard form gives $B_{j-1}(q)=1-qu_j(q)$, and hence  
\[
c_{j,k_{j-1}}
=[q^{k_{j-1}}]u_j(q)
=-[q^{k_{j-1}+1}]B_{j-1}(q).
\]
The first denominator is
\[
B_0(q)=\cron+q^{n+1}
=1+\sum_{h=2}^{n-1}q^h+2q^n.
\]
Since $k_0=n-1$, this gives $c_{1,n-1}=-2$.
The second denominator is $B_1(q)=1$, so $c_{2,0}=0$.
These are the two coefficients in part~(\rmnum{1}).

The $q$-integer denominators in $W_n$ and $U_n$ have the form $[r]_q$ or $[r]_q-q$ with $r\ge2$.  
In the triples of $W_n$ and $U_n$, we have $r=n-i$ or $r=i+2$.
In both cases, the corresponding exponent is $k_{j-1}=r-2$, and the polynomials $u_j(q)$ are given by
\[
\frac{1-[r]_q}{q}=-[r-1]_q, \qquad
\frac{1-([r]_q-q)}{q}=1-[r-1]_q.
\]
The coefficient of $q^{r-2}$ in the first polynomial is $-1$.
In the second, it is $-1$ for $r\ge3$ and $0$ for $r=2$.
Thus the denominator $[i+2]_q-q$ gives $-1+\pi(i=0)$ in part~(\rmnum{3}).
The denominators $1-q$ and $1$ occur with $k_{j-1}=0$ and give the coefficients $1$ and $0$, respectively.
In particular, the last denominator of $U_n$ gives $c_{6n-6,0}=0$, as listed in part~(\rmnum{4}).

The first and third denominators of $V_n$ are both $\cron+q^{n+1}$, with corresponding exponent $n-1$.
Hence $c_{6n-5,n-1}=c_{6n-3,n-1}=-2$.
The second denominator is $\cron=B_0(q)-q^{n+1}$.
Since $k_{6n-5}=n$, we obtain $c_{6n-4,n}=1$.
The final denominator is $1$, so $c_{6n-2,0}=0$.

Finally, since $k_{6n-3}=k_1=0$, the recurrence for $k_j$ starts from the same value when the block $W_nU_nV_n$ repeats.
Thus the pairs $(v_j,k_j)$ repeat with period $6n-4$ from $j=2$ onward. 
The denominators repeat from the same index, so the coefficient formula gives the same periodicity for $c_{j,k_{j-1}}$ from $j=3$ onward. 
The equality $c_{6n-2,0}=c_{2,0}=0$ covers $j=2$ as well.
\end{proof}

\medskip
In Step~2, we use the parameters from Lemma~\ref{met:lem-step1} to obtain the following values of $H_{s_p}\bigl(\Phi_n^{(n+1)}\bigr)$ and $R_p$.

\begin{Lemma}\label{met:lem-step2}
Let $n\ge3$. Then the following assertions hold.
\begin{itemize}
\item[\textnormal{(a)}] For every $p\ge0$, we have
\begin{equation}\label{met:Hankel-Phin+1-pm1}
H_{s_p}\bigl(\Phi_n^{(n+1)}\bigr) \in \{-1, 1\},
\end{equation}
and $H_j\bigl(\Phi_n^{(n+1)}\bigr)=0$ for every $j\notin\mathcal{S}$.

\item[\textnormal{(b)}]
The auxiliary sequence $(R_p)_{p\ge0}$ defined by \eqref{def:Rn} is $(6n-4)$-periodic. For every $p\ge0$, we have
\begin{equation}\label{met:Rp-bound}
R_p\in\{-2,-1,0,1, 2\}.
\end{equation}
\end{itemize}
\end{Lemma}

\begin{proof}
By Lemma~\ref{met:lem-step1}, every $v_i$ belongs to $\{-1, 1\}$  and $v_0=1$.
Thus, for $p\ge1$, the fundamental formula \eqref{eq:HankelDetFundamental} gives
\[
H_{s_p}\bigl(\Phi_n^{(n+1)}\bigr)=(-1)^{\epsilon_p}\,v_0^{s_p}\,v_1^{s_p-s_1}\cdots v_{p-1}^{s_p-s_{p-1}}=(-1)^{\epsilon_p+\nu_p},
\]
where $\nu_p:=\sum_{\substack{0\le i<p\\v_i=-1}}(s_p-s_i)$.
This proves \eqref{met:Hankel-Phin+1-pm1} for $p\ge1$.
The case $p=0$ follows from $s_0=0$ and the empty determinant convention. 
For every nonnegative integer $j\notin\mathcal S$, Theorem~\ref{th:super2}(ii) gives $H_j\bigl(\Phi_n^{(n+1)}\bigr)=0$.
This proves part~\textnormal{(a)}.

\medskip
We now prove part~\textnormal{(b)}.
By \eqref{def:Rn}, we have $R_0=1$, $R_1=c_{1,k_0}=-2$, and, for $p\ge1$,
\begin{equation}\label{met:recR}
R_{p+1}=c_{p+1,k_p}R_p-v_pR_{p-1}.
\end{equation}
Since $c_{2,0}=0$, $k_1=0$, and $v_1=1$ by Lemma~\ref{met:lem-step1}, we obtain $R_2=-1$.

We first consider the triples in $W_n$.
For $0\le i\le n-3$, set $p=3i+2$.
The three successive values of $c_{p+1,k_p}$, as the index runs from $p$ to $p+2$, are $-1$, $-1$, $1$, and we have $v_p=v_{p+1}=v_{p+2}=-1$.
Three applications of \eqref{met:recR} give $(R_{p+2},R_{p+3})=(2R_p-R_{p-1},R_p)$.
Since $(R_1,R_2)=(-2,-1)$, the pairs $(R_{3i+1},R_{3i+2})$ alternate between $(-2,-1)$ and $(0,-1)$.
Hence, for $0\le i\le n-3$, we have
\begin{equation}\label{met:R-Wn-triple}
(R_{3i+2},R_{3i+3},R_{3i+4},R_{3i+5})
=
\begin{cases}
(-1,-1,0,-1) & \text{if $i$ is even},\\
(-1,1,-2,-1) & \text{if $i$ is odd}.
\end{cases}
\end{equation}
Taking $i=n-3$ gives $R_{3n-4}=-1$.
At the final index $p=3n-4$ of $W_n$, we have $c_{3n-3,0}=0$ and $v_{3n-4}=-1$. Therefore,
\[
R_{3n-3}=R_{3n-5}
=
\begin{cases}
0  & \text{if $n$ is odd},\\
-2 & \text{if $n$ is even}.
\end{cases}
\]

We next consider $U_n$. For $0\le i\le n-3$, set $p=3n-3+3i$. 
The recurrence coefficients in each triple are the same as in $W_n$, except that the middle coefficient of the first triple is $c_{3n-1,0}=0$.
For $i=0$, we have $R_{p-1}=-1$, so \eqref{met:recR} gives $R_{p+1}=-1-R_p$ and $(R_{p+2},R_{p+3})=(R_p,-1)$.
By the calculation at the end of $W_n$, this last pair is $(0,-1)$ if $n$ is odd and $(-2,-1)$ if $n$ is even.
The input pairs for any subsequent triples alternate as in $W_n$. 
Hence
\begin{equation}\label{met:R-Un-triples}
(R_p,R_{p+1},R_{p+2},R_{p+3})
=
\begin{cases}
(0,-1,0,-1)  & \text{if $i=0$ and $n$ is odd},\\
(-2,1,-2,-1) & \text{if $i=0$ and $n$ is even},\\
(-1,1,-2,-1) & \text{if $i\ge1$ and $n+i$ is even},\\
(-1,-1,0,-1) & \text{if $i\ge1$ and $n+i$ is odd}.
\end{cases}
\end{equation}
For $i=n-3$, the first case of \eqref{met:R-Un-triples} applies when $n=3$, and the last case applies when $n\ge4$ since $n+i=2n-3$ is odd. 
In both cases, $(R_{6n-10},R_{6n-9})=(0,-1)$.
At the final three indices $p=6n-9,6n-8,6n-7$ of $U_n$, the coefficient pairs $(c_{p+1,k_p},v_p)$ are
\begin{align*}
(c_{6n-8,0},v_{6n-9})&=(-1,-1),
  & (c_{6n-7,n-2},v_{6n-8})&=(-1,-1),\\
(c_{6n-6,0},v_{6n-7})&=(0,-1).
\end{align*}
Applying \eqref{met:recR} with these coefficients gives
\[
(R_{6n-8},R_{6n-7},R_{6n-6})=(1,-2,1).
\]

Finally, for $6n-6\le p\le6n-3$, the coefficient pairs $(c_{p+1,k_p},v_p)$ in $V_n$ are
\begin{align*}
(c_{6n-5,n-1},v_{6n-6})&=(-2,1),
  & (c_{6n-4,n},v_{6n-5})&=(1,-1),\\
(c_{6n-3,n-1},v_{6n-4})&=(-2,-1),
  & (c_{6n-2,0},v_{6n-3})&=(0,1).
\end{align*}
Starting from $(R_{6n-7},R_{6n-6})=(-2,1)$, four applications of \eqref{met:recR} give
\[
(R_{6n-5},R_{6n-4},R_{6n-3},R_{6n-2})=(0,1,-2,-1).
\]

The last three values give
\[
(R_{6n-4},R_{6n-3},R_{6n-2})=(1,-2,-1)=(R_0,R_1,R_2).
\]
By Lemma~\ref{met:lem-step1}, the coefficients in \eqref{met:recR} repeat with period $6n-4$ from $p=2$ onward.
Since $(R_{6n-3},R_{6n-2})=(R_1,R_2)$, induction using \eqref{met:recR} gives $R_{p+6n-4}=R_p$ for every $p\ge1$.
The equality $R_{6n-4}=R_0$ covers $p=0$ as well.
Thus $(R_p)_{p\ge0}$ is $(6n-4)$-periodic.
Since all the values computed over one period belong to  $\{-2,-1,0,1,2\}$, the same holds for every $p\ge0$, proving \eqref{met:Rp-bound}.
\end{proof}

\subsubsection{Computations for $n=2$}\label{met:subsec-n2}

When $n=2$, the range $0\le i\le n-3$ is empty, so the indexed triples in $U_n$ and $W_n$ do not occur.
In this case, $\langle 2\rangle_q=1+2q^2-q^3$, and the three blocks are given by
\begin{align*}
U_2(q)&=\frac{-q^2}{[2]_q}\cfm\frac{-q^2}{[2]_q-q}\cfm\frac{-q^2}{1},\qquad
W_2(q)=\frac{-q^2}{1},\\
V_2(q)&=\frac{q^{3}}{1+2q^2}\cfm\frac{-q^{5}}{1+2q^2-q^3}\cfm\frac{-q^{5}}{1+2q^2}\cfm\frac{q^{3}}{1},
\end{align*}
and the period $6n-4=8$ consists of the three blocks of $U_2$, the four blocks of $V_2$ and the single block of $W_2$. 
Taking $n=2$ in Proposition~\ref{met:prop-HF-n+1}, we obtain
\begin{equation}\label{met:Hfrac-Phi2-3}
\begin{aligned}
\Phi_2^{(3)}(q)=\frac{q}{1+2q^2}\cfm\frac{q^{3}}{1}\cfm\Bigl(&\frac{-q^2}{1}\cfm\frac{-q^2}{1+q}\cfm\frac{-q^2}{1}\cfm\frac{-q^2}{1}\\
&\cfm\frac{q^{3}}{1+2q^2}\cfm\frac{-q^{5}}{1+2q^2-q^3}\\
&\cfm\frac{-q^{5}}{1+2q^2}\cfm\frac{q^{3}}{1}\cfm\Bigr)^*.
\end{aligned}
\end{equation}

We now use \eqref{met:Hfrac-Phi2-3} to carry out Steps~1 and~2.
The following lemma gives the $H$-fraction parameters and the resulting values of $H_{s_p}\bigl(\Phi_2^{(3)}\bigr)$ and $R_p$.

\begin{Lemma}\label{met:lem-n2}
Consider the $H$-fraction \eqref{met:Hfrac-Phi2-3} of $\Phi_2^{(3)}$ in the standard form \eqref{met:Hfrac-Phin-shift}.
For $0\le p\le9$, the following table gives the parameters $(v_p,k_p)$, the coefficients $c_{p+1,k_p}$, the indices $s_p$, and the terms $R_p$ of the sequence defined by \eqref{def:Rn}:
\[
\begin{array}{c|cccccccccc}
p           & 0  & 1  & 2  & 3  & 4  & 5  & 6  & 7  & 8  & 9  \\ \hline
v_p         & 1  & 1  & -1 & -1 & -1 & -1 & 1  & -1 & -1 & 1  \\
k_p         & 1  & 0  & 0  & 0  & 0  & 0  & 1  & 2  & 1  & 0  \\
c_{p+1,k_p} & -2 & 0  & 0  & -1 & 0  & 0  & -2 & 1  & -2 & 0  \\
s_p         & 0  & 2  & 3  & 4  & 5  & 6  & 7  & 9  & 12 & 14 \\
R_p         & 1  & -2 & -1 & -2 & 1  & -2 & 1  & 0  & 1  & -2
\end{array}
\]
From $p=2$ onward, the triples $(v_p,k_p,c_{p+1,k_p})$ repeat with period~$8$. For every $p\ge0$, we have $s_{p+8}=s_p+12$ and $R_{p+8}=R_p$. 
In particular, the conclusions of Lemma~\ref{met:lem-step2} extend to $n=2$.
\end{Lemma}

\begin{proof}
As in the proof of Lemma~\ref{met:lem-step1}, we compare \eqref{met:Hfrac-Phi2-3} with the standard form \eqref{met:Hfrac-Phin-shift}. 
The first numerator is $q$, so $(v_0,k_0)=(1,1)$. For $p\ge1$, the numerator at index $p$ is a monomial with coefficient $v_p$. 
If its degree is $a_p$, then $k_p=a_p-k_{p-1}-2$. This gives the pairs $(v_p,k_p)$
in the table.

Let $B_p(q)$ denote the denominator at index $p$ in \eqref{met:Hfrac-Phi2-3}. For the denominators $1$, $1+q$, $1+2q^2$, and $1+2q^2-q^3$, the corresponding polynomials
$u_{p+1}(q)=(1-B_p(q))/q$ are $0$, $-1$, $-2q$, and $-2q+q^2$, respectively. 
Taking the coefficient of $q^{k_p}$ gives the stated values of $c_{p+1,k_p}$.
Since $k_9=k_1=0$, the recurrence for $k_p$ starts from the same value when the eight-term block repeats. 
Thus the triples $(v_p,k_p,c_{p+1,k_p})$ repeat with period~$8$ from $p=2$ onward.

Starting from $s_0=0$, the relation $s_{p+1}=s_p+k_p+1$ gives the row for $s_p$. Since $k_2+\cdots+k_9=4$, for every $p\ge2$ we have
\[
s_{p+8}-s_p=\sum_{i=p}^{p+7}(k_i+1)=12.
\]
The same identity holds for $p=0,1$, since $s_8-s_0=s_9-s_1=12$.

The last row follows from \eqref{met:recR} with $R_0=1$ and $R_1=c_{1,1}=-2$. One further application gives
\[
R_{10}=c_{10,0}R_9-v_9R_8=-1=R_2.
\]
Thus $(R_9,R_{10})=(R_1,R_2)$. Since the coefficients in \eqref{met:recR} repeat with period~$8$ from $p=2$ onward, induction gives $R_{p+8}=R_p$ for every $p\ge1$.
The equality $R_8=R_0$ covers $p=0$ as well.
All the values in one period belong to $\{-2,-1,0,1\}$, so the same holds for every $R_p$.

Finally, since every $v_p$ belongs to $\{-1,1\}$, the fundamental formula \eqref{eq:HankelDetFundamental} gives $H_{s_p}\bigl(\Phi_2^{(3)}\bigr)\in\{-1,1\}$.
For every nonnegative integer $j\notin\mathcal S$, Theorem~\ref{th:super2}(ii) gives
$H_j\bigl(\Phi_2^{(3)}\bigr)=0$.
\end{proof}

\begin{Remark}
The table agrees with the formulas of Lemma~\ref{met:lem-step1} formally evaluated at $n=2$, except for the coefficients $c_{3,0}$ and $c_{5,0}$, which are $0$ instead of $1$ and $-1$, respectively. 
The first exception arises because the last indexed triple of $U_n$ is absent when $n=2$, and the term at the corresponding index belongs to $W_2$ and has denominator $1$. 
The second arises because the denominator $[n]_q-q$ of the penultimate term of $U_n$ becomes $1$ when $n=2$.
\end{Remark}

\subsubsection{Bounds and (anti-)periodicity}\label{met:subsec-bounds}

To carry out Step~3, we substitute the data obtained in Steps~1 and~2 into \eqref{met:key-formula} and \eqref{met:key-formula-bis}. This yields the following lemma.

\begin{Lemma}\label{met:lem-step3}
For every $n\geq 2$ and every $j\in\{0,1,\ldots,2n(n+1)-1\}$,
\[
H_j^{(n+2)}\in\{-2,-1,0,1,2\}.
\]
\end{Lemma}

\begin{proof}
Fix $n\ge2$ and $0\le j<2n(n+1)$. 
Recall that $\mathcal{S} = \{s_p\mid p\geq 0\}$ and $\mathcal{S}+k = \{s_p+k_p\mid p\geq 0\}$
Since $s_p+k_p=s_p$ when $k_p=0$, the following three cases cover all possible indices $j$.

\noindent\textbf{Case 1.} Suppose that $j=s_p$.
By Lemma~\ref{met:lem-step2} for $n\ge3$ and Lemma~\ref{met:lem-n2} for $n=2$, we have \[
H_{s_p}\bigl(\Phi_n^{(n+1)}\bigr)=\pm1  \qquad \text{and}  \qquad R_p\in\{-2,-1,0,1,2\}.
\]
Plugging these into the key formula~\eqref{met:key-formula}, we obtain
\[
H_{s_p}^{(n+2)}=(-1)^{s_p+p}\,H_{s_p}\bigl(\Phi_n^{(n+1)}\bigr)\,R_p\in\{-2,-1,0,1,2\}.
\]

\noindent\textbf{Case 2.} Suppose that $j=s_p+k_p$ with $k_p\geq 1$.
By Lemma~\ref{met:lem-step1} for $n\ge3$ and Lemma~\ref{met:lem-n2} for $n=2$, every $v_i$ belongs to $\{-1,1\}$.
Hence the prefactor in \eqref{met:key-formula-bis} is $\pm1$, and Case~1 gives
\[
H_{s_p+k_p}^{(n+2)}= (-1)^{\frac{(k_p-1)k_p}{2}} \bigl(v_0v_1\cdots v_p\bigr)^{k_p} H_{s_p}^{(n+2)}\in\{-2,-1,0,1,2\}.
\]

\medskip

\noindent\textbf{Case 3.} Suppose that $j\notin\mathcal{S}\cup(\mathcal{S}+k)$.
By Lemma~\ref{met:lem-reduction} and the vanishing assertion of Theorem~\ref{th:shiftHDet}, we have
\[
H_j^{(n+2)}=H_j^{(1)}\bigl(\Phi_n^{(n+1)}\bigr)=0.
\]
This completes the proof.
\end{proof}

\medskip
We are now ready to prove Theorem~\ref{met:thm-main}.

\begin{proof}[Proof of Theorem~\ref{met:thm-main}]
The case $n=1$ is already known, so let $n\ge2$.
By Lemma~\ref{met:lem-step1} (for $n\ge3$) and Lemma~\ref{met:lem-n2} (for $n=2$),
the triples $(v_p,k_p,c_{p+1,k_p})$ repeat with period $6n-4$ from $p=2$ onward.
The explicit data also give
\[
k_{6n-4}=k_0=n-1,\qquad
k_{6n-3}=k_1=0,\qquad
v_{6n-3}=v_1=1,
\]
and
\[
\sum_{i=2}^{6n-3}(k_i+1)= 2n(n+1),
\qquad
\prod_{i=2}^{6n-3}v_i=(-1)^{6n-6}=1.
\]
Thus, using $s_{p+1}=s_p+k_p+1$ and $V_p=v_0v_1\cdots v_p$, we obtain, for $p \ge 0$,
\begin{equation}\label{met:eq-periodic-params}
s_{p+6n-4}=s_p + 2n(n+1),\quad
k_{p+6n-4}=k_p,\quad
V_{p+6n-4}=V_p.
\end{equation}
Moreover, Lemmas~\ref{met:lem-step2} and~\ref{met:lem-n2} give $R_{p+6n-4}=R_p$ for $p\ge0$.

By definition, $H_j\bigl(\Phi_n^{(n+1)}\bigr)=H_j^{(n+1)}(\Phi_n)$.
The (anti-)periodicity of this sequence established in \cite{Han2025Pedon}, together with $s_{p+6n-4}=s_p+2n(n+1)$, gives
\begin{equation*}
H_{s_{p+6n-4}}\bigl(\Phi_n^{(n+1)}\bigr)
=
(-1)^nH_{s_p}\bigl(\Phi_n^{(n+1)}\bigr)
\end{equation*}
for $p\ge0$. 
Since $6n-4$ and $2n(n+1)$ are even, substituting these identities into \eqref{met:key-formula}, we obtain, for every $p\ge0$,
\[
H_{s_p+2n(n+1)}^{(n+2)}=  (-1)^nH_{s_p}^{(n+2)}.
\]
For the indices $j=s_p+k_p$ with $k_p\ge1$, \eqref{met:key-formula-bis} and \eqref{met:eq-periodic-params} give
\[
H_{s_p+k_p+2n(n+1)}^{(n+2)}  = (-1)^nH_{s_p+k_p}^{(n+2)}.
\]

Finally, if $j\notin\mathcal S\cup(\mathcal S+k)$, then $j=s_p+\ell$ for some $p\ge0$ and $1\le\ell\le k_p-1$.
By \eqref{met:eq-periodic-params}, we have $j+2n(n+1)=s_{p+6n-4}+\ell$ and $ 1\le\ell\le k_{p+6n-4}-1$,
so $j+ 2n(n+1) \notin\mathcal S\cup(\mathcal S+k)$ as well.
Lemma~\ref{met:lem-reduction} and the vanishing assertion of Theorem~\ref{th:shiftHDet} therefore give $H_j^{(n+2)}=H_{j+2n(n+1)}^{(n+2)}=0$.
Thus \eqref{met:eq-antiper-n2} holds for every $j\ge0$.
Together with Lemma~\ref{met:lem-step3}, this gives $H_j^{(n+2)}\in\{-2,-1,0,1,2\}$ for all $j\ge0$ and completes the proof.
\end{proof}

\begin{Remark}
The reason why the determinants $H_j^{(n+2)}$ can take the values $\pm 2$, rather than only $\pm 1$ and $0$ as in the cases $\ell\le n+1$, lies in the auxiliary sequence $(R_p)_{p\ge0}$. 
Indeed, when $\ell\le n+1$, no such auxiliary sequence appears in the determinant formula: each nonzero $H_j^{(\ell)}$ is given, up to a sign, by a product of powers of the numerator coefficients of the corresponding $H$-fraction, all equal to $\pm1$. Hence the determinant values are $0,\pm1$.

For $\ell=n+2$, formula~\eqref{met:key-formula} introduces the extra factor $R_p$. 
For $n\ge2$, the computations using \eqref{def:Rn} show that $R_p$ can equal $-2$, but has no larger absolute value, allowing the determinants to take the values $\pm2$. 
The case $n=1$ is special: the determinant values remain in $\{-1,0,1\}$, as proved in \cite{Ovsienko2025Pedon}.
\end{Remark}


\section{Application: the $(n+3)$-shifted Hankel determinants of $q$-deformed metallic numbers}\label{sec:metallicn3}

We keep the notation of Section~\ref{sec:metallicn2}: $\Phi_n(q)$ is the $q$-deformed metallic number \eqref{met:defPhin}, $H_j^{(\ell)}=H_j^{(\ell)}(\Phi_n)$ its shifted Hankel determinants, and $(v_p,k_p,u_p)_{p\geq 0}$, $(s_p)$, $(\epsilon_p)$, $(R_p)$ the data attached to the $H$-fraction \eqref{met:Hfrac-Phin-shift} of $\Phi_n^{(n+1)}$, which are given explicitly by Lemmas~\ref{met:lem-step1} and~\ref{met:lem-step2} for $n\geq 3$, and in Lemma~\ref{met:lem-n2} for $n=2$.

We now turn to the case $\ell=n+3$ of Conjecture~\ref{met:conj-HKP}\,(b): we prove that the sequence $H^{(n+3)}(\Phi_n)$ is unbounded, in sharp contrast with the bounded value sets $\{-1,0,1\}$ and $\{-2,-1,0,1,2\}$ that occur for $\ell\leq n+2$ (\cite{Han2025Pedon} and Theorem~\ref{met:thm-main}). The strategy parallels that of Section~\ref{met:sec-proof}: a reduction to a second-shift question on $\Phi_n^{(n+1)}$, followed by an application of the second shifted Hankel determinant formula of Theorem~\ref{th:shift2HDet}.

\begin{Theorem}\label{met:thm-n+3}
For every $n\geq 2$, the sequence $H^{(n+3)}(\Phi_n)$ is unbounded. 
More precisely, there exists a positive integer $C_n$ such that, setting
\[
j_m:=s_{6n-3}+m\cdot 2n(n+1)\qquad (m\geq 0),
\] 
where $s_{6n-3}=n(2n+3)$, we have
\begin{equation}\label{met:eq-n+3-growth}
\bigl|H_{j_m}^{(n+3)}(\Phi_n)\bigr|\geq C_n\,m\qquad\text{for every } m\geq 1.
\end{equation}
\end{Theorem}

The case $n=1$, which is not covered by Theorem~\ref{met:thm-n+3}, 
has been proved by E.~Pedon \cite{Pedon2026private}: the sequence 
$H^{(4)}(\Phi_1)$ is also unbounded.

\subsection{The $2$-shift formula and auxiliary sequences}\label{met:subsec-reduction2}

The following lemma expresses the $(n+3)$-shifted Hankel determinants of $\Phi_n$ as the $2$-shifted Hankel determinants of $\Phi_n^{(n+1)}$.

\begin{Lemma}\label{met:lem-reduction2}
For every $n\geq 1$ and every $j\geq 0$,
\begin{equation}\label{met:eq-reduction2}
H_j^{(n+3)}(\Phi_n)=H_j^{(2)}\bigl(\Phi_n^{(n+1)}\bigr).
\end{equation}
\end{Lemma}

\begin{proof}
Write $\Phi_n(q)=\sum_{i\ge0}f_iq^i$.
By definition, the coefficient of $q^i$ in
$\Phi_n^{(n+1)}$ is $f_{i+n+1}$.
Thus, for every $j\ge1$, \eqref{met:HankDet} gives
\[
H_j^{(2)}\bigl(\Phi_n^{(n+1)}\bigr)
=
\det\bigl(f_{a+b+n+3}\bigr)_{0\le a,b\le j-1}
=
H_j^{(n+3)}(\Phi_n).
\]
For $j=0$, both sides equal $1$ by the empty determinant convention.
\end{proof}

By Lemma~\ref{met:lem-reduction2}, proving Theorem~\ref{met:thm-n+3} amounts to establishing the stated linear lower bound for the $2$-shifted Hankel determinants of $\Phi_n^{(n+1)}$. 
We therefore apply Theorem~\ref{th:shift2HDet} to this series.

Throughout the remainder of this section, let $n\ge2$.
We use the parameters $v_p$, $k_p$, $ c_{p+1,k_p}$, and $u_{p+1}(q)$ of the $H$-fraction \eqref{met:Hfrac-Phin-shift} of $\Phi_n^{(n+1)}$, together with the associated sequences $(s_p)$, $(\epsilon_p)$, and $(R_p)$. 
As before, $c_{p+1,r}$ denotes the coefficient of $q^r$ in $u_{p+1}(q)$.
Recall from Theorem~\ref{th:shift2HDet} that $d_{p+1}:=[q^{k_p}]\bigl(1-u_{p+1}(q)q\bigr)$
for $p\ge0$. Equivalently,
\[
d_{p+1}
=
\begin{cases}
1, & \text{if } k_p=0,\\
-c_{p+1,k_p-1}, & \text{if } k_p\ge1.
\end{cases}
\]
By Lemma~\ref{met:lem-step1} for $n\ge3$ and Lemma~\ref{met:lem-n2} for $n=2$, we have $v_p\in\{-1,1\}$ and $0\le k_p\le n$ for every $p\ge0$.

For every $p\ge0$, define
\begin{equation*}
T_p := \sum_{t=0}^{p} d_{t+1}R_t^2\,v_{t+1}v_{t+2}\cdots v_p,
\end{equation*}
where the product is $1$ when $t=p$.
For every $p\ge0$ with $k_p=0$, Lemma~\ref{met:lem-reduction2} and formula~\eqref{eq:H1sj:Hsj:k0} of Theorem~\ref{th:shift2HDet} give
\begin{equation*}
H_{s_p}^{(n+3)}(\Phi_n)=H_{s_p}^{(2)}\bigl(\Phi_n^{(n+1)}\bigr)=H_{s_p}\bigl(\Phi_n^{(n+1)}\bigr)\,T_p.
\end{equation*}
By Lemma~\ref{met:lem-step2} for $n\ge3$ and Lemma~\ref{met:lem-n2} for $n=2$, the factor $H_{s_p}\bigl(\Phi_n^{(n+1)}\bigr)$ belongs to $\{-1,1\}$.
 Hence
\begin{equation}\label{met:eq-modulus}
\bigl|H_{s_p}^{(n+3)}(\Phi_n)\bigr|=|T_p|\qquad\text{whenever } k_p=0,
\end{equation}
Thus, to prove Theorem~\ref{met:thm-n+3}, it suffices to establish the required linear lower bound for $|T_p|$ along the indices $p=6n-3+m(6n-4)$.

We first record the following recurrences, which are valid for all $n\ge2$. Let $p \ge 0$.
The definition of $T_p$, together with $R_0=1$, gives
\begin{equation}\label{met:eq-Tp-rec}
T_0=d_1,\qquad
T_{p+1}=v_{p+1}T_p+d_{p+2}R_{p+1}^2.
\end{equation}
To simplify this recurrence, set
\begin{equation*}
\sigma_p:=v_1v_2\cdots v_p,\qquad
S_p:=\frac{T_p}{\sigma_p},
\end{equation*}
with $\sigma_0=1$.
By Lemma~\ref{met:lem-step1} for $n\ge3$ and Lemma~\ref{met:lem-n2} for $n=2$, every $v_p$ belongs to $\{-1,1\}$.
Thus $\sigma_p\in\{-1,1\}$ and $1/\sigma_p=\sigma_p$.
Dividing \eqref{met:eq-Tp-rec} by $\sigma_{p+1}=v_{p+1}\sigma_p$ therefore gives
\begin{equation}\label{met:eq-Sp-rec}
S_0=d_1,\qquad
S_{p+1}=S_p+d_{p+2}R_{p+1}^2\sigma_{p+1}.
\end{equation}
We denote the successive increments by $S_{p+1}-S_p
=d_{p+2}R_{p+1}^2\sigma_{p+1}.$
Since $\sigma_p\in\{-1,1\}$, we also have $|S_p|=|T_p|$.
We carry out these computations separately for $n\ge3$ and $n=2$.

\subsection{Computations for $n\ge3$}\label{met:subsec-shift-n+3-step1}

In this subsection, we compute the period increment of $S_p$ and the value $S_{6n-3}$ for $n\ge3$. We begin by determining the signs $\sigma_p$ from the parameters $v_p$ in Lemma~\ref{met:lem-step1}.

\begin{Lemma}\label{met:lem-shift2-step1}
Let $n\ge3$. The sign sequence $(\sigma_p)_{p\ge0}$ satisfies $\sigma_0=1$ and $\sigma_p=(-1)^{p-1}$ for $1\le p\le6n-7 $. For the indices corresponding to $V_n$, we have $\sigma_{6n-6}=\sigma_{6n-4}=\sigma_{6n-3}=1$ and $\sigma_{6n-5}=-1$.
Moreover $\sigma_p$ is $(6n-4)$-periodic for $p\geq 1$.
\end{Lemma}

\begin{proof}
By Lemma~\ref{met:lem-step1}, we have $v_1=1$ and $v_j=-1$ for $2\le j\le6n-7$.
Hence $\sigma_0=\sigma_1=1$ and $\sigma_p=(-1)^{p-1}$ for $1\le p\le6n-7$.
The values corresponding to $V_n$ follow from $v_{6n-6}=v_{6n-3}=1$ and $v_{6n-5}=v_{6n-4}=-1$.
Since $(v_j)_{j\ge2}$ is $(6n-4)$-periodic, for every $p\ge1$ we have
\[
\sigma_{p+6n-4}
=
\sigma_p\prod_{j=p+1}^{p+6n-4}v_j
=
\sigma_p\prod_{j=2}^{6n-3}v_j.
\]
The last product contains $3n-5$ negative factors from $W_n$, $3n-3$ from $U_n$, and $2$ from $V_n$.
Thus it equals $(-1)^{6n-6}=1$, proving the required periodicity.
\end{proof}

The coefficients $d_j$ are obtained from the denominators of the $H$-fraction \eqref{met:Hfrac-Phin-shift}. Their values are given in the following lemma.

\begin{Lemma}\label{met:lem-d-values}
Let $n\ge3$. The sequence $(d_j)_{j\ge1}$ takes values in $\{0,1,2\}$ and is $(6n-4)$-periodic from $j=3$ onward. Its initial values and one complete repeating block are given as follows.
\begin{itemize}
\item[(\rmnum{1})] We have $d_1=d_2=1$.

\item[(\rmnum{2})] For the block $W_n$ and $0\le i\le n-3$, we have $d_{3i+4}=d_{3i+5}=1$ and
\[
d_{3i+3}
=
\begin{cases}
0, & \text{if }i=n-3,\\
1, & \text{otherwise}.
\end{cases}
\]
The final term of $W_n$ gives $d_{3n-3}=1$.

\item[(\rmnum{3})] For the block $U_n$ and $0\le i\le n-3$, we have $d_{3n-2+3i}=d_{3n+3i}=1$ and
\[
d_{3n-1+3i}
=
\begin{cases}
0, & \text{if }i=1,\\
1, & \text{otherwise}.
\end{cases}
\]
The final three terms of $U_n$ give $d_{6n-8}=d_{6n-6}=1$ and
\[
d_{6n-7}
=
\begin{cases}
0, & \text{if }n=3,\\
1, & \text{if }n\ge4.
\end{cases}
\]

\item[(\rmnum{4})] The four terms of $V_n$ give $d_{6n-5}=d_{6n-3}=d_{6n-2}=1$ and $d_{6n-4}=2$.
\end{itemize}
\end{Lemma}

\begin{proof}
For every $p\ge0$, let $B_p(q)=1-u_{p+1}(q)q$ denote the denominator at index $p$ in \eqref{met:Hfrac-Phin-shift}.
Then $d_{p+1}=[q^{k_p}]B_p(q)$.
Reading off the relevant coefficients from the denominators listed in the proof of Lemma~\ref{met:lem-step1} gives the stated values.
The initial terms have $k_0=n-1$ and $k_1=0$, with denominators $B_0(q)=\cron+q^{n+1}$ and $B_1(q)=1$, respectively. Since $n\ge3$, we obtain $d_1=[q^{n-1}](\cron+q^{n+1})=1$ and $d_2=[q^0]1=1$.
We detail the cases in which the coefficient is zero.

\medskip
\noindent\textbf{Case 1.}
For the first term of the triple indexed by $i$ in $W_n$, we have $B_{3i+2}(q)=[n-i]_q-q=1+q^2+\cdots+q^{n-i-1}$ and $k_{3i+2}=n-2-i$.
Thus  $d_{3i+3}=[q^{n-2-i}]B_{3i+2}(q)=1$ when $i\le n-4$.
When $i=n-3$, the exponent is $1$ and the denominator has no $q$-term, so $d_{3i+3}=0$.

\medskip
\noindent\textbf{Case 2.}
For the second term of the triple indexed by $i$ in $U_n$, we have $B_{3n-2+3i}(q)=[i+2]_q-q$ and $k_{3n-2+3i}=i$.
The denominator equals $1$ when $i=0$ and $1+q^2+\cdots+q^{i+1}$ when $i\ge1$.
Hence $d_{3n-1+3i}=[q^i]B_{3n-2+3i}(q)=1$ for $i=0$ or $i\ge2$, whereas $d_{3n-1+3i}=0$ for $i=1$.

\medskip
\noindent\textbf{Case 3.}
For the penultimate term of $U_n$, we have $B_{6n-8}(q)=[n]_q-q=1+q^2+\cdots+q^{n-1}$ and $k_{6n-8}=n-2$.
Thus $d_{6n-7}=[q^{n-2}]B_{6n-8}(q)=1$ for $n\ge4$, whereas $d_{6n-7}=[q^1](1+q^2)=0$ for $n=3$.

\medskip

The remaining values follow by the same coefficient extraction.
In particular, $B_{6n-5}(q)=\cron$ and $k_{6n-5}=n$, so $d_{6n-4}=[q^n]\cron=2$.
By Proposition~\ref{met:prop-HF-n+1} and Lemma~\ref{met:lem-step1}, the pairs $(B_p(q),k_p)$ repeat with period $6n-4$ from $p=2$ onward.
Hence $d_{j+6n-4}=d_j$ for every $j\ge3$.
\end{proof}

By Lemmas~\ref{met:lem-step2}, \ref{met:lem-shift2-step1}, and~\ref{met:lem-d-values}, the triples $(\sigma_p,d_{p+1},R_p)$ are $(6n-4)$-periodic for $p\ge2$. 
Since the increments in \eqref{met:eq-Sp-rec} are
$(6n-4)$-periodic from index $2$ onward, the total increment
of $S_p$ over one complete period is independent of the
starting index $p\ge1$. Denoting this increment by $C_n$,
we obtain
\begin{equation}\label{met:def-Cn-sum}
C_n:=S_{p+6n-4}-S_p
=\sum_{q=2}^{6n-3}d_{q+1}R_q^2\sigma_q,
\end{equation}
which is independent of $p\ge1$.

\begin{Lemma}\label{met:lem-shift2-step2}
For every $n\geq 3$, $C_n=4n+2$. In particular $C_n>0$.
\end{Lemma}

\begin{proof}
We evaluate the sum in \eqref{met:def-Cn-sum} by computing the contributions of the blocks $V_n$, $W_n$, and $U_n$.

\medskip
\noindent{\bf{Case 1.}} Consider the block $V_n$.
By Lemma~\ref{met:lem-step1}, the successive coefficient pairs are
\begin{align*}
(c_{6n-5,n-1},v_{6n-6})&=(-2,1),
& (c_{6n-4,n},v_{6n-5})&=(1,-1),\\
(c_{6n-3,n-1},v_{6n-4})&=(-2,-1),
& (c_{6n-2,0},v_{6n-3})&=(0,1).
\end{align*}
The computation for the final terms of $U_n$ in the proof of Lemma~\ref{met:lem-step2} gives $R_{6n-7} = -2$ and $R_{6n-6}=1$.
Applying the recurrence~\eqref{def:Rn}, we obtain
\begin{align*}
R_{6n-5}&=-2R_{6n-6}-R_{6n-7}=0,\\
R_{6n-4}&=R_{6n-5}+R_{6n-6}=1,\\
R_{6n-3}&=-2R_{6n-4}+R_{6n-5}=-2.
\end{align*}
By Lemmas~\ref{met:lem-shift2-step1} and \ref{met:lem-d-values}, we also have
\begin{align*}
    & \sigma_{6n-6}=\sigma_{6n-4}=\sigma_{6n-3}=1, \qquad \sigma_{6n-5}=-1, \\
    & d_{6n-5}=d_{6n-3}=d_{6n-2}=1, \qquad d_{6n-4}=2.
\end{align*}
Consequently,
\begin{equation}\label{met:Vn-contribution}
\sum_{q=6n-6}^{6n-3}d_{q+1}R_q^2\sigma_q
=1+0+1+4=6.
\end{equation}

\medskip

\medskip
\noindent\textbf{Case 2.} Consider the block $W_n$.
For $0\le i\le n-3$, the triple indexed by $i$ occupies the positions $q=3i+2,3i+3,3i+4$.
Combining \eqref{met:R-Wn-triple} with Lemmas~\ref{met:lem-shift2-step1} and~\ref{met:lem-d-values}
gives
\begin{equation*}
\sum_{q=3i+2}^{3i+4}d_{q+1}R_q^2\sigma_q
=
\begin{cases}
0, & \text{if $i$ is even and $i<n-3$},\\
4, & \text{if $i$ is odd and $i<n-3$},\\
1, & \text{if $i=n-3$ and $n$ is odd},\\
3, & \text{if $i=n-3$ and $n$ is even}.
\end{cases}
\end{equation*}
Indeed, for $i<n-3$, the three coefficients $d_{q+1}$ are all $1$, so the contributions are $-1+1+0=0$ when $i$ is even and $1-1+4=4$ when $i$ is odd.
For $i=n-3$, the first coefficient is $0$, giving the last two cases.
Among the indices $0\le i\le n-4$, exactly $\lfloor(n-3)/2\rfloor$ are odd.
Including the last triple, the total contribution of the triples is therefore $2n-5$ for either parity of $n$.
The final term of $W_n$, at $q=3n-4$, has $R_{3n-4}^2=d_{3n-3}=1$ and $\sigma_{3n-4}=(-1)^{n-1}$.
Consequently,
\begin{equation}\label{met:Wn-total}
\sum_{q=2}^{3n-4}d_{q+1}R_q^2\sigma_q
=2n-5+(-1)^{n-1}
=
\begin{cases}
2n-4, & \text{if $n$ is odd},\\
2n-6, & \text{if $n$ is even}.
\end{cases}
\end{equation}

\medskip
\medskip
\noindent{\bf{Case 3.}} Consider the block $U_n$.
This block consists of the triples indexed by $0\le i\le n-3$, followed by the three terms at $q=6n-9,6n-8,6n-7$.
We compute the contributions of these two parts separately.

\smallskip
\noindent{\bf{Subcase 3.1.}} Consider the indexed triples.
By Lemma~\ref{met:lem-shift2-step1}, the signs in the triple indexed by $i$ satisfy $\sigma_{3n-3+3i}=\sigma_{3n-1+3i}=(-1)^{n+i}$ and $\sigma_{3n-2+3i}=(-1)^{n+i+1}$.
By Lemma~\ref{met:lem-d-values}, we have $d_{3n-2+3i}=d_{3n+3i}=1$ and
\[
d_{3n-1+3i}
=
\begin{cases}
0, & \text{if $i=1$},\\
1, & \text{otherwise}.
\end{cases}
\]
We distinguish three possibilities for $i$.
\begin{itemize}
    \item[(1)]  Suppose that $2\le i\le n-3$.
By \eqref{met:R-Un-triples}, the successive values of $R_{3n-3+3i}$, $R_{3n-2+3i}$, and $R_{3n-1+3i}$
are $-1,1,-2$ when $n+i$ is even and $-1,-1,0$ when $n+i$ is odd.
Since all three coefficients $d_{q+1}$ are $1$, we obtain
\begin{equation}\label{met:Un-generic}
\sum_{q=3n-3+3i}^{3n-1+3i}d_{q+1}R_q^2\sigma_q
=
\begin{cases}
4, & \text{if $n+i$ is even},\\
0, & \text{if $n+i$ is odd}.
\end{cases}
\end{equation}

\item[(2)] Suppose that $i=0$.
By \eqref{met:R-Un-triples}, the successive values of $R_{3n-3}$, $R_{3n-2}$, and $R_{3n-1}$ are $0,-1,0$ when $n$ is odd and $-2,1,-2$ when $n$ is even.
Here $d_{3n-2}=d_{3n-1}=d_{3n}=1$, so
\begin{equation}\label{met:Un-i0}
\sum_{q=3n-3}^{3n-1}d_{q+1}R_q^2\sigma_q
=
\begin{cases}
1, & \text{if $n$ is odd},\\
7, & \text{if $n$ is even}.
\end{cases}
\end{equation}

\item[(3)]  Suppose that $i=1$, which occurs only when $n\ge4$.
By \eqref{met:R-Un-triples}, the successive values of $R_{3n}$, $R_{3n+1}$, and $R_{3n+2}$ are $-1,1,-2$ when $n$ is odd and $-1,-1,0$ when $n$ is even.
Here $d_{3n+1}=d_{3n+3}=1$ and $d_{3n+2}=0$.
Thus the middle term contributes $0$, and
\begin{equation}\label{met:Un-i1}
\sum_{q=3n}^{3n+2}d_{q+1}R_q^2\sigma_q
=
\begin{cases}
5, & \text{if $n$ is odd},\\
-1, & \text{if $n$ is even}.
\end{cases}
\end{equation}
\end{itemize}

\smallskip
\noindent{\bf{Subcase 3.2.}} Consider the final three terms of $U_n$.
The computation in the proof of Lemma~\ref{met:lem-step2} gives $R_{6n-9}=-1$, $R_{6n-8}=1$, and $R_{6n-7}=-2$.
By Lemmas~\ref{met:lem-shift2-step1} and~\ref{met:lem-d-values}, we have $\sigma_{6n-9}=\sigma_{6n-7}=1$, $\sigma_{6n-8}=-1$, and $d_{6n-8}=d_{6n-6}=1$.
Moreover, $d_{6n-7}=0$ when $n=3$ and $d_{6n-7}=1$ when $n\ge4$. 
Hence
\begin{equation}\label{met:Un-trailing}
\sum_{q=6n-9}^{6n-7}d_{q+1}R_q^2\sigma_q
=1-d_{6n-7}+4
=
\begin{cases}
5, & \text{if $n=3$},\\
4, & \text{if $n\ge4$}.
\end{cases}
\end{equation}

Combining the above two subcases, we compute the total contribution of $U_n$.
When $n=3$, only the indexed triple $i=0$ occurs.
Among the indices $2\le i\le n-3$, exactly $(n-5)/2$ satisfy $n+i\equiv0\pmod2$ when $n\ge5$ is odd,
and exactly $(n-4)/2$ satisfy this condition when $n\ge4$ is even.
Thus \eqref{met:Un-generic}--\eqref{met:Un-trailing} give
\begin{equation}\label{met:Un-total}
\sum_{q=3n-3}^{6n-7}d_{q+1}R_q^2\sigma_q
=
\begin{cases}
2n, & \text{if $n$ is odd},\\
2n+2, & \text{if $n$ is even}.
\end{cases}
\end{equation}

\smallskip
Adding the contributions \eqref{met:Vn-contribution}, \eqref{met:Wn-total}, and \eqref{met:Un-total} of the three cases, we obtain
\[
C_n=
\begin{cases}
6+(2n-4)+2n=4n+2, & \text{if $n$ is odd},\\
6+(2n-6)+(2n+2)=4n+2, & \text{if $n$ is even}.
\end{cases}
\]
Therefore $C_n=4n+2>0$ for every $n\ge3$. This completes the proof.
\end{proof}

We also need the value of $S_p$ at the initial index $p=6n-3$ of the chosen arithmetic progression.

\begin{Lemma}\label{met:lem-S6n-3}
For every $n\geq 3$, we have $S_{6n-3}=4n+7$. In particular, $S_{6n-3}>0$.
\end{Lemma}

\begin{proof}
Summing the increments in \eqref{met:eq-Sp-rec}
for $0\le p\le6n-4$ and using $S_0=d_1$, we obtain
\[
S_{6n-3}
=d_1+\sum_{q=1}^{6n-3}d_{q+1}R_q^2\sigma_q
=d_1+d_2R_1^2\sigma_1
+\sum_{q=2}^{6n-3}d_{q+1}R_q^2\sigma_q.
\]
By \eqref{met:def-Cn-sum} and
Lemma~\ref{met:lem-shift2-step2}, the last sum equals $C_n=4n+2$.
Since $n\ge3$, Lemmas~\ref{met:lem-d-values}, \ref{met:lem-step2}, and~\ref{met:lem-shift2-step1}
give $d_1=d_2=1$, $R_1=-2$, and $\sigma_1=1$, respectively.
Therefore,
\[
S_{6n-3}
=1+1\cdot(-2)^2\cdot1+(4n+2)
=4n+7>0.
\]
\end{proof}

\subsection{The case $n=2$ and completion of the proof}\label{met:subsec-n3-n2}

Lemmas~\ref{met:lem-shift2-step1}--\ref{met:lem-S6n-3} are written  for $n\geq 3$. 
For $n=2$, the explicit $H$-fraction \eqref{met:Hfrac-Phi2-3} of $\Phi_2^{(3)}$ and the table
in Lemma~\ref{met:lem-n2} yield the following lemma.
Together with the results for $n\ge3$, this will allow us to complete the proof of Theorem~\ref{met:thm-n+3}.

\begin{Lemma}\label{met:lem-n3-n2}
For $n=2$, the numbers $d_{p+1}$ and $\sigma_p=v_1v_2\cdots v_p$ are given, for $0\leq p\leq 9$, by
\[
\begin{array}{c|cccccccccc}
p & 0&1&2&3&4&5&6&7&8&9\\ \hline
d_{p+1} & 0&1&1&1&1&1&0&2&0&1\\
\sigma_p & 1&1&-1&1&-1&1&1&-1&1&1
\end{array}
\]
and $(d_{p+1})_{p\geq 2}$ and $(\sigma_p)_{p\geq 1}$ are $8$-periodic. 
Consequently, the sum $C_2:=S_{p+8}-S_p$ is independent of $p\geq 1$, and
\[
C_2=10=4n+2,\qquad S_9=14.
\]
\end{Lemma}

\begin{proof}
As in the proof of Lemma~\ref{met:lem-d-values}, let $B_p(q)=1-u_{p+1}(q)q$ denote the denominator of the term at index $p$ in \eqref{met:Hfrac-Phi2-3}.
Then $d_{p+1}=[q^{k_p}]B_p(q)$, where the values of $k_p$ are given in Lemma~\ref{met:lem-n2}.
For $p=0,1,\ldots,9$, the successive denominators are
\[
1+2q^2,\quad 1,\quad 1,\quad 1+q,\quad 1,\quad 1,\quad
1+2q^2,\quad 1+2q^2-q^3,\quad 1+2q^2,\quad 1,
\]
and the corresponding exponents are
\[
(k_0,k_1,\ldots,k_9)=(1,0,0,0,0,0,1,2,1,0).
\]
Extracting the indicated coefficients gives the row of $d_{p+1}$-values in the table.
For example, $d_1=[q^1](1+2q^2)=0$ and $d_8=[q^2](1+2q^2-q^3)=2$.
The row of $\sigma_p$-values follows from $\sigma_0=1$ and $\sigma_p=\sigma_{p-1}v_p$, using the values of $v_p$ in Lemma~\ref{met:lem-n2}.

The pairs $(B_p(q),k_p)$ are $8$-periodic for $p\ge2$, so the same is true of $(d_{p+1})_{p\ge2}$.
For $p\ge1$, the periodicity of $(v_p)_{p\ge2}$ gives
\[
\sigma_{p+8}
=\sigma_p\,v_{p+1}\cdots v_{p+8}
=\sigma_p\,v_2\cdots v_9
=\sigma_p,
\]
since six of the eight factors $v_2,\ldots,v_9$ are $-1$ and the other two are $1$.
Together with the $8$-periodicity of $(R_p)$ from Lemma~\ref{met:lem-n2}, these identities and
\eqref{met:eq-Sp-rec} yield
\[
C_2:=S_{p+8}-S_p
=\sum_{q=2}^{9}d_{q+1}R_q^2\sigma_q
\]
for $p\ge1$, which is independent of $p$.

By Lemma~\ref{met:lem-n2}, we have
\[
(R_2,R_3,\ldots,R_9)=(-1,-2,1,-2,1,0,1,-2).
\]
Substituting these values and the table entries into the sum gives the successive contributions $-1, 4, -1, 4, 0, 0, 0, 4$.
The terms at $q=6,8$ are $0$ because $d_7=d_9=0$, and the term at $q=7$ is $0$ because $R_7=0$.
Thus $C_2=10$.
Finally, as in the proof of Lemma~\ref{met:lem-S6n-3},
\[
S_9=d_1+d_2R_1^2\sigma_1+C_2
=0+1\cdot(-2)^2\cdot1+10=14.
\]
\end{proof}

\begin{Remark}
The value $C_2=10$ agrees with the formula $C_n=4n+2$ in Lemma~\ref{met:lem-shift2-step2}.
However, $S_9=14$, rather than the value $15$ obtained by substituting $n=2$ into $4n+7$.
This difference comes from the initial coefficient $d_1$.
For $n\ge3$, the first denominator is $\langle n\rangle_q+q^{n+1} =1+q^2+\cdots+q^{n-1}+2q^n$, whose coefficient of $q^{n-1}$ is $1$.
For $n=2$, this denominator is $1+2q^2$, so $d_1=0$.
Thus the identity $S_{6n-3}=d_1+4+C_n$ gives $S_{6n-3}=4n+7$ for $n\ge3$, as in Lemma~\ref{met:lem-S6n-3}, and $S_9=14$ for $n=2$.
\end{Remark}

\medskip
We are now ready to complete the proof of Theorem~\ref{met:thm-n+3}.

\begin{proof}[Proof of Theorem~\ref{met:thm-n+3}]
Let $n\ge2$.
By \eqref{met:def-Cn-sum} and Lemma~\ref{met:lem-shift2-step2} for $n\ge3$, or Lemma~\ref{met:lem-n3-n2} for $n=2$, we have $S_{p+6n-4}-S_p=C_n$ for every $p\ge1$, where $C_n=4n+2$.
Iterating this identity from $p=6n-3$ gives
\begin{equation}\label{met:eq-Spm}
S_{6n-3+m(6n-4)}
=S_{6n-3}+mC_n
\end{equation}
for $m \ge 0$.
Setting $p_m:=6n-3+m(6n-4)$ and $j_m:=s_{p_m}$.
By Lemma~\ref{met:lem-step1} for $n\ge3$ and Lemma~\ref{met:lem-n2} for $n=2$, we have $k_{6n-3}=0$.
The $(6n-4)$-periodicity of $(k_p)_{p\ge2}$ therefore gives $k_{p_m}=0$ for every $m\ge0$.
Moreover, since $s_{p+1}-s_p=k_p+1$, the same periodicity gives
\[
j_m
=s_{6n-3}
+m\sum_{i=0}^{6n-5}(k_{6n-3+i}+1)
=s_{6n-3}+m\cdot2n(n+1),
\]
where the last equality follows by summing the explicit values over one period.
The initial index is
\[
s_{6n-3}
=\sum_{i=0}^{6n-4}(k_i+1)
=n(2n+3),
\]
again by Lemma~\ref{met:lem-step1} for $n\ge3$ and Lemma~\ref{met:lem-n2} for $n=2$.
Since $k_{p_m}=0$, \eqref{met:eq-modulus}, the identity $|T_p|=|S_p|$, and \eqref{met:eq-Spm} yield
\[
\bigl|H_{j_m}^{(n+3)}(\Phi_n)\bigr|
=|T_{p_m}|
=|S_{p_m}|
=\bigl|S_{6n-3}+mC_n\bigr|.
\]
By Lemma~\ref{met:lem-S6n-3}, we have $S_{6n-3}=4n+7>0$ for $n\ge3$.
For $n=2$, Lemma~\ref{met:lem-n3-n2} gives $S_9=14>0$.
Since $C_n=4n+2>0$, it follows that
\[
\bigl|H_{j_m}^{(n+3)}(\Phi_n)\bigr|
=S_{6n-3}+mC_n
\ge C_nm
\]
for $m\ge1$, which proves \eqref{met:eq-n+3-growth}.
In particular, $\bigl|H_{j_m}^{(n+3)}(\Phi_n)\bigr|\to\infty$ as $m\to\infty$, so the sequence $H^{(n+3)}(\Phi_n)$ is unbounded.
\end{proof}

\begin{Remark}
The proof of Theorem~\ref{met:thm-n+3} gives the exact identity
\[
\bigl|H_{j_m}^{(n+3)}(\Phi_n)\bigr|
=S_{6n-3}+m(4n+2)
\]
for $m\ge0$.
Thus, along this subsequence, we obtain not only unboundedness
but also an explicit linear growth rate. Since
$j_m=s_{6n-3}+2mn(n+1)$, we have
\[
\lim_{m\to\infty}
\frac{\bigl|H_{j_m}^{(n+3)}(\Phi_n)\bigr|}{j_m}
=
\frac{2n+1}{n(n+1)}.
\]
\end{Remark}


\section{Super $1$-fractions, even contraction, and $G$-fractions}\label{sec:gfrac}

Definition \ref{def:super} attaches to each power series $F(x)$,  for each positive integer $\delta$, a unique super $\delta$-fraction. 
For $\delta=2$, this is the $H$-fraction studied in the present paper, which generalizes the classical $J$-fraction.
The case $\delta=1$ similarly generalizes the $S$-fraction.
In \cite{Han2016Adv}, this case is covered by the general existence and uniqueness result recalled in
Theorem~\ref{th:super2}(i). 
Classically, the $S$-fraction and the $J$-fraction of a power series are related by the {\it even contraction} 
\cite[p.~21]{Wall1948}, \cite[pp.~12--13]{Perron1957II}, 
\cite{Stieltjes1894Re}: 
if 
\begin{equation}\label{gf:eq:Sfrac}
F(x)=
  \cFrac{b_0}{1}
  - \cFrac{b_1 x}{1}
  - \cFrac{b_2 x}{1}
  - \cFrac{b_3 x}{1}
  - \cdots,
\end{equation}
then 
\begin{equation}\label{gf:eq:StoJ}
F(x)=
  \cFrac{b_0}{1-b_1x}
  - \cFrac{b_1b_2\, x^2}{1-(b_2+b_3)x}
  - \cFrac{b_3b_4\, x^2}{1-(b_4+b_5)x}
  - \cdots,
\end{equation}
which is the $J$-fraction of $F$, so that, by Heilermann's formula
\eqref{eq:detH}, 
\begin{equation}\label{gf:eq:HnS}
H_n(F)=b_0^{\,n}\prod_{i=1}^{n-1}(b_{2i-1}b_{2i})^{\,n-i}. 
\end{equation}
A power series $F$ admits an infinite $S$-fraction \eqref{gf:eq:Sfrac} with $b_i\ne0$ for every $i\ge0$ if and only if $H_n(F)\ne0$ and $H_n^{(1)}(F)\ne0$ for every $n\ge1$ (see \cite{Wall1948}).
By comparison, it admits an infinite $J$-fraction with all partial numerators nonzero if and only if $H_n(F)\ne0$ for every $n\ge1$.

In this section, we extend these relationships to super $1$-fractions.
Theorem~\ref{gf:th:contraction} gives an explicit construction of the $H$-fraction associated with an arbitrary super $1$-fraction.
Applying Theorem~\ref{th:super2}(ii) to this expansion yields the formulas for all Hankel determinants stated in Theorem~\ref{gf:th:det}.
This $H$-fraction is obtained by even contraction if and only if $\ell_{2i+1}=0$ for every $i\ge0$, where the parameters $\ell_i$ are defined below.
We call a super $1$-fraction satisfying this condition a \emph{$G$-fraction}.
Using Theorem~\ref{th:shiftHDet}, we establish an analogue of the Stieltjes criterion in Theorem~\ref{gf:th:criterion}:
a super $1$-fraction is a $G$-fraction if and only if the $1$-shifted Hankel determinants $H^{(1)}_{s_j}(F)$ at the jump orders $s_j$ of its $H$-fraction are all nonzero.
Finally, applying Theorems~\ref{th:shiftHDet} and~\ref{th:shift2HDet} gives closed-form expressions for the $1$-shifted and $2$-shifted Hankel determinants of a $G$-fraction, as stated in Theorems~\ref{gf:th:shiftG} and~\ref{gf:th:shift2G}.

\smallskip

Throughout this section, $F(x)$ is a power series whose super $1$-fraction is 
\begin{equation}\label{gf:eq:super1}
F(x)=
\cFrac{b_0 x^{\ell_0}}{ 1 - w_1(x) x} 
- \cFrac{ b_1 x^{\ell_0+\ell_1+1}}{ 1 - w_2(x)x}
- \cFrac{ b_2 x^{\ell_1+\ell_2+1}}{ 1 - w_3(x)x}
- {\cdots} 
\end{equation}
where $b_i\neq 0$ are constants, $\ell_i$ are nonnegative integers, and $w_i(x)$ are polynomials of degree at most $\ell_{i-1}-1$.
In particular, $w_{i+1}(x)=0$ whenever $\ell_i=0$.
When all $\ell_i=0$, \eqref{gf:eq:super1} reduces to the $S$-fraction \eqref{gf:eq:Sfrac}.
The letters $b_i$, $\ell_i$, $w_i$ are reserved for the super $1$-fraction of $F$, whereas $(v_j,k_j,u_j)$, $s_j$, $\epsilon_j$, $c_{j,i}$, $V_j$, $R_j$ and $d_j$ denote, as in the rest of the paper, the data attached to the $H$-fraction \eqref{eq:Hfrac} of the same series $F$. 
We assume that \eqref{gf:eq:super1} is nonterminating.
The terminating case can be treated with minor modifications.
For every $i\ge0$, define the tail
\begin{equation}\label{gf:eq:tails}
T_i(x):=
\cFrac{b_i x^{\ell_i}}{1-w_{i+1}(x)x}
-\cFrac{b_{i+1}x^{\ell_i+\ell_{i+1}+1}}{1-w_{i+2}(x)x}
-\cdots.
\end{equation}
Then $T_0=F$, and these tails satisfy
\begin{equation}\label{gf:eq:tailrec}
T_i(x)=
\frac{b_i x^{\ell_i}}
{1-w_{i+1}(x)x-x^{\ell_i+1}T_{i+1}(x)}.
\end{equation}
We call $(v_j,k_j,u_{j+1})$ the $j$-th level of the $H$-fraction \eqref{eq:Hfrac}.
Thus a level $(v,k,u)$ consists of a nonzero numerator coefficient $v$, a nonnegative exponent parameter $k$, and a polynomial $u(x)$ satisfying $\deg u\le k$.
Similarly, the $i$-th level of the super $1$-fraction \eqref{gf:eq:super1} is $(b_i,\ell_i,w_{i+1})$, with $\deg w_{i+1}\le\ell_i-1$.

\subsection{The peeling lemmas}\label{gf:subsec:peel}

The following lemma gives the basic recursive step for $\delta=2$ underlying Theorem~\ref{th:super2}(ii).
It is Lemma~2.2 of \cite{Han2016Adv}, with the sign convention of Remark~\ref{rem:sign}.

\begin{Lemma}\label{gf:lem:peel2}
Let $k\geq0$, let $u(x)$ be a polynomial of degree at most $k$, 
let $G(x)$ be a power series, and let 
\[
F(x)=\frac{x^{k}}{1-u(x)x-x^{k+2}G(x)}.
\]
Then $H_n(F)=0$ for $1\leq n\leq k$, and 
$H_n(F)=(-1)^{k(k+1)/2}H_{n-k-1}(G)$ for $n\geq k+1$.
\end{Lemma}

The $\delta=1$ analogue reads as follows. 
Whereas the $\delta=2$ step of Lemma \ref{gf:lem:peel2} 
sends Hankel determinants to Hankel determinants, 
a $\delta=1$ step exchanges the Hankel determinants and the 
$1$-shifted Hankel determinants.

\begin{Lemma}\label{gf:lem:peel1}
Let $k$ be a nonnegative integer, let $u(x)$ be a polynomial of degree at most $k-1$, and let $G(x)$ be a power series.
When $k=0$, we take $u(x)=0$.
Set
\[
F(x)=\frac{x^{k}}{1-u(x)x-x^{k+1}G(x)}.
\]
Then the following statements hold.
\begin{itemize}
\item[(\rmnum{1})] We have $H_n(F)=0$ for $1\le n\le k$, and for $n\ge k+1$,
\[
H_n(F)=(-1)^{k(k+1)/2}H^{(1)}_{n-k-1}(G).
\]

\item[(\rmnum{2})] We have $H^{(1)}_n(F)=0$ for $1\leq n\leq k-1$, and for $n\geq k$,
\[
H^{(1)}_n(F)=(-1)^{k(k-1)/2}\,H_{n-k}(G).
\]
In particular, $H^{(1)}_k(F)=(-1)^{k(k-1)/2}$.
\end{itemize}
\end{Lemma}

\begin{proof}
 Write $G(x)=g_0+xG_1(x)$, so that $H^{(1)}_m(G)=H_m(G_1)$ for all $m\geq0$. 
 Then $x^{k+1}G=g_0x^{k+1}+x^{k+2}G_1$ and 
\[
F(x)=\frac{x^k}{1-(u(x)+g_0x^k)x-x^{k+2}G_1(x)},
\]
where $u(x)+g_0x^k$ has degree at most $k$. 
Lemma \ref{gf:lem:peel2} applied to this expression gives (\rmnum{1}).

We now prove part~(\rmnum{2}), distinguishing two cases.

\smallskip
\noindent{\bf{Case 1.}} Suppose that $k\ge1$. 
Since $F(0)=0$, we have $H_n^{(1)}(F)=H_n(F/x)$, and
\[
\frac{F(x)}x
=\frac{x^{k-1}}{1-u(x)x-x^{(k-1)+2}G(x)}.
\]
As $\deg u(x) \le k-1$, Lemma~\ref{gf:lem:peel2}, applied with $k-1$ in place of $k$, proves part~(\rmnum{2}) in this case.

\smallskip
\noindent{\bf{Case 2.}} Suppose that $k=0$.
Then $u=0$ and $F=1/(1-xG)$, so $F-1/x=G/(1-xG)$.
Hence $H_n^{(1)}(F)=H_n\bigl(G/(1-xG)\bigr)$, and it remains to prove $H_n\bigl(G/(1-xG)\bigr)=H_n(G)$
for every $n\ge0$.
If $G=0$, the identity is immediate.
Otherwise, let $(v_j,k_j,u_{j+1})_{j\ge0}$ denote, within this proof, the parameters of the $H$-fraction of $G$.
Write $G=v_0x^{k_0}/D$, where
\[
D=1-u_1(x)x
-\cFrac{v_1x^{k_0+k_1+2}}{1-u_2(x)x}
-\cdots.
\]
Then
\[
\frac{G}{1-xG}=\frac{v_0x^{k_0}}{D-v_0x^{k_0+1}}
=\cFrac{v_0x^{k_0}}{1-(u_1(x)+v_0x^{k_0})x}-\cFrac{v_1x^{k_0+k_1+2}}{1-u_2(x)x}-\cdots
\]
Since $\deg(u_1(x)+v_0x^{k_0})\le k_0$, this is a super $2$-fraction and hence, by Theorem~\ref{th:super2}(i), the $H$-fraction of $G/(1-xG)$.
It differs from the $H$-fraction of $G$ only in the first denominator polynomial.
In particular, the parameters $v_j$ and $k_j$ are unchanged.
Theorem~\ref{th:super2}(ii) therefore gives the same Hankel determinants for both series, proving
part~(\rmnum{2}).
\end{proof}

Applying part~(\rmnum{1}) to a level with exponent $k$ and then part~(\rmnum{2}) to the next level, with exponent $m$, we obtain
\[
H_n(F)
=(-1)^{k(k+1)/2+m(m-1)/2}H_{n-k-m-1}(G).
\]
Thus, applying the two reductions successively returns an unshifted Hankel determinant of the remaining tail, just as removing a single $\delta=2$ level does.
The next theorem describes the corresponding transformation at the level of continued fractions.

\subsection{The $H$-fraction of a super $1$-fraction}\label{gf:subsec:contraction}

We first treat one pair of consecutive levels of a super $1$-fraction.

\begin{Lemma}\label{gf:lem:pair}
Let $b,b'$ be nonzero constants, $k,k'$ nonnegative integers, 
$u(x)$ a polynomial of degree at most $k$, 
$u'(x)$ a polynomial of degree at most $k'-1$, 
$T(x)$ a power series, and 
\[
S(x)=\cFrac{b\,x^{k}}{1-u(x)x}-\cFrac{b'x^{k+k'+1}}{1-u'(x)x-x^{k'+1}T(x)}.
\]

\noindent{\rm (a)} If $k'\geq1$, then
\[
S(x)=\cFrac{b\,x^{k}}{1-u(x)x}-\cFrac{b'x^{k+(k'-1)+2}}{1-u'(x)x-x^{(k'-1)+2}T(x)}.
\]

\noindent{\rm (b)} If $k'=0$ (so that $u'=0$), then
\[
S(x)=\frac{b\,x^{k}}{1-\bigl(u(x)+b'x^{k}\bigr)x-x^{k+2}\widetilde T(x)},
\qquad \widetilde T:=\frac{b'T}{1-xT},
\]
and if moreover $T=\dfrac{b''x^{k''}}{1-u''(x)x-x^{k''+1}T'(x)}$ 
for some constant $b''$, integer $k''\geq0$, polynomial $u''$ and 
power series $T'$, then 
\[
\widetilde T=\frac{b'b''x^{k''}}{1-\bigl(u''(x)+b''x^{k''}\bigr)x-x^{k''+1}T'(x)}.
\]
\end{Lemma}

\begin{proof}
Rewriting the exponents as $k+k'+1=k+(k'-1)+2$ and $k'+1=(k'-1)+2$ gives (a).
When $k'=0$, we have $u'=0$, and the identity
\[
\frac{b'x^{k+1}}{1-xT}
=b'x^{k+1}\left(1+\frac{xT}{1-xT}\right)
=b'x^{k+1}+x^{k+2}\frac{b'T}{1-xT}
\]
gives the first formula in part~(b).
Substituting the expression for $T$ into $\widetilde T=b'T/(1-xT)$ yields
\[
\widetilde T
=\frac{b'b''x^{k''}}
{1-u''(x)x-x^{k''+1}T'(x)-b''x^{k''+1}}.
\]
Combining the terms $u''(x)x$ and $b''x^{k''+1}$ in the denominator gives the second formula, completing the proof of (b).
\end{proof}

Note that in case (a), the two levels $(b,k,u)$ and $(b',k',u')$ become the two $H$-fraction levels $(b,k,u)$ and $(b',k'-1,u')$.
In case (b), they merge into the single $H$-fraction level $(b,k,u+b'x^k)$, and $b'$ also multiplies the numerator constant of the next level. 
For a sequence $(\ell_i)_{i\geq0}$ of nonnegative integers, set $\chi_{-1}:=0$ and, for every $t\ge0$, define
\begin{equation*}\label{gf:eq:chi}
\chi_t:=
\begin{cases}
1,&\text{if $\ell_{2t+1}=0$,}\\ 
0,&\text{if $\ell_{2t+1}\geq1$.}
\end{cases}
\end{equation*}

Applying these transformations successively yields the explicit parameters of each level of the corresponding $H$-fraction, as stated in the following theorem.

\begin{Theorem}\label{gf:th:contraction}
Let \eqref{gf:eq:super1} be the super $1$-fraction of a power series $F(x)$.
Then the $H$-fraction of $F$ consists of the following levels, in this order: for $t=0,1,2,\ldots$, the level
\begin{equation}\label{gf:eq:levelA}
\Bigl(\,b_{2t-1}^{\,\chi_{t-1}}\,b_{2t},\ \ \ell_{2t},\ \ 
w_{2t+1}(x)+\bigl(\chi_{t-1}b_{2t}+\chi_t b_{2t+1}\bigr)x^{\ell_{2t}}\Bigr),
\end{equation}
followed, when $\chi_t=0$, by the level
\begin{equation}\label{gf:eq:levelB}
\bigl(\,b_{2t+1},\ \ \ell_{2t+1}-1,\ \ w_{2t+2}(x)\bigr).
\end{equation}
Here we use the convention $b_{-1}^{\,0}=1$ when $t=0$.
\end{Theorem}

\begin{proof}
For every $t\ge0$, set
\[
S_t:=
\frac{\beta_t x^{\ell_{2t}}}
{1-\omega_t(x)x-x^{\ell_{2t}+1}T_{2t+1}(x)},
\]
where $\beta_t:=b_{2t-1}^{\,\chi_{t-1}}b_{2t}$, $\omega_t(x):=w_{2t+1}(x)+\chi_{t-1}b_{2t}x^{\ell_{2t}}$, and $T_i$ are the tails defined in \eqref{gf:eq:tails}.
Then $S_0=T_0=F$ by \eqref{gf:eq:tailrec}, and $\deg\omega_t\le\ell_{2t}$.
We show that extracting the level \eqref{gf:eq:levelA} and, when $\chi_t=0$, the level \eqref{gf:eq:levelB} from $S_t$ leaves the tail $S_{t+1}$.
More precisely, the remaining term in the denominator of the last extracted level is $-x^{K+2}S_{t+1}$, where $K$ is the exponent of that level.
By \eqref{gf:eq:tailrec}, we have
\[
x^{\ell_{2t}+1}T_{2t+1}=\frac{b_{2t+1}x^{\ell_{2t}+\ell_{2t+1}+1}}{1-w_{2t+2}x-x^{\ell_{2t+1}+1}T_{2t+2}},
\] 
so that $S_t$ has the form of the series $S$ of Lemma \ref{gf:lem:pair},  with $(b,k,u)=(\beta_t,\ell_{2t},\omega_t)$,  $(b',k',u')=(b_{2t+1},\ell_{2t+1},w_{2t+2})$ and $T=T_{2t+2}$.
We distinguish two cases according to whether $\ell_{2t+1}$ is zero.

\medskip
\noindent{\bf{Case 1.}} Suppose that $\ell_{2t+1}\ge1$, so that $\chi_t=0$.
Lemma~\ref{gf:lem:pair}(a) gives
\[
S_t
=
\cFrac{\beta_t x^{\ell_{2t}}}{1-\omega_t(x)x}
-
\cFrac{b_{2t+1}x^{\ell_{2t}+(\ell_{2t+1}-1)+2}}
{1-w_{2t+2}(x)x-x^{(\ell_{2t+1}-1)+2}T_{2t+2}}.
\]
Since $\chi_t=0$, we have $\beta_{t+1}=b_{2t+2}$ and $\omega_{t+1}=w_{2t+3}$.
Hence \eqref{gf:eq:tailrec} gives $T_{2t+2}=S_{t+1}$.
The two displayed levels are therefore \eqref{gf:eq:levelA} and \eqref{gf:eq:levelB}, with remaining tail $S_{t+1}$.

\medskip

\noindent{\bf{Case 2.}} Suppose that $\ell_{2t+1}=0$, so that $\chi_t=1$.
Then Lemma~\ref{gf:lem:pair}(b) gives
\[
S_t
=
\frac{\beta_t x^{\ell_{2t}}}
{1-\bigl(\omega_t(x)+b_{2t+1}x^{\ell_{2t}}\bigr)x
-x^{\ell_{2t}+2}\widetilde T},
\]
where, by \eqref{gf:eq:tailrec},
\[
\widetilde T
=
\frac{b_{2t+1}b_{2t+2}x^{\ell_{2t+2}}}
{1-\bigl(w_{2t+3}(x)+b_{2t+2}x^{\ell_{2t+2}}\bigr)x
-x^{\ell_{2t+2}+1}T_{2t+3}}.
\]
Since $\chi_t=1$, we have $ \beta_{t+1}=b_{2t+1}b_{2t+2}$, $ \omega_{t+1} =w_{2t+3}+b_{2t+2}x^{\ell_{2t+2}}$, and therefore $\widetilde T=S_{t+1}$.
Moreover, 
\[
\omega_t(x)+b_{2t+1}x^{\ell_{2t}}
=
w_{2t+1}(x)
+\bigl(\chi_{t-1}b_{2t}+\chi_t b_{2t+1}\bigr)
x^{\ell_{2t}}.
\]
Thus the displayed level is \eqref{gf:eq:levelA}, with remaining tail $S_{t+1}$.

\medskip

Combining the two cases and iterating from $S_0=F$, we obtain a continued fraction whose levels are
exactly those listed in the theorem.
All numerator constants are nonzero, and all exponents are nonnegative.
The polynomial in \eqref{gf:eq:levelA} has degree at most $\ell_{2t}$, while that in \eqref{gf:eq:levelB} has degree at most $\ell_{2t+1}-1$.
Hence this expansion is a super $2$-fraction.
By the uniqueness assertion of Theorem~\ref{th:super2}(i), it is the $H$-fraction of $F$.
\end{proof}

For the $S$-fraction \eqref{gf:eq:Sfrac}, all $\chi_t=1$ and Theorem \ref{gf:th:contraction} is the classical even contraction \eqref{gf:eq:StoJ}. 
To describe the Hankel determinants of the super $1$-fraction \eqref{gf:eq:super1}, we introduce two sequences.
Set $\sigma_0:=0$ and, for every $t\ge0$, define
\begin{equation}\label{gf:eq:sigma}
\sigma_{2t+1}:=\sigma_{2t}+\ell_{2t}+1,
\qquad
\sigma_{2t+2}:=\sigma_{2t+1}+\ell_{2t+1}.
\end{equation}
Equivalently, for every $i\ge0$,  $\sigma_i=\sum_{l<i}\ell_l+\lceil i/2\rceil$.
In particular, $\sigma_{2t+1}=\sigma_{2t+2}$ if and only if $\chi_t=1$.
We also define, for every $i\ge0$,
\begin{equation}\label{gf:eq:e}
e_i
:=\sum_{l=0}^{i-1}
\frac{\ell_l\bigl(\ell_l+(-1)^l\bigr)}{2}
=
\sum_{\substack{0\le l<i\\ l\ \mathrm{even}}}
\binom{\ell_l+1}{2}
+
\sum_{\substack{0\le l<i\\ l\ \mathrm{odd}}}
\binom{\ell_l}{2}.
\end{equation}

The preceding transformation yields the following formulas for the Hankel determinants of a super $1$-fraction.

\begin{Theorem}\label{gf:th:det}
Let \eqref{gf:eq:super1} be the super $1$-fraction of a power series $F(x)$, and let $(\sigma_i)_{i\ge0}$ and $(e_i)_{i\ge0}$ be defined by \eqref{gf:eq:sigma} and \eqref{gf:eq:e}, respectively.
Then the following statements hold.
\begin{itemize}
\item[(\rmnum{1})] We have $H_n(F)=0$ whenever $n\notin\{\sigma_i:i\ge0\}$.

\item[(\rmnum{2})]
For every $i\ge0$, we have
\begin{equation}\label{gf:eq:Hsigma}
H_{\sigma_i}(F)
=
(-1)^{e_i}\prod_{l=0}^{i-1}b_l^{\,\sigma_i-\sigma_l}.
\end{equation}
\end{itemize}
In particular, the jump orders of the $H$-fraction of $F$ are exactly the distinct values of $(\sigma_i)_{i\ge0}$.
\end{Theorem}

\begin{proof}
We first determine the jump orders of the $H$-fraction of $F$.
By Theorem~\ref{gf:th:contraction}, the exponents of the successive levels of the $H$-fraction of $F$ are $\ell_0$, followed, for each $t\ge0$, by $\ell_{2t+1}-1$ if $\chi_t=0$ (this level is omitted if $\chi_t=1$), and then by $\ell_{2t+2}$.

For an $H$-fraction with exponents $(k_j)_{j\ge0}$, the successive jump orders satisfy
$s_{j+1}-s_j=k_j+1$ with $s_0=0$.
Thus the level \eqref{gf:eq:levelA} increases the jump order by $\ell_{2t}+1$, while the level \eqref{gf:eq:levelB}, when present, increases it by $\ell_{2t+1}$.
Starting from $\sigma_0=0$, these increments give
\[
\sigma_{2t}
\longmapsto
\sigma_{2t}+\ell_{2t}+1=\sigma_{2t+1},
\]
followed, when $\chi_t=0$, by
\[
\sigma_{2t+1}
\longmapsto
\sigma_{2t+1}+\ell_{2t+1}=\sigma_{2t+2}.
\]
When $\chi_t=1$, we have $\ell_{2t+1}=0$ and $\sigma_{2t+2}=\sigma_{2t+1}$, so no new jump order is introduced.
Consequently, the jump orders are precisely the distinct values of $(\sigma_i)_{i\ge0}$.
Theorem~\ref{th:super2}(ii) therefore proves part~(\rmnum{1}) and shows that $H_{\sigma_i}(F)\neq0$ for every $i\ge0$.

We now evaluate these nonzero determinants using \eqref{eq:HankelDetFundamental}.
Each $H$-fraction level with exponent $k$ contributes $k(k+1)/2$ to the sign exponent.
The level \eqref{gf:eq:levelA}, whose exponent is $\ell_{2t}$, therefore contributes $\ell_{2t}(\ell_{2t}+1)/2$.
The level \eqref{gf:eq:levelB}, when present, has exponent $\ell_{2t+1}-1$ and contributes $(\ell_{2t+1}-1)\ell_{2t+1}/2$.
When $\chi_t=1$, this level is omitted, but $\ell_{2t+1}=0$, so the same expression is zero.
Hence the total sign exponent at the jump order $\sigma_i$ is
\[
\sum_{\substack{0\le l<i\\ l\ \mathrm{even}}}
\frac{\ell_l(\ell_l+1)}{2}
+
\sum_{\substack{0\le l<i\\ l\ \mathrm{odd}}}
\frac{\ell_l(\ell_l-1)}{2}
=
\sum_{l=0}^{i-1}
\frac{\ell_l\bigl(\ell_l+(-1)^l\bigr)}{2}
=e_i.
\]
The sign factor in $H_{\sigma_i}(F)$ is therefore $(-1)^{e_i}$.

Finally, in \eqref{eq:HankelDetFundamental}, the numerator constant $v$ of a level starting at jump order $s$ contributes $v^{\,\sigma_i-s}$.
The constant $b_{2t}$ occurs in the level \eqref{gf:eq:levelA}, which starts at $\sigma_{2t}$.
The constant $b_{2t+1}$ occurs either in the level \eqref{gf:eq:levelB}, which starts at $\sigma_{2t+1}$ if $\chi_t=0$, or as a factor of the numerator constant in the level \eqref{gf:eq:levelA} for $t+1$, which starts at $\sigma_{2t+2}=\sigma_{2t+1}$ if $\chi_t=1$.
In all cases, $b_l$ contributes $b_l^{\,\sigma_i-\sigma_l}$ for $0\le l<i$, which proves \eqref{gf:eq:Hsigma}.
In particular, when $\sigma_{2t+1}=\sigma_{2t+2}$, the two corresponding expressions agree: we have $e_{2t+2}=e_{2t+1}$, and the additional factor is $b_{2t+1}^{\,0}=1$.
\end{proof}

When all $\ell_i=0$, we have $\sigma_{2t+1}=\sigma_{2t+2}=t+1$,  $e_i=0$, and \eqref{gf:eq:Hsigma} reduces to the classical  formula \eqref{gf:eq:HnS}. 
Theorem \ref{gf:th:det} can also be proved directly by  alternately applying (\rmnum{1}) and (\rmnum{2}) of Lemma \ref{gf:lem:peel1}  to the tails \eqref{gf:eq:tailrec}.
The pairing of levels in \eqref{gf:eq:sigma} and \eqref{gf:eq:e} reflects the alternation between unshifted and $1$-shifted Hankel determinants in these two reductions.

\begin{Remark}\label{gf:rem:shiftgeneral}
Theorem \ref{gf:th:det} concerns the Hankel determinants $H_n(F)$ only. 
The shifted Hankel determinants $H^{(1)}_n(F)$, $H^{(2)}_n(F)$ and $H^{(-1)}_n(F)$ are also determined by the super $1$-fraction data.
Indeed, Theorem~\ref{gf:th:contraction} gives the $H$-fraction of $F$ explicitly.
We can then apply Theorems~\ref{th:shiftHDet}, \ref{th:shift2HDet} and~\ref{th:backshiftHDet}.
However, the levels of the resulting $H$-fraction depend on which odd-indexed exponents $\ell_{2t+1}$ are nonzero.
Consequently, so do the coefficients $c_{j,k_{j-1}}$, $d_j$ and the auxiliary sequences $(R_n)$ and $(E_n)$ appearing in the shifted determinant formulas.
The resulting expressions are more involved than \eqref{gf:eq:Hsigma}, and we do not give explicit formulas in the general case.
The case where all $\ell_{2t+1}=0$, namely that of $G$-fractions, is treated in Section~\ref{gf:subsec:shift}.
\end{Remark}

\subsection{Even contraction and $G$-fractions}\label{gf:subsec:G}

Let us recall the classical even contraction formula. 
Let $(A_i)_{i\geq1}$ and $(B_i)_{i\geq0}$ be power series such that 
$B_i(0)\neq0$ for $i\geq1$ and $A_i(0)=0$ for $i\geq2$, and consider 
the continued fraction 
\begin{equation}\label{gf:eq:GCF}
B_0 +  \cFrac{A_1}{B_1}
+ \cFrac{A_2}{B_2}
+ \cFrac{A_3}{ B_3}
+ {\cdots}
\end{equation}
whose value is a well-defined power series. 
We recall the following classical theorem on even contraction from \cite[p.~21]{Wall1948} and \cite[pp.~12--13]{Perron1957II}.

\begin{Theorem}[{\cite{Perron1957II}}, \cite{Wall1948}]\label{gf:th:evencontraction}
The continued fraction \eqref{gf:eq:GCF} has the same value as
\begin{equation}\label{gf:eq:evenpart}
\begin{aligned}
B_0 +  \cFrac{A_1B_2}{B_1B_2+A_2}
&-  \cFrac{A_2A_3B_4}{B_2B_3B_4+A_4B_2+A_3B_4}\\
&-  \cFrac{A_4A_5B_2B_6}{B_4B_5B_6+A_6B_4+A_5B_6}\\
&-  \cFrac{A_6A_7B_4B_8}{B_6B_7B_8+A_8B_6+A_7B_8}
- {\cdots}
\end{aligned}
\end{equation}
that is, as $B'_0+\cFrac{A'_1}{B'_1}+\cFrac{A'_2}{B'_2}+\cdots$ with
\begin{align*}
	A'_1 &=  A_1 B_2,   \qquad
	A'_2 =  - A_{2}A_{3} { B_4}, \\
	A'_j &=  - A_{2j-2} A_{2j-1} {B_{2j-4}} {  B_{2j}}\quad(j\ge3);\\
	B'_0  &=  B_0, \qquad
	B'_1 = B_1B_2+A_2, \\
	B'_j &=   B_{2j-2} B_{2j-1} B_{2j} + A_{2j} B_{2j-2} +  A_{2j-1} B_{2j}\quad(j\ge2).
\end{align*}
Its $n$-th convergent is the $2n$-th convergent of \eqref{gf:eq:GCF}.
\end{Theorem}

For the $S$-fraction \eqref{gf:eq:Sfrac}, we have $B_0=0$, $B_i=1$ for $i\ge1$, $A_1=b_0$,
and $A_i=-b_{i-1}x$ for $i\ge2$.
Thus \eqref{gf:eq:evenpart} gives the $J$-fraction \eqref{gf:eq:StoJ}.
For the super $1$-fraction \eqref{gf:eq:super1}, we have 
\begin{align*}
B_0&=0, & B_i&=1-w_i(x)x\ \ (i\geq1),\\ 
A_1&=b_0x^{\ell_0},& A_i&=-b_{i-1}x^{\ell_{i-2}+\ell_{i-1}+1}\ \ (i\geq2),
\end{align*}
Substituting these expressions into the formulas of Theorem~\ref{gf:th:evencontraction} gives the following explicit even contraction of a super $1$-fraction.

\begin{Proposition}\label{gf:prop:gec}
For the super $1$-fraction \eqref{gf:eq:super1}, the even contraction \eqref{gf:eq:evenpart} is
\[
\cFrac{A'_1}{B'_1}
+\cFrac{A'_2}{B'_2}
+\cFrac{A'_3}{B'_3}
+\cdots,
\]
where
\begin{align*}
  A'_1 &= b_0x^{\ell_0}\,(1-w_2x),\\
  B'_1 &= (1-w_1x)(1-w_2x)-b_1x^{\ell_0+\ell_1+1},\\
  A'_2 &= -\,b_1b_2\,x^{\ell_0+2\ell_1+\ell_2+2}\,(1-w_4x),\\
  A'_j &= -\,b_{2j-3}b_{2j-2}\,
          x^{\ell_{2j-4}+2\ell_{2j-3}+\ell_{2j-2}+2}\,
          (1-w_{2j-4}x)(1-w_{2j}x)\quad(j\ge3),\\
  B'_j &= (1-w_{2j-2}x)(1-w_{2j-1}x)(1-w_{2j}x)\\
       &\qquad
          -\,b_{2j-1}x^{\ell_{2j-2}+\ell_{2j-1}+1}(1-w_{2j-2}x)\\
       &\qquad
          -\,b_{2j-2}x^{\ell_{2j-3}+\ell_{2j-2}+1}(1-w_{2j}x)
          \qquad(j\ge2).
\end{align*}
\end{Proposition}

In general, the continued fraction of Proposition~\ref{gf:prop:gec} is not a super $2$-fraction for two reasons:
\begin{itemize}
\item[(\rmnum{1})]
The numerators contain polynomial factors $1-w_{2j}(x)x$, rather than being monomials.

\item[(\rmnum{2})]
The denominator terms $b_{2j-1}x^{\ell_{2j-2}+\ell_{2j-1}+1}$ exceed the degree bound $\ell_{2j-2}+1$ when $\ell_{2j-1}>0$.

\end{itemize}
Both obstructions disappear exactly when the exponents of odd index vanish: if $\ell_{2j-1}=0$, then $w_{2j}=0$ since $\deg w_{2j}\le\ell_{2j-1}-1=-1$, and $b_{2j-1}x^{\ell_{2j-2}+\ell_{2j-1}+1}=b_{2j-1}x^{\ell_{2j-2}+1}$ has the admissible degree $\ell_{2j-2}+1$.

\begin{Definition}\label{gf:def:G}
A {\it $G$-fraction} is a super $1$-fraction \eqref{gf:eq:super1} whose exponents of odd index vanish: $\ell_{2j-1}=0$ for all $j\geq1$. Equivalently, it is a continued fraction of the form
\begin{equation}\label{gf:eq:Gfrac}
\begin{aligned}
F(x)={}&\cFrac{b_0 x^{\ell_0}}{ 1 - w_1(x) x} 
- \cFrac{ b_1 x^{\ell_0+1}}{ 1}
- \cFrac{ b_2 x^{\ell_2+1}}{ 1 - w_3(x)x}\\
&- \cFrac{ b_3 x^{\ell_2+1}}{ 1}
- \cFrac{ b_4 x^{\ell_4+1}}{ 1 - w_5(x)x}
- \cFrac{ b_5 x^{\ell_4+1}}{ 1}
- {\cdots} 
\end{aligned}
\end{equation}
where $b_i\neq0$, $\ell_{2j}\geq0$ and $\deg w_{2j+1}\leq\ell_{2j}-1$.
\end{Definition}

The $S$-fraction is the $G$-fraction with all $\ell_{2j}=0$, exactly as the $J$-fraction is the $H$-fraction with all $k_j=0$. 
The even contraction from $G$-fractions to $H$-fractions extends the classical one from $S$-fractions to $J$-fractions, as illustrated below.
\[
  \begin{array}{ccc}
    G\text{-fraction} & \xrightarrow{\ \text{even contraction}\ } & H\text{-fraction}\\[2pt]
    \cup & & \cup\\[2pt]
    S\text{-fraction} & \xrightarrow{\ \text{even contraction}\ } & J\text{-fraction}
  \end{array}
\]
Setting $\chi_t=1$ for every $t\ge0$ in Theorem~\ref{gf:th:contraction} gives the following
explicit description of the even contraction from a $G$-fraction to its $H$-fraction.

\begin{Theorem}\label{gf:th:GtoH}
Let $F(x)$ be a $G$-fraction \eqref{gf:eq:Gfrac}. 
Then the even contraction of $F$ is its $H$-fraction, namely
\begin{equation}\label{gf:eq:GtoH}
  F(x)=\cFrac{b_0x^{\ell_0}}{1-u_1(x)x}
   -\cFrac{b_1b_2\,x^{\ell_0+\ell_2+2}}{1-u_2(x)x}
   -\cFrac{b_3b_4\,x^{\ell_2+\ell_4+2}}{1-u_3(x)x}-\cdots,
\end{equation}
whose data are 
\begin{equation}\label{gf:eq:GtoHdata}
\begin{aligned}
&v_0=b_0,\qquad v_j=b_{2j-1}b_{2j}\quad(j\geq1),\qquad k_j=\ell_{2j}\quad(j\geq0),\\
&u_1(x)=w_1(x)+b_1x^{\ell_0},\\ 
&u_j(x)=w_{2j-1}(x)+(b_{2j-2}+b_{2j-1})\,x^{\ell_{2j-2}}\quad(j\geq2).
\end{aligned}
\end{equation}
In particular $\deg u_j\leq\ell_{2j-2}=k_{j-1}$, and the 
coefficients of $x^{k_{j-1}}$ in the $u_j$ are 
$c_{1,k_0}=b_1$ and $c_{j,k_{j-1}}=b_{2j-2}+b_{2j-1}$ for $j\geq2$.
\end{Theorem}

\begin{proof}
For a $G$-fraction, we have $\ell_{2j-1}=0$ and $w_{2j}=0$ for every $j\ge1$.
Substituting these values into Proposition~\ref{gf:prop:gec} gives
\[
A'_1=b_0x^{\ell_0},
\qquad
B'_1=1-\bigl(w_1(x)+b_1x^{\ell_0}\bigr)x.
\]
For every $j\ge2$, it gives
\begin{align*}
A'_j
&=-b_{2j-3}b_{2j-2}
x^{\ell_{2j-4}+\ell_{2j-2}+2},\\
B'_j
&=1-\Bigl(w_{2j-1}(x)
+(b_{2j-2}+b_{2j-1})x^{\ell_{2j-2}}\Bigr)x.
\end{align*}
Thus the even contraction is \eqref{gf:eq:GtoH}, with the parameters stated in \eqref{gf:eq:GtoHdata}.

Since $\deg w_{2j-1}\le\ell_{2j-2}-1$, we have $\deg u_j\le\ell_{2j-2}=k_{j-1}$.
Moreover, $w_{2j-1}$ has no term of degree $k_{j-1}$, so the coefficient $c_{j,k_{j-1}}$ is as stated.
All numerator constants are nonzero, and all exponents are nonnegative.
Hence \eqref{gf:eq:GtoH} is a super $2$-fraction.
By Theorem~\ref{th:super2}(i), it is the unique $H$-fraction of $F$.
\end{proof}

\begin{Remark}\label{gf:rem:sharp}
The condition $w_{2j}=0$ for every $j\ge1$ is not sufficient for the even contraction to be an $H$-fraction.
Indeed, if $\ell_1\ge1$, then
\[
B'_1=1-w_1(x)x-b_1x^{\ell_0+\ell_1+1}
\]
has degree $\ell_0+\ell_1+1>\ell_0+1$, which exceeds the bound allowed for the first $H$-fraction level.
\end{Remark}

We now characterize $G$-fractions by their $1$-shifted Hankel determinants, in the same way as $S$-fractions are characterized by the non-vanishing of $H_n(F)$ and $H^{(1)}_n(F)$. 
Recall that $s_j$ are the jump orders of the $H$-fraction of $F$ and that $(R_j)$ is the sequence \eqref{def:Rn}.

\begin{Theorem}\label{gf:th:criterion}
Let $F(x)$ be a power series with super $1$-fraction 
\eqref{gf:eq:super1}. The following are equivalent:
\begin{itemize}
\item[(i)] $F$ is a $G$-fraction, that is, $\ell_{2j-1}=0$ for all $j\geq1$;
\item[(ii)] $H^{(1)}_{s_j}(F)\neq0$ for all $j\geq1$;
\item[(iii)] $R_j\neq0$ for all $j\geq1$.
\end{itemize}
Moreover, if $F$ is a $G$-fraction, then, for $j \ge 0$,
\begin{equation}\label{gf:eq:RjG}
R_j=b_1b_3\cdots b_{2j-1},
\end{equation}
and if $F$ is not a $G$-fraction and $t_0$ is the least integer 
such that $\ell_{2t_0+1}\geq1$, then \eqref{gf:eq:RjG} holds for 
$j\leq t_0$, $R_{t_0+1}=0$, and $H^{(1)}_{s_{t_0+1}}(F)=0$ 
with $s_{t_0+1}=\sigma_{2t_0+1}$.
\end{Theorem}

\begin{proof}
By \eqref{eq:H1sj:Hsj}, we have
\[
H^{(1)}_{s_j}(F)
=(-1)^{s_j+j}H_{s_j}(F)R_j,
\]
where $H_{s_j}(F)\neq0$.
Thus conditions \textnormal{(ii)} and \textnormal{(iii)} are equivalent.

Let $t_0$ be the least nonnegative integer such that $\ell_{2t_0+1}\ge1$, and set $t_0=\infty$ if no such integer exists.
Then $\chi_t=1$ for every $t<t_0$.
If $t_0<\infty$, we also have $\chi_{t_0}=0$.
By Theorem~\ref{gf:th:contraction}, for every integer $t$ with $0\le t\le t_0$, the level $(v_t,k_t,u_{t+1})$ of the $H$-fraction of $F$ is given by
\[
\begin{aligned}
 \Bigl(b_{2t-1}^{\,\pi(t\ge1)}b_{2t},\ \ell_{2t},w_{2t+1}(x)
+\bigl(\pi(t\ge1)b_{2t}
+\pi(t<t_0)b_{2t+1}\bigr)x^{\ell_{2t}}\Bigr),
\end{aligned}
\]
where $\pi(P)=1$ if $P$ is true and $0$ otherwise. 
Since $\deg w_{2t+1}\le\ell_{2t}-1$, the coefficient
of $x^{k_t}$ in $u_{t+1}(x)$ is
\[
c_{t+1,k_t}
=\pi(t\ge1)b_{2t}+\pi(t<t_0)b_{2t+1}.
\]
The recurrence \eqref{def:Rn} now gives $R_0=1$, $R_1=c_{1,k_0}=\pi(0<t_0)b_1$, and, for
$1\le t\le t_0$,
\[
R_{t+1}
=\bigl(b_{2t}+\pi(t<t_0)b_{2t+1}\bigr)R_t
-b_{2t-1}b_{2t}R_{t-1}.
\]
Suppose that $R_{t-1}=b_1b_3\cdots b_{2t-3}$ and $R_t=b_1b_3\cdots b_{2t-1}$.
Then $b_{2t}R_t-b_{2t-1}b_{2t}R_{t-1}=0$, so $R_{t+1}=\pi(t<t_0)b_{2t+1}R_t$.
By induction, $R_t=b_1b_3\cdots b_{2t-1}$ for $0\le t\le t_0$, with the empty product equal to $1$, and, if $t_0<\infty$, $R_{t_0+1}=0$.
When $t_0=0$, the latter equality follows directly from $R_1=0$.
This proves \eqref{gf:eq:RjG} when $F$ is a $G$-fraction.
Since all $b_i$ are nonzero, we obtain \textnormal{(i)}$\Rightarrow$\textnormal{(iii)}.
When $F$ is not a $G$-fraction, we have $R_{t_0+1}=0$, which proves \textnormal{(iii)}$\Rightarrow$\textnormal{(i)} by contraposition and gives
$H^{(1)}_{s_{t_0+1}}(F)=0$.
In this case, $s_{t_0+1}=\sigma_{2t_0+1}$ by Theorem~\ref{gf:th:det}, since $\sigma_{2t+1}=\sigma_{2t+2}$ for $t<t_0$.
\end{proof}

\begin{Remark}\label{gf:rem:example}
The exponents of the super $1$-fraction and those of the 
$H$-fraction of the same series $F$ reflect different 
determinants: the latter are governed by the vanishing of the 
Hankel determinants $H_n(F)$ only, whereas the former also 
involve the $1$-shifted Hankel determinants, as Lemma 
\ref{gf:lem:peel1} shows. 
For instance, let $\ell_0=0$, $\ell_1=1$ and $\ell_i=0$ for $i\geq2$, 
so that $w_2(x)=c$ is a constant and 
\[
F(x)=\cFrac{b_0}{1}-\cFrac{b_1x^2}{1-cx}-\cFrac{b_2x^2}{1}-\cFrac{b_3x}{1}-\cFrac{b_4x}{1}-\cdots.
\]
By Theorem \ref{gf:th:contraction} ($\chi_0=0$, $\chi_t=1$ for $t\geq1$), 
the $H$-fraction of $F$ is 
\[
F(x)=\cFrac{b_0}{1}-\cFrac{b_1x^2}{1-cx}-\cFrac{b_2x^2}{1-b_3x}
-\cFrac{b_3b_4x^2}{1-(b_4+b_5)x}-\cdots,
\]
a $J$-fraction: all Hankel determinants $H_n(F)$ are nonzero, 
although the super $1$-fraction of $F$ has a nonzero exponent 
$\ell_1=1$. In accordance with Theorem \ref{gf:th:criterion} 
($t_0=0$), $R_1=c_{1,0}=0$ and $H^{(1)}_1(F)=0$: indeed the 
coefficient of $x$ in $F$ vanishes.
\end{Remark}

\subsection{Shifted Hankel determinants of a $G$-fraction}\label{gf:subsec:shift}

Throughout this subsection, $F$ is a $G$-fraction \eqref{gf:eq:Gfrac}. 
By Theorem \ref{gf:th:GtoH}, the data of its $H$-fraction are given by \eqref{gf:eq:GtoHdata}.
Thus, for every $j\ge0$, its jump orders, sign exponents and products of numerator constants are
\begin{equation}\label{gf:eq:SEP}
s_j=\sum_{t=0}^{j-1}\ell_{2t}+j,\quad 
\epsilon_j=\sum_{t=0}^{j-1}\binom{\ell_{2t}+1}{2},\quad 
V_j=b_0b_1\cdots b_{2j}
\end{equation}
In the notation of \eqref{gf:eq:sigma} and \eqref{gf:eq:e}, we have $s_j=\sigma_{2j-1}=\sigma_{2j}$ for $j\ge1$ and $\epsilon_j=e_{2j}$ for $j\ge0$.
Theorem~\ref{gf:th:det} therefore gives, for every $j\ge0$,
\begin{equation}\label{gf:eq:HsjG}
H_{s_j}(F)=(-1)^{\epsilon_j}\,b_0^{\,s_j}
\prod_{i=1}^{j-1}(b_{2i-1}b_{2i})^{\,s_j-s_i}.
\end{equation}
All other Hankel determinants of $F$ being zero.
Moreover, by \eqref{gf:eq:RjG}, $R_j=b_1b_3\cdots b_{2j-1}$. 

Substituting these data into Theorems~\ref{th:shiftHDet} and~\ref{th:shift2HDet} yields the following two theorems, which express the $1$-shifted and $2$-shifted Hankel determinants of a $G$-fraction directly in terms of its own datas.

\begin{Theorem}\label{gf:th:shiftG}
Let $F(x)$ be a $G$-fraction \eqref{gf:eq:Gfrac}, with $s_j$, 
$\epsilon_j$, $V_j$ given by \eqref{gf:eq:SEP}. 
Then all non-vanishing $1$-shifted Hankel determinants of $F$ 
are given by
\begin{equation}\label{gf:eq:H1sjG}
\begin{aligned}
H^{(1)}_{s_j}(F)&=(-1)^{s_j+j}\,H_{s_j}(F)\,b_1b_3\cdots b_{2j-1}\\
&=(-1)^{\epsilon_j+s_j+j}\,b_0^{\,s_j}\,b_1b_3\cdots b_{2j-1}
\prod_{i=1}^{j-1}(b_{2i-1}b_{2i})^{\,s_j-s_i}
\end{aligned}
\end{equation}
for $j\geq1$ and, for $j \ge 0$,
\begin{equation}\label{gf:eq:H1sjkjG}
H^{(1)}_{s_j+\ell_{2j}}(F)=(-1)^{\ell_{2j}(\ell_{2j}-1)/2}\,
(b_0b_1\cdots b_{2j})^{\ell_{2j}}\,H^{(1)}_{s_j}(F),
\end{equation}
where $H^{(1)}_0(F)=1$. 
In particular $H^{(1)}_{s_j}(F)\neq0$ for all $j$.
\end{Theorem}

\begin{proof}
These are formulas \eqref{eq:H1sj:Hsj} and \eqref{eq:H1sjkj} 
of Theorem \ref{th:shiftHDet}, in which we substitute 
$R_j=b_1b_3\cdots b_{2j-1}$, $k_j=\ell_{2j}$, 
$V_j=b_0b_1\cdots b_{2j}$ and \eqref{gf:eq:HsjG}.
\end{proof}

When all $\ell_i=0$, we have $s_j=j$ and $\epsilon_j=0$.
Thus \eqref{gf:eq:H1sjG} reduces to the classical formula
\[
H^{(1)}_n(F)
=(b_0b_1)^n
\prod_{i=1}^{n-1}(b_{2i}b_{2i+1})^{n-i}
\]
for the $1$-shifted Hankel determinants of the $S$-fraction \eqref{gf:eq:Sfrac}.
This formula also follows from Heilermann's formula \eqref{eq:detH}, applied to the $J$-fraction obtained by odd contraction:
\[
\frac{F-b_0}{x}
=
\cFrac{b_0b_1}{1-(b_1+b_2)x}
-\cFrac{b_2b_3x^2}{1-(b_3+b_4)x}
-\cFrac{b_4b_5x^2}{1-(b_5+b_6)x}
-\cdots.
\]

To evaluate the $2$-shifted Hankel determinants, we also need the coefficients $c_{j,k_{j-1}-1}$ when $k_{j-1}\ge1$.
By \eqref{gf:eq:GtoHdata}, the difference $u_j(x)-w_{2j-1}(x)$ is a scalar multiple of
$x^{\ell_{2j-2}}$.
Hence, for every $j\ge1$ with $\ell_{2j-2}\ge1$,
\begin{equation}\label{gf:eq:Cj}
c_{j,k_{j-1}-1}
=
\bigl[x^{\ell_{2j-2}-1}\bigr]w_{2j-1}(x),
\end{equation}
where $[x^m]w$ denotes the coefficient of $x^m$ in $w$.
This coefficient may be zero.
For every $j\ge1$, the coefficient $d_j$ of $x^{k_{j-1}}$ in $1-u_j(x)x$, as defined in
Theorem~\ref{th:shift2HDet}, is therefore
\begin{equation}\label{gf:eq:Dj}
d_j=
\begin{cases}
1, & \text{if $\ell_{2j-2}=0$},\\
-c_{j,k_{j-1}-1}, & \text{if $\ell_{2j-2}\ge1$}.
\end{cases}
\end{equation}

We now use these coefficients to prove the following formulas for the $2$-shifted Hankel determinants of a $G$-fraction.

\begin{Theorem}\label{gf:th:shift2G}
Let $F(x)$ be a $G$-fraction \eqref{gf:eq:Gfrac}, with $s_j$,  $V_j$ given by \eqref{gf:eq:SEP}, $H_{s_j}(F)$ by \eqref{gf:eq:HsjG},  and $c_{j,k_{j-1}-1}$, $d_j$ by \eqref{gf:eq:Cj}, \eqref{gf:eq:Dj}. 
For $j\geq0$, the $2$-shifted Hankel determinants of $F$  satisfy the following formulas.

\smallskip
\noindent{\bf (1)} If $\ell_{2j}\geq1$, then
\begin{align}
& H^{(2)}_{s_j}(F) =H_{s_j}(F)\,(b_1b_3\cdots b_{2j-1})^2,\label{gf:eq:H2a}\\
& H^{(2)}_{s_j+\ell_{2j}-1}(F) =(-1)^{\frac{(\ell_{2j}-2)(\ell_{2j}-1)}{2}}\,
   V_j^{\,\ell_{2j}-1}\,H^{(2)}_{s_j}(F),\nonumber \\
& H^{(2)}_{s_j+\ell_{2j}}(F)=c_{j+1,k_j-1}\,V_j\,H^{(2)}_{s_j+\ell_{2j}-1}(F)
   +(-1)^{\frac{\ell_{2j}(\ell_{2j}+1)}{2}}\,V_j^{\,\ell_{2j}+1}\,H^{(2)}_{s_j-1}(F), \nonumber
\end{align}
where the last term is absent for $j=0$.

\smallskip
\noindent{\bf (2)} If $\ell_{2j}=0$, then
\begin{equation}\label{gf:eq:H2d}
H^{(2)}_{s_j}(F)=H_{s_j}(F)\sum_{t=0}^{j}d_{t+1}\,(b_1b_3\cdots b_{2t-1})^2
\prod_{i=t+1}^{j}b_{2i-1}b_{2i}.
\end{equation}

\smallskip
\noindent{\bf (3)} All other $2$-shifted Hankel determinants of $F$ vanish.
\end{Theorem}

\begin{proof}
We apply Theorem~\ref{th:shift2HDet} with $k_j=\ell_{2j}$, $R_j=b_1b_3\cdots b_{2j-1}$,
$V_j=b_0b_1\cdots b_{2j}$, and
\[
v_{t+1}v_{t+2}\cdots v_j
=
\prod_{i=t+1}^{j}b_{2i-1}b_{2i}.
\]
The coefficients $c_{j+1,k_j-1}$, when $k_j\ge1$, and $d_{t+1}$ are given by \eqref{gf:eq:Cj}
and \eqref{gf:eq:Dj}, respectively.
For $j\ge1$, the determinant appearing in the term $\Delta_j$ of Theorem~\ref{th:shift2HDet}
has index $s_{j-1}$ if $k_{j-1}=0$, and $s_{j-1}+k_{j-1}$ otherwise.
Since $s_j=s_{j-1}+k_{j-1}+1$, this index is $s_j-1$ in both cases.
The stated formulas and vanishing assertions follow.
\end{proof}

\begin{Remark}\label{gf:rem:DJ}
Comparing \eqref{gf:eq:H2a} with \eqref{gf:eq:H1sjG}, we get
\[
H_{s_j}(F)\,H^{(2)}_{s_j}(F)=\bigl(H^{(1)}_{s_j}(F)\bigr)^2
\qquad\text{whenever $\ell_{2j}\geq1$.}
\] 
This is consistent with the Desnanot--Jacobi identity 
\[
H_{n+1}(F)\,H^{(2)}_{n-1}(F)=H_n(F)\,H^{(2)}_n(F)-\bigl(H^{(1)}_n(F)\bigr)^2,
\] 
since $H_{s_j+1}(F)=0$ when $\ell_{2j}\geq1$. 
When all $\ell_i=0$, only case (2) occurs, all $d_t=1$, and \eqref{gf:eq:H2d} is the classical evaluation of the $2$-shifted Hankel determinants of an $S$-fraction, 
namely \eqref{eq:H2classical} for the $J$-fraction~\eqref{gf:eq:StoJ}.
\end{Remark}

\subsection{The classical cases: $S$-fractions and $J$-fractions}\label{gf:subsec:classical}
We conclude by listing the values of the shifted Hankel determinants of $S$-fractions, which are the case ``all $\ell_i=0$'' of the preceding results, together with the Hankel determinants of $F(x^2)$, which are obtained from those of $F$ by the following 
elementary lemma.
The corresponding formulas for $J$-fractions have already been collected in Corollary~\ref{cor:J}.
These formulas are classical; they go back to Stieltjes \cite{Stieltjes1894Re} and Heilermann \cite{Heilermann1846}, and can be found, with combinatorial proofs by means of 
non-intersecting lattice paths, in \cite{Viennot1983, Krattenthaler1998, Krattenthaler2005}; 
see also \cite{Wall1948, Jones1980T}.

\begin{Lemma}\label{gf:lem:x2}
For every power series $F(x)$ and every $n\geq0$,
\begin{align*}
    H_{2n}\bigl(F(x^2)\bigr) &= H_n(F)\,H^{(1)}_n(F),\\
    H_{2n+1}\bigl(F(x^2)\bigr) &= H_{n+1}(F)\,H^{(1)}_n(F).
\end{align*}
\end{Lemma}

\begin{proof}
The case $n=0$ is immediate, with the convention
$H_0(F)=H_0^{(1)}(F)=1$.
Let $N\geq 1$ and write $F(x)=\sum_{m\geq 0}a_mx^m$.
The coefficients $\tilde a_i$ of $F(x^2)$ satisfy
$\tilde a_{2m}=a_m$ and $\tilde a_{2m+1}=0$.
Thus, in the Hankel matrix
\[
M_N=(\tilde a_{i+j})_{0\leq i,j\leq N-1},
\]
all entries with $i+j$ odd vanish.
For example, when $N=4$, we have
\[
M_4=
\begin{pmatrix}
a_0 & 0   & a_1 & 0   \\
0   & a_1 & 0   & a_2 \\
a_1 & 0   & a_2 & 0   \\
0   & a_2 & 0   & a_3
\end{pmatrix}.
\]
Reorder both the rows and the columns by placing the even
indices before the odd ones, preserving the order within
each group.
Since the same permutation is applied to the rows and the
columns, the determinant remains unchanged.
The resulting matrix is block diagonal.
In the example above, the new index order is $0,2,1,3$,
and the reordered matrix is
\[
\begin{pmatrix}
a_0 & a_1 & 0   & 0   \\
a_1 & a_2 & 0   & 0   \\
0   & 0   & a_1 & a_2 \\
0   & 0   & a_2 & a_3
\end{pmatrix}.
\]
In general, the block corresponding to the even indices is
\[
(\tilde a_{2p+2q})_{0\leq p,q<\lceil N/2\rceil}
=
(a_{p+q})_{0\leq p,q<\lceil N/2\rceil},
\]
whereas the block corresponding to the odd indices is
\[
(\tilde a_{2p+2q+2})_{0\leq p,q<\lfloor N/2\rfloor}
=
(a_{p+q+1})_{0\leq p,q<\lfloor N/2\rfloor}.
\]
The determinant of $M_N$ is therefore the product of the determinants of these two blocks.
\begin{itemize}
    \item If $N=2n$, both blocks have size $n$, giving
\[
H_{2n}\bigl(F(x^2)\bigr)=H_n(F)\,H_n^{(1)}(F).
\]
\item If $N=2n+1$, the even-index block has size $n+1$ and
the odd-index block has size $n$, giving
\[
H_{2n+1}\bigl(F(x^2)\bigr)=H_{n+1}(F)\,H_n^{(1)}(F).
\]
\end{itemize}
This completes the proof.
\end{proof}

Combined with Theorems \ref{th:super2} and \ref{th:shiftHDet}, Lemma \ref{gf:lem:x2} gives all Hankel determinants of $F(x^2)$ in terms of the $H$-fraction of $F$. 
For $S$-fractions, the preceding results yield the explicit formulas below, of which \eqref{gf:eq:S0} and \eqref{gf:eq:S1} go back to Stieltjes \cite{Stieltjes1894Re};
see also \cite{Krattenthaler1998, Krattenthaler2005}.

\begin{Corollary}\label{gf:cor:S}
Let $F(x)$ be the $S$-fraction \eqref{gf:eq:Sfrac}, with all $b_i\neq0$. 
Then, for $n\geq1$,
\begin{align}
H_n(F)&=b_0^{\,n}\,(b_1b_2)^{n-1}(b_3b_4)^{n-2}\cdots(b_{2n-3}b_{2n-2}),\label{gf:eq:S0}\\
H^{(1)}_n(F)&=(b_0b_1)^{n}\,(b_2b_3)^{n-1}(b_4b_5)^{n-2}\cdots(b_{2n-2}b_{2n-1}),\label{gf:eq:S1}\\
H^{(2)}_n(F)&=H_n(F)\sum_{t=0}^{n}(b_1b_3\cdots b_{2t-1})^2
\prod_{i=t+1}^{n}b_{2i-1}b_{2i},\label{gf:eq:S2}\\
H^{(-1)}_n(F)&=-\,b_0\,H_{n-1}(F)\,E_{n-1},\label{gf:eq:Sm1}\\
H_n\bigl(F(x^2)\bigr)&=b_0^{\,n}\,b_1^{\,n-1}\,b_2^{\,n-2}\cdots b_{n-1},\label{gf:eq:Sx2}
\end{align}
where $E_0=0$, $E_1=1$ and 
$E_{m+1}=(b_{2m}+b_{2m+1})E_m-b_{2m-1}b_{2m}E_{m-1}$ for $m\geq1$.
\end{Corollary}

\begin{proof}
Formulas \eqref{gf:eq:S0}--\eqref{gf:eq:S2} follow from \eqref{gf:eq:HsjG}, \eqref{gf:eq:H1sjG}, and \eqref{gf:eq:H2d}, respectively, by taking $\ell_i=0$ for all $i$, so that $s_j=j$, $\epsilon_j=0$, and $d_t=1$.
Formula \eqref{gf:eq:Sm1} follows by applying \eqref{eq:Hm1classical} to the $J$-fraction
\eqref{gf:eq:StoJ}.
Finally, \eqref{gf:eq:Sx2} follows from \eqref{gf:eq:S0}, \eqref{gf:eq:S1}, and
Lemma~\ref{gf:lem:x2}.
Alternatively,
\[
F(x^2)
=\cFrac{b_0}{1}
-\cFrac{b_1x^2}{1}
-\cFrac{b_2x^2}{1}
-\cdots
\]
is a $J$-fraction with $u_j=0$ and $v_j=b_j$, so \eqref{gf:eq:Sx2} also follows directly from
Heilermann's formula \eqref{eq:detH}.
\end{proof}

For $J$-fractions, all the corresponding formulas, including 
those for $H_n(F(x^2))$ obtained from Lemma \ref{gf:lem:x2}, 
have been collected in Corollary \ref{cor:J} of Section \ref{subsec:Jcase}. 
Note that, for an arbitrary $H$-fraction \eqref{eq:Hfrac}, 
the sequences $(R_n)$ and $(E_n)$ of \eqref{def:Rn} and 
\eqref{def:En} are given by the tridiagonal determinants 
\eqref{eq:RdetJ} with $u_j$ replaced by the coefficients of 
highest degree $c_{j,k_{j-1}}$.


\section{Concluding remarks}\label{sec:remarks}

\subsection*{Magnus' P-fractions}\label{subsec:magnus}
The $H$-fractions are, up to an equivalence transformation, the 
P-fractions introduced by Magnus \cite{Magnus1962C, Magnus1962P} 
in connection with the Pad\'e table \cite{Magnus1970, Magnus1974}.
A P-fraction of a power series $L=c_0+c_1x+c_2x^2+\cdots$ is written as 
\[
L=c_0+\cFrac{1}{b_1}
+\cFrac{1}{b_2}
+\cFrac{1}{b_3}
+\cdots.
\]
 where $b_n$ is a polynomial in 
$x^{-1}$ of degree $N_n\ge1$ with leading coefficient $\alpha_n$. 
Setting $F(x)=(L-c_0)/x$, the dictionary with \eqref{eq:Hfrac} is 
\[
N_n=k_{n-1}+1,\qquad v_0=\frac1{\alpha_1},\qquad 
v_n=-\frac1{\alpha_n\alpha_{n+1}}\quad(n\ge1),
\]
and Magnus' determinants $\varphi_m=\det(c_{i+j-1})_{1\le i,j\le m}$ 
are the Hankel determinants $H_m(F)$. 
In \cite{Magnus1962P}, Magnus computes the approximants of the 
P-fraction from determinants of $L$. 
Along the way he proves that $\varphi_m\neq0$ if and only if 
$m=N_1+\cdots+N_n=s_n$ for some $n$, which is the vanishing pattern 
of Theorem \ref{th:super2}(ii), and he establishes two relations 
\cite[(2.9), (2.14)]{Magnus1962P} between $\varphi_{s_n}$, 
$\varphi_{s_{n+1}}$, a bordered determinant $\psi_n$, and the leading 
coefficients $\kappa_n=\alpha_1\cdots\alpha_n$ of the denominators 
of the approximants: 
\begin{equation*}
\begin{aligned}
(-1)^n\varphi_{s_n}&=\kappa_n\kappa_{n+1}\psi_n,\\
\varphi_{s_{n+1}}&=(-1)^{N(N-1)/2}\,\frac{\psi_n^{N}}{\varphi_{s_n}^{N-1}}
\qquad(N=N_{n+1}).
\end{aligned}
\end{equation*}
Eliminating $\psi_n$ gives 
\[
\varphi_{s_{n+1}}=(-1)^{N(N-1)/2}\Bigl(\frac{(-1)^n}{\kappa_n\kappa_{n+1}}\Bigr)^{N}\varphi_{s_n}.
\]
Since $V_n=(-1)^n/(\kappa_n\kappa_{n+1})$ and $N=k_n+1$, this becomes
\[
H_{s_{n+1}}(F)
=(-1)^{k_n(k_n+1)/2}V_n^{k_n+1}H_{s_n}(F).
\]
Iterating this relation yields \eqref{eq:HankelDetFundamental}.
Magnus, however, uses these relations in the opposite direction, 
to compute the P-fraction from the determinants of the series, and 
does not state the evaluation of the Hankel determinants. 

\subsection*{Roblet's P-paths}\label{subsec:roblet}
In his thesis \cite{Roblet1994}, Roblet developed a combinatorial theory
of Pad\'e approximants, based on path models for several families of
continued fractions.
For the P-fractions of Magnus, he introduced the {\it P-paths}
\cite[Chapter 2]{Roblet1994}.
A P-fraction is written as
\[
f(t)=c_0+
\cFrac{\lambda_1t^{N_1}}{P_1(t)}
-\cFrac{\lambda_2t^{N_1+N_2}}{P_2(t)}
-\cFrac{\lambda_3t^{N_2+N_3}}{P_3(t)}
-\cdots.
\]
where $P_k(t)=1-a_{k,1}t-\cdots-a_{k,N_k}t^{N_k}$,
and the dictionary with \eqref{eq:Hfrac} is
\[
F(x)=\frac{f(x)-c_0}{x},\quad
N_k=k_{k-1}+1,\quad \lambda_k=v_{k-1},\quad a_{k,i}=c_{k,i-1}.
\]

\smallskip
\noindent{\bf (1)} 
Our quasi-Dyck paths are very close to the P-paths:
a quasi-Dyck path of length $2n$, representing the coefficient $a_n$ of $F$, is obtained from a P-path of length $2n+2$, representing the coefficient $c_{n+1}$ of $f$, by deleting its first up step and lowering the remaining steps by one unit, the last (long) down step being shortened accordingly; the levels $M_k$ of \eqref{eq:MN} are the levels of the P-paths
minus one.
The two models are therefore essentially equivalent.
We prefer the quasi-Dyck paths because they represent the $H$-fraction $F$ itself, and not the series $c_0+xF$, and hence directly generalize Flajolet's model for $J$-fractions
\cite{Flajolet1980}:
when all $k_j$ are zero, $F$ is a $J$-fraction and the quasi-Dyck paths of length $2n$ are in bijection with Flajolet's weighted Motzkin paths of length $n$, whereas Roblet's Jacobi P-paths represent the coefficients of $f=c_0+xF$, and are mapped onto Motzkin paths by
a compression of their steps \cite[\S2.3.2]{Roblet1994}.

\smallskip
\noindent{\bf (2)} Roblet expressed the determinants of the C-table of $f$ as signed
sums over families of non-intersecting P-paths, by the method of
Gessel and Viennot \cite[\S2.4]{Roblet1994}.
In our notation, the diagonal determinant $C(n/n)$ of \cite{Roblet1994} is $H_n(F)$, and he evaluated it explicitly \cite[Proposition 32]{Roblet1994}.
This evaluation is precisely Theorem~\ref{th:super2}(ii), with essentially the same combinatorial proof as the one given in Section~\ref{sec:lgv}.
The neighboring determinants $C(n+1/n)=H^{(1)}_n(F)$ and
$C(n/n+1)=H_{n+1}(c_0+xF)$ were also interpreted there as signed sums
over (almost) non-intersecting families of P-paths,
so that the model could be used for deriving a formula for
$H^{(1)}_n(F)$.
However, these sums were not evaluated.
They were used to prove the {\it cross identity}
\[
C(p/q)^2=C(p+1/q)C(p-1/q)-C(p/q+1)C(p/q-1)
\]
\cite[Proposition 35]{Roblet1994}, which, for $p=q=n$ and by
\eqref{eq:Hc}, reads
\[
H_n(F)^2=H^{(1)}_n(F)\,H^{(-1)}_n(F)-H^{(1)}_{n-1}(F)\,H^{(-1)}_{n+1}(F),
\]
and which allows one to derive the backward shifted determinants
$H^{(-1)}_n(F)$ from $H_n(F)$ and $H^{(1)}_n(F)$, at least when the
latter do not vanish.
Again, no explicit formula for $H^{(-1)}_n(F)$ was given.
Theorems \ref{th:shiftHDet} and \ref{th:backshiftHDet} provide these
formulas, together with the exact vanishing pattern, for an arbitrary
$H$-fraction.

\subsection*{Buslaev's and Holtz--Tyaglov's formulas}\label{subsec:buslaev}
The evaluation \eqref{eq:HankelDetFundamental} was also obtained by 
Buslaev and by Holtz and Tyaglov, by algebraic means.

\smallskip
\noindent{\bf (1)} Buslaev \cite[Theorem 1]{Buslaev2010} considers a P-fraction
\[
\cFrac{a_1}{b_1(\xi)}
+\cFrac{a_2}{b_2(\xi)}
+\cFrac{a_3}{b_3(\xi)}
+\cdots,
\]
 where $a_k\neq0$ and the 
$b_k(\xi)$ are monic polynomials of degrees $\beta_k$ with 
$\beta_0=0$ and $\beta_{k-1}+\beta_k>0$, and its associated 
(possibly formal) series $\sum_{p\ge1}f_p\xi^{-p}$. 
He proves that, up to a sign,
\[
\det(f_{i+l+1})_{0\le i,l\le n_k-1}=\prod_{j=1}^{k}a_j^{\,n_k-n_{j-1}},
\qquad n_k=\beta_1+\cdots+\beta_k.
\]
When all $\beta_k\geq 1$, we set $\xi=1/x$ and $F(x)=\sum_{p\geq 0}f_{p+1}x^p$.
An equivalence transformation then turns the P-fraction into the $H$-fraction \eqref{eq:Hfrac} of $F$, with $\beta_k=k_{k-1}+1$, $n_k=s_k$, $a_1=v_0$, and $a_{k+1}=-v_k$.
Buslaev's formula is then \eqref{eq:HankelDetFundamental} up to 
the sign, which is not given explicitly in \cite{Buslaev2010}, and 
the vanishing of the other Hankel determinants is not stated there. 
The proof uses column operations based on the denominators of the 
convergents, in the spirit of Magnus' relations and of 
Section \ref{sec:doubleshalg}. 
Buslaev's setting also allows some degrees $\beta_k$ to vanish, 
which covers the $C$-fractions whose exponents satisfy a certain 
alternating condition.
He applies his formula, together with P\'olya's theorem, to the singularities of functions defined by $C$-fractions (Gonchar's theorem on the Leighton conjecture). 
Cigler \cite[Theorem 3.1]{Cigler2013spe} later reproved the evaluation for $C$-fractions by means of the reduction lemma discussed below.

\smallskip
\noindent{\bf (2)} Holtz and Tyaglov \cite[\S1.3]{Holtz2012T} consider a rational 
function 
\[
R(z)=\frac{s_0}{z}+\frac{s_1}{z^2}+\frac{s_2}{z^3}+\cdots
\]
 and its finite continued  fraction
 \[
R(z)=
\cFrac{1}{q_1(z)}
+\cFrac{1}{q_2(z)}
+\cdots
+\cFrac{1}{q_k(z)}.
\]
produced by the Euclidean algorithm, where $q_j(z)=\alpha_jz^{n_j}+\cdots$ with 
$\alpha_j\neq0$ and $n_j\ge1$. 
Their Theorem~1.18 is the one-step reduction lemma, in the 
form (5) below. By iterating it, they obtain an explicit formula, 
with its sign, for the Hankel minors 
$D_{n_1+\cdots+n_j}(R)$, where $D_m(R)=\det(s_{i+l})_{0\le i,l\le m-1}$ 
\cite[(1.62)]{Holtz2012T}, and they note that all the other 
minors vanish \cite[Remark 1.19]{Holtz2012T}. 
For the continued fraction with minus signs 
\[
\cFrac{1}{q_1(z)}
-\cFrac{1}{q_2(z)}
-\cdots
-\cFrac{1}{q_k(z)},
\]
 given by the Sturm algorithm, 
the formula reads \cite[(2.18)]{Holtz2012T}
\[
D_{n_1+\cdots+n_j}(R)=\prod_{i=1}^{j}(-1)^{n_i(n_i-1)/2}\,
\prod_{i=1}^{j}\alpha_i^{-(n_i+2\sum_{\rho=i+1}^{j}n_\rho)}.
\]
With $z=1/x$ and $R(1/x)=xF(x)$, this continued fraction becomes 
a finite $H$-fraction of $F$, with $n_i=k_{i-1}+1$, $v_0=1/\alpha_1$ 
and $v_i=1/(\alpha_i\alpha_{i+1})$ $(i\ge1)$, and $D_m(R)=H_m(F)$. 
Since $n_i+2\sum_{\rho>i}n_\rho=2s_j-s_i-s_{i-1}$ and 
$\sum_i n_i(n_i-1)/2=\epsilon_j$, their formula is exactly 
\eqref{eq:HankelDetFundamental}, sign included. 
Their setting is restricted to rational functions, that is, to 
finite $H$-fractions, but the proof of their Theorem~1.18 only uses 
formal Laurent series, and hence applies to arbitrary power series.

\subsection*{The one-step reduction lemma}\label{subsec:reduction}
The algebraic proofs of Sections \ref{sec:shiftalg} and
\ref{sec:doubleshalg} rest on Lemma \ref{lem:alg:han}: if
\[
F(x)=\frac{v\,x^k}{1-u(x)x-x^{k+2}G(x)},\qquad \deg u\le k,
\]
then $H_n(F)=0$ for $1\le n\le k$ and
$H_n(F)=(-1)^{k(k+1)/2}v^nH_{n-k-1}(G)$ for $n\ge k+1$.
This identity, or equivalent forms of it, appeared several times in
the literature, in different notations and settings.

\smallskip
\noindent{\bf (1)} In \cite[Lemma 2.2]{Han2016Adv} it is stated with $v=1$ and the
denominator $1+u(x)x-x^{k+2}G(x)$, in order to derive
Theorem \ref{th:super2}(ii) by iteration along the $H$-fraction.
The proof is a matrix factorization
$\mathbf F_1\mathbf G=\mathbf F\mathbf B$, where $\mathbf F$ is the
Hankel matrix of $F$, the matrices $\mathbf F_1$ and $\mathbf B$ are
unitriangular, and $\mathbf G$ is block diagonal, with an
anti-triangular block of determinant $(-1)^{k(k+1)/2}$ and the Hankel
matrix of $G$.

\smallskip
\noindent{\bf (2)} Andrews and Wimp \cite[Theorem 1]{Andrews2002Wimp} consider two
sequences with $t_0=1$ and $\sum_{j=0}^Nt_js_{N-j}=s_0\delta_{N,0}$,
that is, $s(x)t(x)=s_0$, and prove that
\[
\det(s_{i+j})_{0\le i,j\le N}=(-1)^Ns_0^{N+1}\det(t_{i+j+2})_{0\le i,j\le N-1},
\]
that is, $
H_{N+1}(s)=(-1)^Ns_0^{N+1}H^{(2)}_N(t).$
Writing $t(x)=1+t_1x+x^2T(x)$, we have
\[
s=\frac{s_0}{1+t_1x-x^2(-T(x))},
\] and Lemma \ref{lem:alg:han} with $k=0$, $v=s_0$ and $G=-T$ gives exactly this identity. 
The sign $(-1)^N$ comes from $H_N(-T)=(-1)^NH_N(T)$.
Thus, \cite[Theorem 1]{Andrews2002Wimp} is the case $k=0$.
Its proof uses row and column operations based on the relation
$s(x)t(x)=s_0$.

\smallskip
\noindent{\bf (3)} Sulanke and Xin \cite[Section 4]{Sulanke2005Xin} study series
defined by a quadratic equation 
\[
F=\frac{x^d}{u(x)+x^kv(x)F(x)},
\]
with $u(0)v(0)\neq0$ and $k\ge1$, in order to evaluate Hankel determinants
by successive transformations.
Their Proposition 4.1 is the factor $v^n$ above.
In their Proposition 4.2(ii) ($k\ge2$, $u(0)=1$), they write
$u=u_L+x^{d+2}u_H$ with $\deg u_L\le d+1$, so that
\[
F=\frac{x^d}{u_L(x)-x^{d+2}G(x)}
\] with $G=-u_H-x^{k-d-2}vF$, and prove
\[
\det H_{n-d-1}(G)=(-1)^{d(d+1)/2}\det H_n(F).\]
This is Lemma \ref{lem:alg:han} with $k=d$ and $v=1$;
the proof, which uses the generating function
$(xF(x)-yF(y))/(x-y)$ of the Hankel matrix and a block-diagonal
matrix, does not use the quadratic equation, and hence holds for an
arbitrary power series $G$.
Their Proposition 4.2(i) ($k=1$) gives, for
\[
F=\frac{x^d}{u_L(x)-x^{d+1}G(x)},
\] the relation
$H^{(1)}_{n-d-1}(G)=(-1)^{d(d+1)/2}H_n(F)$; it reduces to
Lemma \ref{lem:alg:han} by moving $G(0)x^{d+1}$ into $u_L$ and
replacing $G$ by $(G-G(0))/x$, whose Hankel determinants are
the $1$-shifted ones of $G$.

\smallskip
\noindent{\bf (4)} Cigler \cite[Lemma 2.1]{Cigler2013spe} considers
$f(x)=1/(1-ax^pg(x))$ with $g(0)=1$ and $p\ge1$, and evaluates
$\det(f_{i+j-m})_{0\le i,j\le N}$ for $m\ge-1$: these determinants
vanish for $N<m$, are equal to $(-1)^{m(m+1)/2}$ for $N=m$, and
\[
\det(f_{i+j-m})_{0\le i,j\le n+m}
=(-1)^{m(m+1)/2}a^n\det(g_{i+j+m-p+2})_{0\le i,j\le n-1},
\]
with $g_i=0$ for $i<0$.
For $m\ge0$, the left-hand side is $H_{n+m+1}(x^mf)$.
Splitting $ax^pg(x)$ into its terms of degree at most $m+1$ and the
remaining ones, we get
\[
x^mf=\frac{x^m}{1-u(x)x-x^{m+2}aG(x)}
\] with $\deg u\le m$ and
$G(x)=\sum_{j\ge0}g_{j+m+2-p}x^j$, so that the three statements are
Lemma \ref{lem:alg:han} with $k=m$ and $v=1$.
Conversely, every series of the form of Lemma \ref{lem:alg:han} with
$v=1$ can be written as $x^k/(1-ax^pg(x))$, where $ax^p$ is the lowest
term of $u(x)x+x^{k+2}G(x)$ (if it is nonzero), so that the two
statements are equivalent.
The case $m=-1$ reads $H^{(1)}_n(f)=H_n(ax^{p-1}g)$.
Since $(f-1)/x=K/(1-xK)$ with $K=ax^{p-1}g$, this is
Lemma \ref{lem:alg:KxK}.
Cigler uses this lemma to give a proof of Buslaev's evaluation of
the Hankel determinants of $C$-fractions \cite[Theorem 3.1]{Cigler2013spe}.

\smallskip
\noindent{\bf (5)} Holtz and Tyaglov \cite[Theorem 1.18]{Holtz2012T} prove that two 
rational functions $R$ and $G$ vanishing at infinity satisfy 
$R=1/(g+G)$, where $g(z)=\alpha z^m+\cdots$ is a polynomial of degree 
$m\ge1$, if and only if $D_1(R)=\cdots=D_{m-1}(R)=0$ and
\[
D_{m+j}(R)=(-1)^{m(m-1)/2}\,(-1)^j\,\frac{D_j(G)}{\alpha^{m+2j}}
\qquad(j\ge0),
\]
with the notation $D_j$ of the previous subsection. 
With $z=1/x$, write $g(1/x)=\alpha x^{-m}(1-u(x)x)$ with 
$\deg u\le m-1$, $G(1/x)=x\widehat G(x)$ and $R(1/x)=xF(x)$. Then
\[
F=\frac{v\,x^k}{1-u(x)x-x^{k+2}G_1(x)}
\] with $k=m-1$, $v=1/\alpha$ and 
$G_1=-\widehat G/\alpha$, and, since $H_j(G_1)=(-1/\alpha)^jH_j(\widehat G)$, 
the direct part of their theorem is Lemma \ref{lem:alg:han}. 
The proof uses column operations based on the coefficients of $g$.

In summary, the five statements are versions of the same one-step
reduction, covering the case $k=0$ \cite{Andrews2002Wimp}, the
general case in the setting of quadratic equations
\cite{Sulanke2005Xin}, the general case in the setting of
$C$-fractions \cite{Cigler2013spe}, the general case in the
setting of $H$-fractions \cite{Han2016Adv}, and the general case 
for rational functions \cite{Holtz2012T}.
In Section \ref{sec:shiftalg}, the lemma is applied not only to $F$,
but also to $(F-a_0)/x$ and to $xF$, which yields the $1$-shifted and
the backward shifted determinants.

\section*{On the use of artificial intelligence}

An artificial-intelligence assistant based on a large language model 
was used to help find some of the proofs.

\section*{Acknowledgments}
The second author was supported by the China Scholarship Council [No. 202506120133].

\end{document}